%% file: gauss-egregium-arxiv.tex
\documentclass[12pt,twoside,leqno,openany]{amsart}
\usepackage{amssymb,amsbsy,amsmath,amsfonts,amssymb,amscd,times,
graphics,color,xypic,footmisc,fancyhdr,multicol,fancybox,
graphicx,mathrsfs,rotating,ifthen,wasysym}
\usepackage[all]{xy}
\usepackage[french]{babel}
\usepackage[T1]{fontenc} 
\input macros.tex

\let\mathcal\mathscr

\begin{document}

$\:$

\bigskip\bigskip\medskip

\begin{center}

{\large\bf Courbure des surfaces dans l'espace:}

\medskip

{\large\bf le {\sl Theorema Egregium} de Gauss}

\bigskip

Jo\"el {\sc Merker}

\end{center}

\label{gauss-egregium}

\bigskip

\hfill
{\footnotesize\em
Nihil actum reputans si quid superesset agendum.}

\hfill
\REFERENCE{\footnotesize Devise de Gauss}

\Section{1.~Présentation de la {\sl formula egregia}}
\label{presentation-egregia}

\HEAD{1.~Présentation de la {\sl formula egregia}}{
Jo\"el Merker, D\'partement de Math\'matiques d'Orsay}

\medskip
En 1827--28, apr\`es plus de quinze ann\'ees de m\'editation, 
Gauss\footnote{\,
%%%%%%%%%%%%%%%%%%%%%%%-------DEBUT--------%%%%%%%%%%%%%%%%%%%%%%%%%%%
{\sc Gauss}, Carl Friedrich (Bunswick 1777 - G\"ottingen 1851),
math\'ematicien, astronome et physicien allemand, fils d'artisan issu
d'un milieu modeste, remarqu\'e pour la pr\'ecocit\'e de ses talents,
put conduire ses \'etudes gr\^ace \`a la protection et au soutien
du duc de Brunswick. Adolescent, il d\'eveloppe sa virtusosit\'e
technique en calculant des logarithmes et il corrige la table des
nombres premiers de G. von Vega et J. Lambert. \`A 19 ans, il
d\'emontre que le polygone r\'egulier \`a dix-sept c\^ot\'es inscrit
dans un cercle est constructible \`a la r\`egle et au
compas\,\,---\,\,r\'esultat tr\`es inattendu sur un sujet qui n'avait
pas fait de progr\`es substantiels depuis l'Antiquit\'e. En 1799, il
soutient sa th\`ese \`a l'universit\'e d'Helmstedt, donnant les
premi\`eres d\'emonstrations rigoureuses du th\'eor\`eme fondamental
de l'alg\`ebre, d'apr\`es lequel tout polyn\^ome non constant \`a
coefficients complexes poss\`ede au moins une racine r\'eelle ou
imaginaire. En 1801, \`a 24 ans, celui qu'on allait appeler le
<<\,prince des math\'ematiciens\,>> publie un ouvrage majeur de cinq
cents pages intitul\'e {\em Disquisitiones arithmeticae}, qui
incorpore une nouvelle th\'eorie des congruences, une th\'eorie des
formes quadratiques et dont le dernier chapitre est consacr\'e \`a
\'enoncer une condition n\'ecessaire et suffisante pour qu'un polygone
\`a un nombre premier $p$ de c\^ot\'es soit constructible \`a la
r\`egle et au compas. Nomm\'e directeur de l'observatoire de
G\"ottingen en 1807, il consacrera son \'energie \`a la m\'ecanique
c\'eleste, au magn\'etisme, \`a la th\'eorie des erreurs, poursuivant
occasionnellement ses travaux math\'ematiques, qu'il ne publiait pas
et dont il faisait \'etat par lettre et de mani\`ere confidentielle
\`a quelques interlocuteurs choisis.
} %%%%%%%%%%%%%%%%%%%%%%%%-----FIN-----%%%%%%%%%%%%%%%%%%%%%%%%%%%%%%%
publie un m\'emoire intitul\'e {\sl Disquisitiones generales circa
superficies curvas} qui signe l'acte de naissance de la
g\'eom\'etrie diff\'erentielle moderne. Il s'agit l\`a d'un des
articles les plus c\'el\`ebres de l'histoire des math\'ematiques,
article qui frappe encore le lecteur contemporain par la perfection
achev\'ee de son contenu, par la richesse de ses id\'ees novatrices et
par la nettet\'e de son architecture. Le point d'orgue de ce texte est
une formule math\'ematique centrale, complexe et remarquable qui aura
co\^ut\'e plusieurs ann\'ees de recherche au g\'enie calculatoire de
Gauss. Cette formule exprime la courbure d'une surface
bidimensionnelle plong\'ee dans l'espace tridimensionnel \`a l'aide de
donn\'ees m\'etriques purement intrins\`eques, lib\'er\'ees de la
troisi\`eme dimension suppl\'ementaire.

Plus pr\'ecis\'ement\footnote{\,
%%%%%%%%%%%%%%%%%%%%%%%-------DEBUT--------%%%%%%%%%%%%%%%%%%%%%%%%%%%
Les Sections qui suivent sont consacr\'ees \`a pr\'esenter 
la th\'eorie des courbes et des surfaces.
Ici, la notation standard $\frac{ \partial \Phi}{ \partial
x_i}$ désigne la dérivée partielle d'une fonction $\Phi = \Phi ( x_1,
\dots, x_n)$ par rapport à la variable $x_i$.
}, %%%%%%%%%%%%%%%%%%%%%%%%-----FIN-----%%%%%%%%%%%%%%%%%%%%%%%%%%%%%%%
en 1827, Gauss a d\'emontr\'e qu'\'etant donn\'e une surface
bidimensionnelle param\'etr\'ee par deux coordonn\'ees $(u,\, v)$ et
munie d'une m\'etrique infinit\'esimale
\begin{equation}
\label{ds2}
ds^2 
=
E\,du^2 
+ 
2F\,dudv 
+ 
G\,dv^2,
\end{equation}
o\`u $E$, $F$ et $G$ sont des fonctions des deux variables $u$ et de
$v$, il est toujours possible de calculer en tout point de
coordonn\'ees $(u,\, v)$ sa <<\,courbure\,>> $\kappa = \kappa(u,\,
v)$, fonction elle aussi des deux variables $u$ et $v$. La courbure
d'une surface en un point de coordonn\'ees $(u,\, v)$ est une notion
tridimensionnelle con\c cue par Gauss lui-m\^eme qui sera redéfinie en
termes g\'eom\'etriques ci-apr\`es. La formule par laquelle Gauss
calcule la courbure $\kappa$ est spectaculaire en ceci qu'elle ne fait
intervenir que la m\'etrique infinit\'esimale
bidimensionnelle~\thetag{ \ref{ds2}}, sans aucune r\'ef\'erence \`a
l'espace dans lequel est plong\'ee la surface; elle sera
appellée <<\,{\sl formula egregia}\footnote{\,
%%%%%%%%%%%%%%%%%%%%%%%-------DEBUT--------%%%%%%%%%%%%%%%%%%%%%%%%%%%
Cette d\'enomination n'est pas usuelle, mais elle fait \'echo \`a la
d\'enomination <<\,{\sl theorema egregium}\,>> (ou <<\,th\'eor\`eme
remarquable\,>>) que la post\'erit\'e a fix\'ee pour d\'esigner la
cons\'equence principale (et \'el\'ementaire) de la formule~\thetag{
\ref{formula-egregia}}, \`a savoir que la courbure d'une surface
resterait inchangée en chacun de ses points au cours de toute
déformation qui respecterait ses rapports métriques internes~\thetag{
\ref{ds2}}.
}\,>> %%%%%%%%%%%%%%%%%%%%%%%%-----FIN-----%%%%%%%%%%%%%%%%%%%%%%%%%%%%%
ou <<\,formule remarquable\,>>:
\begin{equation}
\label{formula-egregia}
\boxed{
\footnotesize
\aligned
\kappa 
=
& \ 
\frac{1}{4\,(EG-F^2)^2}
\bigg\{\,
E 
\bigg[
\frac{\partial E}{\partial v} 
\cdot
\frac{\partial G}{\partial v}
-
2\,\frac{\partial F}{\partial u} 
\cdot 
\frac{\partial G}{\partial v} 
+
\frac{\partial G}{\partial u}
\cdot
\frac{\partial G}{\partial u}
\bigg] 
+
\\
&
+
F\bigg[
\frac{\partial E}{\partial u}
\cdot\frac{\partial G}{\partial v}
-
\frac{\partial E}{\partial v} 
\cdot
\frac{\partial G}{\partial u} 
-
2\,\frac{\partial E}{\partial v} 
\cdot\frac{\partial F}{\partial v}
+
4\, 
\frac{\partial F}{\partial u} 
\cdot
\frac{\partial F}{\partial v} 
- 
2\,\frac{\partial F}{\partial u} 
\cdot\frac{\partial G}{\partial u}
\bigg] 
+
\\
& 
\ \ \ \ \ \ \ \ \ \ \ \ \ \ \ \ \ \ \ \ \ \ \ \ \ 
+
G\bigg[
\frac{\partial E}{\partial u} 
\cdot
\frac{\partial G}{\partial u} 
- 
2\,\frac{\partial E}{\partial u} 
\cdot 
\frac{\partial F}{\partial v} 
+
\frac{\partial E}{\partial v}
\cdot
\frac{\partial E}{\partial v}
\bigg]
+
\\
&
+2\,
\big(EG-F^2\big) 
\bigg[
-\frac{\partial^2 E}{\partial v^2} 
+
2\,\frac{\partial^2 F}{\partial u \partial v } 
- 
\frac{\partial^2 G}{\partial u^2}
\bigg]
\,\bigg\}.
\endaligned}
\end{equation}
Cette expression math\'ematique, publi\'ee en 1828, a \'et\'e
d\'ecouverte en 1827 grâce à un calcul virtuose.

Dans ce texte, il va être question de
reconstituer les sp\'eculations m\'ethodiques et les calculs formels
qui ont conduit Gauss, apr\`es de nombreuses tentatives infructueuses,
\`a d\'ecouvrir et \`a \'elaborer cette <<\,{\sl formula egregia}\,>>.
Son caract\`ere presque \'esot\'erique et la complexit\'e
quasi-incompressible des calculs formels qui sont n\'ecessaires pour
l'obtenir en ont fait un objet que l'on est parfois tent\'e de
rel\'eguer \`a un rang secondaire dans les manuels contemporains de
g\'eom\'etrie diff\'erentielle. Mais ici, loin d'\'ecarter les
difficult\'es calculatoires, il faudra accepter cette complexit\'e
interne afin de diss\'equer la d\'emonstration originale de Gauss en
la comparant aux simplifications partielles et relatives que la
post\'erit\'e a apport\'ees.

\medskip\noindent{\bf D\'efinition g\'eom\'etrique de la courbure.} 
\label{definition-courbure-courbes}
Le c\^ot\'e paradoxal et remarquable de
cette formule se dévoilera en 
se rapportant à la d\'efinition premi\`ere de
la courbure, due \`a Gauss lui-m\^eme. La courbure $\kappa$ d'une
surface math\'ematique plong\'ee dans l'espace est une notion
quantitative pr\'ecise qui t\'emoigne des ondulations
qu'une surface physique r\'eelle pr\'esente \`a l'{\oe}il dans
l'espace. En diff\'erents points d'une surface, la courbure est en
g\'en\'eral variable.

Pour la calculer, on introduit d'abord une sph\`ere fixe $\Sigma$ de
rayon $1$ (dans un syst\`eme pr\'ed\'efini d'unit\'es, {\em e.g.} de
rayon $1$ m\`etre, s'il s'agit de l'unit\'e de longueur qui est
standard en physique) centr\'ee en un point de l'espace, et appel\'ee
par Gauss {\sl sph\`ere auxiliaire}. Soit $S$ une surface quelconque,
par exemple celle qui est dessin\'ee sur la figure ci-dessous,
et soit $p$ l'un de ses points.

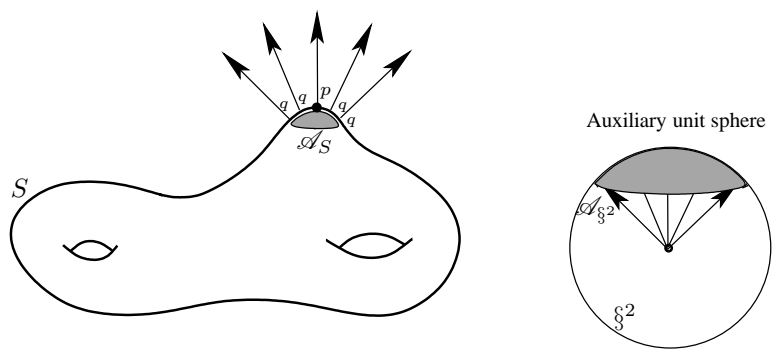
\begin{center}
\input application-de-Gauss.pstex_t
\end{center}

\noindent
On d\'elimite une petite r\'egion autour du point $p$ sur la surface
$S$\,\,---\,\,par exemple la r\'egion appel\'ee $\mathcal{ A}_S$ que
l'on a dessin\'ee en gris\'e sur la figure. En tout point $q$
de cette r\'egion $\mathcal{ A}_S$, on consid\`ere le vecteur
orthogonal \`a la surface \footnote{\,
%%%%%%%%%%%%%%%%%%%%%%%-------DEBUT--------%%%%%%%%%%%%%%%%%%%%%%%%%%%
Pour cela, on sous-entend que la surface a \'et\'e orient\'ee, au
moins localement autour du point $p$. Sur la figure,
l'orientation choisie est telle que l'ext\'erieur <<\,physique\,>> et
intuitif de la surface co\"{\i}ncide avec l'ext\'erieur
<<\,math\'ematique\,>>. En fait, la courbure reste inchangée si la
surface est munie de l'orientation opposée.
} %%%%%%%%%%%%%%%%%%%%%%%%-----FIN-----%%%%%%%%%%%%%%%%%%%%%%%%%%%%%%%
de longueur $1$, et on le d\'eplace parall\`element jusqu'\`a ce que
son origine co\"{\i}ncide avec le centre de la sph\`ere
auxiliaire. Les extr\'emit\'es de ces vecteurs d\'ecrivent alors une
r\'egion que l'on appelle $\mathcal{ A}_\Sigma$ sur la sph\`ere
auxiliaire $\Sigma$ et que l'on dessine en gris\'e sur la 
figure. Cette r\'egion $\mathcal{ A}_\Sigma$ d\'epend de la
r\'egion $\mathcal{ A}_S$; lorsque la r\'egion $\mathcal{ A}_S$ se
r\'etr\'ecit autour du point $p$, la r\'egion correspondante
$\mathcal{ A}_\Sigma$ diminue aussi autour du point de la sph\`ere
auxiliaire qui se situe \`a l'extr\'emit\'e du vecteur orthogonal \`a
la surface en $p$.

Pour fixer rigoureusement la terminologie, on appellera {\sl aire}
d'une surface ou d'une portion limit\'ee de surface le {\sl nombre qui
mesure son extension}, calcul\'ee dans un syst\`eme pr\'ed\'efini
d'unit\'es, {\em e.g.} en <<\,m\`etres-carr\'e\,>> (${\rm m}^2$) s'il
s'agit d'une mesure physique standard.

D'apr\`es la d\'efinition premi\`ere due \`a Gauss, la courbure
$\kappa (p)$ d'une surface $S$ en l'un de ses points $p$ est \'egale
\`a la limite\,\,---\,\,quand la r\'egion $\mathcal{ A}_S$ se
r\'etr\'ecit autour de $p$, ce qui sera noté <<\,$\mathcal{ A}_S \to
p$\,>>\,\,---\,\,d'un quotient dans lequel l'aire de la r\'egion
$\mathcal{ A}_S$ appara\^{\i}t au d\'enominateur, et dans lequel
l'aire de la r\'egion $\mathcal{ A}_\Sigma$ appara\^{\i}t au
num\'erateur ({\em cf.} à nouveau la figure ci-dessus):
\begin{equation}
\label{quotient-aires}
\kappa(p) 
=
\lim_{\mathcal{A}_S\,\to\,p}
\frac{\text{ \rm 
aire de la r\'egion}\
\mathcal{A}_\Sigma\
\text{\rm sur la sph\`ere auxiliaire}
}{
\text{\rm aire de la r\'egion}\ 
\mathcal{A}_S\ \text{\rm sur la surface}}.
\end{equation}
La sph\`ere unit\'e $\Sigma$ \'etant consid\'er\'ee comme la surface
canonique de r\'ef\'erence, cette formule exprime visiblement un {\em
principe de comparaison des aires <<\,courb\'ees\,>>}, valable dans
l'infiniment petit.

Sur la figure en question, la r\'egion $\mathcal{ A}_\Sigma$ est
plus grande que la r\'egion $\mathcal{ A}_S$, y compris lorsque la
r\'egion $\mathcal{ A}_S$ r\'etr\'ecit autour du point $p$. Il en
d\'ecoule que le quotient des aires et sa limite sont sup\'erieurs \`a
$1$: la surface est en effet plut\^ot <<\,assez courb\'ee\,>> au point
$p$. Comme on pourra le constater en exer\c cant son intuition, plus
la surface est <<\,courb\'ee\,>> ou <<\,pointue\,>> au voisinage du
point $p$, plus l'aire $\mathcal{ A}_\Sigma$ est grande, par rapport
\`a l'aire $\mathcal{ A}_S$. La plupart des surfaces\,\,---\,\,par
exemple, celle de la France\,\,---\,\,ont une courbure non nulle,
positive ou n\'egative, qui varie beaucoup d'un point \`a
l'autre\,\,---\,\,notamment en montagne\footnote{\,
%%%%%%%%%%%%%%%%%%%%%%%-------DEBUT--------%%%%%%%%%%%%%%%%%%%%%%%%%%%
On ne pense jamais assez souvent \`a la surface terrestre comme
r\'eservoir tr\`es riche d'exemples transcendant le caract\`ere
parfois exag\'er\'ement sobre de l'intuition math\'ematique.
}. %%%%%%%%%%%%%%%%%%%%%%%%-----FIN-----%%%%%%%%%%%%%%%%%%%%%%%%%%%%%%%
Au contraire, si la surface $S$ est un plan, tous les vecteurs qui lui
sont orthogonaux sont parall\`eles entre eux; la r\'egion $\mathcal{
A}_\Sigma$ est alors r\'eduite \`a un point, donc son aire est nulle
et le num\'erateur de la formule~\thetag{ \ref{quotient-aires}} est
toujours nul; par cons\'equent: {\em la courbure d'un plan dans
l'espace est nulle en tout point $p$}\,\,---\,\ fort heureusement!

\medskip\noindent{\bf Paradoxe remarquable.}
\label{paradoxe-remarquable}
De mani\`ere cruciale, pour d\'efinir la courbure d'une surface en
l'un de ses points, on utilise l'application de Gauss qui transporte
les vecteurs orthogonaux \`a la surface jusqu'\`a une sph\`ere
auxiliaire de rayon égal à $1$. De m\^eme qu'il suffit seulement de
deux coordonn\'ees pour rep\'erer un point quelconque du plan
euclidien, il suffit seulement de deux coordonn\'ees $(u,\, v)$ pour
rep\'erer un point quelconque $q$ de la surface proche de $p$.
Toutefois, il faut disposer aussi d'une troisi\`eme coordonn\'ee pour
repr\'esenter les vecteurs orthogonaux \`a la surface et pour
d\'efinir la r\'egion $\mathcal{ A}_\Sigma$ qui entre dans la
d\'efinition géométrique de la courbure. Il semble alors impossible de
parler de la courbure en se limitant aux deux seuls degr\'es de
libert\'e internes dont jouit une surface; la troisi\`eme dimension et
le plongement de la surface dans l'espace semblent absolument
indispensables; une saisie externe semble n\'ecessaire tant pour
l'appr\'ehension intuitive de la courbure que pour sa
d\'efinition. Ainsi, au premier abord, la courbure est une notion
<<\,extrins\`eque\,>>, c'est-\`a-dire ext\'erieure
\`a la surface, n'appartenant pas \`a son essence interne et semblant
d\'ependre fondamentalement de sa forme dans l'espace.

\CITATION{
Gauss reconnaissait lui-même que sa supériorité venait de son
entêtement.
\REFERENCE{Félix~{\sc Klein}}
}

Mais la perspicacit\'e de Gauss aura deviné la nature
<<\,intrins\`eque\,>> cach\'ee de cette notion de courbure: en effet,
gr\^ace \`a la {\em formula egregia}~\thetag{ \ref{formula-egregia}},
la courbure de la surface $S$\,\,---\,\,dont les points arbitraires
sont rep\'er\'es par les deux coordonn\'ees $(u,\, v)$\,\,---\,\,peut
\^etre saisie de mani\`ere interne et bidimensionnelle, sans sph\`ere
auxiliaire, sans vecteurs normaux \`a la surface, sans troisi\`eme
dimension: tel est le paradoxe remarquable. Il fallait donc toute la
pers\'ev\'erance de celui qu'on a appel\'e le <<\,{\em princeps
mathematicorum}\,>> pour parachever la d\'ecouverte de la courbure des
surfaces par l'\'elaboration d'une formule aussi complexe que~\thetag{
\ref{formula-egregia}}, gr\^ace \`a laquelle le concept de courbure
est envisageable sans aucun recours \`a une spatialit\'e externe.
\`A la fois pour l'exégète historien et pour le métaphysicien 
mathématicien, le myst\`ere est donc de comprendre, ne serait-ce
qu'{\em a posteriori} et pour l'`honneur' de la philosophie des
mathématiques, comment Gauss a su, par des encha\^{\i}nements
intuitifs et des méditations complexes, entrevoir l'existence de la
{\em formula egregia}, et comment, aussi, il a pu achever le tour de force
calculatoire qui exprime {\em intrinsèquement} l'extrinsèque
illusoire. 

Entendu comme anticipation des théories d'\'Elie Cartan qui traitent
plus systématiquement des problèmes d'équivalence entre structures
géométrico-différentielles, le {\sl Theorema egregium} montre que la
formation des concepts mathématiques tente de répondre à des tensions
organisationnelles dans des labyrinthes de calculs formels.

L'objectif analytique principal de cette ouverture sera
donc\,\,---\,\,sous toutes ses facettes\,\,---:

\smallskip
\centerline{\fbox{\ la gen\`ese conceptuelle et \underline{calculatoire} 
de la {\em formula egregia}~\thetag{\ref{formula-egregia}}.}}

\medskip\noindent{\bf Méthodologie expositionnelle.}
\label{methodologie-expositionnelle}
Les analyses écrites\footnote{\,
%%%%%%%%%%%%%%%%%%%%%%%-------DEBUT--------%%%%%%%%%%%%%%%%%%%%%%%%%%%
{\em Cf.}, pour l'intervention orale, la
Journ\'ee <<\,{\sf Philosophie et math\'ematiques}\,>>, organis\'ee par le
{\sl Centre international d'\'etude de la philosophie fran\c caise
contemporaine}\, le Samedi 24 Mai 2003 \`a l'\'Ecole Normale
Sup\'erieure de Paris, en Salle des Actes.
Programme:
{\sf 10h:} Ivahn {\sc Smadja} (ENS-Ulm): 
{\em Sch\`emes et structures: r\'eflexion sur le
savoir math\'ematique}; 
{\sf 11h:} Jo\"el {\sc Merker} (CNRS): 
{\em M\'etaphysique de l'ouverture math\'ematique};
{\sf 14h:} Quentin {\sc Meillassoux} (ENS-Ulm):
{\em Sur la port\'ee sp\'eculative des math\'ematiques:
le probl\`eme du fossile et le probl\`eme de Hume};
{\sf 15h:} Jean-Jacques {\sc Szczeciniarz} (Paris 7):
{\em Philosophie et m\'etaphysique de la finitude};
{\sf 16h.:} Alain {\sc Badiou} (ENS-Ulm): 
{\em Concepts \`a la lisi\`ere des math\'ematiques et de la
philosophie. Un exemple: l'Ouvert}.
} %%%%%%%%%%%%%%%%%%%%%%%%-----FIN-----%%%%%%%%%%%%%%%%%%%%%%%%%%%%%%% 
seront organisées de mani\`ere \`a rendre
intuitivement accessibles tous les concepts math\'ematiques
utilisés.

Except\'e une formation math\'ematique honnête et une pratique
minimale du calcul formel, nulle autre connaissance sp\'ecifique n'est
requise pour aborder la lecture de ce texte. Un style
<<\,litt\'eraire-spéculatif\,>> r\'esumera l'essentiel des
commentaires les plus techniques, en liaison fréquente avec une {\em
pratique spéculative plus large du calcul}. Certains paragraphes
techniques pourront \^etre mis entre parenth\`eses afin d'offrir une
lecture plus rapide, le Calcul étant par ailleurs régulièrement
analysé et replacé {\em au centre} de la pensée mathématique.

Ainsi, plusieurs types de discours seront associ\'es, assembl\'es et
coordonn\'es. L'expos\'e didactique, volontairement d\'evelopp\'e,
structur\'e tr\`es progressivement et agr\'ement\'e d'illustrations
g\'eom\'etriques \'epur\'ees donnera un aspect <<\,vulgarisation
scientifique\,>> aux analyses techniques. Toujours, la pr\'ecision et
la rigueur intuitive l'emporteront sur le discours math\'ematique
standard, souvent peu accessible \`a cause d'un exc\`es d'abstraction
et d'id\'ealisation. Toutefois, le caract\`ere descriptif
<<\,agr\'eable\,>> de certains passages ne devra pas faire oublier
qu'il existe des r\'esultats math\'ematiques extr\^emement pr\'ecis et
<<\,impeccables\,>> qui constituent le v\'eritable soubassement des
r\'esultats cit\'es. L'originalit\'e d'une telle approche r\'eside
probablement dans le fait que ce seront les calculs de Gauss {\em en
tant que calculs} qui seront analysés, sans céder à la 
tentation\,\,---\,\,récurrente dans l'`après-Hilbert'\,\,---\,\,d'un 
certain finalisme conceptuel. 

Dans ces analyses, le discours philosophique viendra
s'enchâsser\,\,---\,\,presque implicitement\,\,---\,\,\`a la mani\`ere
d'un <<\,r\'eveil\,>> des m\'etaphysiques internes et d'une reformulation
des causalit\'es conceptuelles, trop souvent masqu\'ees et
sous-entendues dans le discours math\'ematique usuel ou formel. Ainsi
seront privilégiées en premier lieu: la clart\'e intuitive,
l'accessibilit\'e des concepts et la lisibilit\'e des mouvements
dialectiques.

En somme, des r\'eflexions philosophico-math\'ematiques {\em
sp\'ecialis\'ees} seront élaborées afin de reformuler en langue
naturelle quelques-unes des innombrables dialectiques implicites dans
lesquelles circule le math\'ematicien en action. Des tr\'esors de
sp\'eculation aux formes insoup\c conn\'ees pars\`ement l'{\oe}uvre
publi\'ee de Gauss ainsi
que sa correspondance, et ces pensées qui sont
extrêmement précieuses pour accéder à la réalité mathématique
adéquate sont nées dans le berceau arithmétiques
de ses {\em calculs} virtuoses.

\CITATION{
Gauss est attiré par l'art du calcul et compte
sans se lasser avec une énergie infatigable. 
Et c'est grâce à ces exercices incessants de manipulation
des nombres (notamment des nombres décimaux
avec un nombre incroyable de chiffres
après la virgule) qu'il acquiert cette faculté
étonnante de tout calculer, que sa mémoire
devient capable d'emmagasiner un nombre de données
extraordinaire et de jongler avec les chiffres
comme personne ne le fit jamais ni ne le fera par la suite.
\`A force de manier les chiffres, il découvre
<<\,expérimentalement\,>> les principales lois; cette
méthode qui va à l'encontre des principes professés aujourd'hui
en mathématiques s'était répandue au XVIII\textsuperscript{ème}
siècle, on la rencontre par exemple chez Euler. 
\REFERENCE{Félix~{\sc Klein}.}}

\Section{2.~Courbes mathématiques dans le plan}
\label{courbes-dans-le-plan}

\HEAD{2.~Courbes mathématiques dans le plan}{
Jo\"el Merker, D\'partement de Math\'matiques d'Orsay}

\medskip\noindent{\bf Filaments et trajectoires.}
\label{filaments-trajectoires}
Deux cat\'egories d'objets math\'ematiques sugg\'er\'es et transmis
par l'exp\'erience physique poss\`edent une caract\'eristique que l'on
appelle commun\'ement la <<\,courbure\,>>: ce sont les {\sl lignes
courbes} et les {\sl surfaces courbes}. Cette caract\'eristique
s'oppose
\`a la <<\,droiture\,>> parfaite qui se manifeste
ext\'erieurement dans l'espace, \`a la vue des lignes droites et des
surfaces planes dont regorge l'architecture classique et
contemporaine.

Si l'on fixe un plan math\'ematique abstrait dans lequel les courbes
sont trac\'ees, on parle alors de {\sl courbes planes}. Dessin\'ees
\`a l'int\'erieur d'une telle surface absolument droite et sans aucune
\'epaisseur, elles ont néanmoins toute la libert\'e de se courber, de
s'enrouler et de s'entortiller avec souplesse.

Gr\^ace \`a un long processus d'id\'ealisation conceptuelle, les courbes
{\em math\'ematiques}\, planes parviennent \`a incarner
abstraitement la notion intuitive, sugg\'er\'ee par l'exp\'erience, de
{\em ficelle continue et infimement fine}, d\'epos\'ee par exemple sur
une table. Elles parviennent aussi \`a incarner abstraitement
la trajectoire d'un point mat\'eriel\,\,---\,\,poussi\`ere, mol\'ecule
aromatique ou atome ionis\'e\,\,---\,\,se d\'epla\c cant dans un
plan (ou dans l'espace),
autre exemple physique simple de courbe ayant une \'epaisseur
infime.

Les math\'ematiques fondamentales sont donc parvenues \`a id\'ealiser
de tels exemples; elles ont \'elabor\'e en leur sein tout un arsenal
de concepts abstraits pour se libérer des intuitions physiques et
briller de tous leurs feux en apportant des r\'eponses parfaitement
rigoureuses \`a des probl\'ematiques autonomes. 

Toutefois, la justesse partielle d'une telle <<\,image d'\'Epinal\,>>
ne doit jamais faire oublier que les myst\`eres historiques et
philosophiques abondent dans toute th\'eorie achev\'ee et qu'il existe
encore bien des aspects intuitifs, transmis imparfaitement par
l'exp\'erience physique du monde, qui attendent leur conceptualisation
math\'ematique.

\medskip\noindent{\bf Identification du continu et coordonnées.} 
\label{identification-continu}
En g\'eom\'etrie analytique contemporaine,
toute ligne droite id\'eale est {\em identifi\'ee} au {\sl continu
fondamental unidimensionnel}, c'est-\`a-dire \`a l'ensemble, not\'e
$\R$, de tous les nombres $x$ possibles, dits {\sl nombres r\'eels}.
Ces nombres sont tous ceux que l'on peut \'ecrire en base 10 avec une
infinit\'e de d\'ecimales. Ainsi l'essence g\'eom\'etrique d'une
droite est-elle identifi\'ee \`a l'essence num\'erique du continu. Ce
dernier est ordonn\'e lin\'eairement, il est complet et il n'a aucun
<<\,trou\,>> en lui-m\^eme, exactement comme pourrait l'\^etre une
ligne droite id\'eale dont la conception aurait \'et\'e approfondie
par l'intuition.

En v\'erit\'e, l'un des plus profonds ab\^{\i}mes \'epist\'emologiques
aux jalons historiques multiples s'ouvre au c{\oe}ur d'une telle
identification. Il s'agit l\`a d'un v\'eritable <<\,tour de
passe-passe\,>> m\'etaphysique qui magn\'etise encore l'attention de
la philosophie des math\'ematiques. Mais pour simplifier, sans
\'evoquer les difficult\'es sp\'eculatives auxquelles conduit cette
identification, on considérera les nombres r\'eels comme connus, et on
admettra que toute droite est en correspondance biunivoque avec la
succession des nombres r\'eels.

Par d\'efinition, le plan math\'ematique id\'eal est identifi\'e \`a
la collection de tous les couples $(x,\, y)$ de {\em deux}\, nombres
r\'eels quelconques $x$ et $y$. Chaque point du plan est rep\'er\'e par
un unique couple $(x,\, y)$ de nombres r\'eels appel\'es {\sl
coordonn\'ees} de ce point. Dans ce plan num\'erique, un point 
particulier et deux droites sp\'eciales m\'eritent attention:

\begin{itemize}

\smallskip\item[$\bullet$]
l'origine, point central par excellence, dont les deux coordonn\'ees
$(0,\, 0)$ sont nulles;

\smallskip\item[$\bullet$]
la droite qui est constitu\'ee de tous les points de coordonn\'ees
$(x,\, 0)$, o\`u $x\in \R$ est un nombre r\'eel quelconque; cette
premi\`ere droite est <<\,{\sl horizontale}\,>> sur la figure; on
l'appelle {\sl axe des $x$};

\smallskip\item[$\bullet$]
la droite qui est constitu\'ee de tous les points de coordonn\'ees
$(0,\, y)$, o\`u $y\in \R$ est un nombre r\'eel quelconque; cette
deuxi\`eme droite est <<\,{\sl verticale}\,>> sur la figure; on
l'appelle {\sl axe des $y$}.

\end{itemize}\smallskip

\begin{center}
\input plan.pstex_t
\end{center}

On oriente ces deux droites dans le sens des $x$ et des $y$
croissants, d'o\`u les fl\`eches sur la figure. G\'eom\'etriquement
parlant, les deux coordonn\'ees $x$ et $y$ d'un point $(x,\, y)$ sont
obtenues par projection sur la droite horizontale et sur la droite
verticale, respectivement. Ainsi, travaille-t-on dans ce plan
math\'ematique id\'eal: on peut y tracer des droites vari\'ees et des
courbes diverses. Or les deux nombres r\'eels $x$ et $y$ peuvent avoir
un nombre infini de d\'ecimales ils sont toujours sans aucune
<<\,\'epaisseur\,>> physique ou math\'ematique. Donc comme l'intuition
semble effectivement le deviner sur une figure, les courbes id\'eales
seront d'une finesse infinie dans ce plan cart\'esien abstrait. Mais
alors, {\em comment d\'efinit-on math\'ematiquement une courbe}?

\medskip\noindent{\bf Trois saisies analytiques
des courbes dans le plan.} 
\label{trois-saisies-courbes}
Il existe essentiellement trois fa\c cons
de d\'efinir math\'ematiquement une courbe, comme l'illustreront les
trois figures qui vont suivre. La numérotation {\bf I}, {\bf II} et
{\bf III} de ces trois d\'efinitions est celle que Gauss a choisie
dans~\cite{ ga1828} pour pr\'esenter la th\'eorie des surfaces;
ultérieurement, on constatera en effet qu'elles se g\'en\'eralisent
aisément aux surfaces sises dans l'espace des $x, y, z$.

La repr\'esentation {\bf III} d'une courbe sous forme de graphe $y =
f(x)$ est celle que l'on apprend le plus couramment dans
l'enseignement secondaire, mais on commencera par la repr\'esentation
{\bf I} sous la forme d'une \'equation implicite, qui est
la plus immédiatement riche sur le plan fonctionnel.

\medskip\noindent{\small\sf Repr\'esentation implicite I.} 
Dans ce premier cas, on représente une courbe comme le lieu de tous
les points $(x_p, y_p)$ du plan en lesquels s'annule une certaine
fonction $W(x,\, y)$ des deux variables $(x,\, y)$, mais toutes les
fonctions ne conviennent pas, il est nécessaire d'en exclure quelques
unes qui sont trop dégénérées ou pas assez régulières. Du point de vue
de la théorie des ensembles, la courbe est donc l'ensemble des points
$(x,\, y)$ tels que $W(x,\, y) = 0$. Par exemple, l'axe des $x$ est
repr\'esent\'e par l'\'equation $y=0$; de m\^eme, l'axe des $y$ est
repr\'esent\'e par l'\'equation $x=0$; un cercle de rayon $1$ centr\'e
\`a l'origine est tout simplement repr\'esent\'e par l'\'equation $x^2
+ y^2 - 1=0$; enfin, une {\sl ellipse g\'en\'erale} est
repr\'esent\'ee par l'\'equation
\[
\frac{(x-x_1)^2}{a^2}
+
\frac{(y-y_1)^2}{b^2} 
- 
1 
=
0, 
\]
o\`u $(x_1,\, y_1)$ est un point quelconque du plan et o\`u $a >0$ et
$b>0$ sont des nombres r\'eels arbitraires\footnote{\,
%%%%%%%%%%%%%%%%%%%%%%%-------DEBUT--------%%%%%%%%%%%%%%%%%%%%%%%%%%%
Sur la figure, le point $(x_1,\, y_1)$ est au centre de
l'ellipse; autour d'un point quelconque $q$ de la courbe, on a
dessin\'e un petit carr\'e $\Delta_q$ pour signifier qu'on {\sl
localise} la courbe dans ce petit voisinage. Cet {\sl acte de
localisation} qui possède un statut à la fois mental et conceptuel est
fondamental en géométrie, parce qu'il permet de focaliser l'attention
de l'esprit sur une fraction de la courbe en {\em oubliant} sa forme
globale.
}. %%%%%%%%%%%%%%%%%%%%%%%%-----FIN-----%%%%%%%%%%%%%%%%%%%%%%%%%%%%%%%

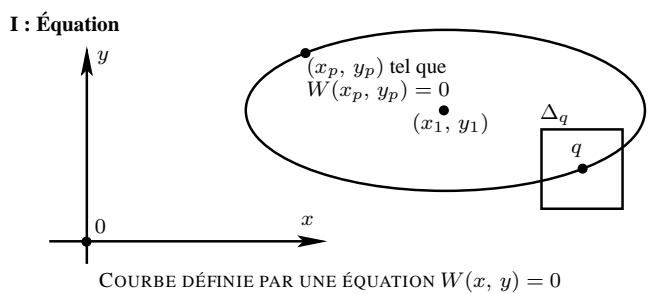
\smallskip
\begin{center}
\input equation.pstex_t
\end{center}

Dans les traités de géométrie différentielle, on dit que la
repr\'esentation par l'\'equation $W(x,\, y) = 0$ est {\sl implicite},
parce que ni $y$ n'est exprim\'e explicitement en fonction de $x$, ni
$x$ n'est exprim\'e explicitement en fonction de $y$. Mais en
v\'erit\'e, dans l'\'equation $W (x,\, y)= 0$, il existe une
d\'ependance implicite cach\'ee entre $x$ et $y$ gr\^ace \`a laquelle
on a effectivement affaire (en général) \`a un objet unidimensionnel,
c'est-à-dire \`a une vraie courbe.

\medskip\noindent{\bf Dialectique de l'ontologie imprécise.}
\label{dialectique-ontologie-imprecise}
En tout cas, il importe avant tout d'observer que {\em toute la
diversit\'e et la variabilit\'e des courbes ainsi définies repose sur
le choix de la fonction $W(x,\, y)$}. Et il va sans dire qu'il existe
un très grand nombre de telles fonctions de deux variables $x$ et $y$,
puisque par exemble tous les polynômes sont concernés, aussi bien que
toutes les fonctions analytiques, toutes les fonctions
différentiables, et même en un certain sens toutes les fonctions
fractales qui satisfont certaines conditions de non-dégénérescence.
Il est vrai aussi qu'il est essentiellement impossible de se
représenter par la pensée dans une généralité pure ce que seraient ces
fonctions $W ( x, y)$, polynomiales, analytiques, différentiables ou
fractales, afin d'embrasser véritablement toute l'extension du concept
de courbe. Mais par ailleurs et du point de vue de l'intuition
géométrique en recherche de pensée\,\,---\,\,laquelle tente de se
représenter les courbes possibles en exerçant de délicats
efforts afin de variabiliser mentalement la mobilité idéale de
filaments sans épaisseur dans le plan\,\,---, la richesse ontologique
de l'Idée mathématique de courbe plane semble se dessiner sans
parvenir à se réaliser pleinement par une saisie définitive. Toujours
duelle, la géométrie fonctionnelle, ou `fonctionnellité' géométrique,
ne parvient jamais à mobiliser
la mobilité absolue des êtres idéels.

\medskip\noindent{\bf Thèse fondamentale.}
{\sf [Dialectique de l'ontologie imprécise]}
{\em Ni le fonctionnel {\em pur}, ni le géométral {\em pur}
ne se départissent du caractère éminemment {\em potentiel} 
de la saisie qu'ils proposent.}

\medskip
Autrement dit et à un niveau plus
élevé, l'extension ontologique, en mathématiques,
n'est jamais donnée comme un monde
(erreur du <<\,cela est là\,>>), ni postulable par pétition 
de principe (erreur du positivisme symbolique), 
mais elle demeure éternellement
contrainte de respecter l'universalité
catégorique du potentiel, 
forcée par les labyrinthes à ne se déployer que partiellement
dans des imperfections non prévisibles
qui lui deviennent propres.

\medskip\noindent{\small\sf Repr\'esentation graph\'ee III.}
Contrairement à 
la représentation implicite, 
la repr\'esentation graph\'ee consiste \`a exprimer
l'une des deux variables, $y$ par exemple, comme une fonction {\em
explicite} de l'autre variable, ici $x$. La courbe est alors le lieu
des points $(x,\, f(x))$ du plan, o\`u $x$ varie et o\`u $f(x)$ est
une fonction d\'etermin\'ee.

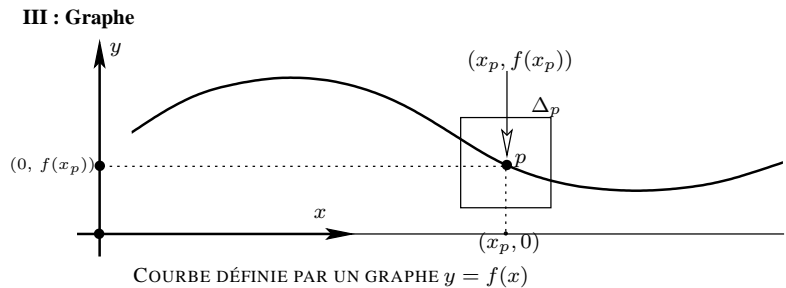
\medskip

\begin{center}
\input graphe.pstex_t
\end{center}

Cette troisième saisie (présentée ici dans un second temps)
pré-oriente donc {\em par la pensée} la représentation {\em
unidimensionnelle} de la courbe afin de {\em briser immédiatement}
l'ambiguïté inorientée de la liaison du $x$ et du $y$ à travers
l'équation $W ( x, y) = 0$. D'une certaine façon, la représentation
graphée {\em refuse} la généralité initiale du concept et se {\em projette
d'emblée dans l'unidimensionnalité décidée qu'elle rend expressément
visible}.
 
Sur la figure ci-dessus, on voit aisément que le nombre $f(x)$
mesure la <<\,hauteur\,>> du point de la courbe qui se situe \`a la
verticale du point $x$ de l'axe horizontal. Lorsqu'on se d\'eplace
horizontalement sur l'axe des $x$, la hauteur $f(x)$ varie tandis que
le filament infiniment fin se d\'eforme contin\^ument et ondule
librement.

\`A nouveau, {\em toute la diversit\'e et toute la variabilit\'e des
courbes reposent sur le choix de la fonction $f(x)$}, mais cette
fois-ci, au lieu 
de considérer une fonction de deux variables $W ( x, y)$, on a
affaire à une fonction d'une seule variable $f(x)$ et il s'avère alors
que l'extension ontologique du concept (au moins local) de courbe dans
le plan s'identifie mieux avec la collection de toutes les fonctions
possibles d'une seule variable.

\Scholie
Cette dernière affirmation est à nuancer car les identifications ne
sont pas absolument exactes sans faire d'hypothèses mathématiques
légèrement restrictives. En tout état de cause, si l'on ne choisit
pas d'effectuer de telles hypothèses, l'ontologie se ramifie et se
diversifie (pour le plus grand bien des mathématiques), en tant que le
définitionnel multiple ouvre vers des objets de type 
<<\,courbe\,>> non
équivalents entre eux et non univoquement saisissables.
<<\,Comment saisir?\,>> demeure une question mathématique
permanente.
\qed

\medskip
En v\'erit\'e, \`a partir
de l'\'equation $W(x,\, y) =0$, il est possible de reconstituer une
certaine d\'ependance cach\'ee entre $x$ et $y$ gr\^ace au
th\'eor\`eme dit <<\,des fonctions implicites\,>>, dont l'\'enonc\'e
est le suivant. \`A dessein, ce th\'eor\`eme ne sera pas énoncé ici
formellement, c'est-à-dire avec toute la rigueur qui si\'erait \`a un
texte math\'ematique technique, puisque cela
pourrait freiner l'expression des intuitions
métaphysiques archaïques\,\,---\,\,lesquelles demeurent
nécessairement coprésentes dans toute approche qui exige
d'elle-même un 
cadencement par la rhétorique de la rigueur. 

Soit donc $p$ un point de coordonn\'ees $(x_p,\, y_p)$ appartenant
\`a la courbe. Si l'on suppose qu'une 
(au moins) des deux d\'eriv\'ees
partielles\footnote{\,
%%%%%%%%%%%%%%%%%%%%%%%-------DEBUT--------%%%%%%%%%%%%%%%%%%%%%%%%%%%
Classiquement, la d\'eriv\'ee $\frac{ df}{ dx}(x_p)$ d'une fonction
$f(x)$ d'une variable $x$ en un point $x_p$ est la limite, quand
$\varepsilon$ tend vers z\'ero, du quotient $\frac{ f(x_p+
\varepsilon) - f(x_p)}{ \varepsilon}$, si cette limite existe. Pour
les fonctions $W(x,\, y)$ de deux variables, la d\'efinition est
analogue: la d\'eriv\'ee partielle $\frac{ \partial W}{\partial y}
(x_p,\, y_p)$ est la limite, quand $\varepsilon$ tend vers z\'ero, du
quotient $\frac{ W( x_p,\, y_p+ \varepsilon) - W(x_p,\, y_p)}{
\varepsilon}$, si cette limite existe. De m\^eme et bien évidemment
aussi, pour calculer $\frac{ \partial W}{ \partial x} (x_p,\, y_p)$,
on \'etudie l'existence de la limite du quotient $\frac{ W(x_p+
\varepsilon,\, y_p) - W(x_p,\, y_p)}{ \varepsilon}$ lorsque
$\varepsilon$ tend vers zéro.
} %%%%%%%%%%%%%%%%%%%%%%%%-----FIN-----%%%%%%%%%%%%%%%%%%%%%%%%%%%%%%%
suivantes:
\[
\frac{\partial W}{\partial y}(x_p,\,y_p)\ \ \
\text{\scriptsize\sf [premier cas]}
\ \ \ \ \ \ \
\text{\rm ou}
\ \ \ \ \ \ \
\frac{\partial W}{\partial x}(x_p,\,y_p)\ \ \ 
\text{\scriptsize\sf [deuxi\`eme cas]}
\]
ne s'annule pas, alors dans un voisinage suffisamment petit du point
$p$, l'ensemble des points $(x,\, y)$ tels que $W(x,\, y)= 0$ peut
\^etre repr\'esent\'e explicitement soit 
par une \'equation graph\'ee le long de l'axe des $x$ de la forme $y=
f(x)$ [premier cas] soit par une \'equation graph\'ee le long de l'axe
des $y$ de la forme $x = g(y)$ [deuxi\`eme cas]. Lorsque les deux
d\'eriv\'ees partielles en question ci-dessus ne s'annulent pas (ce
qui se produit très souvent), les deux repr\'esentations graph\'ees
$y= f(x)$ et $x= g(y)$ sont possibles. Sur la figure,
un tel petit voisinage de $p$ a été noté
$\Delta_p$.

Dans les tous les cas, les fonctions $g(y)$ et $f(x)$ d\'ependent de
la fonction $W(x,\, y)$. Il existe un th\'eor\`eme math\'ematique
et des algorithmes num\'eriques impl\'ement\'es sur ordinateur, qui
permettent de construire de telles fonctions $f(x)$ et $g(y)$ 
\`a partir de la donn\'ee de base qu'est la fonction $W(x,\, y)$.

\Scholie
La n\'ecessit\'e de tels algorithmes provient du fait qu'en
g\'en\'eral, la liste des fonctions <<\,connues\,>> ou
<<\,classiques\,>> ne contient pas toutes les fonctions $g(y)$ et
$f(x)$ que l'on peut obtenir de la sorte \`a partir des fonctions
$W(x,\, y)$ <<\,classiques\,>>. Par exemple, la fonction $y = f(x)$
implicitement d\'efinie par $\sin (y) + y^2 = x$ au voisinage de
$(x_p,\, y_p) = (0,\, 0)$ n'est pas <<\,classique\,>> et ne poss\`ede
pas de d\'enomination particuli\`ere, bien que la fonction $W(x,\, y)
:= \sin (y)+ y^2 - x$ soit une fonction <<\,classique\,>>.

En fait, une multitude innombrable de fonctions sans nom apparaissent
lorsque l'on veut r\'esoudre $y$ (ou $x$) en fonction de $x$ (ou de
$y$) dans les diverses \'equations implicites $W(x,\, y)= 0$. Tout
domaine de fonctions <<\,classiques\,>>, <<\,explicites\,>> ou
<<\,concr\`etes\,>>, tel que l'univers des fonctions alg\'ebriques et
des fonctions circulaires, est imm\'ediatement \'epuis\'e et
d\'epass\'e par l'application de ce th\'eor\`eme. C'est pourquoi le
th\'eor\`eme des fonctions implicites poss\`ede un caract\`ere
m\'etaphysique sp\'ecial, en tant qu'il impose par nature le
d\'epassement de toute donation premi\`ere, pens\'ee comme explicite
ou concr\`ete. On comprend qu'il ait fallu attendre le vingti\`eme
si\`ecle pour que ce th\'eor\`eme soit formul\'e comme un \'enonc\'e
math\'ematique pourvu d'une d\'emonstration rigoureuse (théorème du
point fixe {\em via} des itérations <<\,à la Picard\,>>) qui cèle en
fait l'incroyable extension ontologique que le procédé impose. En
effet, il fallait accepter un nouvel univers d'existence abstrait que
l'on appelle maintenant la classe des {\sl fonctions
diff\'erentiables}, et accepter des moyens de travail concrets mais
limités sur le plan de l'explicitation des fonctions impliqués. Avec
l'apparition du théorème des fonctions implicites dans des contextes
non résolubles par des formules finies, il fallait en quelque sorte
accepter que le monde des êtres se réduise à la plus simple expression
de leur différentiabilité, sans référence à des caractéristiques
individuantes, structurales, internes. Quant à une justification
fondationnelle, on peut s'accorder pour dire que de tels
êtres\,\,---\,\,les fonctions que l'on peut extraire d'une relation
implicite donnée explicitement, voire même aussi implicitement\,\,
sont en dernier recours de nature purement num\'erique, à savoir des
correspondances univaluées entre nombres decimaux. 

Ainsi, pour que les consid\'erations pr\'ec\'edentes, ainsi que
celles qui suivront, aient un sens, certaines hypoth\`eses de
r\'egularit\'e doivent \^etre satisfaites par toutes les fonctions
impliqu\'ees dans la discussion. 
Mais afin de rester strictement dans le
sujet du calcul, 
on laissera de c\^ot\'e les consid\'erations historiques et
philosophiques soulev\'ees par la recherche de ces hypoth\`eses, et
on supposera toujours que les fonctions utilisées 
sont diff\'erentiables jusqu'\`a un ordre suffisant.
\qed

\smallskip
En conclusion, dans la repr\'esentation implicite et dans la
repr\'esentation graph\'ee, il faut retenir l'id\'ee qu'il {\em existe
une d\'ependance fonctionnelle\,\,---\,\,implicite ou
explicite\,\,---\,\,entre les deux variables $x$ et $y$}. Cette
d\'ependance ne laisse plus qu'un seul degr\'e de libert\'e pour la
mobilit\'e et le d\'eplacement interne le long de la courbe. Une telle
contrainte confirme l'intuition fondamentale d'apr\`es laquelle les
courbes sont par nature unidimensionnelles.

\medskip\noindent{\bf Caractère intrinsèque
de la repr\'esentation param\'etr\'ee II.}
\label{caractere-intrinseque-parametree}
Enfin, la dernière possibilit\'e consiste \`a faire d\'ependre
chacune des deux variables $x$ et $y$ d'une variable $t$ auxiliaire
unique. Cette repr\'esentation correspond \`a l'intuition d'une
courbe param\'etr\'ee par un param\`etre temporel: le point de
coordonn\'ees $(x(t),\, y(t))$ se d\'eplace alors dans le plan <<\,sur
la courbe\,>> lorsque le <<\,temps\,>> $t$ s'\'ecoule. Cette
deuxi\`eme mani\`ere de voir inscrit d'embl\'ee l'unidimensionnalit\'e
au c{\oe}ur de la saisie de la courbe.

Il est vrai que l'unidimensionnalit\'e est aussi inscrite au c{\oe}ur
de la repr\'esentation graph\'ee $y= f(x)$, puisque seule la variable
$x$ varie, mais cette repr\'esentation pr\'esente l'inconv\'enient de
privil\'egier, sans raison suffisante, l'une des deux coordonn\'ees,
en l'occurence ici, la premi\`ere coordonn\'ee. Dans certains cas, il
faut plut\^ot privil\'egier la coordonn\'ee $y$. \`A cet \'egard
d'ailleurs, par rapport \`a la repr\'esentation graph\'ee, la
repr\'esentation implicite poss\'edait l'avantage de ne privil\'egier
aucune des deux coordonn\'ees, bien qu'elle e\^ut le d\'efaut de
laisser dans l'ombre la relation de d\'ependance pr\'ecise entre $x$
et $y$.

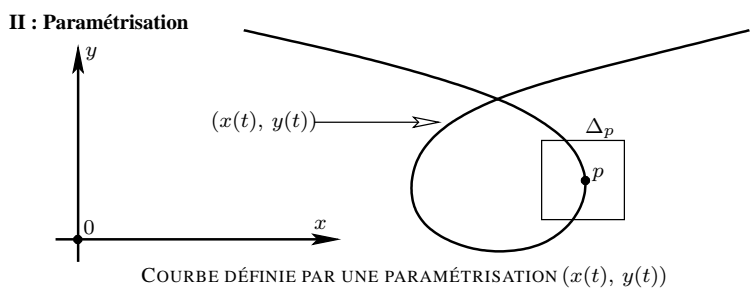
\begin{center}
\input parametrisation.pstex_t
\end{center}

{\em La repr\'esentation param\'etr\'ee {\bf
II} est la plus ad\'equate des trois}. En effet, elle ne poss\`ede
aucun des deux d\'efauts sus-mentionn\'es: premi\`erement, dans
l'expression $(x(t),\, y(t))$, l'unidimensionnalit\'e est clairement
exprim\'ee; et deuxi\`emement, les deux variables sont trait\'ees sur
le m\^eme plan sans que l'une d'entre elles soit privil\'egi\'ee.
Ici, la variable $t$ est un param\`etre {\em interne}\, pour la
courbe: {\em l'intrins\`eque s'inscrit d'embl\'ee dans la
repr\'esentation param\'etr\'ee}.

Dans les paragraphes consacr\'es aux surfaces sera formulée
une d\'efinition analogue\,\,---\,\,aussi num\'erot\'ee
{\bf II}, d'apr\`es la d\'enomination originale de
Gauss\,\,---\,\,dans laquelle ce seront deux param\`etres {\em
internes} $(u,\, v)$ qui permettront de repr\'esenter la surface.

\medskip\noindent{\bf \'Equivalence locale 
entre les trois repr\'esentations.} 
\label{equivalence-locale-trois-representations}
Toute courbe peut \^etre {\sl localis\'ee} au voisinage de l'un de ses
points $p$. Cette op\'eration consiste \`a choisir un voisinage
suffisamment $\Delta_p$ de $p$\,\,---\,\,le petit carr\'e sur les
figures précédentes\,\,---\,\,et \`a restreindre la
consid\'eration de la courbe dans ce petit voisinage.
Math\'ematiquement parlant, on d\'emontre l'\'enonc\'e suivant.

\smallskip\noindent{\bf Théorème.} 
{\em Du point de vue local, les trois repr\'esentations {\bf I}, {\bf
II} et {\bf III} sont \'equivalentes.
}\medskip

Ces \'equivalences seront utiles, puisque dans toute la suite des
analyses, les consid\'erations seront exclusivement locales.
Autrement dit, {\em dans un voisinage suffisamment petit de l'un de ses
points, une courbe peut \^etre repr\'esent\'ee aussi bien par une
\'equation, que par une param\'etrisation ou par un graphe}. 

La figure suivante diagrammatise ce th\'eor\`eme, au moyen de
trois grandes fl\`eches doubles dessinant un triangle; chacune de ces
fl\`eches correspond \`a un raisonnement math\'ematique rigoureux;
par exemple, la fl\`eche oblique situ\'ee \`a l'extr\^eme gauche de la
figure correspond au th\'eor\`eme des fonctions implicites.

\begin{center}
\input univers-courbes.pstex_t
\end{center}

Les trois doubles fl\`eches creuses situ\'ees au centre du diagramme
indiquent trois <<\,transsubstantiations\,>> math\'ematiques
remarquables et \'equivalentes entre elles: l'id\'ee g\'eom\'etrique
de courbe se transmue, de trois mani\`eres diff\'erentes, en l'id\'ee
de fonction math\'ematique; on passe de l'univers g\'eom\'etrique \`a
l'univers fonctionnel, de telle sorte que du point de vue technique,
ce ne sont plus que les fonctions qui existent. La suite des
\'ev\'enements montrera que les questions qui \'etaient
g\'eom\'etriques au d\'epart se transforment rapidement en des
probl\`emes purement alg\'ebriques ou analytiques qui portent sur un
ensemble de fonctions ainsi que sur la collection de leurs d\'eriv\'ees.

\Section{3.~Préliminaire sur les surfaces.}
\label{preliminaire-surfaces}

\HEAD{3.~Préliminaire sur les surfaces}{
Jo\"el Merker, D\'partement de Math\'matiques d'Orsay}

\medskip\noindent{\bf Illustrations.}
\label{illustrations}
Dans l'espace euclidien standard de dimension trois, tel que l'espace
<<\,physique\,>> de la m\'ecanique newtonienne, les surfaces peuvent
revêtir des formes innombrables, tant le monde physique en regorge. Du
point de vue de l'idéalisation mathématique cependant, seules les
structures générales comptent, et les morphologies se limitent à
quelques échantillons de surfaces-types, intéressantes à titre de
modèles ou en vertu de théorèmes abstraits de classification.

\begin{center}
\input echantillon-surfaces.pstex_t
\end{center}

Or il n'est question dans ce texte que des propriétés locales de
courbure (gaussienne) qu'ont les surfaces, et dans ces conditions, le
deuxième schéma géométrique (archaïque) ci-dessus, à savoir un
<<\,petit morceau\,>> de surface (de tissu, de feuille, de peau)
suffira amplement pour diagrammatiser le concept de surface locale.
Bien entendu, tout petit morceau découpé par la pensée dans l'une
quelconque des cinq autres surfaces dessinées conviendrait aussi. La
surface peut alors être imaginée comme le bord d'un solide plein,
comme la limite entre deux profondeurs-épaisseurs (l'air et le
solide), ou comme l'interface sans épaisseur entre deux milieux
physiques intrinsèquement cohésifs. Elle devient donc pour la pensée
mathématique <<\,instinctive\,>> une généralisation naturelle et
presque évidente du concept de courbe, puisqu'il suffit d'admettre que
la dimensionnalité d'extension de l'objet initial en lui-même (la
courbe) augmente d'une unité par l'effet d'une mobilité douée d'un
degré de liberté, et que simultanément aussi, la dimensionnalité
d'extension dans lequel est plongé l'objet s'élargit aussi par
spatialisation. Du point de vue de la métaphysique archaïque
des mathématiques, l'<<\,un-deux\,>> en s'étoffe en un <<\,deux-trois\,>>.

\medskip\noindent{\bf Trois saisies analytiques des surfaces 
dans l'espace.}
\label{trois-saisies-surfaces}
Aux trois représentations possibles d'une courbe dans le plan
correspondent alors de manière complètement analogue
exactement trois possibilités de représenter une surface
et on démontre, comme dans le cas des courbes, que ces trois
manières de représenter les surfaces sont équivalentes entre
elles, d'un point de vue mathématique rigoureux.

\begin{itemize}

\smallskip\item[{\bf I:}] {\sf Représentation implicite}.

\smallskip\item[{\bf II:}] {\sf Représentation paramétrée}. 

\smallskip\item[{\bf III:}] {\sf Représentation graphée}. 

\end{itemize}\smallskip

\CITATION{
Il y a deux méthodes générales pour l'étude des propriétés d'une
surface courbe. Dans la {\em première}, on se sert de l'équation entre
les coordonnées $x$, $y$, $z$ que l'on supposera ramenée à la forme
$W = 0$, où $W$ sera une fonction des indéterminées $x$, $y$,
$z$. [\dots] Dans la {\em seconde} méthode, on exprime les coordonnées
en forme de fonctions des deux variables $p$ et $q$. [\dots] \`A ces
deux méthodes générales se rattache une {\em troisième} méthode, dans
laquelle l'une des coordonnées, par exemple $z$, est exprimée en
fonction des deux autres, $x$, $y$; cette méthode n'est évidemment
autre chose qu'un cas particulier, soit de la première méthode, soit
de la seconde. \REFERENCE{{\sc Gauss},~\cite{ga1828},~10--13}}

Soient donc $x, y, z$ trois coordonnées qui permettent de repérer tous
les points possibles d'un espace à trois dimensions de type euclidien
ou cartésien tel que conceptualisé par la mécanique classique; $z$
décrit la troisième dimension. 

Dans la représentation implicite {\bf I}, on devine alors aisément
qu'une surface devra être définie comme le lieu de tous les points
$(x_p, y_p, z_p)$ en lesquels s'annule une certaine fonction $W ( x,
y, z)$ des trois variables $(x, y, z)$. En y réfléchissant plus
avant, on se convainc effectivement qu'une unique équation du genre $W
( x, y, z) = 0$ prive les points d'un seul degré de liberté parmi les
trois qu'ils possèdent au départ, et il ne leur en reste alors plus
que deux, et c'est bien là ce qu'on attend d'une surface\footnote{\,
%%%%%%%%%%%%%%%%%%%%%%%-------DEBUT--------%%%%%%%%%%%%%%%%%%%%%%%%%%%
En toute rigueur, comme dans le cas des courbes, il faut effectuer une
hypothèse sur la fonction $W$, et la plus naturelle est de demander
qu'en tout point de l'ensemble $\{ W = 0\}$, 
au moins une des trois dérivées partielles $\frac{ \partial W}{
\partial x}$, ou $\frac{ \partial W}{ \partial y}$, ou $\frac{
\partial W}{ \partial z}$ ne s'annule pas.
}. %%%%%%%%%%%%%%%%%%%%%%%%-----FIN-----%%%%%%%%%%%%%%%%%%%%%%%%%%%%%%%
Par exemple, un ellipsoïde quelconque dans l'espace
(généralisation de la notion d'ellipse dans le plan)
est défini par une équation du type général: 
\[
\frac{(x-x_1)^2}{a^2}
+
\frac{(y-y_1)^2}{b^2}
+
\frac{(z-z_1)^2}{c^2}
-
1
=
0,
\]
où le centre $p_1$ de coordonnées 
$(x_1, y_1, z_1)$ est un point quelconque dans l'espace et où 
$a > 0$, $b > 0$ et $c > 0$ sont des nombres réels positifs
arbitraires. 

\begin{center}
\input ellipsoide.pstex_t
\end{center}

\noindent
Un tel ellipsoïde provient d'une sphère ronde (cas $a = b = c$)
par déformation, élongation ou aplatissement, le long des
axes de coordonnées. 

Dans la représentation graphée {\bf III}, la plus explicite et la plus
concrète des trois, l'une des trois coordonnées, par exemple $z$, est
fonction des deux autres, par exemple fonction de $x$ et de $y$, de
telle sorte que les points de la surface sont exactement tous les
points $(x, y, z)$ de l'espace tels que:
\[
z
=
f(x,y),
\]
pour une certaine fonction $f$ qui dépend de la surface\footnote{\,
%%%%%%%%%%%%%%%%%%%%%%%-------DEBUT--------%%%%%%%%%%%%%%%%%%%%%%%%%%%
---\,\,et qui dépend aussi des coordonnées $(x, y, z)$ choisies 
au départ sur l'espace.
}. %%%%%%%%%%%%%%%%%%%%%%%%-----FIN-----%%%%%%%%%%%%%%%%%%%%%%%%%%%%%%%
Cette représentation graphée peut d'ailleurs être envisagée comme un
cas particulier de la représentation implicite, si l'on introduit 
simplement la fonction:
\[
W(x,y,z)
:=
z
-
f(x,y),
\]
mais en toute rigueur, elle exprime davantage que la représentation
implicite, puisqu'elle explicite la coordonnée $z$ comme fonction des
deux autres coordonnées, lesquelles désignent
sans ambiguïté l'horizontalité
bidimensionnelle relative\footnote{\,
%%%%%%%%%%%%%%%%%%%%%%%-------DEBUT--------%%%%%%%%%%%%%%%%%%%%%%%%%%%
Presque immédiatement, on devine les généralisations de ces concepts
analytiques dans un espace réel à un nombre arbitraire $n\geqslant 2$
(au lieu de $n = 2$ ou $n = 3$) de dimensions, muni de
coordonnées $(x_1, x_2, \dots, x_{ n-1}, x_n)$: la représentation
d'une {\em hyper}surface (telle est la terminologie consacrée)
s'effectue ou bien par une équation implicite $W( x_1, \dots, x_n) =
0$, ou bien par une équation graphée $x_n = f ( x_1, \dots, x_{
n-1})$, ou bien par une une paramétrisation:
\[
(u_1,\dots,u_{ n-1})
\longmapsto
\big(
x_1(u_1,\dots,u_{n-1}),
\,\dots\dots,\,
x_n(u_1,\dots,u_{n-1})
\big).
\]
Ce qui est surprenant dans ces formulations par le langage
mathématique, c'est que toute la difficulté réelle et légitime que
l'on a à se représenter mentalement des objets géométrique possédant
un grand nombre de dimensions {\em s'efface} derrière la simplicité
troublante du <<\,$n$ quelconque\,>> et des morphèmes syntaxiques. 
Il faut en déduire que la pensée mathématique peut constamment être
porteuse de ce qu'elle ne sait pas mentaliser complètement, 
surtout lorsqu'il s'agit des représentations géométriques.
}. %%%%%%%%%%%%%%%%%%%%%%%%-----FIN-----%%%%%%%%%%%%%%%%%%%%%%%%%%%%%%%

Enfin, dans la représentation paramétrée {\bf II}, la bidimensionnalité
de la surface est rendue visible
d'emblée, puisqu'on se donne une application:
\[
(u,v)
\longmapsto
\big(x(u,v),\,y(u,v),\,z(u,v)\big)
\] 
d'un espace de deux variables $(u, v)$ à valeurs dans la surface qui
en représente ses trois coordonnées $x$, $y$ et $z$ comme certaines
fonctions spécifiques $x ( u, v)$, $y( u, v)$ et $z ( u, v)$ de ces
deux variables. Autrement dit, un petit morceau de surface est
constitu\'e de l'ensemble des points de coordonn\'ees $(x(u,\, v),\,
y(u,\, v),\, z(u,\, v))$ dans l'espace, lorsque le point de
coordonn\'ees $(u,\, v)$ varie dans un plan auxiliaire
bidimensionnel. 

\medskip

\begin{center}
\input u-v-x-y-z.pstex_t
\end{center}

\noindent
Intuitivement, on peut s'imaginer que l'on d\'eforme en pens\'ee le
plan <<\,droit\,>> des $(u,\, v)$ dans l'espace gr\^ace \`a la
libert\'e que fournit la troisi\`eme dimension, avec une libert\'e
illimit\'ee, comme s'il s'agissait d'une membrane id\'eale infiniment
plus \'elastique et plus mobile que toute membrane existant dans le
monde physique. Dans ces conditions, deux réseaux orthogonaux
de courbes parallèles deviennent deux réseaux de courbes
(pas forcément orthogonales) et approximativement parallèles
sur la surface, si tant est qu'on puisse parler de courbes
parallèles.

Bien entendu, il est n\'ecessaire de pr\'eciser la r\'egularit\'e de
l'application qui associe aux deux coordonn\'ees $(u,\, v)$, le point
ayant les trois coordonn\'ees $(x(u,\, v),\, y(u,\, v),\, z(u,\, v))$
dans l'espace tridimensionnel: on supposera que cette application est
{\sl continue}, c'est-\`a-dire que des petites variations des
coordonn\'ees $(u,\, v)$ induisent des petites variations des
coordonn\'ees $(x(u,\, v),\, y(u,\, v),\, z(u,\, v))$\,\,---\,\,la
membrane ne se d\'echire pas\,\,---\,\,et de plus que cette
application est {\sl diff\'erentiable}, c'est-\`a-dire que les
d\'eriv\'ees partielles de $x$, de $y$ et de $z$ par rapport \`a $u$
et $v$ existent, au moins jusqu'\`a l'ordre deux\,\,---\,\,la
membrane ne pr\'esente pas de lignes d'angle et elle est bien
<<\,lisse\,>>.

\Section{4.~Courbure des courbes planaires}
\label{courbure-courbes-planaires}

\HEAD{4.~Courbure des courbes planaires}{
Jo\"el Merker, D\'partement de Math\'matiques d'Orsay}

\medskip\noindent{\bf Longueur d'une courbe.} 
\label{longueur-courbe}
D'après le th\'eor\`eme de Pythagore,
un segment rectiligne dans le plan dont les projections horizontale
et verticale sont \'egales \`a $a$ et \`a $b$, respectivement,
poss\`ede une longueur \'egale \`a $\sqrt{ a^2 + b^2}$.

\medskip

\begin{center}
\input theoreme-pythagore.pstex_t
\end{center}

\noindent
Pour calculer la longueur d'une courbe graph\'ee $\{ y= f(x) \}$ entre
deux bornes $x_1$ et $x_2 > x_1$, on la décompose en une infinit\'e de
segments rectilignes infinit\'esimaux. Puisque le segment
infinit\'esimal situ\'e au point $x$ poss\`ede une projection
horizontale \'egale \`a $dx$ et une projection verticale \'egale \`a
$dy = f'(x)\, dx$, en appliquant Pythagore dans l'infinitésimal, on en
déduit que la longueur du segment rectiligne infinit\'esimal est
\'egale \`a:
\[
\aligned
\sqrt{dx^2+(f'(x)\, dx)^2 
} 
=
& \
\sqrt{dx^2[1+ f'(x)^2]} 
\\
=
& \
\vert dx\vert\,
\sqrt{1+f'(x)^2}.
\endaligned
\]
Si l'on identifie alors la courbe \`a une cha\^{\i}ne infinie de tels
segments infinit\'esimaux plac\'es les uns \`a la suite des autres, sa
longueur est bien \'evidemment \'egale \`a la somme des longueurs de
ces segments infinit\'esimaux. Les fondements du calcul intégral 
permettent alors d'obtenir la formule suivante. 

\smallskip\noindent{\bf Théorème.}
{\em La longueur d'une courbe graph\'ee continûment
différentiable $y = f(x)$ entre
deux bornes $x_1$ et $x_2 > x_1$ est \'egale \`a:}
\[
\int_{x_1}^{x_2}\, 
\sqrt{1+f'(x)^2}\, 
dx.
\]

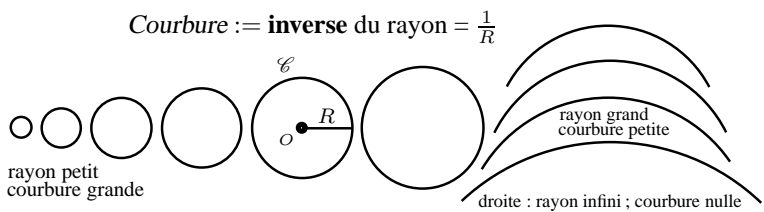
\medskip\noindent{\bf Courbure des cercles et des droites.}
\label{courbure-cercles-droites}
Dans le plan euclidien, apr\`es les droites, ce sont les cercles qui
sont les courbes les plus simples. Tout cercle $\mathcal{ C}$ est
caract\'eris\'e par son centre $O$ et par son rayon $R$. Les petits
cercles sont tr\`es <<\,courb\'es\,>>, {\em i.e.} leur
<<\,courbure\,>> est tr\`es accentu\'ee; les grands cercles sont peu
<<\,courb\'es\,>>, {\em i.e.} ils sont presque <<\,droits\,>>, comme
la surface de la mer \`a l'\'echelle des cartes g\'eographique. Sur
l'autoroute, la <<\,courbure\,>> des virages est beaucoup moins
sensible\,\,---\,\,{\em via} la force centrifuge\,\,---\,\,que sur une
route de montagne en \'epingle \`a cheveux.

Quelle quantit\'e peut-on alors proposer pour parler de la {\sl
courbure math\'ematique} d'un cercle? Pour r\'epondre \`a cette
question, il faudra tenir compte de trois faits \'el\'ementaires:

\begin{itemize}

\smallskip\item[{\bf (1)}]
plus le
cercle est petit, plus il est <<\,courb\'e\,>>, intuitivement parlant;

\smallskip\item[{\bf (2)}]
un cercle \'etant invariant par rotation autour de son centre ({\em
cf.} l'invention de la roue), sa courbure doit
\^etre la m\^eme en tout point;

\smallskip\item[{\bf (3)}]
les droites n'ont pas de <<\,courbure\,>> et
peuvent \^etre identifi\'ees \`a des cercles de rayon infini.

\end{itemize}\medskip

\begin{center}
\input courbure-cercle.pstex_t
\end{center}
\medskip

\noindent
L'id\'ee la plus simple\footnote{\,
%%%%%%%%%%%%%%%%%%%%%%%-------DEBUT--------%%%%%%%%%%%%%%%%%%%%%%%%%%%
Toute autre fonction $\kappa( R)$, d\'efinie sur l'ensemble des
rayons, d\'ecroissante et satisfaisant $\lim_{ R \to 0} \, \kappa (R)
= + \infty$ ainsi que $\lim_{ R\to \infty} \, \kappa (R) = 0$,
pourrait aussi convenir.
} %%%%%%%%%%%%%%%%%%%%%%%%-----FIN-----%%%%%%%%%%%%%%%%%%%%%%%%%%%%%%%
consiste alors \`a d\'efinir {\sl
courbure math\'ematique} d'un cercle (ou d'une droite) comme l'{\em
inverse} de son rayon: 
\[
\frac{1}{R} 
:=
\text{{\sl courbure} d'un cercle de rayon } 
R.
\]
En bref, la courbure est inversement proportionnelle au rayon. C'est
cette d\'efinition simple que les math\'ematiques ont retenue. Aux
droites dont le rayon est infini correspond alors une courbure nulle:
$\frac{ 1}{ \infty} = 0$, ce qui est en accord avec l'intuition
élémentaire.

\medskip\noindent{\bf Cha\^{\i}ne de segments rectilignes et
courbure.} 
\label{chaine-segments-courbure}
Cependant, en identifiant une courbe \`a une collection infinie de
segments de droite pour en calculer la longueur à l'aide du calcul
intégral, ne perd-on pas irr\'em\'ediablement toute représentation de
sa courbure? La r\'eponse \`a ce paradoxe se devine ais\'ement sur une
figure macroscopique.

\medskip

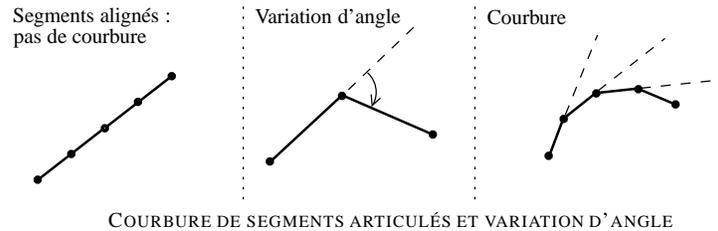
\begin{center}
\input angle-courbure.pstex_t
\end{center}

\noindent
Lorsque les segments sont align\'es, la courbure est \'evidemment
nulle; au contraire, lorsque les segments ne sont pas align\'es, la
variation d'angle entre deux segments adjacents devrait permettre de
mesurer une <<\,courbure approximative\,>> de la cha\^{\i}ne. En effet,
si la cha\^{\i}ne se courbe dans une certaine direction, la variation
des angles se fait toujours dans le m\^eme sens. Pour \'etudier cette
variation des angles, il faut introduire la notion de 
{\sl vecteur unitaire
normal} \`a une courbe.

Soit donc 
$C$ une courbe quelconque et soit $p$ l'un de ses points. Dans la
suite, on appellera {\sl unitaire} tout vecteur de longueur \'egale
\`a $1$. Sur la tangente en $p$, on trace un vecteur unitaire
d'origine $p$. Deux possibilit\'es se pr\'esentent (vers la droite ou
vers la gauche; vers le haut ou vers le bas), mais on sous-entendra
que la courbe est orient\'ee, {\em i.e.} qu'on a convenu \`a l'avance
d'un sens de parcours sur la courbe\,\,---\,\,et alors il existe un
seul vecteur tangent unitaire qui se dirige dans le m\^eme sens que la
courbe. Soit $T(p)$ ce vecteur tangent unitaire orienté, la lettre
<<\,$T$\,>> \'etant l'initiale de l'adjectif <<\,tangent\,>>. Par une
rotation de $90$ degr\'es, chaque vecteur $T(p)$ se transforme en un
vecteur unitaire $N(p)$ orthogonal \`a la droite tangente en $p$. 
On appelle $N(p)$ le {\sl vecteur normal unitaire \`a la
courbe en $p$}, la lettre <<\,$N$\,>> \'etant l'initiale de l'adjectif
<<\,normal\,>>.

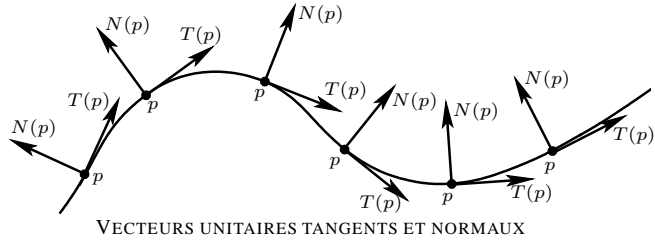
\medskip

\begin{center}
\input normales-tangentes.pstex_t
\end{center}

\noindent
Sur la figure, par souci de simplicit\'e, on donne le
m\^eme nom <<\,$p$\,>> aux six points choisis sur la courbe
pour tracer les vecteurs unitaires tangents et normaux. On obtient
ainsi deux champs de vecteurs, trac\'es sur la courbe $C$.

\medskip\noindent{\bf Courbure des cercles via l'application de Gauss.}
\label{cercles-application-Gauss} Maintenant, 
sur une telle courbe quelconque $C$, puisque le vecteur normal
unitaire $N(p)$ se d\'eduit par une rotation d'angle fixe du vecteur
tangent unitaire $T(p)$, la variation d'angle des vecteurs $T(p)$ est
\'evidemment égale \`a la variation d'angle des vecteurs $N(p)$,
lorsque le point $p$ se déplace sur la courbe. On se concentrera donc
sur le vecteur normal unitaire $N(p)$, puisqu'il est plus facile de
repr\'esenter sa variation sur les diagrammes. Tout d'abord, est-il
possible de comprendre et de voir sous un autre angle la courbure des
cercles?

\Scholie
Question universelle d'ordre proprement métaphysique, et question
indéfiniment reproductible: étant donné un concept mathématique,
défini ou appréhendé d'une certaine manière, n'est-il pas possible
qu'il existe d'{\em autres} manières de l'appréhender ou de le
définir, sachant que de telles {\em autres} manières ont, à cause d'un
passé mathématique sédimenté, la capacité {\em métaphysique} de se
poser dans l'attente de fécondations intuitives rétroactives. Principe
de raison inductive, aussi, basé sur la constatation
<<\,aristotélicienne\,>> que la multiplicité ramifiée du dire de
l'être\,\,---\,\,l'être mathématique en l'occurrence\,\,---\,\,ne
s'interrompt jamais. Enfin, au-delà d'une démultiplication de l'Un
conceptuel par la pensée, l'examen du caractère problématique de
l'Idée pousse inlassablement aux genèses en généralisation.
\qed

\smallskip

{\small 

La {\sl courbure de Menger} pour les compacts du plan complexe $\C$,
par exemple, étend la conceptualisation de l'Idée de courbure planaire
bien au-delà des objets unidimensionnels, à des compacts arbitraires
du plan complexe $\C$ relativement à une mesure de Radon, et ce, comme
suit. Tout d'abord, si $z_1, z_2, z_3$ sont trois points de $\C$
distincts deux à deux, et si $S ( z_1, z_2, z_3)$ désigne l'aire
(absolue) du triangle $z_1z_2z_3$, on vérifie que la quantité:
\[
c(z_1,z_2,z_3)
:=
\frac{S(z_1,z_2,z_3)}{
\vert z_1-z_2\vert\,\vert z_1-z_3\vert\,\vert z_2-z_3\vert}
=
\frac{1}{R(z_1,z_2,z_3)}
\]
est égale à l'inverse du rayon $R(z_1, z_2, z_3)$ du cercle
circonscrit à ce triangle, quantité qui est nulle lorsque, et
seulement lorsque ces trois points sont alignés. Ensuite, si $K
\Subset \C$ est un sous-ensemble compact arbitraire, à savoir un
sous-ensemble qui est fermé pour la topologie induite de $\C$ et qui
est aussi borné,
\qed

}

\smallskip
Dans le plan euclidien, soient alors par
exemple deux cercles, un <<\,petit\,>> cercle
de rayon $r$ et de courbure $\frac{ 1}{r}$, ainsi qu'un <<\,grand\,>>
cercle de rayon $R$ et de courbure $\frac{ 1}{ R}$.

\medskip

\begin{center}
\input cercle-auxiliaire.pstex_t
\end{center}

\`A c\^ot\'e d'eux, on trace aussi un cercle de rayon $1$ que l'on 
appellera {\sl cercle auxiliaire}; la longueur de sa
circonf\'erence est \'egale \`a $2\pi$. Dans un tel cercle
math\'ematique id\'eal, les angles sont mesur\'es non pas en degr\'es,
mais en {\sl radians}\,\,---\,\,la plus simple des unit\'es de mesure:
par d\'efinition, la {\sl mesure en radians} d'un angle est \'egale
\`a la longueur de l'arc de cercle qu'il d\'ecoupe sur le cercle
auxiliaire. Par exemple, le tour complet vaut $2\pi$
radians\,\,---\,\,au lieu de 360 degr\'es\,\,---\,\, et le quart de
tour vaut $\frac{ \pi}{ 2}$ radians\,\,---\,\,au lieu de 90
degr\'es. Pour simplifier, les angles considérés ne seront pas
orient\'es.

Soient donc deux angles infinit\'esimaux $d \theta$ et $d \varphi$
trac\'es sur le cercle auxiliaire. En admettant sans s'y attarder que
la métaphysique du calcul infinitésimal est bien fondée, on dessinera
tous les infinit\'esimaux avec une certaine extension physique;
toutefois, sous une autre perspective, l'intuition id\'eale doit les
<<\,embrasser\,>> comme de vrais infiniment petits. 

\Scholie
Troublante puissance de l'intuition mathématique experte!
Capable à la fois d'emboîter le pas au sensible et de maintenir
vivante l'exigence de traductibilité théorique rigoureuse, elle
parcourt les espaces de la pensée à tous leurs étages.
\qed

\smallskip
Ainsi les
deux angles $d \theta$ et $d \varphi$ co\"{\i}ncident chacun avec une
longueur d'arc de cercle infinit\'esimal. On translate alors l'angle
$d\theta$ jusqu'au petit cercle et on translate l'angle $d\varphi$
jusqu'au grand cercle. La longueur des deux arcs de cercle ainsi
d\'ecoup\'es est \'egale\,\,---\,\,respectivement\,\,---\,\,\`a
$rd\theta$ et \`a $Rd\varphi$.

R\'eciproquement, si l'on se donne la longueur d'un arc
infinit\'esimal situ\'e sur l'un des deux cercles (par exemple le
grand), c'est-\`a-dire, si l'on se donne la longueur d'un arc situ\'e
entre deux points $p$ et $q$ infiniment proches, on peut construire un
arc de cercle infinit\'esimal {\em associ\'e}, situ\'e sur le cercle
auxiliaire: il suffit de translater le vecteur normal {\em unitaire}\,
$N(p)$ jusqu'au centre du cercle auxiliaire, de translater aussi le
vecteur $N(q)$ et de mesurer l'angle qui est d\'elimit\'e sur la
sph\`ere auxiliaire par les extr\'emit\'es de ces deux vecteurs. Ce
proc\'ed\'e de translation était en fait déjà illustr\'e sur la
figure précédente.

Or on constate une identit\'e alg\'ebrique triviale: la courbure
de chaque cercle s'identifie au quotient de ces deux longueurs
infinit\'esimales\footnote{\,
%%%%%%%%%%%%%%%%%%%%%%%-------DEBUT--------%%%%%%%%%%%%%%%%%%%%%%%%%%%
Gr\^ace \`a l'homog\'en\'eit\'e de la courbure d'un cercle, les
m\^emes propri\'et\'es sont satisfaites par des arcs de longueur
finie, non infinit\'esimale. Au contraire, pour une courbe arbitraire,
la considération
des variations d'angle infinit\'esimales\,\,---\,\,qui
sera généralisée {\em infra}\,\,---\,\,est
incontournable.
}: %%%%%%%%%%%%%%%%%%%%%%%%-----FIN-----%%%%%%%%%%%%%%%%%%%%%%%%%%%%%%%
\[
\frac{ 1}{ r} = \frac{ d\theta}{r d\theta}
\ \ \ \ \ \ \
{\rm et}
\ \ \ \ \ \ \
\frac{ 1}{ R} = \frac{ d\varphi}{ R d \varphi}.
\]
En résumé, 
la courbure de tout cercle s'identifie au quotient:
\[
\frac{\text{\rm longueur correspondante sur le cercle auxiliaire}}{
\text{\rm longueur d'un arc infinit\'esimal sur la circonf\'erence}}.
\]

\medskip\noindent{\bf Définition de la courbure des courbes 
quelconques dans le plan.}
\label{courbure-courbes-quelconques}
Soit maintenant $C$ une courbe arbitraire trac\'ee dans le plan
euclidien et soit à nouveau le cercle auxiliaire de rayon unité, placé
à droite de la figure explicative ci-dessous. 
Soit $p$ un point de la courbe,
soit $p+dp$ un point infinit\'esimalement proche de $p$ et soit
$d\ell$ la distance infinit\'esimale qui s\'epare les deux
points. On translate alors\footnote{\,
%%%%%%%%%%%%%%%%%%%%%%%-------DEBUT--------%%%%%%%%%%%%%%%%%%%%%%%%%%%
On appellera {\sl application de Gauss} cette op\'eration de
translation des vecteurs unitaires normaux.
} %%%%%%%%%%%%%%%%%%%%%%%%-----FIN-----%%%%%%%%%%%%%%%%%%%%%%%%%%%%%%%
les deux vecteurs unitaires normaux $N(p)$ et $N(p+ dp)$ jusqu'au
cercle auxiliaire et et on calcule la distance infinit\'esimale $d\theta$
qui est d\'elimit\'ee par ces deux vecteurs sur le cercle auxiliaire.

\medskip

\begin{center}
\input cercle-application-Gauss.pstex_t
\end{center}

\smallskip\noindent{\bf Définition.}
La courbure d'une courbe $C$ quelconque en l'un de ses
points $p$ est le quotient infinit\'esimal:
\[
\frac{ d\theta}{d\ell}
=
\frac{\text{\rm longueur correspondante sur le cercle auxiliaire}}{
\text{\rm longueur d'un arc infinit\'esimal sur la courbe}},
\]
lequel peut s'interpréter rigoureusement comme une limite
qui est bien définie lorsque la courbe est au moins deux fois
continûment différentiable (on dit alors
que la courbe est {\sl de classe $\mathcal{ C}^2$}).

\Scholie
Compréhension et cohérence de la définition par variation d'angle dans
un cas simple (les cercles); généralisation-idéalisation 
de la même définition au cas de courbes
de classe $\mathcal{ C}^2$ absolument 
quelconques; saisie intuitive de la mesure de la courbure par
<<\,déroulement\,>> sur un cercle unité auxiliaire; 
vérification méditative de l'adéquation de cette 
désarticulation mesurante;
inspiration implicite qu'un certain 
<<\,blotissement de cercle\,>> épouse
et dirige la courbure de la courbe en chaque point.

\qed

\smallskip\noindent{\bf Cercle osculateur \`a une courbe en un point.}
\label{cercle-osculateur-courbe}
Soit à nouveau $C$ une courbe arbitraire du plan et soit $p$ l'un de
ses points. Si $q$ et $r$ sont deux points distincts situ\'es sur $C$
de part et d'autre de $p$, il existe un unique cercle passant par les
trois points $p$, $q$ et $r$, 
cercle se réduit à une droite lorsque ces
points sont alignés. Soit donc $\mathcal{ C}(p,\, q,\, r)$ ce cercle
(ou droite), qui d\'epend de la position de $q$ et de $r$, comme sur
la figure ci-dessous.

\medskip

\begin{center}
\input cercle-secant.pstex_t
\end{center}

\noindent
On d\'emontre math\'ematiquement\footnote{\,
%%%%%%%%%%%%%%%%%%%%%%%-------DEBUT--------%%%%%%%%%%%%%%%%%%%%%%%%%%%
---\,\,en supposant que la courbe est
(au moins) deux fois contin\^ument diff\'erentiable, 
{\em i.e.} de classe $\mathcal{ C}^2$\,\,---\,\,
} %%%%%%%%%%%%%%%%%%%%%%%%-----FIN-----%%%%%%%%%%%%%%%%%%%%%%%%%%%%%%%
que lorsque les deux points $q$ et $r$ se rapprochent ind\'efiniment
du point $p$, le cercle $\mathcal{ C}(p,\, q,\, r)$ se rapproche
ind\'efiniment d'un unique cercle passant par le point $p$. Ce cercle
est appel\'e {\sl cercle osculateur\footnote{\,
%%%%%%%%%%%%%%%%%%%%%%%-------DEBUT--------%%%%%%%%%%%%%%%%%%%%%%%%%%%
Du latin {\em osculari} <<\,embrasser\,>>, pour exprimer que le cercle
a un ordre de contact \'elev\'e avec la courbe, plus \'elev\'e que sa
droite tangente.
} %%%%%%%%%%%%%%%%%%%%%%%%-----FIN-----%%%%%%%%%%%%%%%%%%%%%%%%%%%%%%%
\`a la courbe $C$ au point $p$}.

Soit $\mathcal{ C}_p$ ce cercle, soit $O_p$ son centre et 
soit $R_p$ son
rayon. Intuitivement parlant, ce cercle <<\,osculateur\,>> est le cercle
qui approxime le mieux\footnote{\,
%%%%%%%%%%%%%%%%%%%%%%%-------DEBUT--------%%%%%%%%%%%%%%%%%%%%%%%%%%%
La propri\'et\'e de meilleure approximation du cercle osculateur
raffine et am\'eliore la propri\'et\'e d'approximation de la droite
tangente. En effet, on d\'emontre math\'ematiquement que lorsqu'un
point $r$ varie sur le cercle osculateur $\mathcal{ C}_p$ et se
rapproche de $p$, la distance de $r$ \`a son projet\'e orthogonal
$q_r$ sur la courbe $C$ devient infiniment petite par rapport au {\em
carr\'e}\, de la distance de $r$ \`a $p$, c'est-à-dire que
le quotient:
\[
\frac{\text{\rm distance}\,(r,\,q_r)}{
\left[\text{\rm distance}\,(r,\,p)\right]^2}
\]
tend vers z\'ero quand le point $r$ tend vers le point $p$. 
} %%%%%%%%%%%%%%%%%%%%%%%%-----FIN-----%%%%%%%%%%%%%%%%%%%%%%%%%%%%%%%
la courbe dans un voisinage infime du point $p$. Aussi pourrait-on
d\'efinir aussi la {\sl courbure math\'ematique de la courbe $C$ au
point $p$} comme
\'etant l'inverse du rayon de son cercle osculateur en $p$:
\[
\frac{1}{R_p} 
=
\text{{\sl courbure} de la courbe $C$ au point $p$},
\]
et on admettra, sans le redémontrer, le résultat
élémentaire suivant.

\smallskip\noindent{\bf Théorème.}
{\em Les deux définitions de la courbure, par variation infinitésimale
d'angle du vecteur normal unitaire ou par 
inversion du rayon du cercle osculateur, coïncident.}\medskip

\Scholie
Convergence et coïncidence de deux multiples définitionnels qui sont
<<\,naturels\,>> chacun de leur côté: il faut y voir la manifestation,
à un niveau supérieur, d'une cohérence {\em métaphysique} des
mathématiques. La tangente à une courbe offre un principe de
comparaison avec (ou d'approximation par) une droite; l'osculation par
des cercles raffine le principe par simple enrichissement d'une
panoplie-modèle.

Bien entendu, l'ontologie génétique intuitive commande immédiatement
d'étendre ces principes de comparaison (ou d'approximation) à des
cadres plus généraux, voire les plus généraux possibles. Formules de
Taylor-Young à divers ordres, ou variétés de jets en dimension
supérieure ouvrent la voie vers des univers calculatoirement
extrêmement riches, qui naissent en partie ici, dans le <<\,berceau
simplissime\,>> des courbes. En théorie de Cartan, les courbures que
possèdent les déformations d'espaces homogènes quelconques recèlent
une {\em ontologie indéfinie, imprévisible et inépuisable}, non
seulement parce que les algèbres (groupes) de Lie non semi-simples se
ramifient à l'infini de manière non totalement classifiable, mais
encore parce que la construction des structures de parallélisme
absolu à la Cartan-Tanaka, virtuellement possibles dans tous les cas
de figure, exige en fait dans chaque situation spécifique de très
longs calculs d'élimination qui sont très emboîtés les uns dans les
autres et presque jamais effectués(ables) en totalité! Arbres de
Hilbert-Gödel (incomplétude) et arbres de calculs (indéfinitude),
telle est l'incroyable {\em situation
universelle d'ouverture} à laquelle sont exposées
les mathématiques.
\qed

\smallskip
\Scholie
Toutefois, les considérations géométriques précédentes n'apportent
aucune information quantitative. On soupçonne en effet facilement
qu'il devrait exister une {\em formule} qui exprime {\em
quantitativement} la courbure d'une courbe en fonction des éléments
qui la définissent. Autrement dit, si l'être d'une courbe est
saisissable non pas seulement en tant qu'être géométrique, mais en
tant qu'être fonctionnel effectif à cause du mystérieux et indubitable

\smallskip
\centerline{\fbox{\sl principe de transsubstantiation du géométral en le
fonctionnel-calculatoire pur}\,,}

\smallskip\noindent
alors puisque la courbure géométrique
dépend directement des caractères de la courbe, et donc puisque l'on
ne sait pas encore de quelle manière la dépendance a lieu,
c'est-à-dire puisque la courbure géométrique peut {\em a priori}
dépendre sans restriction de {\em tous} les caractères potentiellement
étudiables de la courbe, il doit être possible de donner la parole à
{\em certains} caractères adéquats de cette courbe afin d'en extraire,
{\em par le raisonnement et par le calcul pur}, une formule 
qui exprime la courbure sous une forme explicite satisfaisante.
Pétition de principe métaphysique\,\,---\,\,objecterait Heidegger\,\,---,
et domination du principe de raison\,\,---\,\,ajouterait-il\,\,---, 
mais toutefois, accepter une telle 
<<\,règle du jeu supérieure\,>> comme le font les mathématiciens
ne limite en rien la composante proprement {\em méditative et
spéculative} de l'exploration technique.
\qed

\smallskip\Scholie
Encore une fois, chaque élément de spéculation
métaphysico-mathématique compte. Ici, il faut affirmer fermement que
le calcul n'est jamais <<\,un calcul\,>> placé en position secondaire
ou latérale par rapport à des raisonnements conceptuels abstraits ou
géométriques purs. Seul le calcul effectif pousse la recherche de
vérité dans les retranchements de son inexpression initiale et de son
inexplicitation embryonnaire. Poursuivre des raisonnements en
direction d'une complétude exige qu'on accepte, tant qu'une
question n'est pas résolue, qu'elle n'est pas encore résolue.
\qed

\smallskip
Heureusement, il se trouve qu'il est relativement facile de {\em
calculer} l'expression explicite de la courbure d'une courbe graph\'ee
$y= f(x)$ en chacun de ses points. Ce calcul constitue un
pr\'eliminaire important aux calculs plus \'elabor\'es que Gauss a
conduits pour les surfaces.

\Section{5.~Expression analytique de la courbure des courbes planaires}
\label{expression-courbure-planaire}

\HEAD{5.~Expression analytique de la courbure des courbes planaires}{
Jo\"el Merker, D\'partement de Math\'matiques d'Orsay}

\medskip\noindent{\bf Vecteurs tangents et vecteurs normaux.}
\label{tangents-normaux}
En un point arbitraire de coordonn\'ees $(x, f(x))$, le vecteur
infinit\'esimal tangent a pour composantes $(dx, f'(x)\, dx)$.
Apr\`es division par $dx$, le vecteur macroscopique (non
infinit\'esimal) de coordonn\'ees $(1,\, f'(x))$ est tangent \`a la
courbe. La norme de ce vecteur \'etant \'egale \`a $\sqrt{ 1+
f'(x)^2}$, on en d\'eduit que le vecteur unitaire tangent $T(x)$
\`a la courbe en ce point a pour coordonn\'ees:
\[
T(x) = \left(
\frac{ 1}{
\sqrt{ 1 + f'(x)^2}}, \ \
\frac{ f'(x)}{
\sqrt{ 1 + f'(x)^2}}
\right).
\]

\medskip

\begin{center}
\input tp-np-projection.pstex_t
\end{center}

\noindent
Apr\`es une rotation d'angle $\frac{ \pi}{2}$, le vecteur normal
unitaire $N(x)$ a donc pour coordonn\'ees:
\[
N(x) 
= 
\left(
-
\frac{f'(x)}{
\sqrt{1+f'(x)^2}}, \ \
\frac{1}{
\sqrt{1+f'(x)^2}}
\right).
\]
La figure aide \`a deviner\,\,---\,\,ou 
\`a se souvenir\,\,---\,\,pourquoi il y a un
signe <<\,$-$\,>> devant la premi\`ere composante de $N(x)$. D'après
sa définition en termes de variation de l'angle du vecteur normal, la
courbure $c(x)$ de la courbe au point de coordonn\'ees $\big(x,\, 
f(x) \big)$
est donn\'ee par le quotient:
\[
c(x) 
= 
\frac{\text{\rm norme de}\,
dN(x)}{\text{\rm longueur de}\ 
x\ \text{\rm \`a}\ x+dx}.
\]
Or la longueur infinit\'esimale de la courbe entre le point rep\'er\'e
par $x$ et le point rep\'er\'e par $x+ dx$ est
\'egale \`a $\sqrt{ 1+ f'(x)^2} \, dx$: c'est 
le d\'enominateur. Pour calculer le numérateur, il faut prendre la
norme des deux membres de la différentielle $dN(x) = N'(x) \cdot dx$,
ce qui donne:
\[
\text{\rm norme de}\ 
dN(x) 
=
\vert dx\vert\cdot
\text{\rm norme de}\,N'(x),
\]
la notation classique $\vert a \vert$ d\'esignant la {\sl valeur
absolue} d'un nombre r\'eel $a$. Maintenant, pour calculer la
d\'eriv\'ee $N'(x)$, puis sa norme, il faut appliquer des formules
classiques:

\begin{itemize}

\smallskip\item[{\bf (1)}]
la d\'eriv\'ee de la fonction $\frac{ 1}{\sqrt{ x}}$ est \'egale \`a
$-\frac{ 1}{2\, x \, \sqrt{ x}}$;

\smallskip\item[{\bf (2)}]
la d\'eriv\'ee d'une fonction compos\'ee $\frac{ 1}{\sqrt{ g(x)}}$ 
est \'egale \`a $-\frac{ g'(x)}{2 \, g(x) \, \sqrt{ g(x)}}$;

\smallskip\item[{\bf (3)}]
la d\'eriv\'ee d'une fonction d\'eriv\'ee $f'(x)$ est
appel\'ee {\sl d\'eriv\'ee seconde} et not\'ee $f''(x)$;

\smallskip\item[{\bf (4)}]
la d\'eriv\'ee de la fonction 
$g(x):= 1+ f'(x)^2$ est \'egale \`a $2\, f'(x) \, f''(x)$;

\smallskip\item[{\bf (5)}]
la d\'eriv\'ee d'un produit de fonctions 
$f'(x) \, h(x)$ est \'egale \`a 
$f''(x) \, h(x)+ f'(x) \, h'(x)$.

\end{itemize}\smallskip

\noindent
En appliquant toutes ces formules,
on d\'eduit l'expression de la
d\'eriv\'ee des deux composantes du vecteur normal unitaire; pour
all\'eger l'\'ecriture, l'argument $(x)$ des
fonctions $f'$ et $f''$ sera sous-entendu:
\[
\small
\aligned
N'(x)
& \
=
\left(
\frac{ - f''}{\sqrt{ 1+ f'{}^2}}+ 
\frac{f'{}^2 \, f''}{(1+ f'{}^2) \, \sqrt{ 1+ f'{}^2}}, \ \ 
\frac{- f'{} \, f''}{(1+ f'{}^2) \, \sqrt{ 1+ f'{}^2}}\right) \\
& \
=
\left(
\frac{ - f''}{(1+ f'{}^2) \, \sqrt{ 1+ f'{}^2}}, \ \ 
\frac{- f'{} \, f''}{(1+ f'{}^2) \, \sqrt{ 1+ f'{}^2}}\right), 
\endaligned
\]
d'o\`u il découle que la norme de $N'(x)$, not\'ee $\norm
N'(x) \norm$, est \'egale \`a:
\[
\aligned
\norm N'(x)\norm 
& \
= 
\sqrt{
\frac{f''{}^2 }{(1+ f'{}^2)^3 }+ 
\frac{
f'{}^2 \, f''{}^2 }{(1+ f'{}^2)^3 }
} \\
& \
=
\frac{ \vert f'' \vert}{ 1+ f'{}^2}.
\endaligned
\]
Pour terminer, il suffit de remplacer les valeurs du
num\'erateur et du d\'enominateur qui viennent
d'être calculées:
\[
\aligned
c(x) 
& \
= \frac{
\norm N'\norm\cdot dx}{
\sqrt{1+f'{}^2}\cdot dx} 
& \
= 
\frac{\vert f'' \vert}{(1+f'{}^2)
\sqrt{ 1 + f'{}^2}}.
\endaligned
\]
Jusqu'\`a pr\'esent, pour simplifier la pr\'esentation des calculs, il
n'a pas \'et\'e question de l'{\sl orientation} de la courbe $C$, {\em
i.e.} de son sens de parcours. Il existe toujours deux orientations
oppos\'ees l'une \`a l'autre. Si l'on choisit une orientation, il est
possible de d\'efinir la {\sl courbure orient\'ee}, concept qui permet
de d\'eterminer de quel c\^ot\'e de la courbe se trouve son cercle
osculateur. Cette courbure orient\'ee est donc positive ou n\'egative,
sa valeur absolue restant bien entendu \'egale \`a l'inverse du rayon
du cercle osculateur. Pour l'obtenir, il suffit d'\'eliminer le signe
de valeur absolue dans le numérateur $\vert f'' \vert$ et de choisir
l'un des signes <<\,$+$\,>> ou <<\,$-$\,>> devant l'expression
compl\`ete, par exemple le signe <<\,$+$\,>>.
En définitive, on a le:

\smallskip\noindent{\bf Théorème.}
{\em La courbure {\rm orient\'ee} 
$c(x)$ d'une courbe graph\'ee $y = f(x)$ en un point
rep\'er\'e par son abscisse $x$ s'exprime quantitativement par
la formule:}
\[
c(x) 
= 
\frac{f''(x)}{ 
\big(1+f'(x)^2\big) 
\sqrt{1+f'(x)^2}}.
\]

\Scholie
Au premier ordre, le raisonnement
infinit\'esimal a permis d'identifier une courbe \`a une
cha\^{\i}ne infinie de segments rectilignes infinit\'esimaux, comme si
la courbure disparaissait totalement \`a une telle \'echelle.
L'orientation de ces segments rectilignes est repr\'esent\'ee par le
vecteur unitaire tangent \`a la courbe $T(p)$, ou encore\,\,---\,\,\`a
une rotation d'angle $\frac{ \pi}{2}$ pr\`es\,\,---, par le vecteur
unitaire normal $N(p)$, tous deux obtenus par diff\'erentiation
infinit\'esimale. C'est donc en \'etudiant les {\em diff\'erences
infinit\'esimales entre ces vecteurs obtenus par diff\'erentiation}\,
que l'on peut appr\'ehender la courbure. Autrement dit, la courbure
se per\c coit apr\`es deux diff\'erentiations, ce qui est confirm\'e
par la pr\'esence de la d\'eriv\'ee seconde $f''(x)$ dans l'expression
quantitative compl\`ete de $c(x)$. On dit que la courbure
est une quantit\'e {\sl du second ordre}.
\qed

\medskip\noindent{\bf Expression de la courbure en représentation
paramétrée.}\label{courbure-courbes-parametres}
Par des calculs analogues \`a ceux qui viennent d'être conduits, on
trouve l'expression explicite de la courbure d'une courbe, lorsqu'elle
est repr\'esent\'ee sous une forme paramétrique:
\[
t
\longmapsto
\big(x(t),\,y(t)\big),
\]
où la variable $t$ appartient à $\R$. 

\smallskip\noindent{\bf Théorème.}
{\em Dans une
telle repr\'esentation param\'etr\'ee, la courbure {\rm orient\'ee}
d'une courbe en un point de coordonn\'ees $\big( x(t),\, y(t) \big)$ 
s'exprime quantitativement au moyen de la formule:}
\[
c(t) = 
\frac{ x'(t) \, y''(t) - y'(t) \, x''(t)}{ 
\big[
(x'(t))^2 + (y'(t))^2
\big]^{ 3/2}}.
\]

\proof 
Le vecteur normal unitaire $N(t)$ \`a la courbe param\'etr\'ee
$(x(t),\, y(t))$ en un point rep\'er\'e par $t$ a pour coordonn\'ees:
\[
N(t) = 
\left(
- \frac{ y'(t) }{ \sqrt{ x'(t)^2+ 
y'(t)^2}}, \ \ 
\frac{ x'(t) }{ \sqrt{ x'(t)^2+ 
y'(t)^2}}
\right).
\]
La longueur infinit\'esimale de la courbe entre le point rep\'er\'e
par $t$ et le point rep\'er\'e par $t+ dt$ est \'egale \`a $\sqrt{
x'(t)^2+ y'(t)^2} \cdot dt$. La courbure (non orient\'ee) est donn\'ee
par:
\[
c(t) 
=
\frac{\text{\rm norme de}\ N'(t)}{ 
\sqrt{x'(t)^2+ y'(t)^2}}.
\]
Il suffit donc de calculer les deux composantes de $N'(t)$:
\[
N'(t) = 
\left(
\frac{ y' \, x' \, x'' - y'' \, x' \, x' }{(x'{}^2+ 
y'{}^2 ) \sqrt{ x'{}^2+ y'{}^2}}, \ \ 
\frac{ x''\, y'\, y'- x' \, y' \, y'' }{(x'{}^2+ 
y'{}^2) \sqrt{ x'{}^2+ 
y'{}^2}}
\right),
\]
puis la norme de $N'(t)$, pour trouver l'expression annoncée.
\endproof 

\medskip\noindent{\bf Expression de la courbure en représentation
implicite.}\label{courbure-courbes-implicite}
Enfin, on peut supposer que la courbe est donnée sous forme implicite:
\[
0
=
W(x,y),
\]
où $W$ est une fonction de classe au moins $\mathcal{ C}^2$ par
rapport à ses deux variables $(x, y)$, et où l'on suppose qu'en tout
point du lieu des zéros $\{ (x_p, y_p) \in \R^2 \colon W(x_p, y_p) =
0\}$, l'une au moins des deux dérivées partielles $W_x (x_p, y_p)$ ou
$W_y (x_p, y_p)$ ne s'annule pas, 
condition qui assure que la courbe est
géométriquement lisse. Cette condition s'accorde avec le fait
que l'on doit diviser par $W_x^2 + W_y^2$ pour calculer
quantitativement la courbure.

\smallskip\noindent{\bf Théorème.}
{\em Dans la repr\'esentation implicite, en un point quelconque de
coordonn\'ees $(x,\, y)$ appartenant \`a la courbe d'\'equation
$W(x,\, y) = 0$, la courbure
s'exprime quantitativement au moyen de la formule:}
\[
c(x,\, y) 
=
\frac{W_{xx}\,(W_y)^2 
- 
2\,W_{ xy}\,W_x\,W_y 
+ 
W_{ yy}\,(W_x)^2}{
\big[(W_x)^2+(W_y)^2\big]^{ 3/2}}(x,y).
\]

\proof
Deux démonstrations s'offrent: la première, indirecte car basée sur le
précédent théorème, est peu coûteuse en calcul; la seconde est directe
et indépendante, mais elle requiert des calculs d'élimination
différentiels moins élémentaires. Il sera néanmoins utile 
d'effectuer de
tels calculs afin de mieux préparer la présentation de la courbure
dans le cas beaucoup plus délicat des surfaces.

Pour la première démonstration (indirecte), localement au voisinage
d'un point $p = (x_p, y_p)$ quelconque de la courbe, on peut, par
symétrie des coordonnées, supposer que $W_y ( x, y) \neq 0$ ne
s'annule pas. Grâce au théorème des fonctions implicites, la courbe
est alors en fait représentable comme le graphe $y = f ( x)$ d'une
certaine fonction de classe au moins $\mathcal{ C}^2$, définie dans le
voisinage de $p$ et satisfaisant $y_p = f ( x_p)$. Par définition, ce
graphe entier est contenu dans l'ensemble des zéros de $f$, d'où il
découle que l'équation:
\[
0
\equiv
W\big(x,f(x)\big)
\]
est satisfaite identiquement\footnote{\,
%%%%%%%%%%%%%%%%%%%%%%%-------DEBUT--------%%%%%%%%%%%%%%%%%%%%%%%%%%%
Dans les {\oe}uvres de Sophus Lie, des raisonnements par équations
identiquement satisfaites s'avèrent être très fréquemment nécessaires
et constamment utiles, et ils sont aussi poussés, presque toujours
avec un nombre quelconque de variables, à un niveau de généralité qui
dépasse de mille coudées le cas des courbes.
} %%%%%%%%%%%%%%%%%%%%%%%%-----FIN-----%%%%%%%%%%%%%%%%%%%%%%%%%%%%%%%
pour tout $x$ dans un certain voisinage ouvert non vide de $x_p$. En
tout état de cause, cette équation doit dans l'instant être
différentiée par rapport à $x$, ce qui donne une autre équation
identiquement satisfaite:
\[
0
=
W_x+f'\,W_y
\]
{\em i.e.} avec les arguments:
\[
0
\equiv
W_x\big(x,f(x)\big)+f'(x)\,W_y\big(x,f(x)\big),
\]
d'où l'on tire la valeur, en un point arbitraire $\big(x, f(x) \big)$ 
de la courbe proche de $p$, de la dérivée de la fonction
graphante:
\[
f'
=
-\,
{\textstyle{\frac{W_x}{W_y}}}
\ \ \ \ \ 
\text{\rm {\em i.e.} avec les arguments:}
\ \ \ \ \
f'(x)
=
-\,
{\textstyle{\frac{W_x(x,f(x))}{W_y(x,f(x))}}},
\]
comme simple quotient des dérivées partielles d'ordre $1$ de 
la fonction implicitante $W$. Maintenant, dans l'expression
de la courbure: 
\[
c(x)
=
\frac{f''(x)}{\big(1+f'(x)^2\big)^{3/2}}
\]
en termes de cette fonction graphante $f = f (x)$, non seulement la
dérivée d'ordre $1$, mais encore la dérivée d'ordre $2$ de $f$
interviennent. Rien de plus immédiat que la suggestion de dériver
encore une fois l'équation identiquement satisfaite, ce qui donne, en
n'oubliant aucun terme dans le calcul des dérivées partielles, une
nouvelle équation identiquement satisfaite:
\[
0
\equiv
W_{xx}+W_{xy}\,f'+f''\,W_y+f'\,W_{xy}+(f')^2\,W_{yy},
\] 
écrite maintenant sans les arguments $\big(x, f(x) \big)$ des dérivées
partielles de la fonction $W$. Résoudre $f''$ est aisé:
\[
f''\,W_y
=
-\,W_{xx}
+
2\,W_{xy}\,
{\textstyle{\frac{W_x}{W_y}}}
-
W_{yy}\,
{\textstyle{\frac{W_x^2}{W_y^2}}},
\]
et lorsqu'on remplace $f'$ et $f''$ par les valeurs ainsi obtenues en
fonction de $W_x$, $W_y$, $W_{ xx}$, 
$W_{ xy}$ et $W_{ yy}$, on obtient effectivement, après
harmonisation, l'expression annoncée pour
la courbure:
\[
c
=
\frac{W_{xx}\,(W_y)^2 
- 
2\,W_{ xy}\,W_x\,W_y 
+ 
W_{ yy}\,(W_x)^2}{
\big[(W_x)^2+(W_y)^2\big]^{ 3/2}}.
\]

\medskip\noindent{\bf Deuxième preuve directe systématique.}
\label{preuve-directe-systematique}
La deuxième preuve ne supposera pas qu'une expression de la courbure
ait déjà été obtenue par une autre voie, elle calculera
{\em directement} la variation d'angle du vecteur normal
à la courbe en se rapportant à l'équation implicite.
Mais tout d'abord, comment s'écrit le vecteur normal à une
courbe en représentation implicite? 

Si deux points infiniment proches l'un de l'autre $(x + dx, y + dy)$
et $(x, y)$ appartiennent à la courbe, à savoir si:
\[
\aligned
0
=
W(x,y)
&
=
W(x+dx,y+dy)
\\
&
=
\zero{W(x,y)}
+
W_x(x,y)\,dx
+
W_y(x,y)\,dy,
\endaligned
\]
alors le vecteur infinitésimal $(dx, dy)$ est
évidemment tangent à la courbe
et la condition ainsi obtenue:
\[
0
=
W_x\,dx+W_y\,dy
\]
exprime visiblement que
$(dx, dy)$ est orthogonal au vecteur $(W_x, W_y)$ pour
le produit scalaire euclidien standard sur $\R^2$. De manière
immédiatement équivalente, ce vecteur $(W_x, W_y)$ est orthogonal à la
droite tangente, et donc il suffit de le diviser par
sa norme pour obtenir le vecteur normal unitaire:
\[
N
=
N(x,y)
=
\Big(
{\textstyle{\frac{W_x}{(W_x^2+W_y^2)^{1/2}}}},\,\,
{\textstyle{\frac{W_y}{(W_x^2+W_y^2)^{1/2}}}}
\Big)
\]
à la courbe au point $(x, y)$.

Maintenant par définition, la courbure s'exprime comme la variation
infinitésimale du vecteur normal divisée par la longueur du segment
infinitésimal de $(x, y)$ à $(x + dx, y + dy)$, c'est-à-dire:
\[
c(x,y)
=
\frac{\bignorm N(x+dx,y+dy)-N(x,y)\bignorm}{\sqrt{dx^2+dy^2}}.
\]
Il est donc nécessaire de calculer la différentielle:
\[
\aligned
dN
&
=
N(x+dx,y+dy)-N(x,y)
\\
&
=
N_x(x,y)\,dx
+
N_y(x,y)\,dy
\endaligned
\]
du vecteur normal en fonction des dérivées partielles de la fonction
$W$, et c'est ici-même que les calculs commencent à se déployer.

En effet, une application soignée des règles classiques de
différentiation de fonctions composées, produit ou quotient, donne
les différentielles des deux composantes du vecteur
normal $N$:
\[
\footnotesize
\aligned
dN
&
=
\bigg(
\frac{W_{xx}\,dx+W_{xy}\,dy}{(W_x^2+W_y^2)^{1/2}}
-
\frac{
[W_x^2\,W_{xx}+W_x\,W_y\,W_{xy}]\,dx
+
[W_x^2\,W_{xy}+W_x\,W_y\,W_{yy}]\,dy}{
(W_x^2+W_y^2)^{3/2}},\,\,
\\
&
\ \ \ \ \ ,\ \
\frac{W_{xy}\,dx+W_{yy}\,dy}{(W_x^2+W_y^2)^{1/2}}
-
\frac{
[W_y\,W_x\,W_{xx}+W_y^2\,W_{xy}]\,dx
+
[W_y\,W_x\,W_{xy}+W_y^2\,W_{yy}]\,dy}{
(W_x^2+W_y^2)^{3/2}}
\bigg),
\endaligned
\]
et l'on va constater, après réduction (nécessaire) au même
dénominateur:
\[
\tiny
\!\!\!\!\!\!\!\!\!\!\!\!\!\!\!\!\!\!\!\!\!\!\!\!\!\!\!\!\!\!
\aligned
dN
&
=
\bigg(
\frac{
W_{xx}(\zero{W_x^2}+W_y^2)\,dx
+
W_{xy}(\zerozero{W_x^2}\!+W_y^2)\,dy
-
[\zero{W_x^2\,W_{xx}}\!+W_x\,W_y\,W_{xy}]\,dx
-
[\zerozero{W_x^2\,W_{xy}}\!+W_x\,W_y\,W_{yy}]\,dy}{
(W_x^2+W_y^2)^{3/2}},\,\,
\\
&
\ \ \ \ \ ,\ \
\frac{
W_{xy}(W_x^2+\zero{W_y^2})\,dx
+
W_{yy}(W_x^2+\zerozero{W_y^2})\,dy
-
[W_y\,W_x\,W_{xx}+\zero{W_y^2\,W_{xy}}]\,dx
-
[W_y\,W_x\,W_{xy}+\zerozero{W_y^2\,W_{yy}}]\,dy}{
(W_x^2+W_y^2)^{3/2}}
\bigg),
\endaligned
\]
que pour chacune des deux composantes de
ce vecteur infinitésimal $dN \in \R^2$, deux paires de
termes, soulignées, s'annihilent, d'où il reste, après ces
simplifications agréables:
\[
\small
\aligned
dN
&
=
\frac{1}{(W_x^2+W_y^2)^{3/2}}
\bigg(
[W_y^2\,W_{xx}-W_x\,W_y\,W_{xy}]\,dx
+
[W_y^2\,W_{xy}-W_x\,W_y\,W_{yy}]\,dy,
\\
&
\ \ \ \ \ \ \ \ \ \ \ \ \ \ \ \ \ \ \ \ \ \ \ \ \ \ \ \ \ \ \ \ \ \ 
,\ \
[W_x^2\,W_{xy}-W_x\,W_y\,W_{xx}]\,dx
+
[W_x^2\,W_{yy}-W_y\,W_x\,W_{xy}]\,dy
\bigg).
\endaligned
\]
Il est maintenant possible de calculer la norme au carré de la
différentielle, après placement à gauche du dénominateur:
\[
\scriptsize
\aligned
&
(W_x^2+W_y^2)^3\,
\bignorm dN\bignorm^2
=
dx^2
\Big[
W_y^2(W_y\,W_{xx}-W_xW_{xy})^2
+
W_x^2(W_x\,W_{xy}-W_y\,W_{xx})^2
\Big]
+
\\
&
+
2dxdy
\Big[
W_y^2(W_yW_{xx}-W_xW_{xy})(W_yW_{xy}-W_xW_{yy})
+
W_x^2(W_xW_{xy}-W_yW_{xx})(W_xW_{yy}-W_yW_{xy})
\Big]
+
\\
&
\ \ \ \ \ \ \ \ \ \ \ \ \ \ \ \ \ \ \ \ \ \ \ \ \ \ \ \ \ \ \ \ \ \
\ \ \ \ 
+dy^2
\Big[
W_y^2(W_y\,W_{xy}-W_x\,W_{yy})^2
+
W_x^2(W_x\,W_{yy}-W_y\,W_{xy})^2
\Big].
\endaligned
\]

\Scholie
Rappel sur l'intention du calcul: en principe, la courbure au carré
$c(x,y)^2$ de la courbe en un point $(x, y)$ situé sur le lieu
implicite des zéros $\{ W ( x, y) = 0\}$ s'obtient en {\em divisant}
cette norme au carré $\norm dN \norm^2$ par $dx^2 + dy^2$. Or, il
n'est nullement évident ici que le membre de droite soit effectivement
divisible par $dx^2 + dy^2$! Outre l'inévitable erreur de calcul
(évitée ici dans l'{\em a posteriori} de la constitution textuelle
d'un manuscrit élaboré et relu de nombreuses fois), l'intuition
chercheuse doit {\em interroger} ce membre de droite et lui demander
instamment s'il ne serait pas possible de le réécrire comme $dx^2 +
dy^2$ que multiplie une certaine expression
algébrico-différentielle. Outre aussi l'inévitable facilité que
procure la connaissance, offerte auparavant, de l'expression à
obtenir\,\,---\,\,à savoir justement, le carré du numérateur déjà
présenté dans le théorème ci-dessus\,\,---, l'intuition chercheuse
doit examiner tous les germes possibles de transformations symboliques
signifiantes qui demeurent peut-être déposés au sein de la dernière
équation écrite. Aboutissement, suspension, rebondissement: telle est
la dynamique, fût-elle itérée de nombreuses fois au cours des
analyses. En tout cas, le carré caché est toujours invisible, tel est
le suspens.

On se rapproche ici de l'idée tentatrice qu'un certain {\em
créationnisme de l'algébrique et du conceptuel} se développe au
contact même et à l'intérieur même des calculs ardus, comme une {\em
rencontre de réalités mathématiques potentielles} le long de ces
chemins tracés dans l'inconnu que sont les intentions (ou conjectures)
{\em provocatrices} de mouvements de recherche.
\qed

\Scholie
Quelles que soient les tentatives, il est en effet remarquable ici que
le chemin direct vers la courbure, en partant de l'équation implicite,
engendre d'imprévisibles obstacles d'illisibilité dans le flot continu
et interne des calculs. En effet, contre l'attente précédente, 
puisque les cinq quantités $W_x$, $W_y$, $W_{
xx}$, $W_{ xy}$ et $W_{ yy}$ peuvent recevoir des valeurs réelles
arbitraires en tout point $(x, y)$ de la courbe, il est en fait
impossible, en général, que le membre de droite 
ci-dessus soit
multiple de $dx^2 + dy^2$: errances de l'intuition, fausses pistes
dans le labyrinthe partiel d'un Inconnu pré-exploré. Or pourquoi?
Pourquoi en est-il ainsi? Pourquoi cet obstacle qui obstrue et qui
bloque?

Dans toute exploration mathématique, seuls les {\em retours en
arrière} et les {\em reprises générales} permettent de sauver
dialectiquement les phénomènes: remise en cause d'options ou
d'hypothèses, relèvement aux points de bifurcations incidemment
ignorés et défiance à l'égard de pseudo-métaphysiques non supportées
par des spéculations adéquates.

Ici, c'est la mémoire qui a failli: oubli regrettable que les
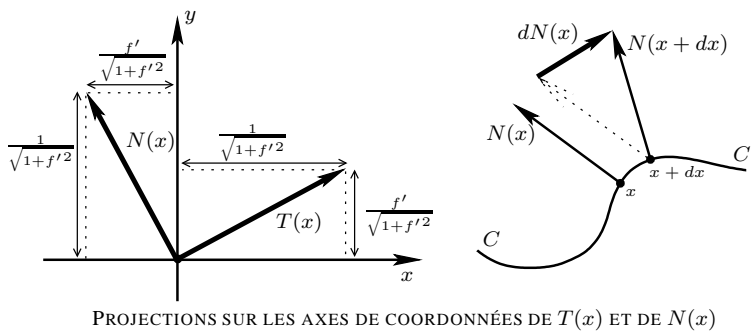
projections $dx$ et $dy$ sur les axes de coordonnées de l'élément
différentiel infinitésimal étaient en fait {\em liées depuis le début}
par une relation fondamentale:
\[
0
=
W_x\,dx+W_y\,dy,
\]
dont il fallait fatalement tenir compte dans tous les calculs. Mais à
présent, l'espoir renaît de {\em faire aboutir les calculs}.
L'erreur conscientisée est en fait salvatrice.
\qed

\medskip
Ainsi, on supposera que $W_y \neq 0$ (ce qui est justifié par
interchangeabilité des variables $x$ et $y$), et on remplacera $dy$
partout par $- \frac{ W_x}{ W_y}\, dx$. Mais à l'instant même de la
relecture de la précédente équation dans laquelle 
on s'apprêtait à effectuer cette
substitution de $dy$, une observation imprévue car insoupçonnée et
agréable se révèle: {\em les trois différentielles quadratiques
$dx^2$, $2dx dy$ et $dy^2$ sont toutes 
multiples de $W_x^2 + W_y^2$}, et donc
l'on peut {\em simplifier à gauche et à droite} par ce facteur commun
avant d'effectuer le remplacement, ce qui donne pour
commencer:
\[
\small
\aligned
(W_x^2+W_y^2)^2\,
\bignorm dN\bignorm^2
&
=
dx^2
(W_y\,W_{xx}-W_xW_{xy})^2
+
\\
&\ \ \ \ \
+
2dxdy
(W_yW_{xx}-W_xW_{xy})(W_yW_{xy}-W_xW_{yy})
+
\\
&
\ \ \ \ \
+
dy^2
(W_y\,W_{xy}-W_x\,W_{yy})^2.
\endaligned
\]

\Scholie
Ce sera seulement quand les machines à calcul symbolique posséderont
une panoplie d'actes dialectiques suffisamment riche pour {\em réagir
dialectiquement} dans l'imprévisibilité bifurcationnelle du calcul
irréversiblement synthétique-orienté que l'on pourra prétendre que le
calcul est véritablement automatisable. 
\qed\medskip

Une fois cette simplification
acquise, on peut remplacer $dy$ comme convenu par $- \frac{ W_x}{
W_y}\, dx$ d'abord dans le dénominateur $dx^2 + dy^2$ de la courbure
au carré:
\[
c(x,y)^2
=
\frac{\norm dN\norm^2}{dx^2+dy^2}
=
\frac{1}{W_x^2+W_y^2}\,
\frac{W_y^2\,\norm dN\norm^2}{dx^2},
\]
puis effectuer ce même remplacement dans $W_y^2\, \norm dN \norm^2$,
en constatant que le facteur $W_y^2$ élimine les dénominateurs, ce qui
donne au total:
\[
\footnotesize
\aligned
c^2
=
\frac{1}{(W_x^2+W_y^2)^3}\,
\bigg\{\,
&
\frac{
\zero{dx^2}\,W_y^2(W_y\,W_{xx}-W_xW_{xy})^2
}{\zero{dx^2}}
-
\\
&\
-
\frac{
2\zero{dx^2}W_xW_y(W_yW_{xx}-W_xW_{xy})(W_yW_{xy}-W_xW_{yy})
}{\zero{dx^2}}
+
\\
&\
+
\frac{
\zero{dx^2}W_x^2(W_y\,W_{xy}-W_x\,W_{yy})^2
}{\zero{dx^2}}
\,\bigg\},
\endaligned
\]
et après la simplification évidente de $dx^2$ au numérateur et au
dénominateur, il reste une expression entre crochets:
\[
\footnotesize
\aligned
c^2
=
\frac{1}{(W_x^2+W_y^2)^3}\,
\bigg\{\,
&
W_y^2(W_y\,W_{xx}-W_xW_{xy})^2
-
\\
&\
-
2\,W_xW_y(W_yW_{xx}-W_xW_{xy})(W_yW_{xy}-W_xW_{yy})
\\
&\
+
W_x^2(W_y\,W_{xy}-W_x\,W_{yy})^2
\,\bigg\},
\endaligned
\]
que l'{\oe}il averti du calculateur reconnaîtra comme étant de la forme
générale $A^2 - 2AB + B^2 = (A-B)^2$:
\[
\footnotesize
\aligned
c^2
=
\frac{1}{(W_x^2+W_y^2)^3}\,
\bigg\{\,
&
\Big(
W_y(W_y\,W_{xx}-W_xW_{xy})
-
W_x(W_y\,W_{xy}-W_x\,W_{yy})
\Big)^2
\,\bigg\},
\endaligned
\]
ce qui se simplifie pour terminer en:
\[
c(x,y)^2
=
\frac{W_x^2\,W_{yy}-2W_xW_y\,W_{xy}+W_y^2W_{xx}}
{(W_x^2+W_y^2)^3},
\]
et qui conclut finalement la deuxième démonstration du théorème.
\endproof

\medskip\noindent{\bf \'Evanouissement intrinsèque de la courbure
des courbes.} \label{evanouissement-intrinseque-courbes}
Toute courbe suffisamment régulière (de classe $\mathcal{ C}^1$
convient) s'identifie essentiellement
\`a une succession de segments
rectilignes infinit\'esimaux. Si l'on imagine que chacun de ces
segments est reli\'e aussi bien \`a celui qui le pr\'ec\`ede et \`a
celui qui le suit par un pivot autour duquel la rotation est libre, de
telle sorte que la courbe est articul\'ee comme une r\`egle de
charpentier, il est alors clair pour des raisons m\'ecaniques
évidentes, que l'on peut <<\,d\'e-courber\,>> la courbe et la
transformer en une cha\^{\i}ne absolument rectiligne en tirant
simplement avec ses mains sur les deux extr\'emit\'es.

\medskip

\begin{center}
\input articulation-chaine.pstex_t
\end{center}

Cette op\'eration conserve non seulement la longueur de la courbe, 
mais elle conserve aussi la longueur de chacun de ses
segments infinit\'esimaux, qui sont repr\'esent\'es par
des segments de longueur finie sur la figure.

\smallskip\noindent{\bf Observation.}
{\em 
Toute courbe dans le plan peut \^etre transform\'ee en un segment de
droite par une transformation qui conserve la longueur de tous ses
\'el\'ements infinit\'esimaux.}\medskip

Plus g\'en\'eralement, on peut transformer toute courbe quelconque en
toute autre courbe quelconque de telle sorte que leurs longueurs
infinit\'esimales se correspondent. Autrement dit, on peut {\em
supprimer intrinsèquement la courbure de toutes les courbes}.

\Section{6.~Opposition dialectique impr\'evisible du dimensionnel}
\label{opposition-imprevisible-dimensionnel}

\HEAD{6.~Opposition dialectique impr\'evisible du dimensionnel}{
Jo\"el Merker, D\'partement de Math\'matiques d'Orsay}

\medskip\noindent{\bf Anticipation du Theorema egregium.}
\label{anticipation-theorema-egregium}
Au contraire, la plupart des surfaces ne peuvent {\em pas}\, \^etre
transform\'ees en un plan par une transformation qui conserve les
longueurs de leurs \'el\'ements infinit\'esimaux. Par exemple, pour
aplatir (écraser) l'\'ecorce d'une demi-orange pos\'ee sur une table,
il est n\'ecessaire qu'elle se d\'echire en plusieurs endroits. C'est
Gauss qui le premier a su \'elaborer un concept ad\'equat de {\sl
courbure} d'une surface et a pu <<\,expliquer\,>> math\'ematiquement
l'in\'equivalence entre une surface plane et une surface
sph\'erique. Avant de pr\'esenter la th\'eorie g\'en\'erale, 
on d\'eveloppe l'\'enonc\'e difficile et fondamental que
Gauss a obtenu comme une cons\'equence directe de la {\em formula
egregia}~\thetag{ \ref{formula-egregia}}.

\smallskip\noindent{\bf Theorema egregium.}
{\em Si l'on transforme une surface $S$ en une autre surface
$S'$ de telle sorte que les longueurs infinit\'esimales de
toutes les courbes trac\'ees sur les surfaces soient conserv\'ees,
alors la courbure de Gauss en un point $p$ de la premi\`ere surface
$S$ est \'egale \`a la courbure de Gauss en le point $p'$
qui lui correspond sur la seconde surface $S'$.}\medskip

Ainsi na\^{\i}t une {\em opposition dialectique impr\'evisible}\,
entre la dimension $1$ (les courbes) et la dimension $2$ (les
surfaces):

\begin{itemize}

\smallskip\item[$\bullet$]
la courbure des courbes planes est modifiable \`a volont\'e en
conservant les longueurs infinit\'esimales: il suffit de
se représenter une ficelle souple mais {\em inextensible};

\smallskip\item[$\bullet$]
au contraire,
la courbure de Gauss des surfaces dans l'espace demeure {\em invariante}
à travers
toutes les transformations qui préservent infnitésimalement 
les longueurs: pour modifier la courbure d'une surface, 
seules des surfaces qui sont
à la fois {\em souples} et {\em extensibles} conviendraient.

\end{itemize}\smallskip

\Scholie
Dans son r\'eservoir infini et in\'epuisable de complexit\'e et de
nouveaut\'e\,\,---\,\,dont le th\'eor\`eme d'incompl\'etude de G\"odel
fournit une image partielle\,\,---\,\,, le monde math\'ematique
r\'eserve des surprises qui an\'eantissent les espoirs trop simplistes
de l'induction hâtive. Incomparabilit\'e des th\'eories \`a une
variable et des th\'eories \`a deux variables (ou plus): ce
ph\'enom\`ene g\'en\'eral traverse toutes les math\'ematiques
r\'ecentes, {\em cf.} le {\sl phénomène de Hartogs} (1906), 
fondateur de l'analyse moderne à plusieurs variables
complexes.
\qed

\Section{7.~Courbure des surfaces dans l'espace}
\label{courbure-surfaces-espace}

\HEAD{7.~Courbure des surfaces dans l'espace}{
Jo\"el Merker, D\'partement de Math\'matiques d'Orsay}

\medskip\noindent{\bf Détermination du plan tangent à une surface.}
\label{plan-tangent-surface}
\`A la suite de travaux précurseurs
de Parent au début du dix-huitième
siècle\footnote{\,
%%%%%%%%%%%%%%%%%%%%%%%-------DEBUT--------%%%%%%%%%%%%%%%%%%%%%%%%%%%
{\em Cf.} Liberman~\cite{li1978}, pp.~181--183
et Kreyszig~\cite{kr1994}.
}, %%%%%%%%%%%%%%%%%%%%%%%%-----FIN-----%%%%%%%%%%%%%%%%%%%%%%%%%%%%%%%
Clairaut, Euler et Monge utilisent la notion géométrique et
différentielle de plan tangent à une surface, en admettant que la
lissité est génériquement satisfaite lorsque l'équation de la
surface s'exprime au moyen de fonctions élémentaires générales. Le
fameux théorème de Clairaut\,\,---\,\,qu'il 
obtint à l'âge de seize ans\,\,---\,\,énonce
que le long d'une géodésique $C$ tracée sur une surface $S$ de
révolution, on a:
\[
r(p)\sin\theta(p)
=
{\rm const.},
\]
où $r(p)$ est le rayon du cercle parallèle passant par un point
quelconque $p \in S$ de la surface, et où $\theta (p)$ est l'angle que
la géodésique $C$ fait avec le méridien passant par $p$.

En 1803, Lacroix interprète la condition de tangence d'un plan à une
surface en termes de proximité maximale possible, c'est-à-dire à
l'ordre deux, précisément
comme pour le cas plus
simple des droites tangentes à une
courbe. Dupin en 1813 et Cauchy en 1826 établissent que
la collection des droites tangentes à toutes les courbes tracées
sur une surface passant par un point fixe engendrent\,\,---\,\,et 
sont contenues dans\,\,---\,\,un 
unique {\sl plan tangent} à la surface, et ils en
donnent l'équation affine:
\[
0
=
(\xi-x)\,W_x
+
(\upsilon-y)\,W_y
+
(\zeta-z)\,W_z,
\]
lorsque la surface est fournie en représentation
implicite par une équation générale du type:
\[
0
=
W(x,y,z),
\]
pour une certaine fonction $W = W ( x, y, z)$ 
définie au voisinage de la surface dans $\R^3$.

\medskip\noindent{\bf Vecteurs unitaires normaux à une surface.}
\label{vecteurs-unitaires-normaux}
Soit maintenant $S \subset \R^3$ une surface locale ou
globale, saisie dans l'une des trois représentations possibles:
graphée, implicite, ou paramétrique\,\,---\,\,comme dans la
Section
p.~\pageref{preliminaire-surfaces} ci-dessus. 
Soit $p \in S$ un point
quelconque. Le plan tangent $T_p S$ à la surface en $p$ étant de
codimension $1$ dans $\R^3$, il existe une unique droite {\em
orthogonale} à ce plan\,\,---\,\,soit donc $D_p$ cette droite passant par
$p$. Les considérations de courbure qui vont suivre étant locales, on
peut supposer que la surface est munie d'une certaine orientation;
lorsque $S$ est globale, on supposera implicitement qu'elle est
orientable, 
et qu'elle est munie d'une orientation choisie parmi les deux
orientations opposées qui sont possibles.

Ainsi, des deux directions de sortie que définit la droite $D_p$,
l'une d'elles est privilégiée, à savoir celle qui forme un
trièdre direct dans $\R^3$ avec tout dièdre de $T_pS$ direct pour
l'orientation de la surface. Par conséquent, sur la surface, on peut
définir une certaine collection $p \mapsto {\bf n}(p)$ de vecteurs
dirigés par les droites $D(p)$ dans la direction privilégiée en
question, donc {\em normaux} ({\em i.e.}
orthogonaux à la surface), et qui sont de plus {\em unitaires},
c'est-à-dire de norme $1$, 
deux conditions qui les déterminent en tout point de manière
unique et sans ambiguïté. On a vu que dans le cas des courbes, y
compris celles qui seraient tracées sur les surfaces, c'est la
variation infinitésimale du vecteur unitaire normal qui permet
d'introduire un concept adéquat de courbure. Comment cette vision
doit-elle se généraliser en dimension supérieure?

\Scholie
Les liens entre les dimensions d'étude s'expriment par analogie et
généralisation. Tout concept acquis en petite dimension appelle son
extension dans le cas des dimensions immédiatement supérieures, puis
dans le cas le plus général possible de la dimension arbitraire. Dans
les travaux ultérieurs de Lie, l'étude de la dimension quelconque $n$
avec un nombre arbitraire $r$ de paramètres de mouvement constituera
un objectif décidé et systématique de la théorie des groupes de
transformation.

En vérité, l'exigence de généralisation exprime une tension
universelle et immédiatement éprouvable, tension
à laquelle doivent se
soumettre les réalisations mathématiques. La philosophie des
mathématiques quant à elle se trouve en quelque sorte {\em déchirée}
par l'incapacité actuelle de l'esprit humain 
qui aurait à embrasser
l'extension métaphysique de cette exigence universelle, puisque à la
fois le général est du ressort de la philosophie, et à la fois le
général n'est pas un général abstrait et indifférencié dans les
mathématiques, puisqu'il se réalise dans tous les domaines de
l'Algèbre, de l'Analyse, de la Géométrie en s'individuant de manière
imprévisible et extrêmement diverse. Quel peut être alors le statut
épistémologique de l'exigence mathématique de généralité? Faut-il y
voir un aspect transcendantal? Ou bien seulement voir en elle un
principe simple, accessible et reproductible qui signale
indéfiniment l'ouverture en tant que {\em recherche}?
En tout cas, la tension en direction de la généralité est omniprésente
dans les mathématiques, toujours germe fertile de travaux nouveaux et
originaux, et c'est donc, à tout le moins, un {\em principe moteur}
qu'on doit voir en elle dans la pratique.
\qed

\smallskip
Maintenant, lorsque la surface
dans l'espace $\R^3$ est représentée sous la forme
d'un graphe $z = z (x, y)$ pour lequel la troisième 
coordonnée est une certaine fonction
définie mais quelconque des deux premières
coordonnées, 
un petit vecteur infinitésimal de coordonnées
$(dx, dy, dz)$ placé en un point $\big(x, y, z(x, y) \big)$
de la surface est {\em tangent} à la surface
en ce point si et seulement si
la variation infinitésimale de la
troisième coordonnée 
$z$ suit, par simple prise
de différentielle, 
la contrainte de son égalité à la fonction $z ( x, y)$:
\[
dz
=
z_x\,dx+z_y\,dy.
\]
Ainsi, les deux vecteurs de coordonnées:
\[
{\bf t}_x
:=
(1,0,z_x)
\ \ \ \ \ \ \ \ \
\text{\rm et}
\ \ \ \ \ \ \ \ \
{\bf t}_y
:=
(0,1,z_y)
\]
appartiennent à ce plan tangent et ils l'engendrent
linéairement. Quelles sont alors
les trois coordonnées:
\[
\big(
X(x,y),\,Y(x,y),\,Z(x,y)
\big)
=
{\bf n}(x,y)
\]
du vecteur normal unitaire à la surface en un point quelconque
repéré par ses deux coordonnées horizontales
$(x, y)$? La réponse est simple: l'orthogonalité
aux deux vecteurs tangents ci-dessus: 
\[
0
=
{\bf t}_x\cdot{\bf n}
\ \ \ \ \ \ \ \ \
\text{\rm et}
\ \ \ \ \ \ \ \ \
0
=
{\bf t}_y\cdot{\bf n}
\]
s'exprime de manière équivalente
comme système de deux équations linéaires:
\[
\aligned
0
&
=
X+0+z_x\,Z
\\
0
&
=
0+Y+z_y\,Z,
\endaligned
\]
mais alors la condition que le vecteur ${\bf n}$ est de norme unité:
\[
X^2+Y^2+Z^2
=
1
\]
permet de résoudre uniquement\,\,à un signe
$\pm$ global près\,\,---\,\,ledit système:
\[
\footnotesize
\aligned
\pm X
=
\frac{z_x}{(1+z_x^2+z_y^2)^{1/2}},
\ \ \ \ \ \ \ \
\pm Y
=
\frac{z_y}{(1+z_x^2+z_y^2)^{1/2}},
\ \ \ \ \ \ \ \
\pm Z
=
\frac{-1}{(1+z_x^2+z_y^2)^{1/2}}.
\endaligned
\]

\medskip\noindent{\bf Transfert sur la sphère auxiliaire et orientation.}
\label{transfert-auxiliaire}
Sans réveiller la métaphysique du passage à la dimension
supérieure, ni même expliciter les motivations astronomiques de
repèrage de la voûte céleste sur une sphère projective orientée idéale
située à l'infini du monde observable, Gauss introduit d'emblée dans
son mémoire la {\sl mesure de courbure}\,\,---\,\,{\em mensura
curvatur{\ae}}\,\,---\,\,d'une surface en termes purement
géométriques, d'abord non analytiques, de
{\em représentation imagée} sur le
`miroir-modèle' de la sphère.

\CITATION{
De même qu'en imaginant par le centre de notre sphère
auxiliaire des droites respectivement parallèles à chacune des
normales d'une surface courbe, à chaque point déterminé de la deuxième
surface vient correspondre un point déterminé de la première; de la
même manière, toute ligne ou toute figure tracée sur la surface courbe
sera représentée par une ligne ou une figure tracée sur la surface
sphérique. Dans la comparaison des deux figures qui se correspondent
ainsi, et dont l'une sera comme l'image de l'autre, on peut se placer
à deux points de vue: on peut avoir égard seulement aux quantités; ou
bien ne s'occuper que des relations de position, abstraction faite des
relations de quantité.
\REFERENCE{{\sc Gauss},~\cite{ga1828},~13--14}}

\Scholie
Ici, le langage est visiblement inspiré par
les catégories de la philosophie: quantité, qualité, 
deux points de vue possibles et complémentaires sur la {\em
correspondance} des figures, de la surface, vers la sphère auxiliaire,
{\em via} l'application de Gauss qui transporte
tous les vecteurs unitaires
normaux à l'origine du système de coordonnées. Une mentalisation
difficile du transfert, notamment du transfert des courbes, est
évoquée, puisque Gauss pense à des figures quelconques tracées sur la
surface, mais cet aspect de la conception géométrique disparaît
essentiellement complètement lorsqu'il est ressaisi symboliquement,
dans un traité moderne, en termes d'{\em application ensembliste}
d'une surface $S$ vers la sphère auxiliaire. En effet, le problème de
concevoir une {\em transformation} entre deux espaces bidimensionnels,
qui est un problème de pensée spécifique à l'étude des surfaces
beaucoup plus délicat que dans le cas unidimensionnel\,\,---\,\,car
les graphes d'application de $\R^2$ dans $\R^2$ vivent dans un espace
à quatre dimensions alors que les graphes de fonctions de $\R$ dans
$\R$ vivent simplement dans le plan\,\,---\,\,ce problème est en
quelque sorte noyé et éliminé par l'option
théorique\,\,---\,\,auquel les
enseignements et traités se cantonnent depuis des
décennies\,\,---\,\,d'architecturer `ensemblistement' les mathématiques.
On définit la notion d'application ou de fonction entre deux ensembles
abstraits quelconques, et on relègue en ce lieu préliminaire l'effort
de conception mentale sous forme de quelques exercices élémentaires
sur des ensembles finis ou sur de simples fonctions de $\R$ dans $\R$.
Autrement dit, la structuration par inclusions théoriques prend en
quelque sorte le risque de faire {\em disparaître}\,\,---\,\,au moins
localement au niveau des abstractions\,\,---\,\,les conceptualisations
régionales qui sont néanmoins nécessaires dans les contextes
déterminés. Il n'est donc pas avisé de croire que la conceptualisation
axiomatico-formelle puisse fixer l'extension des idées par
déclarations formelles, puisque les questions que l'on se pose parfois
inopinément pour {\em comprendre et méditer} un concept peuvent
toujours éventuellement {\em ouvrir} vers des domaines nouveaux des
mathématiques, comme par exemple la dynamique des applications entre
{\em ensembles} fractals, absolument invisibles pour l'intuition des
ordinaux et des cardinaux formée au contact de la théorie abstraite
des ensembles. \qed\smallskip

\Scholie
Plusieurs niveaux de pensée s'expriment en géométrie différentielle.
L'ambition de tout géomètre synthétique vise à 
mettre entre parenthèses, si ce n'est éliminer,
toute trace d'exposition symbolique. Mais la puissance d'engendrement
de la pensée analytico-algébrique est telle que, 
pour d'autres fibrés principaux que celui
de groupe ${\rm SO}(2)$, de nombreux invariants
dans la théorie de Cartan du problème d'équivalence entre structures
géométriques s'introduisent {\em sans incarnation géométrique
évidente}.
\qed

\smallskip
Plus loin, Gauss discute du {\em signe} qu'il faut
affecter à la mesure de courbure d'une surface
(redéfinie dans le prochain
paragraphe),
et il montre qu'un tel signe peut être déterminé
en examinant si la disposition relative
d'une paire de segments est préservée, 
ou inversée.

\CITATION{
La position de la figure tracée sur la surface sphérique peut
être ou semblable ou opposée (inverse) à celle
de la figure qui lui correspond sur la surface
courbe; le premier cas a lieu lorsque deux lignes
sur la surface courbe partant du même point et dans
deux directions différentes, mais non opposées,
sont représentées sur la surface sphérique
par deux lignes placées semblablement, c'est-à-dire
lorsque l'image de la ligne située vers la droite
est aussi à droite; le second cas, lorsque c'est
le contraire qui a lieu.
\REFERENCE{{\sc Gauss},~\cite{ga1828},~14--15}}

Mais comme en général la courbure des surfaces est variable, lorsqu'on
veut estimer quantitativement la mesure de courbure contenue dans une
certaine région finie de la surface, un problème imprévu de topologie
de l'espace s'interpose, analogue
mais géométriquement plus complexe, au problème
du signe d'une fonction réelle. 

\CITATION{
Le signe positif ou négatif dont nous
affectons la {\em mesure} de la courbure
d'une figure infiniment petite d'après la position
de cette figure, nous l'étendons aussi à la courbure
intégrale d'une figure finie sur la surface
courbe. Toutefois, si nous voulions embrasser ce sujet
dans toute sa généralité\footnotemark~certains 
éclaircissements seraient nécessaires.
\REFERENCE{{\sc Gauss},~\cite{ga1828},~15}}

\footnotetext{\,
%%%%%%%%%%%%%%%%%%%%%%%-------DEBUT--------%%%%%%%%%%%%%%%%%%%%%%%%%%%
Expression explicite d'une exigence mathématique
universelle et reproductible.
} %%%%%%%%%%%%%%%%%%%%%%%%-----FIN-----%%%%%%%%%%%%%%%%%%%%%%%%%%%%%%%

\Scholie
\`A cause d'une question aussi anodine que le découpage d'une région
bidimensionnelle en sous-régions dans lesquelles une certaine
quantité\,\,---\,\,ici la courbure de Gauss\,\,---\,\,est ou bien
strictement positive, ou bien strictement négative, ou bien nulle, on
entre incidemment dans le domaine de recherche absolument autonome des
ensembles de niveaux de fonctions absolument {\em quelconques}. La
thèse d'{\em ouverture} pour les mathématiques, maintes fois défendue
ici, s'illustre alors comme éclatement de questionnements multiples
coprésents. Une véritable phénoménologie, non pas seulement des
pratiques, mais aussi des actes d'interrogation éprouvés par le sujet
mathématicien, non pas seulement d'un point de vue cognitif, mais
aussi du point de vue des structures de la pensée, serait à inventer,
en tant qu'une telle phénoménologie pourrait être encadrée par
quelques principes universels du questionnement mathématique.
\qed

\smallskip
En tout cas ici, au début de ses {\em Disquisitiones generales circa
superficies curvas}, Gauss dépense une certaine énergie textuelle à
exprimer ses efforts de pensée pour comprendre la manière dont le
positif et le négatif de la courbure s'entrelacent dans la surface
à travers le miroir de la sphère auxiliaire.

\CITATION{
Lorsqu'une figure tracée sur une surface courbe est telle qu'à 
chacun des points qu'elle comprend il
correspond, sur la sphère auxiliaire, des points
{\em différents}, alors nulle
ambiguïté. Mais si cette condition n'est pas remplie, il
sera nécessaire de faire entrer deux ou plusieurs
fois en ligne de compte certaines portions de la
surface sphérique, et de là suivant que la
similitude sera directe ou inverse, des termes qui 
s'ajouteront ensemble ou se détruiront
partiellement. 
\REFERENCE{{\sc Gauss},~\cite{ga1828},~15}}

\CITATION{
Ce qu'il y aura de plus simple, en pareil
cas, sera d'imaginer qu'on ait divisé la figure
tracée sur la surface courbe en parties telles, que
chacune d'elles, considérée isolément, satisfasse à la condition
énoncée tout à l'heure, d'attribuer à chaque partie
la courbure qui lui convient, courbure
dont la grandeur sera donnée par l'aire de la figure qui 
lui correspond sur la surface sphérique
et dont le signe dépendra de la position même de la figure,
et enfin de prendre pour la courbure totale de la figure
entière la quantité qu'on obtiendra en ajoutant ensemble
les courbures intégrales correspondant à chacune des parties
de la figure.
\REFERENCE{{\sc Gauss},~\cite{ga1828},~15--16}}

Au total, la somme de la mesure de courbure comprise dans une certaine
région $V \subset S$ de la surface sera égale à l'{\em intégrale de
surface}: 
\[
\int_V\,\kappa\cdot 
d\sigma, 
\]
étendue à $V$, de la courbure ponctuelle que multiplie l'élément
infinitésimal d'aire $d\sigma$ induit sur $S$ par la métrique
euclidienne de $\R^3$. Ainsi la notion d'intégrale double pour une
fonction dont le signe varie sans restriction
élimine-t-elle\,\,---\,\,grâce à l'emploi sous-jacent de l'infini pour
la sommation\,\,---\,\,toute la difficulté qu'il y avait à penser la
complexité des zones de positivité et de négativité de la fonction.
En géométrie algébrique complexe, la positivité partielle et la
négativité partielle, au sens de Green-Griffiths, des fibrés
vectoriels holomorphes de rang $\geqslant 2$, est relativement mal
comprise quant à la cohomologie et quant aux directions dans l'espace
cotangent. Une stratégie analogue d'élimination provisoire des
difficultés par élévation et submersion consiste à tensoriser par une
très grande puissance d'un fibré en droites ample, et à faire tendre
cette puissance vers l'infini de manière à ce que sa positivité
auxiliaire {\em absorbe} toutes les traces de négativité non comprises
dans le fibré vectoriel initial. Autrement dit et d'un point de vue
plus général, les difficultés, en mathématiques, sont réellement
omniprésentes, et dans des contextes variés, certaines {\sl stratégies
d'évitement} sont déployées afin d'apporter au moins quelques réponses
partielles à des questions qui ouvrent sur des problèmes d'une trop
grande complexité à une époque donnée.

Contrairement à ce qu'on observe chez Riemann, l'ouverture chez
Gauss se manifeste rarement, mais dans ce passage, 
Gauss discute de la manière dont on doit penser
le {\em bord} des régions sur lesquelles on
intègre la courbure\,\,---\,\,bord constitué d'un
nombre fini de courbes lisses
dans les cas simples. 

\CITATION{ 
Le périmère de la figure tracée sur une surface courbe correspondra
toujours, sur la surface sphérique auxiliaire, à une ligne fermée. Que
si cette ligne ne se coupe elle-même en aucun point, elle divisera la
surface sphérique en deux parties; la courbure intégrale de la figure
sera donnée par l'aire de cette partie, cette aire étant positive ou
négative suivant que, par rapport à son périmètre, elle aura une
position semblable ou inverse à celle que la figure a elle-même par
rapport à son propre périmètre. Mais lorsque cette ligne se coupera
elle-même une ou plusieurs fois, elle donnera une figure compliquée, à
laquelle cependant on peut attribuer légitimement une aire déterminée,
comme s'il s'agissait d'une figure sans n{\oe}uds; et cette aire,
convenablement entendue, sera toujours la valeur exacte de la courbure
intégrale.
\REFERENCE{{\sc Gauss},~\cite{ga1828},~15--16}}

En fait, au moment où il finalise son mémoire et où il rédige ce
passage, Gauss pense certainement à un résultat qu'il exposera plus
loin et qui devint l'un de ses plus célèbres théorèmes:

\CITATION{
La somme des angles d'un triangle formé par des lignes
géodésiques sur une surface quelconque, 
est supérieure à $180^{\sf o}$ si cette
surface est concavo-concave, et inférieure à 
$180^{\sf o}$ si cette surface est concavo-convexe, d'une
quantité qui a pour mesure l'aire du triangle
sphérique qui lui
correspond, d'après les directions des normales, 
en comptant la surface totale de la sphère pour
$720^{\sf o}$.
\REFERENCE{{\sc Gauss},~\cite{ga1828},~15--16}}

Autrement dit, sur une surface dont la courbure est ou bien partout
positive, ou bien partout négative, si $\widehat{ A}$, $\widehat{ B}$,
$\widehat{ C}$ sont les mesures des angles d'un triangle curviligne
$\mathcal{ T} = ABC$ découpé sur la surface dont les trois côtés sont
des géodésiques, alors on a:
\[
\widehat{A}+\widehat{B}+\widehat{C}
=
\pi
+
\int\!\!\int_{\mathcal{T}}\,\kappa\,d\sigma.
\]
On peut aussi deviner dans ce passage
qui évoque une ouverture le germe des travaux ultérieurs de 
Bonnet\footnote{\,
%%%%%%%%%%%%%%%%%%%%%%%-------DEBUT--------%%%%%%%%%%%%%%%%%%%%%%%%%%%
En 1848, Ossian Bonnet a généralisé un théorème de Gauss relatif à
l'aire d'un triangle géodésique ({\em cf.} Libermann, \cite{ li1978}).
Pour une courbe $C$ tracée sur une surface $S$, il définit en un point
$p$ de $C$ la {\sl courbure géodésique} $\frac{ 1}{ \rho_g}$ de $C$
comme étant la courbure, au point $p$, de la projection orthogonale de
$C$ sur le plan tangent à $S$ en $p$; pour une géodésique, on a donc
automatiquement $\frac{ 1}{ \rho_g} = 0$ en tout point. Si 
$\mathcal{ U}$ est
une portion de $S$ limitée par une courbe {\em fermée} $C \subset S$,
Bonnet établit la formule:
\[
\int_C\,d\omega
-
\int_C\,\frac{ds}{\rho_g}
=
\int\!\!\int_{\mathcal{U}}\,\kappa\,d\sigma,
\]
où $\kappa$ est la courbure de Gauss, où $d\sigma$ est
l'élément d'aire de $S$ et où 
$\omega$ est l'angle que fait la tangente à $C$ 
avec une direction qui reste fixe par transport parallèle.
} %%%%%%%%%%%%%%%%%%%%%%%%-----FIN-----%%%%%%%%%%%%%%%%%%%%%%%%%%%%%%%
sur la relation entre la topologie globale des surfaces à bord et la
courbure intégrée\footnote{\,
%%%%%%%%%%%%%%%%%%%%%%%-------DEBUT--------%%%%%%%%%%%%%%%%%%%%%%%%%%%
L'importance de la formule `locale' de Bonnet provient du 
fait qu'elle a effectivement été généralisée au vingtième siècle en 
une formule globale
({\em cf.} à nouveau Libermann, \cite{ li1978}). 
Soit $S$ une surface orientable, 
compacte et sand bord. 
En appliquant la formule de
Bonnet à chaque triangle d'une triangulation 
de $S$, on obtient ce qu'on appelle maintenant
le {\sl Théorème de Gauss-Bonnet}:
\[
\int\!\!\int_S\,
\kappa\,d\sigma
=
2\pi\,(n_2-n_1+n_0),
\]
où $n_2$, $n_1$ et $n_0$ sont, respectivement, le
nombre de triangles, d'arêtes et de sommets
de la triangulation. Dans cette relation d'égalité, le premier
membre étant visiblement indépendant de la
triangulation, il en découle que le second l'est aussi.
En fait, le nombre entier $\chi ( S) := 
n_2 - n_1 + n_0$, appelé aujourd'hui {\sl caractéristique
d'Euler-Poincaré} de $S$, ne dépend que de la
topologie de $S$, et non de sa métrique
infinitésimale. Une telle relation fondamentale entre
la géométrie différentielle et la topologie
algébrique a reçu ultérieurement de très nombreuses extensions
et généralisations.
}. %%%%%%%%%%%%%%%%%%%%%%%%-----FIN-----%%%%%%%%%%%%%%%%%%%%%%%%%%%%%%%

Toutefois, même s'il mentionne des ouvertures
coprésentes au champ interrogatif de la théorie des surfaces, Gauss
doit choisir de s'écarter de ces questions qui se rapportent en fait à
d'autres branches futures de la topologie et de la géométrie en
général.

\CITATION{
Au surplus, nous croyons devoir réserver pour une autre
occasion des explications plus amplement développées, concernant les
figures envisagées au point de vue le plus général.
\REFERENCE{{\sc Gauss},~\cite{ga1828},~16}}

Exigence, donc, proprement gaussienne, de penser les domaines
découpés sur les surfaces {\em de la manière la plus générale
possible}, la question reste d'actualité en tant que
l'approfondissement de la notion de lieu concerne
l'histoire de la Topologie.

\medskip\noindent{\bf Définition géométrique de la courbure de Gauss.}
\label{definition-geometrique-courbure}
Si donc ${\bf n} (q)$ est le vecteur unitaire normal à la surface $S$
en l'un de ses points quelconque $q$, Gauss le {\em transporte}
systématiquement à l'origine du système des coordonnées: telle est la
fameuse {\sl application de Gauss}, qui décalque sur une
sphère auxiliaire pure $\Sigma$ la
variation normale au second ordre de toute surface $S$.

Pour une région quelconque $V_\varepsilon ( p)$ délimitée autour du
point $p$ dans la surface qui se rétrécit autour de $p$ lorsque
$\varepsilon >0$ tend vers $0$, l'ensemble des extrémités des vecteurs
${\bf n} (q)$ pour tous les points $q \in V_\varepsilon (p)$ décrit
alors une certaine région ${\bf n} \big( V_\varepsilon (p)\big)$ sur
la sphère auxiliaire $\Sigma$. Avec ces deux régions $V_\varepsilon
(p)$ et ${\bf n} \big( V_\varepsilon (p) \big)$, la {\sl courbure de
Gauss} de la surface au point $p$ est définie comme étant la limite,
quant $\varepsilon$ tend vers $0$, du quotient de l'aire
infinitésimale image par l'aire infinitésimale initiale:
\[
\kappa(p)
:=
\lim_{V_\varepsilon(p)\,\to\,p}\,
\frac{\text{\sf aire}\,\big[{\bf n}\big(V_\varepsilon(p)\big)\big]}{
\text{\sf aire}\,V_\varepsilon(p)},
\]
de manière exactement analogue à la définition connue pour les
courbes. On peut prendre par exemple comme région $V_\varepsilon (p)$
simplement l'intersection $V_\varepsilon (p) := \B_\varepsilon (p)
\cap S$ de la surface avec la boule dans $\R^3$ de rayon $\varepsilon$
et de centre $p$:
\[
\B_\varepsilon(p)
:=
\big\{
(x,y,z)\in\R^3\colon\,\,
x^2+y^2+z^2<\varepsilon^2\big\}.
\]
En vérité, la limite\,\,---\,\,quand elle existe, par exemple lorsque
la surface est de classe au moins $\mathcal{ C}^2$\,\,---, ne dépend
pas de la forme de la région $V_\varepsilon (p)$ qui rétrécit autour
de $p$, fait qui peut et qui doit être démontré. Bien entendu, Gauss
raisonne directement en termes d'infinitésimaux, ayant à l'esprit que
les différences de forme entre éléments infimes de surface
s'évanouissent dans l'infiniment petit grâce à la présence du quotient.
De plus incidemment, Gauss sait qu'une formule analytique quantitative
va réaliser sous peu cette définition géométrique, à savoir d'un point
de vue plus métaphysique: le premier moment géométrique de la
définition sait par avance qu'il va être dépassé dialectiquement par
un calcul d'explicitation qui rendra immédiatement visible
une indépendance vis-à-vis de la forme.

\Scholie
Bifurcations possibles, ici, pour la pensée structurante: ou bien on
choisit de spécifier une forme simple de voisinages $V_\varepsilon$ en
des termes qui peuvent être définis géométriquement, et on démontre
ensuite éventuellement aussi qu'une et une seule quantité de courbure
serait obtenue avec d'autres régions de forme quelconque, ou bien
on choisit d'emblée la généralité dans la définition, mais on doit
alors établir rapidement que le résultat existe indépendamment de la
forme de la région. En tout état de cause, la limite s'avère plus
complexe qu'en dimension $1$, où le $dx$ horizontal, simple segment
infinitésimal, a un bord constitué de deux points seulement. Donc
d'un point de vue général, la pensée structurante est-elle en
permanence confrontée à des choix de présentation qui ne sont pas
métaphysiquement équivalents, 
bien que techniquement
interchangeables; c'est là l'un des aspects les
plus labyrinthiques
des mathématiques qui contribue
régulièrement à parsemer d'embûches supplémentaires
l'appropriation intuitive
et la compréhension des théories.
\qed

\Section{8.~Deux expressions analytiques extrins\`eques 
\\
de la courbure des surfaces}
\label{deux-expressions-extrinseques}

\HEAD{8.~Deux expressions analytiques extrins\`eques de 
la courbure des surfaces}{
Jo\"el Merker, D\'partement de Math\'matiques d'Orsay}

\medskip\noindent{\bf Courbure des surfaces graphées.} 
\label{courbure-surfaces-graphees}
On cherche
maintenant\,\,---\,\,poursuit Gauss\,\,---\,\,une
formule propre à exprimer la mesure de
la courbure en chaque point d'une surface courbe. Il s'agit de
calculer, en termes des éléments analytiques différentiels de la
surface, un quotient de deux aires infinitésimales, et donc
évidemment, il est
nécessaire de commencer par calculer effectivement ces aires. Parmi
les deux représentations possibles: graphée ou implicite, d'une
surface dans l'espace, c'est la représentation graphée qui fournit la
formule renfermant le nombre minimal d'éléments, 
et donc, c'est par cette
formule que l'on commencera.

L'élément de surface
$V_\varepsilon (p)$ choisi par Gauss sur la surface est un
simple triangle infinitésimal non dégénéré de sommets $p$, $p + dp$ et
$p + \delta p$ sur la surface, où $dp$ et $\delta p$ sont deux
vecteurs infinitésimaux tangents en $p$ qui ne sont pas colinéaires.
En projection horizontale sur le plan des $(x, y)$, 
ces trois points forment un triangle
non dégénéré dont les
coordonnées seront notées:
\[
p^\pi
:=
(x,\,y),
\ \ \ \ \
p^\pi+dp^\pi
:=
(x+dx,\,y+dy),
\ \ \ \ \ 
p^\pi+\delta p^\pi
:=
(x+\delta x,\,y+\delta y).
\]
Comme le montrent des raisonnements
pythagoriciens élémentaires,
l'aire (infinitésimale) du triangle infinitésimal
de sommets $p$, $p + dp$ et
$p + \delta p$ est dans une relation de proportionalité
simple avec l'aire (infinitésimale) 
du triangle projeté de sommets
$p^\pi$, $p^\pi + d p^\pi$ et
$p^\pi + \delta p^\pi$, 
le coefficient de proportionalité étant
simplement égal à la troisième composante
$Z$ du vecteur normal.

\CITATION{
En appelant $d\sigma$ l'aire d'un élément
de cette surface, $Z\, d\sigma$ sera l'aire de la projection
de cet élément sur le plan des coordonnées $x$ et $y$;
et de même, si $d\Sigma$ est l'aire de l'élément correspondant
sur la surface sphérique auxiliaire, $Z\, d\Sigma$ sera l'aire
de la projection de cet élément sphérique sur le même plan; 
et il est manifeste que ces projections auront entre elles
les mêmes relations de grandeur et
de position que les éléments eux-mêmes.
\REFERENCE{{\sc Gauss},~\cite{ga1828},~16--17}}

Autrement dit, pour la surface infinitésimale initiale et pour sa
surface image sur la sphère auxiliaire, {\em les deux coefficients de
proportionalité entre l'aire et son aire projetée:}
\[
\frac{d\Sigma^\pi}{d\Sigma}
=
\frac{d\sigma^\pi}{d\sigma}
=
Z
\]
{\em sont égaux entre eux et égaux à la troisième composante $Z$ de
leur vecteur normal commun}. Cette circonstance est très favorable,
car pour calculer la courbure définie
par Gauss comme quotient d'aires infinitésimales:
\[
\kappa
=
\frac{d\Sigma}{d\sigma}
=
\frac{d\Sigma^\pi}{d\sigma^\pi},
\] 
on peut se ramener, grâce à une identité algébrico-géométrique
dont la connaissance intuitive ou pratique
remonte vraisemblablement jusqu'aux mathématiques de la 
Préhistoire\footnote{\,
%%%%%%%%%%%%%%%%%%%%%%%-------DEBUT--------%%%%%%%%%%%%%%%%%%%%%%%%%%%
Telle est en fin de compte implicitement ici, l'une des thèses
défendue dans ce mémoire: chaque principe élémentaire des
mathématiques renvoie à une métaphysique propre que l'on peut
réveiller à tout instant dans la pensée, comme s'il existait une
phénoménologie de la mobilisation des aspects historico-philosophiques
qui enrichirait constamment les raisonnements techniques. En effet,
l'échange de numérateur et de dénominateur à travers l'égalité $\frac{
a}{ b} = \frac{ c}{ d}$ est un acte archaïque, très ancien dans
l'histoire humaine des mathématiques, présent peut-être sans traces
écrites et sans formalisation dans une tradition orale disparue, et ce
mouvement atomique de `calcul' n'en possède pas moins une importance
cruciale pour accéder à la formule explicite $\kappa = \frac{ z_{ xx}
\, z_{ yy} - z_{ xy}^2}{ (1 + z_x^2 + z_y^2)^2}$ qui exprime la courbure
d'une surface graphée $z = z ( x, y)$, {\em cf.} ce qui
va suivre. Autrement dit, l'élémentarité la plus
antique rayonne dans les mathématiques les plus sophistiquées.
}: %%%%%%%%%%%%%%%%%%%%%%%%-----FIN-----%%%%%%%%%%%%%%%%%%%%%%%%%%%%%%%
\[
{\textstyle{\frac{a}{b}}}
=
{\textstyle{\frac{c}{d}}}
\Longrightarrow
{\textstyle{\frac{a}{c}}}
=
{\textstyle{\frac{b}{d}}},
\]
à calculer seulement les quotients des aires infinitésimales
{\em projetées:}
\[
\kappa
=
\frac{d\Sigma^\pi}{d\sigma^\pi}.
\]

\`A un signe près qui dépend seulement d'un choix d'orientation,
la formule qui
exprime l'aire, dans le plan des $(x, y)$, de tout triangle
formé par trois points est bien connue elle aussi.

\CITATION{
Considérons maintenant un élément triangulaire
de la surface courbe, et supposons que les coordonnées
des trois points qui forment la projection de cet
élément sont:
\[
\begin{array}{cc}
x, & y,
\\
x+dx, & y+dy,
\\
z+\delta z, & z+\delta z;
\end{array}
\]
le double de l'aire de ce triangle sera alors exprimé
par la formule:
\[
dx\cdot\delta y
-
dy\cdot\delta x,
\]
expression positive ou négative, suivant que la position
du côté qui se dirige du premier point au
troisième, comparée à celle du côté qui
se dirige du premier au second, est semblable
ou inverse à la position
de l'axe coordonné $y$ par rapport à l'axe coordonné
$x$.
\REFERENCE{{\sc Gauss},~\cite{ga1828},~17}}

La {\em même} formule s'appliquera alors aussi immédiatement
pour calculer le double de l'aire 
du triangle formé par les trois points projetés
sur l'espace des $(X, Y)$
du triangle infinitésimal correspondant:
\[
\begin{array}{cc}
X, & Y,
\\
X+dX, & Y+dY,
\\
Z+\delta Z, & Z+\delta Z,
\end{array}
\] 
tracé sur
la sphère auxilaire, à savoir:
\[
dX\cdot\delta Y
-
dY\cdot\delta X.
\]

\CITATION{
La mesure de la courbure sera donc, en ce point de la surface
courbe:
\[
\kappa
=
\frac{dX\cdot\delta Y-dY\cdot\delta X}{dx\cdot\delta y-dy\cdot\delta x}.
\]
\REFERENCE{{\sc Gauss},~\cite{ga1828},~17}}

\medskip\noindent{\bf Prémices d'élimination différentielle 
systématique.} \label{premices-elimination-differentielle}
Cette dernière formule est encore loin d'expliciter la
courbure en termes de la fonction graphante $z ( x, y)$, puisqu'elle
incorpore encore l'élément d'arbitraire que constitue l'aire d'un
triangle infinitésimal quelconque dirigé par les deux vecteurs
$dp^\pi$ et $\delta p^\pi$. Si donc on admet que la courbure ne dépend
pas de la région infinitésimale choisie
sur la surface, on peut prévoir que le numérateur
$dX \cdot \delta Y - dY \cdot \delta X$
devrait s'avérer divisible par le dénominateur $dx \cdot \delta y - dy
\cdot \delta x$, le facteur de proportionalité 
{\em rémanent} se réduisant à
être une {\em fonction} déterminée des deux variables $x$ et $y$,
celle que l'on aimerait connaître. Cet argument anticipateur montre
donc qu'un principe métaphysique de simplification et d'annihilation
devrait s'exercer ici. La dynamique du calcul est donc dominée
et dirigée
par les principes de la pensée spéculative.

En effet, les deux différentielles $d$ et $\delta$ appliquées aux deux
composantes horizontales $X$ et $Y$ du vecteur normal ${\bf n}$
donnent, si on les écrit systématiquement toutes les quatre:
\[
\aligned
dX
&
=
\bigg(
\frac{\partial X}{\partial x}
\bigg)\,dx
+
\bigg(
\frac{\partial X}{\partial y}
\bigg)\,dy,
\\
dY
&
=
\bigg(
\frac{\partial Y}{\partial x}
\bigg)\,dx
+
\bigg(
\frac{\partial Y}{\partial y}
\bigg)\,dy,
\\
\delta X
&
=
\bigg(
\frac{\partial X}{\partial x}
\bigg)\,\delta x
+
\bigg(
\frac{\partial X}{\partial y}
\bigg)\,\delta y,
\\
\delta Y
&
=
\bigg(
\frac{\partial Y}{\partial x}
\bigg)\,\delta x
+
\bigg(
\frac{\partial Y}{\partial y}
\bigg)\,\delta y.
\\
\endaligned
\]

\Scholie
Dans toutes ses approches de la réalité mathématique du calcul, Gauss
procède en effet d'une manière systématique et par {\em complétude
exploratoire}. On peut dire que le principe d'une systématicité
décidé {\em a priori} est une saine {\em méthode} de recherche
en mathématique, si ce n'est que la complétude des
parcours demande toujours un travail considérable.
\qed

\smallskip
Quand on calcule alors le numérateur de la courbure:
\[
\aligned
dX\cdot\delta Y
-
dY\cdot\delta X
&
=
\big(
X_x\,dx+X_y\,dy
\big)\cdot
\big(
Y_x\,\delta x+Y_y\,\delta y
\big)
-
\\
&
\ \ \ \ \
-
\big(
Y_x\,dx+Y_y\,dy
\big)\cdot
\big(
X_x\,\delta x+X_y\,\delta y
\big),
\endaligned
\]
deux produits de deux termes sont soustraits l'un à l'autre;
rien que de très élémentaire pour tout calculateur expérimenté, 
mais il est toutefois nécessaire d'interroger les symétries
formelles du résultat afin de savoir si l'intuition 
d'anticipation va réellement se confirmer, 
{\em i.e.} afin de savoir si ce numérateur
est effectivement multiple de l'aire infinitésimale
$dx \cdot \delta y - dy \cdot \delta x$
lorsqu'on développe lesdits produits: 
\[
\footnotesize
\aligned
dX\cdot\delta Y
-
dY\cdot\delta X
&
=
\zero{X_x\,Y_x\,dx\delta x}
+
X_x\,Y_y\,dx\delta y
+
X_y\,Y_x\,dy\delta x
+
\zerozero{X_y\,Y_y\,dy\delta y}
-
\\
&
\ \ \ \ \
-
\zero{Y_x\,X_x\,dx\delta x}
-
Y_x\,X_y\,dx\delta y
-
Y_y\,X_x\,dy\delta x
-
\zerozero{Y_y\,X_y\,dy\delta y}
\\
&
=
\big(
X_x\,Y_y-X_y\,Y_x
\big)\cdot
\big(
dx\,\delta y-dy\,\delta x
\big).
\endaligned
\]
Ainsi comme on aura
su l'anticiper, l'expression de la courbure se 
change-t-elle\,\,---\,\,{\em irréversiblement pour ce qui concerne
l'acquisition corrélative de connaissance 
mathématique}\,\,---\,\,en une expression maintenant
purement fonctionnelle et parfaitement déterminée:
\[
\kappa
=
X_x\cdot Y_y
-
X_y\cdot Y_x,
\]
dans laquelle a disparu tout l'arbitraire 
du triangle infinitésimal. Cette fois-ci en effet, on
ne voit plus que des dérivées partielles d'ordre $1$ des deux
composantes $X$ et $Y$ du vecteur normal {\bf n}. Cette expression de
la courbure est donc potentiellement définitive. Toutefois, 
en mathématiques, la
question se pose en permanence de savoir si un
résultat obtenu est complet et définitif. 

\smallskip
\centerline{\fbox{\em Nihil actum reputans si quid superesset agendum},}

\smallskip\noindent
telle était la devise {\em métaphysique} de Gauss, 
devise qui exprime notamment l'exigence permanente de s'{\em interroger 
en pensée} sur l'achèvement des travaux de calcul. 

Or puisque $X$ et $Y$
s'exprimaient en fait en termes de la fonction graphante $z ( x, y)$,
il est clair que le calcul n'est pas achevé, 
et par conséquent, il reste encore à insérer
ces formules et aussi, à simplifier, à symétriser et à harmoniser les
formules qu'on sera susceptible d'obtenir. Mais Gauss ne procède pas
de manière directe pour effectuer ce calcul. Pour faire face à une
multiplicité de dérivées partielles d'ordre $2$, il introduit comme Euler
et Lagrange le faisaient régulièrement des {\em lettres} avec
lesquelles il va raisonner algébriquement et organiser les calculs de
manière synoptique: il faut y voir les prémices du calcul
d'élimination beaucoup plus considérable qu'il présentera
dans la suite de son mémoire afin d'établir que la courbure s'exprime
en fonction de la métrique infinitésimale
{\em intrinsèque}. 

\CITATION{
En posant comme ci-dessus:
\[
\frac{dz}{dx}
=
t,
\ \ \ \ \ \ \
\frac{dz}{dy}
=
u,
\]
et en outre:
\[
\frac{d^2z}{dx^2}
=
T,
\ \ \ \ \
\frac{d^2z}{dxdy}
=
U,
\ \ \ \ \
\frac{d^2z}{dy^2}
=
V,
\]
ce qui équivaut à:
\[
dt
=
T\,dx+U\,dy,
\ \ \ \ \ \ \ \ \
du
=
U\,dx+V\,dy,
\]
nous aurons, d'après des formules données précédemment:
\[
X
=
-t\,Z,
\ \ \ \ \
Y
=
-u\,Z,
\ \ \ \ \
(1+t^2+u^2)\,Z^2
=
1,
\]
et par suite:}
\CITATION{
\[
\aligned
&
dX
=
-Z\,dt-t\,dZ,
\\
&
dY
=
-Z\,du-u\,dZ,
\\
&
(1+t^2+u^2)\,dZ
+
Z\,(t\,dt+u\,du)
=
0,
\endaligned
\]
ou bien:
\[
\aligned
dZ
&
=
-Z^3\,(t\,dt+u\,du),
\\
dX
&
=
-Z^3\,(1+u^2)\,dt
+
Z^3\,tu\,du,
\\
dY
&
=
Z^3\,tu\,dt
-
Z^3\,(1+t^2)\,du,
\endaligned
\]
et de là on tire:}
\CITATION{
\[
\aligned
\frac{dX}{dx}
&
=
Z^3\big[-(1+u^2)\,T+tu\,U\big],
\\
\frac{dX}{dy}
&
=
Z^3\big[-(1+u^2)\,U+tu\,V\big],
\\
\frac{dY}{dx}
&
=
Z^3\,\big[tu\,T-(1+t^2)\,U\big],
\\
\frac{dY}{dy}
&
=
Z^3\,\big[tu\,U-(1+t^2)\,V\big].
\endaligned
\]
En substituant ces valeurs dans l'expression précédente, 
il vient:
\[
\kappa
=
Z^6(TV-U^2)(1+t^2+u^2)
=
Z^4(TV-U^2)
=
\frac{TV-U^2}{(1+t^2+u^2)}.
\]
\REFERENCE{{\sc Gauss},~\cite{ga1828},~18}}

\Scholie
Fluidité enchaînée d'une longue phrase de calculs en synthèse: le
style de présentation frappe de netteté quant aux
opérations sous-entendues
qui sont laissées à la compréhension du lecteur. Or
excepté pour le résultat final,
Gauss {\em élimine ici intentionnellement toute 
fraction rationnelle et toute trace de division}. En effet, 
au lieu de calculer
directement les dérivées partielles par rapport à $x$ et à $y$ des
deux fractions:
\[
\pm X
=
\frac{z_x}{(1+z_x^2+z_y^2)^{1/2}}
\ \ \ \ \ \ \ \
\text{\rm et}
\ \ \ \ \ \ \ \
\pm Y
=
\frac{z_y}{(1+z_x^2+z_y^2)^{1/2}},
\] 
pour en soustraire ensuite les paires croisées multipliées entre
elles\,\,---\,\,ce qui conduirait à des expressions fort 
développées\footnote{\,
%%%%%%%%%%%%%%%%%%%%%%%-------DEBUT--------%%%%%%%%%%%%%%%%%%%%%%%%%%%
Un calcul direct donnerait en effet:
\[
\aligned
\pm X_x
=
\frac{z_{xx}+z_{xx}\,z_y^2-z_x\,z_y\,z_{xy}}{
\Big(\sqrt{1+z_x^2+z_y^2}\Big)^3},
\ \ \ \ \ \ \ \ \ \ 
\pm X_y
=
\frac{z_{xy}+z_{xy}\,z_y^2-z_x\,z_y\,z_{yy}}{
\Big(\sqrt{1+z_x^2+z_y^2}\Big)^3},
\\
\pm Y_x
=
\frac{z_{xy}+z_{xy}\,z_x^2-z_x\,z_y\,z_{xx}}{
\Big(\sqrt{1+z_x^2+z_y^2}\Big)^3},
\ \ \ \ \ \ \ \ \ \ 
\pm Y_y
=
\frac{z_{yy}+z_{yy}\,z_x^2-z_x\,z_y\,z_{xy}}{
\Big(\sqrt{1+z_x^2+z_y^2}\Big)^3},
\endaligned
\]
quantités qu'il faudrait ensuite multiplier par paires
puis soustraire pour simplifier finalement $X_x\, Y_y - X_y\, Y_x$.
Il se trouve que l'on devrait réaliser que
le numérateur est divisible par $1 + z_x^2 + z_y^2$, 
ce que Gauss réussit à faire voir de manière directe.
}%%%%%%%%%%%%%%%%%%%%%%%%-----FIN-----%%%%%%%%%%%%%%%%%%%%%%%%%%%%%%%
et quelque peu inélégantes en raison de la présence d'une racine
carrée au dénominateur\,\,---, Gauss travaille en éliminant ledit
dénominateur grâce à l'introduction d'une indéterminée auxiliaire $Z$
satisfaisant:
\[
(1+t^2+u^2)\,Z^2
=
1,
\]
et qui apparaîtra au final élevée
à une puissance {\em paire}, de telle
sorte que toute racine carrée présente
implicitement dans les calculs intermédiaires disparaîtra. 
Non seulement le calcul est
plus économique ainsi, mais il s'articule aussi d'une 
manière exclusivement polynomiale, donc plus satisfaisante
pour l'esprit du calculateur.
On notera aussi au passage que l'ambiguïté du signe
commun $\pm$
de $X$ et de $Y$ causée initialement
par le choix d'une orientation de la normale
d'un côté ou de l'autre de la surface s'évanouit dans 
la formule de départ
$\kappa = X_x \, Y_y - X_y \, Y_x$, puisque
l'on a $\pm \cdot \pm = +$ où 
la multiplication s'effectue ligne à ligne.
\qed

\medskip\noindent{\bf Courbure des surfaces en représentation
implicite.} \label{courbure-en-representation-implicite}
Dans un second moment, Gauss cherche à déduire de cette
première formule une deuxième formule pour la courbure, 
valable lorsque la
surface est représentée par une équation implicite de la forme:
\[
0
=
W(x,y,z),
\]
pour une certaine fonction $W = W ( x, y, z)$. 

\CITATION{
La formule donnée précédemment pour la mesure de courbure
est la plus simple de toutes les formules générales,
en ce qu'elle ne renferme que cinq éléments; 
nous arriverons à une formule plus compliquée, renfermant
neuf éléments, si nous voulons employer la première
des méthodes que nous avons dit être propre à étudier
les caractères des surfaces.
\REFERENCE{{\sc Gauss},~\cite{ga1828},~21}}

Ainsi, dans les notations de Gauss qui consistent à introduire
des lettres spécifiques et différenciées pour désigner
les dérivées partielles des fonctions et 
à travailler {\em algébriquement} avec
ces lettres, la première différentielle 
$dW =
W_x\, dx + W_y\, dy + W_z\, dz$ de la fonction définissante $W$ 
sera écrite:
\[
dW
=
P\,dx+Q\,dy+R\,dz.
\]
En poursuivant ces notations, il est alors judicieux
d'introduire la collection de toutes les dérivées
partielles du second ordre de la fonction $W$, puisqu'il
semble clair, du point
de vue de la pensée anticipante, que la formule de la courbure en termes
de $W$ va les faire intervenir:
\[
\aligned
\frac{\partial^2W}{\partial x^2}
&
=
P',
\ \ \ \ \ \ \ \
\frac{\partial^2W}{\partial y^2}
=
Q',
\ \ \ \ \ \ \ \
\frac{\partial^2W}{\partial z^2}
=
R',
\\
\frac{\partial^2W}{\partial y\partial z}
&
=
P'',
\ \ \ \ \ \ \ \
\frac{\partial^2W}{\partial x\partial z}
=
Q'',
\ \ \ \ \ \ \ \
\frac{\partial^2W}{\partial x\partial y}
=
R''.
\endaligned
\]
\'Evidemment, ces dérivées partielles interviennent comme coefficients
des trois $1$-formes $dP$, $dQ$, $dR$, à savoir on a:
\[
\aligned
dP
&
=
P'\,dx+R''\,dy+Q''\,dz,
\\
dQ
&
=
R''\,dx+Q'\,dy+P''\,dz,
\\
dR
&
=
Q''\,dx+P''\,dy+R'\,dz,
\endaligned
\]
collection tabulaire de neuf termes dont les éléments diagonaux
n'incorporent qu'un seul `prime', les éléments extra-diagonaux en
comportant deux et étant égaux un à un symétriquement par rapport à la
diagonale principale. 

\Scholie
Ainsi s'explique le fait que la lettre $P''$ a été choisie
intentionnellement pour correspondre à la dérivée partielle seconde de
$W$ par rapport à $y$ et $z$, et non à celle par rapport à $x$ et $y$,
comme on aurait pu le croire si l'ordre lexicographique commun choisi
pour les fonctions $P'$, $Q'$, $R'$ et pour les variables $x$, $y$,
$z$ avait été à nouveau appliqué pour les fonctions $P''$, $Q''$,
$R''$: chez un calculateur prodige, la micro-organisation algébrique
fourmille de choix subtils non notifiés au lecteur par des
explications langagières; car le jeu avec des calculs délicats exige
une pensée autonome riche de principes non écrits.

En outre, tout calcul exprime à chaque étape un dynamisme de mobilité
potentielle qui doit chercher à se {\em manifester} et à se {\em
déclencher} à travers une organisation synoptique, structurée et
harmonieuse de ses éléments atomiques. Les choix de notation chez
Gauss sont donc gouvernés par certaines micro-structures formelles
qu'il aura par exemple découvertes lors d'explorations manuscrites.
\qed

\smallskip
L'objectif est maintenant de {\em traduire} en termes
de la fonction $W$ la formule
pour la courbure obtenue
dans le cas d'un graphe, 
à savoir la formule: $\kappa =
\frac{ z_{ xx}\, z_{ yy} - z_{ xy}^2}{ (1 + z_x^2 + z_y^2)^2}$. 
De la relation différentielle déduite de
l'équation $0 = W$ par différentiation:
\[
0
=
P\,dx+Q\,dy+R\,dz,
\]
on tire par identification que si la surface est imaginée comme
représentée aussi 
par une équation graphée de la forme: $z = z(x,y)$
de telle sorte que l'équation
suivante: 
\[
0
\equiv
W\big(x,y,z(x,y)\big)
\]
soit satisfaite identiquement par rapport à $x$ et $y$, alors 
en dérivant cette identité par rapport
à $x$ et à $y$, on voit que les deux
dérivées partielles du second ordre de la fonction $z$ sont données
par:
\[
z_x
=
t
=
-\,{\textstyle{\frac{P}{R}}}
=
-\,{\textstyle{\frac{W_x}{W_z}}}
\ \ \ \ \ \
\text{\rm et}
\ \ \ \ \ \
z_y
=
u
=
-\,{\textstyle{\frac{Q}{R}}}
=
-\,{\textstyle{\frac{W_y}{W_z}}},
\]
lorsque l'argument $(x, y, z)$ appartient à la surface.
 
\CITATION{
Maintenant, puisqu'on a $t = - \frac{ P}{ R}$, nous obtenons, 
par la différentiation:
\[
\aligned
R^2\,dt
&
=
-R\,dP+P\,dR
\\
&
=
(PQ''-RP')\,dx
+
(PP''-RR'')\,dy
+(PR'-RQ'')\,dz,
\endaligned
\]
ou bien, en éliminant $dz$ à l'aide de l'équation
$P\, dx + Q\, dy + R\, dz = 0$:
\[
\aligned
R^3\,dt
&
=
(-R^2\,P'+2\,PRQ''-P^2\,R')\,dx
+
\\
&
\ \ \ \ \
(PRP''+QRQ''-PQR'-R^2R')\,dy.
\endaligned
\]
\REFERENCE{{\sc Gauss},~\cite{ga1828},~22}}

En effet, sur la surface graphée, seules $dx$ et $dy$
sont des différentielles internes, la troisième
$dz$ s'exprimant nécessairement en fonction d'elles;
aussi cette étape de calcul est-elle {\em nécessaire}, 
nullement entachée d'indécision. 

\Scholie
Micro-organisation gaussienne non explicitée: l'ordre alphabétique
entre les trois lettres $P$, $Q$, $R$ passe en seconde position par
rapport à l'ordre concernant le nombre de `primes'.
\qed 

\CITATION{
On a de même:
\[
\aligned
R^3\,du
&
=
(PRP''+QRQ''-PQR'-R^2R'')\,dx
+
\\
&
\ \ \ \ \
+
(-R^2Q'+2\,QRP''-Q^2R')\,dy.
\endaligned
\]
Et de là nous concluons:
\[
\aligned
R^3T
&
=
-R^2P'+2\,PRQ''-P^2R',
\\
R^3U
&
=
PRP''+QRQ''-PQR'-R^2R'',
\\
R^3V
&
=
-R^2Q'+2\,QRP''-Q^2R.
\endaligned
\]
\REFERENCE{{\sc Gauss},~\cite{ga1828},~22}}

\Scholie
Ici, l'expression <<\,on a de même\,>> sous-entend qu'il y a des
calculs entièrement analogues lorsqu'on considère $d z_y$ au lieu de
$dz_x$. Du point de vue du mathématicien professionnel, accepter la
non explicitation d'un ou de plusieurs calculs fait partie des règles
du jeu communément admises, surtout lorsque le calcul sous-entendu
présente des similitudes directes et évidentes avec un calcul que
l'auteur a déjà explicité en détail. Cette <<\,règle du jeu\,>>
appelle de nombreux commentaires.

Dans son mémoire, Gauss ne cache essentiellement aucun calcul,
contrairement à de nombreuses autres sources modernes des
mathématiques, dans lesquelles il est de bon ton de faire admettre à
son lecteur que certains calculs intermédiaires sont sans intérêt et
en quelque sorte indignes d'être explicités ou plus encore publiés,
notamment après l'avènement du formalisme de la géométrie
différentielle globale depuis les années 1950--60. Ce fait pose une
véritable question philosophique quant au {\em statut} littéral
accordé au calcul dans les mathématiques, puisque 
qu'il prend le contrepied des 
disciplines littéraires pour lesquelles {\em tout le parcours de la
pensée} doit s'écrire dans la continuité langagière d'une {\em
complétude}. Les actes mathématiques seraient-ils trop
complexes et trop nombreux pour que les pratiques de rédaction aient
convergé, depuis des décennies, vers une abréviation opaque?

D'un côté, la ré-effectuation complète et le re-parcours total des
calculs est absolument nécessaire afin de {\em vérifier} la véracité
d'un résultat et de {\em contrôler} qu'aucune erreur ne se
cache. Souvent, de nombreux calculs intermédiaires non explicités sont
soumis à des tensions formelles qui forcent les résultats finaux à
être <<\,habillés\,>> d'harmonies internes qui sont {\em suffisantes}
comme outils de contrôle de véracité, et dans de tels cas, le
lecteur sait qu'il peut se dispenser d'avoir à ré-effectuer tout le
parcours de l'auteur.

Mais dans d'autres cas, certains arguments cruciaux qui ne sont
visibles ni {\em a priori} ni {\em a posteriori}, se jouent dans la
conduite des calculs et se placent en des endroits spécifiques
qu'il est difficile de localiser au moment de la recherche et des
premières découvertes\,\,---\,\,c'est-à-dire {\em avant} qu'un
irréversible de connaissance synthétique partagé par une communauté
mathématique ne scelle définitivement l'acceptation d'un
résultat\,\,---, et alors dans de telles circonstances,
tous les calculs doivent être rigoureusement re-parcourables.

Ces paradoxes et ces difficultés soulèvent en définitive la question
de savoir quel {\em statut} on devrait philosophiquement accorder au
calcul en tant que sa non-effectuation partielle demeure
conventionnelle,
ce qui montre aussi qu'il y a une {\em localité spatio-temporelle}
de nombreux théorèmes démontrés par les mathématiciens
contemporains.
\qed

\CITATION{
En substituant ces valeurs dans la formule du $\S VII$\footnotemark, 
nous obtenons pour la mesure de la courbure $\kappa$ l'expression
symétrique suivante:
\[
\aligned
(P^2+Q^2+R^2)\,\kappa
&
=
P^2\,(Q'R'-{P''}^2)
+
Q^2(P'R'-{Q''}^2)
+
\\
&
\ \ \ \ \
+R^2(P'Q'-{R''}^2)
+
2\,QR\,(Q''R''-P'P'')
+
\\
&
\ \ \ \ \
+
2\,PR\,(P''R''-Q'Q'')
+
2\,PQ\,(P''Q''-R'R'').
\endaligned
\]
\REFERENCE{{\sc Gauss},~\cite{ga1828},~22}}

\footnotetext{\,
%%%%%%%%%%%%%%%%%%%%%%%-------DEBUT--------%%%%%%%%%%%%%%%%%%%%%%%%%%%
---\,\,à savoir dans la formule:
\[
\kappa
=
\frac{TV-U^2}{(1+t^2+u^2)^2}
=
\frac{z_{xx}\,z_{yy}-z_{xy}^2}{(1+z_x^2+z_y^2)^2}\text{\rm \,\,---}
\]
} %%%%%%%%%%%%%%%%%%%%%%%%-----FIN-----%%%%%%%%%%%%%%%%%%%%%%%%%%%%%%%

\Section{9.~Gen\`ese des {\em Disquisitiones generales
circa superficies curvas}}
\label{genese-disquisitiones}

\HEAD{9.~Gen\`ese des {\em Disquisitiones generales
circa superficies curvas}}{
Jo\"el Merker, D\'partement de Math\'matiques d'Orsay}

\smallskip\noindent{\bf Cinq concepts novateurs.}
\label{cinq-concepts-novateurs}
En ce qui concerne tant la maîtrise des calculs que la
richesse des idées conceptuelles,
Stäckel suggère dans~\cite{
stackel1890} que le mémoire des {\em Disquisitiones Generales circa
superficies curvas} long seulement d'une cinquantaine de pages est
comparable aux 470 pages des {\em Disquisitiones arithmetic{\ae}}.
En effet, les recherches de Gauss en géométrie l'ont conduit à concevoir,
à introduire et à ouvrir un champ entièrement nouveau, à savoir
l'étude intrinsèque des sous-variétés de l'espace, et outre une
dizaine de théorèmes profonds et mûrs, on trouve dans son mémoire cinq
concepts novateurs\footnote{\,
%%%%%%%%%%%%%%%%%%%%%%%-------DEBUT--------%%%%%%%%%%%%%%%%%%%%%%%%%%%
{\em Cf.}~Dombrowski, \cite{ do1979}.
}: %%%%%%%%%%%%%%%%%%%%%%%%-----FIN-----%%%%%%%%%%%%%%%%%%%%%%%%%%%%%%%

\begin{itemize}

\smallskip\item[{\bf 1)}]
l'application de Gauss, qui transporte les vecteurs unitaires normaux
à une surface au centre d'une sphère auxiliaire de rayon $1$;

\smallskip\item[{\bf 2)}]
la notion géométrique et analytique de {\sl mesure de courbure},
définie comme quotient des aires infinitésimales (orientées) sur la
surface et sur une sphère auxiliaire;

\smallskip\item[{\bf 3)}]
la courbure totale, intégrale de surface de la
courbure ponctuelle sur n'importe quelle région
de la surface; 

\smallskip\item[{\bf 4)}]
la variation angulaire $\Theta$ d'une paramétrisation de la surface en
coordonnées $(u, v)$, à savoir la $1$-forme différentielle:
\[
\Theta
:=
\frac{1}{2\,\sqrt{EG-F^2}}\,
\bigg(
\frac{F}{E}\,dE
+
E_v\,du
-
G_u\,dv
-
2\,F_u\,du
\bigg);
\]

\smallskip\item[{\bf 5)}]
les coordonnées géodésiques normales dans lesquelles existent deux
familles parallèles de géodésiques qui sont mutuellement
orthogonales en tout point, ce qui s'exprime, en termes des
coefficients métriques, par les normalisations suivantes:
\[
E
\equiv
1,\ \ \ \ \ \ 
F
\equiv
0,\ \ \ \ \ \
\text{\rm d'où:}\ \ \ \ \
ds^2
=
du^2+G\,dv^2.
\] 

\end{itemize}\smallskip

En tout cas, l'historiographie contemporaine s'est beaucoup interrogée
sur le {\sl temps long} de la genèse gaussienne
et sur les circonstances de l'apparition
de cette théorie entièrement nouvelle. 
L'exigence\,\,---\,\,grandissante avec l'âge\,\,---\,\,d'achèvement et
de maturité dans l'obtention de résultats que Gauss s'imposait, l'a
conduit à différer pendant près
d'une quinzaine d'année toute publication intermédiaire. 
En reprenant des éléments transmis notamment par
Stäckel, Dombrovski et Keyszig, 
ce sera l'occasion d'expliciter les tensions mathématiques
et les causalités philosophiques qui ont dirigé
Gauss vers l'obtention de son {\em Theorema Egregium}. 
Comme convenu, l'objectif des analyses
se focalisera ensuite principalement sur l'examen des
{\em calculs} qui ont été conduits par Gauss.

\smallskip\noindent{\bf Lien avec un mémoire d'Euler.}
\label{lien-memoire-Euler}
On trouve, dans le \S~VIII des {\em Disquisitiones generales
circa superficies curvas} un théorème
classique 
qui fournit une autre interprétation géométrique 
de la courbure de Gauss. 

\CITATION{
{\em La mesure de la courbure en chaque point d'une surface
est égale à une fraction dont le numérateur est l'unité,
et dont le dénominateur est le produit des deux
rayons de courbure extrêmes dans les
sections par des plans normaux}.
\REFERENCE{{\sc Gauss},~\cite{ga1828},~22}}

\noindent
En fait, on sait que dès
1813, Gauss était déjà en possession de cette formule:
\[
\kappa
=
\frac{1}{r_{\sf min}\cdot r_{\sf max}}
\]
qui égale la mesure de la courbure d'une surface 
en l'un de ses points à l'inverse du
produit entre les deux rayons de courbure principaux
en ce point, extrema 
découverts,
introduits et étudiés auparavant par Euler.

\CITATION{
{\sc Problème.} {\em Une surface dont la nature est connue étant coupée
par un plan quelconque, déterminer la courbure de la section,
qui en est formée}. [$\cdots$] 
Supposons donc qu'on tire par la différentiation\footnotemark
$dz = p\, dx + q \, dy$, de sorte que
$p = \big( \frac{ dz}{ dx} \big)$ et
$q = \big( \frac{ dz}{ dy} \big)$. [$\cdots$] 
Soit $z = \alpha \, y - \beta \, x + \gamma$ l'équation
qui détermine le plan. [$\cdots$] Alors le rayon osculateur
de la section sera exprimé en sorte:
\[
-\,
\frac{(\alpha\alpha+\beta\beta-2\alpha q+2\beta p
+(\alpha p+\beta q)^2+pp+qq)^{\frac{3}{2}}}{
\big(
(\alpha-q)^2\big(\frac{dp}{dx}\big)
+
(\beta+p)^2\big(\frac{dq}{dy}\big)
+
2(\alpha-q)(\beta+p)\big(\frac{dp}{dy}\big)
\big)\,
\sqrt{1+\alpha\alpha+\beta\beta}}.
\]
[$\cdots$] 
{\em Réflexion.} Dès qu'on connoit les rayons osculateurs pour trois
sections différentes, ceux pour toutes les autres en sont parfaitement
connues. [$\cdots$] {\em Réflexion.} Quelle que soit la courbure d'un
élément, les deux sections, dont l'une contient la plus grande
courbure et l'autre la plus petite, sont toujours normales entr'elles.
\REFERENCE{{\sc Euler},~\cite{eule1760},~2--4;~20--21}}

\footnotetext{\,
%%%%%%%%%%%%%%%%%%%%%%%-------DEBUT--------%%%%%%%%%%%%%%%%%%%%%%%%%%%
Euler considère la surface sous la forme d'un graphe $z = z ( x, y)$.
Pour être concis dans la présentation de la formule, les extraits
ciblés sont très légèrement abrégés.
} %%%%%%%%%%%%%%%%%%%%%%%%-----FIN-----%%%%%%%%%%%%%%%%%%%%%%%%%%%%%%%

Les recherches géométriques d'Euler étaient motivées 
par ses travaux sur la mécanique du solide
et par des problèmes issus du calcul des 
variations\footnote{\,
%%%%%%%%%%%%%%%%%%%%%%%-------DEBUT--------%%%%%%%%%%%%%%%%%%%%%%%%%%%
{\em Cf.} Kreyszig, \cite{kr1994}.
}. %%%%%%%%%%%%%%%%%%%%%%%%-----FIN-----%%%%%%%%%%%%%%%%%%%%%%%%%%%%%%%
Ses publications consacrées
aux courbes géodésiques tracées sur les surfaces
paraissent entre 1728 et 1732, mais c'est surtout le mémoire~\cite{
eule1760} paru à l'académie des scienes de Berlin que la postérité
paradigmatique a
conservé
en mémoire. 

Dans les années 1730, on trouve en effet une série d'articles sur les
courbes planes dans lesquels Euler développe l'idée de relier
systématiquement les deux coordonnées cartésiennes $x(s)$ et $y
(s)$\,\,---\,\,où $s$ est une coordonnée par longueur
d'arc sur la courbe\,\,---\,\,à 
une équation `naturelle' $r = r (s)$ qui relie
le rayon de courbure $r$ au paramètre
$s$. C'est à cette occasion qu'Euler
démontre qu'une particule-test massive libre de champ de force sur une
surface se déplace toujours le long d'une géodésique.
En 1744, il découvrit en appliquant
son calcul des variations que la {\sl caténoïde}, à savoir la 
surface de révolution qu'on obtient par rotation
autour de son axe de 
la chaînette d'équation polaire:
\[
\rho
=
{\rm const}\cdot \cosh\theta,
\]
est une {\sl surface minimale}, théorie dont les aspects globaux sont
encore très actifs aujourd'hui. Ce n'est que plus tard en 1776
(travail publié en 1785) que Meusnier, un élève de Monge, établit que
l'annulation de la {\sl courbure moyenne}:
\[
\frac{1}{2}\,\bigg(\frac{1}{r_{\sf min}}
+
\frac{1}{r_{\sf max}}\bigg)
\]
est une propriété nécessaire pour qu'une surface de classe au moins
$\mathcal{ C}^2$ soit d'aire minimale, localement dans l'espace de ses
déformations $\mathcal{ C}^2$. Ensuite, c'est 
dans son mémoire de 1760 qu'Euler montra que
le rayon de courbure $r_\theta$ d'une section normale quelconque à une
surface en un point $p$ qui fait un angle $\theta$ avec la direction
principale minimale est donné par la formule:
\[
\frac{1}{r_\theta}
=
\frac{\cos^2\theta}{r_{\sf min}}
+
\frac{\sin^2\theta}{r_{\sf max}}.
\]

\CITATION{
Ces conclusions renferment à peu près tout ce que l'illustre
Euler a le premier enseigné sur la courbure des surfaces.
\REFERENCE{{\sc Gauss},~\cite{ga1828},~20}}

Cette formule classique due
à Euler est habituellement considérée comme l'une des
premières contributions à une théorie générale des surfaces, en tant
qu'elle propose une mesure de la courbure des surfaces relativement à
toutes ses sections par des plans orthogonaux. 

\CITATION{
{\sc Réflexion.}
Pour comparer donc les courbures de deux élémens entr'elles,
on n'a qu'à chercher pour chacun les sections qui donnent
le plus grand et le plus petit rayon osculateur, et si
l'on trouve ces deux rayons les mêmes dans l'une
et l'autre, on peut prononcer hardiment que ces deux
élémens sont doués de la même courbure.
Et partant, pour connoitre la véritable courbure d'un 
élément quelconque de surface, il suffit d'en
chercher le plus grand et le plus petit
rayon osculateur: puisque ceux de toutes les autres
sections en sont déterminés parfaitement, en sorte qu'aucune
variété n'y sauroit plus avoir lieu.
\REFERENCE{{\sc Euler},~\cite{eule1760},~21}}

\Scholie
On comprend maintenant mieux pourquoi Gauss prononçait cette appréciation
quelque peu sévère <<\,Ces conclusions renferment à peu près tout ce
que l'illustre Euler a le premier enseigné sur la courbure des
surfaces\,>>. Deux raisons à ce `jugement
mathématique': 1) la formule $\frac{ 1}{
r_\theta} = \frac{ \cos^2 \theta}{ r_{\sf min}} + \frac{ \sin^2
\theta}{ r_{\sf max}}$ est si élémentaire, en comparaison de celles
que l'on trouve dans la {\em vraie théorie proprement
bidimensionnnelle des surfaces}, que Gauss est à même, en a peine
moins de deux pages (\S~VIII des
{\em Disquistiones generales circa superficies curvas}), 
de la re-démontrer {\em presque sans calcul},
d'en énoncer les conséquences, et de la rapporter au {\em vrai}
concept de la mesure de courbure qu'il vient d'introduire; 2) Euler, qui
avoue presque une
petite réticence de s'autoriser sans plus ample justification à 
<<\,prononcer hardiment\,>> que la {\em courbure} d'une
surface est parfaitement déterminée seulement par la courbure de
paires de courbes orthogonales, 
n'a en vérité pas rigoureusement respecté la
véritable métaphysique de l'{\em ouverture mathématique conceptuelle},
laquelle aurait commandé, d'après Gauss, de {\em maintenir fermement
ouverte} la question de savoir s'il n'existe pas un concept plus
adéquat de courbure qui serait proprement {\em bidimensionnel}. Il
est vrai que la courbure de Gauss des surfaces plongées dans l'espace
euclidien $\R^3$ est connue dès que les rayons de courbure extrémaux
{\em dans $\R^3$} de ses courbes sont connus en chaque point, et donc
d'une certaine manière, il est vrai que la métaphysique implicite de
la pensée d'Euler, qui consistait à saisir la courbure d'un objet
bidimensionnel par la collection de ses sous-objets unidimensionnels
librement mobiles,
possède un sens mathématique substantiel. N'est-il pas vrai
aussi que l'étude des courbes holomorphes entières unidimensionnelles
contenues dans des variétés projectives de type général 
de grand degré et de grande
dimension sert à étudier, grâce à ces {\em courbes algébriques-tests},
la structure interne
de ces variétés projectives, et éventuellement aussi, à déterminer
la présence de points à
coordonnées entières en relation avec l'arithmétique des équations
diophantiennes? N'est-il pas vrai aussi que les travaux de Lie sur
les groupes continus de transformations commencent par étudier à
fond la théorie des groupes à un terme \deutsch{einglidrige Gruppe},
{\em parce que} la structure d'un groupe quelconque
se reflète dans la collection de tous ses 
sous-groupes à un 
paramètre\,\,---\,\,toujours notoirement
plus nombreux que les sous-groupes à un nombre $r \geqslant 2$ de
termes, puisque nulle condition de fermeture algébrique par crochets
n'est exigible lorsque $r = 1$\,\,---, 
{\em cf.} notamment l'importance des sous-algèbres de Cartan
dans la présentation et l'organisation des travaux de classification
par Killing des algèbres de Lie complexes semi-simples.

Autrement dit, la pensée qui inspire Euler procède d'une métaphysique
mathématique légitime et toujours actuellement présente dans le
développement des mathématiques, métaphysique qui déclare
et décide une {\em
option naturelle 
d'étude} face à des objets qui
sont réellement complexes, notamment
à cause de leur dimension élevée. Mais il n'en
reste pas moins que Gauss\,\,---\,\,quelque peu implicitement puisque
les pensées qu'il a certainement développées dans ses {\em méditations
intérieures en recherche}, ne sont nullement 
(d)écrites\,\,---\,\,reproche en
quelque sorte à Euler sans vraiment 
le dire de ne pas avoir su maintenir plus
{\em ouverte} la question de la conceptualisation
de la courbure, question qui devait se poser fermement
dans un champ aussi neuf que la théorie des surfaces.
Et en vérité, Gauss avait si fondamentalement raison
de s'imposer une telle rigueur dans l'exigence de pensée
que sans cela, le caractère {\em intrinsèque} de la courbure
en dimension $\geqslant 2$ serait probablement resté invisible 
encore pendant quelques décennies, puisque la courbure
des courbes n'a de sens qu'en relation avec une dimensionalité
ambiante supplémentaire\,\,---\,\,fin 
de ce commentaire spéculatif étendu.
\qed

\smallskip
Euler découvrit donc
l'orthogonalité des directions principales dans lesquelles le rayon de
courbure est minimal et maximal, hors des points qui sont des {\sl
ombilics}\,\,---\,\,terminologie introduite par Monge et passée dans la
langue anglo-saxonne contemporaine de la géométrie
différentielle\,\,---\,\, {\em i.e.} qui satisfont $r_{\sf min} = r_{\sf
max}$. Ultérieurement, Meusnier a étudié la courbure des sections
obliques, et obtenu ce qu'on appelle aujourd'hui classiquement le {\sl
Théorème de Meusnier}: {\em la courbure en un point $p$ de la section
(courbe) $S \cap \Pi_{ \theta, \omega}$ d'une surface $S$ par un plan
$\Pi_{ \theta, \omega}$ faisant un angle $\omega$, où $\vert \omega
\vert < \frac{ \pi}{ 2}$, avec le plan normal $\Pi_\theta$ repéré par
l'angle $\theta$ relatif la direction principale minimale et contenant
la même droite tangente que $\Pi_\theta$:
\[
\Pi_\theta\cap T_pS
=
\Pi_{\theta,\omega}\cap T_pS,
\]
est égale au quotient par $\sin \omega$ de la
courbure de la section par ce plan normal:}
\[
\frac{1}{r_{\theta,\omega}}
=
\frac{1}{\sin\omega}\,
\frac{1}{r_\theta}
=
\frac{1}{\sin\omega}
\bigg(
\frac{\cos^2\theta}{r_{\sf min}}
+
\frac{\sin^2\theta}{r_{\sf max}}
\bigg).
\]
Cette formule et d'autres établies par Meusnier lui ont permis
de démontrer que les seules surfaces dont les deux rayons
de courbure sont égaux en tout point sont les plans et
les sphères.

\medskip\noindent{\bf Chronologie sommaire de la genèse.}
\label{chronologie-sommaire}
Pour commencer, quelques dates-clé dans l'émergence\footnote{\,
%%%%%%%%%%%%%%%%%%%%%%%-------DEBUT--------%%%%%%%%%%%%%%%%%%%%%%%%%%%
Gauss publiait en latin, et relativement 
peu, et attendait que les fruits de ses recherches soient vraiment
parfaits\,--\,``{\em Pauca, sed matura}'', disait-il. Ainsi, pour la
plupart des datations historiques de ses d\'ecouvertes, les historiens
se sont bas\'es sur deux sources d'information importantes~: sa
correspondance, et son c\'el\`ebre {\em Journal math\'ematique}. Ce
Journal, l'{\oe}uvre r\'esum\'ee d'une vie, long d'une vingtaine de pages
seulement, contient 146 \'enonc\'es extr\^emement brefs et dat\'es
pr\'ecis\'ement, de tous les r\'esultats que Gauss a d\'emontr\'es dans
sa vie et qu'il jugeait importants.
} %%%%%%%%%%%%%%%%%%%%%%%%-----FIN-----%%%%%%%%%%%%%%%%%%%%%%%%%%%%%%%
des idées de Gauss sur la géométrie différentielle\footnote{\,
%%%%%%%%%%%%%%%%%%%%%%%-------DEBUT--------%%%%%%%%%%%%%%%%%%%%%%%%%%%
{\em Cf.} Stäckel~\cite{stackel1890}, Dombrowski~\cite{ do1979}
et Kreyszig~\cite{ kr1994} pour ce qui va suivre.
}. %%%%%%%%%%%%%%%%%%%%%%%%-----FIN-----%%%%%%%%%%%%%%%%%%%%%%%%%%%%%%%

\smallskip\noindent$\bullet$
{\sf 1794:} premières réflexions sur les géométries 
non-euclidiennnes\footnote{\,
%%%%%%%%%%%%%%%%%%%%%%%-------DEBUT--------%%%%%%%%%%%%%%%%%%%%%%%%%%%
Dans une lettre adressées à Gerling le 10 octobre 1846, Gauss
affirme:

\CITATION{
Le théorème que M. Schweikart vous a signalé, d'après lequel
pour toute géométrie, la somme de tous les angles extérieurs
à un polygone diffère de $360^{\sf o}$ degrés d'une quantité
qui est proportionnelle à l'aire de sa surface, est le
premier théorème au seuil de cette théorie, 
un théorème dont j'ai reconnu la nécessité déjà en 1794.
\REFERENCE{\cite{do1979},~126}}

}. %%%%%%%%%%%%%%%%%%%%%%%%-----FIN-----%%%%%%%%%%%%%%%%%%%%%%%%%%%%%%%

\smallskip\noindent$\bullet$
{\sf 1810--13:} application de Gauss; 
concept de mesure de courbure; formule
$\kappa = \frac{ 1}{ r_{\sf min} \cdot r_{\sf max}}$. 

\smallskip\noindent$\bullet$
{\sf 1816:} \deutschplain{Schöne Theorem:} invariance (extrinsèque)
de la courbure par isométries extrinsèques, {\em cf.} la prochaine
Section
p.~\pageref{principe-raison-suffisante}
ci-dessous.

\smallskip\noindent$\bullet$
{\sf 1822:} \deutschplain{Copenhagen Preisschrift}.
Théorèmes de comparaison pour les angles d'un 
triangle géodésique. Coordonnées isothermes. 
Formules intrinsèques pour la courbure.

\smallskip\noindent$\bullet$
{\sf 1825:} \deutschplain{Neue Untersuchungen}. 
Somme des angles dans un triangle géodésique. 
Lemme d'orthogonalité. 
\'Equations de Gauss. Expression intrinsèque de la
courbure dans des 
coordonnées polaires géodésiques.

\smallskip\noindent$\bullet$
{\sf 1826--27:} Conception du plan final et
rédaction des
{\em Disquisitiones generales circa superficies curvas}. 

\medskip\noindent{\bf Débuts des Disquisitiones.} 
\label{debuts-Disquisitiones}
Il s'agit
maintenant d'explorer comment Gauss développa ses idées et conçut la
rédaction finale des {\em Disquisitiones}.

De même que le développement du calcul infinitésimal doit beaucoup aux
problèmes de la mécanique newtonienne classique, la théorie des
surfaces s'est principalement développée en relation avec des
problèmes appliqués de géodésie, de cartographie terrestre et
d'astronomie, trois facteurs de motivation extrêmement importants qui
ont su mobiliser les recherches générales de Gauss à 
un niveau supérieur d'abstraction.

Eu égard au fait que les calculs numériques ou formels de géodésie et
de représentation conforme s'avéraient rapidement assez considérables,
il convient de rappeler que Gauss entama 
dès son adolescence d'intenses calculs, d'abord
numériques\,\,---\,\,inventant `spontanément' la méthode
des moindres carrés à l'âge de dix-sept ans\,\,---, puis 
d'une nature plus
théorique. De 1795 à 1801, date de la publication finale des {\em
Disquisitiones arithmetic{\ae}}, Gauss débuta, entreprit et acheva ses
travaux en arithmétique supérieure. Ses découvertes en algèbre, sur
la cyclotomie\footnote{\,
%%%%%%%%%%%%%%%%%%%%%%%-------DEBUT--------%%%%%%%%%%%%%%%%%%%%%%%%%%%
Il s'agit des extensions de corps li\'ees
\`a l'\'equation alg\'ebrique $x^n-1=0$, dite 
{\em cyclotomique}, du grec {\em kyklos}, cercle 
et {\em tom\^e}, coupure. Dans le dernier
chapitre des {\em
Disquisitiones arithmetic{\ae}}, Gauss élabore
une th\'eorie nouvelle 
qui lui permet de g\'en\'eraliser
consid\'erablement son r\'esultat de 1796 sur la division du cercle en
dix-sept parties \'egales, et d'obtenir une {\em condition
n\'ecessaire et suffisante} pour qu'un polygone r\'egulier \`a un
nombre entier $n$ de c\^ot\'e soit constructible \`a la r\`egle et au
compas. Gauss d\'emontre en effet que
ce polygone r\'egulier \`a $n$ c\^ot\'es est constructible \`a la
r\`egle et au compas {\em si et seulement si} $n=2^{\mu}p_1p_2\cdots
p_k$ est le produit par une puissance de $2$ de nombres {\em premiers}
de la forme $p_j=2^{2^j}+1$, o\`u $j$ est un entier, $j=0,
1,2,\ldots$. Un nombre $p$ est dit {\em premier} s'il n'est pas
divisible par un nombre $q$ strictement inf\'erieur \`a $p$.  Ces
nombres $p_j$ sont appel\'es {\em nombres de Fermat}, car Pierre de
Fermat (1601-1665), conseiller au Parlement de Toulouse, et renomm\'e
pour ses recherches en arithm\'etique, avait conjectur\'e que tous les
nombres $p_j$ sont premiers. Les nombres $p_0=3$, $p_1=5$, $p_2=17$,
$p_3=257$, $p_4=65537$, sont premiers.  Mais le math\'ematicien suisse
Euler \'etablit en 1732 que $p_5=2^{32}+1 =4294967297$ n'est pas
premier, et qu'il est divisible par 641~; Legendre en 1780 montra que
$p_6$ est divisible par 274 177 et {\em on ne conna\^{\i}t
explicitement aucun $p_j$ premier pour $j\geq 5$}~! Une telle remarque
limite donc de mani\`ere inattendue la port\'ee du th\'eor\`eme de
Gauss, dans lequel tous les $p_j$ apparaissant dans
$n=2^{\mu}p_1p_2\cdots p_k$ doivent, en plus, \^etre des nombres
premiers. Ainsi, une conjecture, m\^eme c\'el\`ebre, peut tr\`es bien
s'av\'erer \^etre compl\`etement fausse!
} %%%%%%%%%%%%%%%%%%%%%%%%-----FIN-----%%%%%%%%%%%%%%%%%%%%%%%%%%%%%%%
et sur la théorie des fonctions elliptiques ont marqué
ces années les plus créatives de sa vie
mathématique, période durant laquelle la
géométrie ne joua presqu'aucun rôle. Qui plus est, vu de l'extérieur,
c'est-à-dire en se basant seulement sur les publications effectives,
les travaux ultérieurs de Gauss en astronomie\,\,---\,\,n'impliquant
que la géométrie élémentaire et la géométrie
sphérique\,\,---\,\,n'indiquent aucune réorientation en direction de
la géométrie abstraite ou supérieure. Il n'en reste pas moins que les
lettres et manuscrits non publiés découverts dans son
\deutschplain{Nachlass} ont montré que lorsque Gauss n'était 
pas encore officiellement astronome, il
créait déjà régulièrement des concepts géométriques fondamentaux, et
l'examen de ses notes a convaincu les exégètes qu'il a découvert
progressivement, au cours de ces années 1810 à 1816, tous les
théorèmes importants de sa théorie des surfaces, même si les {\em
preuves} que Gauss continuait à améliorer ont
évolué entretemps, c'est-à-dire de 1816 à 1827, toujours vers une plus
haute {\em causalité-en-vérité-interne}. C'est donc {\em avant}
de débuter ses travaux sur la géodésie (après 1816), que Gauss a
entrepris ses recherches sur les surfaces, en relation avec les
fondements de la géométrie et avec les géométries non euclidiennes. En
effet, en 1816, l'année où il fut nommé directeur du nouvel
observatoire de Göttingen, Gauss avait déjà finalisé l'idée de
représentations paramétriques {\em intrinsèques} pour les surfaces; il
avait aussi finalisé la notion d'application conforme entre surfaces,
le transfert des vecteurs unitaires normaux sur une sphère auxiliaire,
la notion d'isométrie infinitésimale, la notion de mesure de courbure
$\kappa$, la notion de courbure intégrée, l'invariance de la courbure
à travers les isométries infinitésimales, les propriétés des
géodésiques sur les sphéroïdes généraux ou elliptiques, et le théorème
généralisé de Legendre pour les surfaces arbitraires. Cette liste de
réalisations et d'accomplissements s'avère même encore plus
impressionante si l'on rappelle que Gauss s'occupait aussi durant la
même période de fonctions hypergéométriques, qu'il obtint deux
nouvelles preuves du théorème fondamental de l'algèbre, qu'il trouva
de nouvelles approches de la loi de réciprocité quadratique, qu'il
effectua de nombreuses observations astronomiques, et 
bien entendu aussi, qu'il conduisit de très nombreux
calculs numériques et formels.

\medskip\noindent{\bf Le Preisschrift de 1822.} 
\label{Preisschrift-1822}
La période de
géodésie pour Gauss commence avec ses travaux théoriques au sujet des
géodésiques tracées sur les sphéroïdes\,\,---\,\,le modèle étant la
surface de la Terre\,\,---, et du point de vue réellement pratique,
cette période comprend aussi tous ses travaux topographiques sur la
géographie du royaume de Hanovre entre 1821 et 1825. En 1822, Gauss
prépare son \deutschplain{Copenhagen Preisschrift} (\cite{
gauss1822}), sur les applications conformes qui lui valut le prix de
l'Académie de Copenhague. La {\sl conformalité}, un terme suggéré par
Gauss, est une condition naturelle pour les applications, puisque par
exemple la projection stéréographique et la projection de Mercator
sont toutes deux conformes, et des formules pour les applications
conformes générales de la sphère sur un plan ont été obtenues 
dans les années 1770 par Lambert,
Lagrange et Euler. Le progrès important de Gauss par rapport à ces
idées\footnote{\,
%%%%%%%%%%%%%%%%%%%%%%%-------DEBUT--------%%%%%%%%%%%%%%%%%%%%%%%%%%%
L'exigence inflexible en direction de la généralité s'exprime 
à ce sujet dans une lettre envoyée le 5 juillet 1816 à 
Schumacher\footnote{\,
%%%%%%%%%%%%%%%%%%%%%%%-------DEBUT--------%%%%%%%%%%%%%%%%%%%%%%%%%%%
Heinrich Christian {\sc Schumacher} (1780--1850),
docteur en droit puis
professeur en astronomie à Copenhague à 
partir de 1810, fut ami de Gauss de 1808 jusqu'à
sa mort.
}: %%%%%%%%%%%%%%%%%%%%%%%%-----FIN-----%%%%%%%%%%%%%%%%%%%%%%%%%%%%%%%

\CITATION{
J'ai aussi conversé avec Lindenau au sujet d'une question offerte
à la compétition \deutsch{Preisfage} qui devait être posée
dans le même journal. J'avais déjà pensé à un problème
intéressant: 
{\em Dans le cas général, projeter (appliquer) une surface
donnée sur une autre surface donnée, de telle sorte que
la surface image et la surface initiale soient
en similitude dans l'infinitésimal}.
Un cas spécial se produit lorsque la première surface
est une sphère, et la seconde un plan. Dans ce cas, la projection
stéréographique et la projection de Mercator sont
deux solutions particulières. Cependant, on veut la solution
générale pour tous les types de surfaces, 
de telle sorte qu'elle contienne tous ces cas particuliers.
\REFERENCE{\cite{do1979},~127}}

}, %%%%%%%%%%%%%%%%%%%%%%%%-----FIN-----%%%%%%%%%%%%%%%%%%%%%%%%%%%%%%%
c'est qu'il traite des applications conformes entre
surfaces quelconques munies d'une métrique infinitésimale
arbitraire:
\[
ds^2
=
E\,dp^2+2\,F\,dpdq+G\,dq^2,
\]
où $p$ et $q$ sont des paramètres internes de la surface. Dans le cas
où toutes les données sont analytiques réelles, Gauss a établi dans ce
mémoire qu'il existe toujours un changement de coordonnées locales:
\[
(p,q)
\longmapsto
(u,v)
=
\big(u(p,q),\,v(p,q)\big)
\]
qui transforme ce $ds^2$ dans de nouvelles coordonnées $(u, v)$ dites
{\sl isothermes} en un nouveau:
\[
ds^2
=
\lambda^2\,\big(du^2+dv^2)
\]
normalisé de telle sorte que 
le coefficient central $F$ s'annule identiquement
et de telle sorte aussi que les deux autres
coefficients $E$ et $G$ soient égaux à une certaine fonction
strictement positive $\lambda^2 = \lambda^2 ( u,v)$. Il est
intéressant de noter que la notion d'{\sl isométrie infinitésimale}
apparaît à cette occasion, lorsque Gauss considère deux surfaces $S$
et $\overline{ S}$ munies, respectivement, d'un $ds^2$ et d'un $d
\overline{ s}^2$, qui sont reliées par une relation {\em conforme} de
la forme:
\[
E\,du^2+2\,F\,dudv+G\,dv^2
=
\lambda^2\big(
\overline{E}\,d\overline{u}^2
+
2\,\overline{F}\,d\overline{u}d\overline{v}
+
\overline{G}\,d\overline{v}^2\big),
\]
sachant que lorsque $\lambda^2 = 1$, 
la première surface est `développable' 
\deutsch{abbildbar} sur la seconde.

Dans la préface de son \deutschplain{Preisschrift}, 
Gauss écrivait:

\smallskip
\centerline{\fbox{{\em Ab his via sternitur ad majora}},}

\smallskip\noindent
en s'inspirant visiblement de la conclusion aphoristique: 

\smallskip
\centerline{{\em Et his principiis via ad major sternitur},}

\smallskip\noindent
que Newton avait énoncée à la fin de son article fondateur sur le
calcul infinitésimal. Très probablement, les `accomplissements
majeurs' auxquels Gauss pensait à travers cet aphorisme référaient à
l'idée de {\sl géométrie intrinsèque} qui résulte de la relation
d'équivalence\,\,---\,\,nouvellement découverte\,\,---\,\,par 
isométries infinitésimales
entre surfaces. Ceci montre donc que Gauss considérait son {\em
Preisschrift} comme étant une {\em percée 
angulaire vers sa théorie générale des
surfaces à venir}.

On a aussi trouvé dans son \deutschplain{Nachlass} un manuscrit non
publié intitulé <<\,\deutschplain{Stand meiner Untersuchungen über die
Umformung der Flächen}\footnote{\,
%%%%%%%%%%%%%%%%%%%%%%%-------DEBUT--------%%%%%%%%%%%%%%%%%%%%%%%%%%%
<<\,{\em
\'Etat de mes recherches sur la transformation des surfaces}\,>>, 
\cite{ gauss1981}, VIII, 374--384.
}\,>> %%%%%%%%%%%%%%%%%%%%%%%%-----FIN-----%%%%%%%%%%%%%%%%%%%%%%%%%%%%%
dans lequel il consigne par écrit tous les
résultats obtenus ou compris à l'époque
de la publication de son \deutschplain{Preisschrift}. 
Notamment dans cette courte note
à usage strictement personnel, il mentionne qu'en coordonnées isothermes
$ds^2 = \lambda^2 ( du^2 + dv^2)$, la mesure de courbure
d'une surface se calcule
grâce à une formule simple, brève et symétrique:
\[
\kappa
=
-\,\frac{1}{\lambda^2}\,
\bigg(
\frac{\partial^2\log\lambda}{\partial u^2}
+
\frac{\partial^2\log\lambda}{\partial v^2}
\bigg).
\]

\Scholie
Eu égard aux objectifs philosophiques de ces analyses, il importe de
faire remarquer le caractère d'ores et déjà {\em définitif} des
arguments prouvant le caractère intrinsèque de la courbure, 
plusieurs années
{\em avant} que Gauss ne soit en mesure d'établir la fameuse et
complexe {\em formula egregia} valable dans tout système de
coordonnées curvilignes. 
En effet, on a dit plus haut que Gauss avait déjà établi
l'existence de coordonnées isothermes dans toutes les circonstances
dignes de considération à son époque\,\,---\,\,c'est-à-dire sous
l'hypothèse quelque peu implicite que toutes les données géométriques
et fonctionnelles soient toutes développables en série entière
convergente\,\,---, et donc par conséquent, la formule précédente
montre {\em dans toutes les circonstances possibles} que la mesure de
courbure est indépendante du plongement des surfaces dans l'espace
$\R^3$, et qu'elle dépend seulement de la métrique infinitésimale {\em
intrinsèque} à la surface.
\qed

\Scholie
Telle est l'un des aspects du (long) chemin vers les vérités
adéquates dans le labyrinthe de la recherche mathématique:
plusieurs moments diachroniques de conquêtes
progressives sur des territoires inconnus signalent des 
vérités eu égard à des parcours indirects et
néanmoins complètement convaincants quant à la certitude
mathématique; il n'en reste pas moins que 
l'exigence de comprendre en profondeur les causes et le pourquoi 
des choses demande de {\em découvrir d'autres chemins plus
satisfaisants qui mèneraient plus adéquatement aux mêmes
vérités}.
\qed

\Scholie
Ainsi Gauss ne sera-t-il parfaitement satisfait des arguments
qui prouvent le caractère intrinsèque de la courbure
que lorsqu'il pourra montrer par une formule explicite
{\em directe} que la courbure s'exprime en fonction 
d'un $ds^2$ absolument quelconque. 
\qed

\medskip\noindent{\bf Les Neue Untersuchungen de 1825.} 
\label{neue-1825}
Après le
\deutschplain{Preisschrift}, l'étape suivante en direction d'une
théorie générale des surfaces avant la réalisation finale des {\em
Disquisitiones} est la rédaction des \deutschplain{Neue
Untersuchungen} (\cite{ gauss1825}), écrites vers la fin de l'année
1825, et conçues comme une extension de la théorie des géodésiques à
un niveau supérieur et plus abstrait\footnote{\,
%%%%%%%%%%%%%%%%%%%%%%%-------DEBUT--------%%%%%%%%%%%%%%%%%%%%%%%%%%%
On doit à Stäckel d'avoir collecté avec grand soin
et édité ces manuscrits de Gauss.
}. %%%%%%%%%%%%%%%%%%%%%%%%-----FIN-----%%%%%%%%%%%%%%%%%%%%%%%%%%%%%%%
Cet avant-projet diffère des
{\em Disquisitiones} non seulement par une plus petite quantité de
matériel traité, mais aussi par la méthode générale des démonstrations,
qui repose sur l'utilisation fréquente des coordonnées géodésiques au
lieu des coordonnées paramétriques intrinsèques générales et
quelconques. Mis à part les considérations sur la courbure
géodésique, tous les thèmes traités dans l'avant-projet ont été inclus
par Gauss dans son mémoire final.

Qu'en est-il alors, à ce dernier stade d'évolution, des aspects reliés
à l'invariance de la mesure de courbure $\kappa$ à travers les
isométries infinitésimales? Malheureusement, la source exégétique du
{\em Theorema egregium} reste essentiellement inconnue, parce
qu'aucune note relative à sa dérivation n'a été découverte dans le
\deutschplain{Nachlass} pour la période qui précède l'année 1816.
Stäckel (\cite{ stackel1890}, p.~95) pense que c'est l'étude des
géodésiques qui a ouvert la voie, et on y reviendra dans la prochaine
Section p.~\pageref{principe-raison-suffisante} ci-dessous 
qui sera consacrée à
l'action du principe de raison suffisante. De plus, Gauss était
parfaitement en mesure dès 1822 de démontrer l'équivalent du {\em Theorema
egregium} pour des coordonnées isothermes, et aussi à nouveau
en 1825 pour des coordonnées polaires géodésiques dans lesquelles la
métrique infinitésimale est de la forme simplifiée:
\[
ds^2
=
dp^2+G(p,q)\,dq^2.
\]
Dans ces coordonnées spéciales, les deux familles de courbes $\{ p =
{\rm const.} \}$ et $\{ q = {\rm const.} \}$, 
dont la seconde est constituée
de géodésiques, s'intersectent
orthogonalement,
et Gauss y 
exprime en effet
la mesure de courbure $\kappa$ en termes de $G$ et de ses
dérivées partielles
du second ordre comme suit:
\[
\kappa
=
-\,\frac{1}{\sqrt{G}}\,
\frac{\partial^2}{\partial u^2}
\big(\sqrt{G}\big).
\]
Pour ces raisons qui démontrent
complètement l'invariance de la courbure
à travers les isométries 
infinitésimales, Gauss était intimement convaincu
qu'il devrait exister une formule analogue 
mais plus générale qui permettrait
de calculer la 
courbure dans un système de coordonnées
$(u, v)$ quelconques en termes des coefficients
métriques $E$, $F$, $G$ et de leurs dérivées
partielles d'ordre $1$ et $2$. 
Mais lorsqu'il se rendit compte que les preuves
valables pour les systèmes normalisés de coordonnées
ne se généralisaient pas aisément, il 
abandonna son travail sur les 
\deutschplain{Neue Untersuchungen} à la fin
de l'année 1825, écrivant le 19 février
1826 à son proche ami Olbers, astronome à Bremen qui
avait usé de son influence pour soutenir
sa nomination à Göttingen:

\CITATION{
Je ne connais presqu'aucune période de ma vie durant
laquelle j'aie obtenu un gain relativement si négligeable
\deutsch{so wenig reinen Gewinn}
d'un travail aussi acharné que celui conduit cet hiver.
J'ai découvert de {\em nombreux} beaux résultats, tandis
que mes efforts concernant d'autres idées restaient souvent
infructueux pendant des {\em mois}.
\REFERENCE{{\sc Gauss},~\cite{kr1994},~103}}

Ainsi en 1825, Gauss s'était-il fixé comme objectif de développer
et d'exposer une théorie générale des surfaces
seulement après être parvenu, au prix d'efforts
extraordinaires en calcul, à vaincre l'obstacle de la
{\em formula egregia}, ce qu'il réalisa enfin à la fin 
de l'année 1826.

\Section{10.~Action du principe de raison suffisante}
\label{principe-raison-suffisante}

\HEAD{10.~Action du principe de raison suffisante}{
Jo\"el Merker, D\'partement de Math\'matiques d'Orsay}

\medskip\noindent{\bf Triangles géodésiques.}
\label{triangles-geodesiques}
L'obstacle d'insatisfaction mathématique aura été d'autant plus
considérable que Gauss était en fait déjà en possession d'une preuve
indirecte du caractère purement intrinsèque de la mesure de courbure
{\em dès l'année 1816}\,\,---\,\,c'est-à-dire dix ans
avant la rédaction
des {\em Disquisitiones}\,\,---\,\,période
à laquelle il maîtrisait la {\sl formule d'écart
angulaire} relativement à $\pi$ de la somme des angles d'un triangle
géodésique $ABC$ tracé sur une surface:
\[
\widehat{A}+\widehat{B}+\widehat{C}
=
\pi
+
\int\!\!\int_{ABC}\,\kappa\,d\sigma,
\]
au moins dans les circonstances où
la courbure ponctuelle $p \mapsto \kappa (p)$ est de signe localement
constant sur 
la surface, quoique de magnitude {\em variable}.

\CITATION{
La somme des angles d'un (petit) triangle géodésique $\Delta$
sur une surface courbe dans $\E^3$ est égale à la somme
de $\pi$ et de l'aire orientée de l'image sphérique
de $\Delta$, où l'aire orientée est considérée comme
positive ou négative, suivant que le bord de l'image
sphérique de $\Delta$ tourne autour de l'image
dans la même direction ou dans la direction opposée, en 
comparaison à la manière dont le bord de $\Delta$ tourne 
autour $\Delta$.
\REFERENCE{{\sc Gauss},~\cite{gauss1981},~VIII,~435}}

Cette formule {\em directement inspirée} de la géométrie
sphérique était aussi largement connue dans le cas des surfaces
développables. Gauss possédait par ailleurs depuis 1794 un résultat
analogue valable dans le cadre de la géométrie hyperbolique, et
par voie de conséquence, il maîtrisait en vérité ladite formule pour
essentiellement toutes les surfaces de courbure constante\footnote{\,
%%%%%%%%%%%%%%%%%%%%%%%-------DEBUT--------%%%%%%%%%%%%%%%%%%%%%%%%%%%
{\em Cf.} Dombrowski~\cite{ do1979} pour ce qui va suivre.
}, %%%%%%%%%%%%%%%%%%%%%%%%-----FIN-----%%%%%%%%%%%%%%%%%%%%%%%%%%%%%%%
ce qu'il énonça notamment dans une lettre à Gerling en 
1819\footnote{\,
%%%%%%%%%%%%%%%%%%%%%%%-------DEBUT--------%%%%%%%%%%%%%%%%%%%%%%%%%%%
D'autres passages de cette lettre ont incité la tradition
commentatrice à penser que les mesures terrestres ou astronomiques
d'angle ont été effectuées par Gauss afin de `tester' la validité ou
la non-validité de la géométrie euclidienne dans l'espace physique.
}. %%%%%%%%%%%%%%%%%%%%%%%%-----FIN-----%%%%%%%%%%%%%%%%%%%%%%%%%%%%%%%

\CITATION{
Non seulement le défaut de la somme des angles dans
un triangle par rapport à $180^{\sf o}$ degrés
croît lorsque sa surface augmente, mais encore il
lui est exactement proportionnel.
\REFERENCE{{\sc Gauss},~\cite{gauss1981},~VIII,~10}}
 
\noindent
Il est donc tout
à fait légitime de s'imaginer que Gauss ait pu être conduit à deviner, 
dès 1812,
que la validité de cette formule devait être ensuite testée dans le
cas des surfaces quelconques plongées dans l'espace, quelle que soit
la variabilité de leur courbure, sachant
qu'il était en outre motivé par l'étude des 
géodésiques sur des sphéroïdes {\em non}
sphériques, et aussi par
l'étude générale des applications conformes.

\Scholie
Déformer un objet mathématique par extension en généralité, c'est
l'exposer à une ouverture hypothétique qui consiste à {\em lever} une
ou plusieurs conditions spécifiques auxquelles sa définition initiale
le soumettait. La dialectique entre la constance et la
non-constance\,\,---\,\,ici de la courbure\,\,---\,\,est toujours
universellement susceptible de s'incarner en rapport avec toute
fonction qui est construite en relation interne avec une catégorie
d'objets mathématiques, quel que soit le niveau de sophistication de
la construction en question, et la non-constance opposée à la
constance d'une donation initiale est {\em fréquemment} susceptible de
{\em créer} une ontologie nouvelle dont l'existence, si elle est
attestée, confirme le nouveau type d'être à étudier. D'un seul coup
alors, une multiplicité automatique de questions mathématiques
surgissent, dont la trame métaphysique se réduit à la conservation, ou
à la disparition, en totalité, ou en partie, des vérités connues dans
des circonstances déjà étudiées. Aussi du point de vue de la
philosophie générale des mathématiques, la structure interne du
questionnement mathématique se réduit-elle souvent à la plus simple
interrogation au sujet d'une permanence ou d'une disparition de
propriétés établies. \`A un niveau plus élevé, ce fait soulève
la méta-question extrêmement délicate du statut des {\em questions}
dans le développement des mathématiques, et il
est vraisemblablement impossible d'apporter une
`réponse' à cette méta-question. 
\qed

\Scholie
Spécialement ici, le concept\,\,---\,\,ici, de courbure\,\,---\,\,se
retourne en lui-même dans la réflexivité de son expression. En effet,
la formule quantitative pour l'excès angulaire $\widehat{ A} +
\widehat{ B} + \widehat{ C} - \pi$ de tout triangle sur une sphère
bordé par des grands cercles s'interprétait comme multiple (constant)
de l'aire du triangle, donc comme dérivant d'un principe de
distance à l'homogénéité, et c'est parce que la notion de courbure
exprime le rapport infinitésimal de toute surface à une sphère qu'une
formule {\em intégrée} pour la différence moyenne était
vraisemblablement susceptible
d'exister. Autrement dit, la variabilité se {\em reflète}
infinitésimalement dans l'homogénéité, et son rapport 
en généralité par rapport à l'homogénéité {\em cause en quelque
sorte métaphysiquement} l'extension de la
formule d'écart angulaire grâce 
à l'intégration qui rend fini l'infinitésimal.
\qed

\smallskip\noindent{\bf Déduction du caractère intrinsèque
de la courbure.} \label{deduction-intrinseque} C'est en 1816 que Gauss
a découvert le caractère proprement intrinsèque de la courbure dans le
cas général des surfaces quelconques, grâce à la formule d'écart
angulaire rapellée à l'instant. Tel est le contenu du résultat qu'il
appela `joli théorème' \deutsch{Schöne Theorem}, 
et que l'on peut restituer comme suit:

\smallskip 
{\em Si une courbe sur laquelle une figure est fixée prend différentes
formes dans l'espace, alors l'aire de l'image sphérique de cette
figure sur la sphère auxiliaire est toujours la même, pour toutes les
formes que la surface peut prendre par flexions sans extension, sans
duplicature, sans déchirure}.

\smallskip
L'attrait de ce théorème, c'est la découverte d'une invariance, 
d'une permanence de propriété lors de tout mouvement isométrique
d'une surface. Effort de mentalisation, questionnements, 
visions, interrogations: la véritable 
{\em genèse causale}, chez
Gauss, de la 
{\em rigidité bidimensionnelle} de la courbure ne pourra
être 
connue en totalité, 
puisque cette genèse impliquait des jours, des mois, 
des années de réflexions et de méditations, riches
d'aspects métaphysiques que tout texte mathématique
formel se refuse à expliciter. 

Le 5 juillet 1816, Gauss fait part dans une lettre à son ancien
étudiant et ami (confident?) Schumacher de son insatisfaction
concernant ce `joli théorème'\footnote{\,
%%%%%%%%%%%%%%%%%%%%%%%-------DEBUT--------%%%%%%%%%%%%%%%%%%%%%%%%%%%
La suggestion qu'avait faite Gauss en 1816 aux
éditeurs des {\em Astronomische Abhandlundgen} de mettre
à la compétition le problème général des applications
conformes au sujet duquel il possédait
déjà de nombreuses idées n'avait pas été entendue. C'est pourquoi 
Schumacher usa en 1821 
de la première opportunité qui
lui fut offerte pour suggérer à la 
{\em Société scientifique de Copenhague}
de mettre ce problème au concours; toutefois, Gauss n'en fut pas
immédiatement informé. Alors en 1822, 
le problème n'ayant toujours pas reçu
de solution, il fut remis au concours, et Schumacher en informa Gauss
le 4 juin. Dès le 10 juin, Gauss répondit: <<\,Je suis désolé
d'apprendre seulement maintenant le renouveau de cette compétition
[$\cdots$]\,>>. Le 25 novembre, Gauss demanda à Schumacher quelle est
la date limite pour la soumission des propositions: fin de l'année
1822, lui répondit-il. Le 11 décembre, sous
la pression de ce délai, Gauss soumit son premier mémoire~\cite{
gauss1822} sur la théorie des surface qui répondait à
l'une de ses propres questions, et 
qui ne sera publié qu'en 1825 aux 
{\em Astronomischen
Abhandlungen}, après avoir été récipiendaire du prix.
}, %%%%%%%%%%%%%%%%%%%%%%%%-----FIN-----%%%%%%%%%%%%%%%%%%%%%%%%%%%%%%%
et on reviendra plus bas sur les aspects fondamentaux
de l'exigence gaussienne qui légitimaient ses insatisfactions récurrentes.

\CITATION{
Ce résultat est sans doute tout à fait correct, mais il ne
suffit pas à résoudre la question générale posée [$\cdots$].}

D'un point de vue géométrique moderne toutefois, l'argument de 1816 
qui déduit l'invariance de la courbure à partir
de la formule d'excès angulaire est correct, complet
et convaincant\,\,---\,\,le voici donc. 
Soit comme supposé:
\[
\phi\colon\,\,\,\,
S\to S',
\ \ \ \ \
(u,v)
\longmapsto
(u',v')
=
\big(u'(u,v),\,v'(u,v)\big)
\]
une application (analytique réelle) inversible
entre deux surfaces munies, respectivement,
de métriques quadratiques
infinitésimales $ds^2$ et $d{s'}^2$, qui est une
{\em isométrie}\,\,---\,\,ni flexion, ni duplicature, 
ni déchirure\,\,---, {\em i.e.} qui satisfait 
$\phi^* ( d{ s'}^2 ) = ds^2$, c'est-à-dire
plus précisément:
\[
E'\,{du'}^2
+
2\,F'\,du'dv'
+
G'\,{dv'}^2
\Big\vert_{u'=u'(u,v),\,v'=v'(u,v)}
=
E\,du^2
+
2\,F\,dudv
+
G\,dv^2.
\]
Puisque l'application 
$\phi$ respecte les distances dans l'infiniment petit,
elle laisse aussi inchangée la longueur des
courbes macroscopiques lors du passage de
$S$ à $S'$, et {\em vice versa}. 
Comme les courbes géodésiques minimisent par
définition localement les longueurs, 
$\phi$ transforme les géodésiques en géodésiques.
Pour ce qui est des éléments d'aire infinitésimale 
associés par $ds^2$ à $S$ et par
$d{ s'}^2$ à $S'$, à savoir:
\[
d\sigma
=
\sqrt{EG-F^2}\,dudv
\ \ \ \ \ \ \ \
\text{\rm et}
\ \ \ \ \ \ \ \
d\sigma'
=
\sqrt{E'G'-{F'}^2}\,du'dv',
\]
il est clair aussi que $\phi$ les transforme aussi l'un vers l'autre:
\[
\phi^*(d\sigma')
=
d\sigma,
\]
puisque toute aire de rectangle infinitésimal orthogonal
s'exprime comme produit de sa longueur par sa largeur.

Par conséquent, $\phi$ envoie tout petit voisinage
$V_\varepsilon ( p)$
d'un point $p$ bordé par trois ou plus 
géodésiques\footnote{\,
%%%%%%%%%%%%%%%%%%%%%%%-------DEBUT--------%%%%%%%%%%%%%%%%%%%%%%%%%%%
Du théorème sur les triangles, Gauss 
déduit immédiatement par découpage un 
théorème analogue plus général sur l'excès
angulaire de tout polygone géodésique: 

\CITATION{
La somme de tous les angles d'un petit polygone géodésique
$\Pi$ à $n$ côtés tracé sur une surface courbe $S$ dans
$\E^3$ est égale à la quantité $(n-2)\, \pi$ à laquelle
il faut additionner l'aire de surface orientée
de l'image sphérique de $\Pi$ obtenue
{\em via} l'application 
vers la sphère auxiliaire.
\REFERENCE{{\sc Gauss},~\cite{gauss1981},~VIII,~435}}
} %%%%%%%%%%%%%%%%%%%%%%%%-----FIN-----%%%%%%%%%%%%%%%%%%%%%%%%%%%%%%%
vers un autre voisinage: 
\[
V_\varepsilon'(p')
:=
\phi\big(V_\varepsilon(p)\big)
\]
du point image $p' = \phi (p)$ qui est lui aussi bordé par trois ou 
plus géodésiques et qui possède en outre la {\em même} aire:
\[
{\sf aire}\big(V_\varepsilon'(p')\big)
=
{\sf aire}\big(V_\varepsilon(p)\big).
\]

Par ailleurs, en supposant pour
fixer les idées que ledit voisinage est un triangle géodésique:
\[
V_\varepsilon (p) 
=
A_\varepsilon B_\varepsilon C_\varepsilon,
\]
l'application $\phi$ produisant une isométrie 
dans l'infinitésimal, elle préserve tous les angles entre
courbes ou droites tangentes, et donc on a des
égalités par paires pour les angles aux sommets entre
le triangle et son triangle image:
\[
\widehat{A}_\varepsilon'
=
\widehat{A}_\varepsilon,
\ \ \ \ \ \ \ 
\widehat{B}_\varepsilon'
=
\widehat{B}_\varepsilon,
\ \ \ \ \ \ \ 
\widehat{B}_\varepsilon'
=
\widehat{B}_\varepsilon.
\]
Par conséquent, si l'on compare ces trois égalités d'angles
aux deux formules respectives d'excès angulaire:
\[
\aligned
\widehat{A}_\varepsilon+\widehat{B}_\varepsilon+\widehat{C}_\varepsilon
&
=
\pi
+
\int\!\!\int_{A_\varepsilon B_\varepsilon C_\varepsilon}\,
\kappa\,d\sigma,
\\
\widehat{A}_\varepsilon'+\widehat{B}_\varepsilon'+\widehat{C}_\varepsilon'
&
=
\pi
+
\int\!\!\int_{A_\varepsilon'B_\varepsilon'C_\varepsilon'}\,
\kappa'\,d\sigma',
\endaligned
\]
on déduit immédiatement l'égalité entre
intégrales de courbures:
\[
\int\!\!\int_{A_\varepsilon B_\varepsilon C_\varepsilon}\,
\kappa\,d\sigma
=
\int\!\!\int_{A_\varepsilon'B_\varepsilon'C_\varepsilon'}\,
\kappa'\,d\sigma',
\]
d'où on déduit
aussi en faisant tendre $\varepsilon$ vers $0$\,\,---\,\,puisque l'on
peut naturellement supposer que l'aire des triangles est un infiniment
petit d'ordre deux en $\varepsilon$\,\,---\,\,:
\[
\kappa(p)
\cdot
\zero{{\sf aire}
\big(A_\varepsilon B_\varepsilon C_\varepsilon
\big)}
+
{\rm O}(\varepsilon^3)
=
\kappa'(p')
\cdot
\zero{{\sf aire}
\big(A_\varepsilon'B_\varepsilon'C_\varepsilon'
\big)}
+
{\rm O}(\varepsilon^3),
\]
ce qui donne en simplifiant de part et d'autre par les coefficients
d'aires {\em qui sont égaux}\footnote{\,
%%%%%%%%%%%%%%%%%%%%%%%-------DEBUT--------%%%%%%%%%%%%%%%%%%%%%%%%%%%
Une autre manière de déduire l'égalité ponctuelle des courbures aurait
été d'écrire:
\[
\kappa(p)
=
\lim_{V_\varepsilon(p)\,\to\,p}\,
\frac{{\sf aire}\,\big[{\bf n}\big(V_\varepsilon(p)\big)\big]}{
{\sf aire}\,V_\varepsilon(p)}
=
\lim_{\phi(V_\varepsilon(p))\,\to\,p'}\,
\frac{{\sf aire}\,\big[{\bf n}\big(V_\varepsilon'(p')\big)\big]}{
{\sf aire}\,V_\varepsilon'(p')}
=
\kappa'(p').
\]
}: %%%%%%%%%%%%%%%%%%%%%%%%-----FIN-----%%%%%%%%%%%%%%%%%%%%%%%%%%%%%%%
\[
\boxed{
\kappa'(p')
=
\kappa(p)
}\,.
\] 
En conclusion, la courbure de Gauss au point $p$ et celle au
point-image $p' = \phi (p)$ sont égales entre elles à travers
l'isométrie $\phi$ entre les deux surfaces $S$ et $S'$: telle est la
démonstration géométrique du {\em Theorema egregium} obtenue par Gauss
dès 1816.

\medskip\noindent{\bf Insatisfaction gaussienne.} 
\label{insatisfaction-gaussienne}
Toutefois, Gauss ne publia jamais cette preuve.  L'une des raisons
tient à ce qu'elle en reste au niveau d'une présentation informelle,
sans préciser rigoureusement la notion d'aire orientée de surface sur
la sphère auxiliaire, sans réellement traiter le cas général des
surfaces dont
la courbure est alternativement positive ou négative\,\,---\,\,il
faut attendre 1825 pour que tel soit le cas.  Gauss jugeait durement
les démonstrations non rigoureuses, et s'appliquait à lui-même une
exigence mathématique dont il informait ses contemporains dans sa
correspondance pour se justifier de publier
relativement peu\footnote{\,
%%%%%%%%%%%%%%%%%%%%%%%-------DEBUT--------%%%%%%%%%%%%%%%%%%%%%%%%%%%
La personnalit\'e de Gauss et son choix de vie apparaissent
habituellement tr\`es fascinants aux yeux des math\'ematiciens, et son
exp\'erience t\'emoigne d'une possible {\em autonomie compl\`ete} de
l'activit\'e math\'ematique. Gauss s'\'etait en effet
convaincu, \`a partir d'exp\'erien\-ce v\'ecues, qu'il n'aurait que
peu de choses \`a apprendre \`a vouloir communiquer et \`a \'echanger
avec les autres math\'emati\-ciens; aussi pr\'ef\'era-t-il
s'isoler presque compl\`etement du champ des influences de
l'activit\'e math\'ematique de l'\'epoque.
}. %%%%%%%%%%%%%%%%%%%%%%%%-----FIN-----%%%%%%%%%%%%%%%%%%%%%%%%%%%%%%%
Refuser de publier comme le fit constamment Gauss dans sa maturité,
c'est refuser de se libérer artificiellement de ce qui n'est pas
encore compris.

\Scholie
D'autres raisons plus profondes peuvent être considérées comme
véritables causes d'une auto-critique récurrente chez Gauss, et donc
aussi chez tout mathématicien conscient de l'imperfection ubiquitaire
des mathématiques (<<\,platonisme des terrains vagues\,>>). Ces
causes ont trait à l'{\em essence-en-recherche} des vérités
mathématiques, c'est-à-dire à ce qui fait que la non-accession à des
parcours de preuves réellement adéquats oblige à {\em maintenir
intacte} l'exigence de tracer d'{\em autres} chemins démonstratifs qui
soient plus convaincants, plus systématiques, plus
inter-expressifs. Certes, l'informel, l'intuitif et l'approximatif offrent
un premier contact avec les vérités mathématiques, mais dans des
moments ultérieurs, on doit chercher à comprendre synthétiquement
l'ensemble d'une manière unifiée tout en contrôlant mieux l'extension
ontologique des objets concernés\footnote{\,
%%%%%%%%%%%%%%%%%%%%%%%-------DEBUT--------%%%%%%%%%%%%%%%%%%%%%%%%%%%
Pour ce qui est des applications conformes, Gauss a dû en élargir
considérablement le domaine ontologique, passant de transformations
{\em géométrico-algébriques} simples et de bas degré à un univers
fonctionnel arbitraire satisfaisant {\em seulement} des équations aux
dérivées partielles d'ordre $1$.  Une extension ontologique analogue
mais postérieure sera effectuée dans la thèse de Riemann, où les
fonctions holomorphes seront définies {\em seulement} comme annihilées
par l'opérateur $\frac{ \partial}{ \partial \overline{ z}}$.
}. %%%%%%%%%%%%%%%%%%%%%%%%-----FIN-----%%%%%%%%%%%%%%%%%%%%%%%%%%%%%%%
Autrement dit, les causes des vérités mathématiques exigent une
maturation que seuls peuvent procurer des parcours démonstratifs {\em
pluriels}, et la pensée du mathématicien se doit de comparer et
d'évaluer les niveaux de satisfaction abstraite que procure telle ou
telle preuve. C'est bien cela qui fascine dans la pensée
épistolière de Gauss: la présence d'une philosophie implicite des
mathématiques dont les aspects saillants demeurent universels et
trans-historiques.
\qed 

\smallskip
Dans le manuscrit pré-finalisatoire de 1825 (\cite{ gauss1825}), on
trouve une organisation plus élaborée des arguments de preuve, dans
laquelle la démonstration géométrique de l'invariance de la courbure
exposée ci-dessus est suivie du lemme d'orthogonalité aujourd'hui dit
<<\,de Gauss\,>>, puis d'une expression explicite $\kappa = - \frac{
1}{ \sqrt{ G}}\, \frac{ \partial^2}{ \partial u^2} \big( \sqrt{ G}\big)$
de la courbure en coordonnées polaires géodésiques. Dombrowski (\cite{
do1979}, p.~135) suggère alors que la raison pour laquelle ce
manuscrit s'achève abruptement à ce moment-là est vraisemblablement
due au fait suivant: Gauss aurait réalisé qu'une {\em autre et
meilleure preuve} du {\em Theorema egregium} découlait directement de
cette formule $\kappa = - \frac{ 1}{ \sqrt{ G}}\,
\frac{ \partial^2}{ \partial u^2} \big( \sqrt{ G}\big)$,
puisque toute quantité géométrique telle qu'un $ds^2$ se transforme
d'une manière déterminée à travers tout changement de coordonnées, et
donc cette formule pour la courbure devrait se transformer en une
formule valable dans un système de coordonnées quelconque, quoique le
calcul n'ait rien d'élémentaire ou d'aisé.

\CITATION{
Eu égard à l'article 21 des {\em Disquisitiones generales}, on
pourrait même être conduit à conclure qu'à cette époque, Gauss perçut
la possibilité\,\,---\,\,en principe\,\,---\,\,d'exprimer les
coefficients $E$, $F$, $G$ de la première forme fondamentale dans une
carte locale de la surface en termes des coefficients $E'$, $F'$, $G'$
d'une autre carte. Par conséquent, il aurait pu être possible
d'exprimer la courbure $\kappa = - \frac{ 1}{ \sqrt{ G}}\,
\frac{ \partial^2}{ \partial u^2} \big( \sqrt{ G}\big)$
obtenue d'abord pour le cas spécial des coordonnées géodésiques
polaires dans lesquelles $ds^2 = dp^2 + G\, dq^2$ sous la forme de
l'équation de Gauss [{\em formula egregia}]. De cette manière, Gauss
pourrait avoir découvert la formule explicite~\thetag{
\ref{formula-egregia}} p.~\pageref{formula-egregia} dans des
coordonnés paramétriques quelconques sur la surface. Mais il est
presque certain que Gauss ne trouva pas cette expression avant l'année
1826.
\REFERENCE{\cite{do1979},~135}}

Il est donc vraisemblable toutefois que Gauss ait été en mesure
d'effectuer un tel calcul, qu'il ait découvert ainsi la {\em formula
egregia}, mais qu'à nouveau, il ait été insatisfait du caractère {\em
toujours indirect} de l'accès à ladite formule.

\CITATION{
Ainsi, il semble que Gauss considéra d'une façon conclusive que sa
première démonstration géométrique de l'invariance de la courbure
était `périmée' et `détrônée'.
\REFERENCE{\cite{do1979},~135}}
 
\Scholie
\`A cet instant précis, une spéculation cruciale
sur la philosophie générale des mathématiques s'impose. L'{\em
exploration-en-recherche} de la `réalité' mathématique n'est
aucunement prédéterminée par l'insertion du sujet dans des structures
prédéfinies, fussent-elles transcendantales. Au contraire, la {\em
pensée-en-recherche} se déploie comme un {\em mobile (dés-)articulé de
questionnements multiples}, chaque décision pouvant être soumise à
révision. L'irréversible-synthétique est une quête à travers laquelle
il faut interroger, s'orienter, absorber, tel sera aussi le point de
vue de Lie, systématique dans la ramification du possible qu'il
déploiera sans relâche. Gauss réoriente inlassablement sa vision de la
théorie des surfaces jusqu'à ce que tous les rôles y soient distribués
réflexivement. {\em La vérité mathématique n'est pas indépendante des
preuves, elle n'est pas non plus langagièrement démultipliable par
identification à ses variations démonstratives, mais même acquise ou
finalisée, {\sf elle continue à s'\underline{interroger}
inlassablement et à s'évaluer réflexivement dans ses propres preuves
métamorphosables.}} Aussi Gauss n'a-t-il multiplié les
preuves\,\,---\,\,du {\em Theorema egregium}, du théorème fondamental
de l'algèbre, de la loi de réciprocité quadratique\,\,---\,\,que parce
que les mathématiques éprouvent toujours leurs réalisations dans la
tension d'un {\em Multiple-du-parcours} qui est leur essence même.
\qed

\medskip\noindent{\bf Raison suffisante leibnizienne.} 
\label{raison-suffisante-leibnizienne}
D'après
Leibniz, {\sl vérités nécessaires} et {\sl vérités de fait} doivent
être soigneusement distinguées. En effet, les {\sl vérités
nécessaires}, puisqu'elles sont nécessaires en elles-mêmes, ne
pourraient pas ne pas être, et elles reposent en dernier recours sur
des vérités primitives, et notamment
sur le principe de contradiction d'après lequel tout ce qui est
contradictoire est faux, et tout ce qui est logiquement
non-contradictoire a des chances d'être nécessaire. 

\CITATION{
33. Il y a aussi deux sortes de vérités, celles de Raisonnement et
celles de Fait. Les vérités de Raisonnement sont nécessaires et leur
opposé est impossible, et celles de Fait sont contingentes et leur
opposé est possible. Quand une vérité est nécessaire, on en peut
trouver la raison par l'analyse, la résolvant en idées et en vérités
plus simples, jusqu'à ce qu'on vienne aux primitives.
\REFERENCE{\cite{leibniz}}}

En particulier, les vérités obtenues par le raisonnement mathématique
sont nécessaires, puisqu'on y parvient, d'après une
affirmation leibnizienne\,\,---\,\,approfondie au vingtième
siècle
avec les théories de la démonstration\,\,---\,\,par des 
raisonnements analytiques qui ramènent en
dernier recours tous les énoncés à la certitude incontournable et
absolue de la non-contradiction.

En revanche, toujours d'après Leibniz, les {\sl vérités de fait}
concernent des étants et des événements qui sont intrinsèquement
contingents, comme l'existence de tels arbres ou de tels êtres humains,
qui existent ici et là, mais qui auraient tout à fait pu ne pas exister,
à cause d'un affouage, à 
cause d'une construction ou à cause d'une guerre, de
la maladie. Mais ce n'est pas pour autant que l'existence contingente
est sans raison, puisque chaque arbre existe parce qu'une graine a
germé, et puisque chaque homme existe parce que ses parents se sont 
mutuellement fécondés.

Leibniz choisit alors de distinguer ces deux types de vérités à l'aide
d'un critère de finitude ou d'indéfinitude dans la régression des
raisons: en effet, la différence entre nécessité analytique et
contingence mondaine, c'est que dans le premier cas, la Raison qui
raisonne peut conduire en un nombre fini d'étapes l'analyse des
raisons jusqu'à un terme, tandis que dans le second cas, la Raison ne
peut jamais comprendre et embrasser tout le détail {\em
potentiellement indéfini} des choses de l'univers afin de rendre
complètement raison d'un fait contingent spécifique donné.

\CITATION{
36. Mais la raison suffisante se doit trouver aussi dans les vérités
contingentes ou de fait, c'est-à-dire, dans la suite des choses
répandues par l'univers des créatures ; où la résolution en raisons
particulières pourrait aller à un détail sans bornes, à cause de la
variété immense des choses de la Nature et de la division des corps
à l'infini. Il y a une infinité de figures et de mouvements présents
et passés qui entrent dans la cause efficiente de mon écriture
présente; et il y a une infinité de petites inclinations et
dispositions de mon âme, présentes et passées, qui entrent dans la
cause finale.
\REFERENCE{\cite{leibniz}}}

\medskip\noindent{\bf Racine métaphysique du principe.}
\label{racine-metaphysique-principe}
Malgré cette pluralité et cette indéfinitude des raisons contingentes
et particulières, il reste néanmoins\,\,---\,\,en toute circonstance
interrogative\,\,---\,\,de la nature, de l'{\em essence} de la {\em
Raison} de penser que ce qui n'a aucune {\em raison} d'exister devrait
n'exister pas, et partant aussi dans le sens inverse-réciproque, que
tout étant qui existe, que ce soit de manière nécessaire ou de manière
contingente, parce qu'il s'insère dans un certain tissu voilé de
relations spatio-temporelles ou abstraites en vertu duquel il
entretient des liens souvent mystérieux avec d'autres êtres
nécessaires ou contingents, tout étant qui 
existe, donc, devrait trouver en lui-même
et autour de lui-même certaines {\sl raisons d'être} qui lui sont
propres dans la localité métaphysique partiellement invisible de
causes internes ou externes à son 
existence. C'est donc à cet instant-là
qu'apparaît le si profond {\sl Principe de raison suffisante} attribué
à Leibniz:

\CITATION{
\!\!\!\!\!\!\!\!\!\!\!\!
principe en vertu
duquel nous considérons qu'aucun fait ne saurait se trouver vrai ou
existant, aucune énonciation véritable, sans qu'il y ait une raison
suffisante pourquoi il en soit ainsi et non pas autrement, quoique ces
raisons le plus souvent ne puissent point nous être connues,
\REFERENCE{\cite{leibniz}}}

\smallskip\noindent
à savoir au moment où la méditation métaphysique déclare au sujet du
monde de tous les étants\,\,---\,\,qu'ils soient nécessaires ou
contingents, corporels
ou incorporels\,\,---\,\,que des raisons profondes doivent toujours agir
en réponse au moins
partielle à des interrogations spontanées et irrépressibles de la
Raison. 

Enfin, d'après Leibniz, le principe de raison suffisante est un des
<<\,deux grands principes de nos raisonnements\,>>, avec le principe
de non-contradiction\footnote{\,
%%%%%%%%%%%%%%%%%%%%%%%-------DEBUT--------%%%%%%%%%%%%%%%%%%%%%%%%%%%
Le principe de raison suffisante ne peut être réduit, chez Leibniz, au
principe de raison nécessaire. Le principe de raison suffisante est
lié au principe selon lequel tout prédicat est inhérent au sujet
({\em Praedicatum inest subjecto}). Il découlerait même de celui-ci, car
s'il y avait une vérité sans raison, alors, nous aurions une
proposition dont le sujet ne contiendrait pas le prédicat, ce qui est
absurde. Mais il va sans dire que le
concept d'{\sl irréversible-synthétique} en mathématique
annule complètement cette illusion que la connaissance
mathématique soit analytique.
}. %%%%%%%%%%%%%%%%%%%%%%%%-----FIN-----%%%%%%%%%%%%%%%%%%%%%%%%%%%%%%%
Il peut se résumer de manière aphoristique en
l'expression latine: {\em nihil est sine ratione}, à savoir: <<\,rien
n'est sans raison\,>>. Christian Wolff, qui reprend et étend le
rationalisme de Leibniz, l'affirme à son tour : <<\,{\em Rien n'existe
sans qu'il y ait une raison pour qu'il en soit ainsi et non
autrement}.\,>>

\medskip\noindent{\bf Mathématique du principe de raison.} 
\label{mathematique-du-principe-de-raison}
Or du
point de vue des mathématiques, il serait nul et non avenu de
discuter du bien fondé de ce principe à caractère
métaphysique\,\,---\,\,analysé par ailleurs et contesté d'un point de
vue purement philosophique par un ouvrage éponyme de
Heidegger\,\,---\,\,puisque la pensée mathématique est le domaine par
excellence où la Raison se donne les moyens d'élaborer et
d'approfondir les arguments les plus indubitables quant à des chaînes
finies de raisons qui concernent le nombre, l'étendue, l'infini. Le
principe de raison suffisante affirme donc en tant que principe moteur
pour l'activité mathématique que l'on doit constamment croire qu'il y
a des liens de genèse entre des propositions vraies et établies par
ceraines voies, fussent-elles imparfaites et transitoires, et
l'existence d'autres propositions et d'autres vérités mathématiques
établies par d'autres voies qui sont posées dans l'attente d'une
réalisation ou d'une découverte: 
{\em tout simplement}!

\medskip\noindent{\bf Exemple paradigmatique: la formula egregia.} 
\label{exemple-paradigmatique-formula-egregia}
La recherche par Gauss et l'élaboration sur un temps long de la {\em
formula egregia} constitue alors un exemple paradigmatique de l'action
du principe de raison suffisante dans la genèse des contenus
mathématiques, puisque c'est bien parce qu'une première démonstration
{\em extrinsèque} du caractère intrinsèque
et invariant par isométrie de la courbure
lui était connue (\deutschplain{Schöne Theorem}, 1816), que Gauss a
pensé qu'il devait y avoir une {\em démonstration supérieure et purement
intrinsèque qui établît au mieux cette invariance}.

\smallskip
\centerline{\fbox{\em \ Le principe de raison
mobilise constamment la création de réalité mathématique.}}\smallskip

\Section{11.~Caract\'erisation diff\'erentielle des surfaces
de courbure nulle}
\label{caracterisation-surfaces-de-courbure-nulle}

\HEAD{11.~Caract\'erisation diff\'erentielle des surfaces
de courbure nulle}{
Jo\"el Merker, D\'partement de Math\'matiques d'Orsay}

\medskip\noindent{\bf \'Equivalence à la métrique pythagoricienne.}
\label{equivalence-metrique-pythagoricienne}
\`A la fin de l'ann\'ee 1826, toujours en qu\^ete d'une formule qui
eût exprimé la mesure de courbure $\kappa$ explicitement en fonction
des trois coefficients métriques $E$, $F$, $G$, Gauss con\c cut
l'id\'ee de caract\'eriser les surfaces localement isom\'etriques \`a
un plan, c'est-\`a-dire celles pour lesquelles il existe un changement
de coordonn\'ees:
\[
(u,v)
\longmapsto 
\big(\overline{u},\,\overline{v}\big)
=
\big(\overline{u}(u,v),\,\overline{v}(u,v)\big)
\]
qui transforme leur m\'etrique en la m\'etrique euclidienne standard;
cela revient à supposer que l'on a:
\[
d\overline{u}^2+d\overline{v}^2 
=
E\,du^2+2F\,du dv+G\,dv^2,
\]
apr\`es remplacement des nouvelles variables $\overline{ u}$ et
$\overline{v}$ dans le membre de gauche par leurs expressions en
fonction de $u$ et $v$. Si cette condition est satisfaite, on dira que
la surface est {\sl localement isométrique au plan euclidien}. Cette
derni\`ere \'etape fut cruciale dans la gen\`ese de la {\em computatio
egregia}, puisque Gauss d\'ecouvrit \`a cette occasion le {\em
num\'erateur}\, de l'expression de la courbure donn\'ee par la {\em
formula egregia}~\thetag{ \ref{formula-egregia}}
p.~\pageref{formula-egregia}. 

\medskip\noindent{\bf Numérateur de la formula egregia.}
\label{numerateur-formula-egregia}
Autrement dit, dans des notes
personnelles non publi\'ees, Gauss 
parvint à \'etablir la caractérisation différentielle
suivante:

\smallskip\noindent{\bf Théorème.}
{\em Une surface paramétrée 
en coordonnées $(u, v)$ et 
munie d'une m\'etrique infinit\'esimale: 
\[
ds^2
=
E\,du^2+2F\,du dv+G\,dv^2
\]
est localement euclidienne {\em si et seulement si} les trois
coefficients $E$, $F$, $G$ satisfont l'\'equation diff\'erentielle non
lin\'eaire du second ordre suivante:}
\[
\aligned
0
&
= 
E\big[
E_v\,G_v-2\,F_u\,G_v+G_u\,G_u
\big]
+
\\
&\ \ \ \ \
+
F\big[
E_u\,G_v-E_v\,G_u-2\,E_v\,F_v
+
4\,F_u\,F_v-2\,F_u\,G_u
\big]
+
\\
&\ \ \ \ \
+
G\big[
E_u\,G_u-2\,E_u\,F_v+E_v\,E_v
\big]
+
\\
&\ \ \ \ \
+
2\,(EG-F^2)
\big[
-E_{vv}+2\,F_{uv}-G_{uu}
\big].
\endaligned
\]

\`A un changement de notation pr\`es, on reconnaît effectivement
le num\'erateur de la {\em formula egregia}, co\"{\i}ncidence
remarquable, int\'eressante et avantageuse. 
Par un heureux effet d'économie, 
les calculs qui ont conduit Gauss
\`a la d\'emonstration de ce th\'eor\`eme sont notablement plus
ais\'es et plus courts que les calculs de la {\sl computation egregia}
de 1827. D'après Kreyszig (\cite{ kr1994}, p.~106), Gauss aurait alors pu
s'imaginer (avant de l'avoir vraiment établie) que la {\em formula
egregia}\, ne diff\'erait de cette expression complexe que par un
facteur inessentiel, à savoir le dénominateur $(EG-F^2)^2$. C'est donc
gr\^ace \`a un tel travail pr\'eliminaire que Gauss aurait \'et\'e en
mesure d'anticiper et de r\'ealiser d\'efinitivement le calcul qui le
conduisit \`a la {\em formula egregia}.

\medskip\noindent{\bf Factorisation par complexification.}
\label{factorisation-par-complexification}
Soit donc une métrique gaussienne quelconque:
\[
ds^2 
= 
E\,du^2+2F\,dudv+G\,dv^2
\]
telle qu'il existe une transformation $(u, v)
\longmapsto (\overline{ u}, \overline{ v})$ qui la transforme en la
métrique pythagoricienne standard 
$d \overline{ u}^2 + d \overline{ v}^2$. Il
s'agit d'établir que $E$, $F$ et $G$ satisfont le systèmes d'équations
aux dérivées partielles du premier ordre exhibé;
la réciproque, simple application
du théorème de Clebsch-Frobenius, sera laissée de côté.

Si l'on s'autorise à utiliser le nombre imaginaire $i := \sqrt{ -1}$,
on peut d\'ecomposer le carr\'e $(d \overline{ u})^2 +
(d \overline{ v})^2$ en facteurs linéaires:
\[
(d\overline{u})^2+(d\overline{v})^2 
= 
(d\overline{u}+ id\overline{v})(d\overline{u} - id\overline{v}) 
=: 
d\lambda \, d\mu,
\]
en posant simplement: 
\[
\lambda
:= 
\overline{u}+i\overline{v}
\ \ \ \ \
\text{\rm et} 
\ \ \ \ \ 
\mu 
:= 
\overline{u}-i\overline{v}. 
\]
Le fait que la fonction $\mu = \mu ( u,v)$ à valeurs complexes soit la
conjuguée de la fonction $\lambda = \lambda ( u, v)$ ne sera pas
utilisé dans les raisonnements qui vont suivre, et à partir de
maintenant, on va {\em oublier} les deux fonctions réelles initiales
$\overline{ u} ( u, v)$ et $\overline{ v} ( u, v)$ pour ne considérer
que $\lambda$ et $\mu$ en tant que {\em nouvelles} fonctions
significatives de $(u, v)$.

Soient ensuite les deux différentielles: 
\[
d\lambda 
=
\lambda_u\,du+\lambda_v\,dv 
\ \ \ \ \ \ \
\text{\rm et} 
\ \ \ \ \ \ \
d\mu
=
\mu_u\,du+\mu_v\,dv 
\]
de ces deux
nouvelles fonctions de $(u, v)$. Dans la nouvelle \'equation qui
exprime la condition d'isométrie infinitésimale à la métrique
pythagoricienne:
\[
\aligned
E\,du^2 
+ 
2F\,dudv 
+ 
G\,dv^2
&
=
d\overline{u}^2
+
d\overline{v}^2
=
\\
&
=
d\lambda\,d\mu
=
\big(\lambda_u\,du+\lambda_v\,dv\big) 
\big(\mu_u\,du+\mu_v\,dv\big),
\endaligned
\] 
si l'on identifie alors les trois coefficients de $du^2$, de $du dv$
et de $dv^2$ entre les deux membres situés
aux extrémités gauche et droite, on
obtient que $E$, $F$ et $G$ doivent nécessairement être
liés aux dérivées partielles du premier ordre de ces fonctions par
les formules quadratiques suivantes:
\[
\aligned
E 
& \
= 
\lambda_u \, \mu_u, \\
2F 
& \
=
\lambda_u \, \mu_v+ \lambda_v \, \mu_u, \\
G 
& \
=
\lambda_v \, \mu_v.
\endaligned
\]

\medskip\noindent{\bf Relations différentielles systématiques.}
\label{relations-differentielles-systematiques}
C'est ici que la combinatoire diff\'erentielle entre en sc\`ene. De
la sorte, deux fonctions $\lambda$ et $\mu$ imposent des conditions
aux trois fonctions $E$, $F$ et $G$. Puisqu'il y a quatre d\'eriv\'ees
partielles d'ordre $1$, à savoir: $\lambda_u$, $\lambda_v$, $\mu_u$ et
$\mu_v$, et seulement trois \'equations, ce syst\`eme ne permet pas de
d\'eterminer immédiatement les inconnues $\lambda$ et $\mu$. C'est
pourquoi il faut diff\'erentier par rapport \`a $u$ et par rapport \`a
$v$ chacune de ces trois
\'equations, ce qui donne six \'equations aux d\'eriv\'ees
partielles écrites simplement \`a la suite les unes des autres:
\begin{equation}
\label{syst-lin-E-F-G-lambda-mu}
\aligned
E_u 
& \
=
\lambda_{uu} \, \mu_u + \lambda_u \, \mu_{uu}, \\
E_v 
& \
= 
\lambda_{uv} \, \mu_u + \lambda_u\, \mu_{uv}, \\
2 F_u 
& \
=
\lambda_{ uu} \, \mu_v + 
\lambda_u \, \mu_{vu} +
\lambda_{vu} \, \mu_{ u} + 
\lambda_v\, \mu_{uu}, \\
2F_v
& \
=
\lambda_{uv} \, \mu_v + 
\lambda_u \, \mu_{ vv} + 
\lambda_{vv} \, \mu_u + 
\lambda_v \, \mu_{ uv}, \\
G_u 
& \
= 
\lambda_{vu} \,\mu_v +
\lambda_v \, \mu_{vu}, \\
G_v 
& \
=
\lambda_{vv} \, \mu_v + 
\lambda_v \, \mu_{vv}.
\endaligned
\end{equation}
Dans ce syst\`eme de six \'equations, il importe d'observer
qu'il y a autant de d\'eriv\'ees partielles d'ordre deux des deux
fonctions inconnues $\lambda$ et $\mu$ que de d\'eriv\'ees partielles
d'ordre un des fonctions donn\'ees $E$, $F$ et $G$. En effet, ce
nombre commun est \'egal \`a six:
\[
\aligned
\lambda_{uu}, \ \ & \
\lambda_{uv}, \ \ & \
\lambda_{vv}, \ \ & \
\mu_{uu}, \ \ & \
\mu_{uv}, \ \ & \
\mu_{vv}, \\
E_u, \ \ & \
E_v, \ \ & \
F_u, \ \ & \
F_v, \ \ & \
G_u, \ \ & \
G_v,
\endaligned
\]
en tenant bien s\^ur compte de la commutativit\'e des d\'eriv\'ees
partielles, qui assure que $\lambda_{uv} = \lambda_{ vu}$ et que
$\mu_{uv} = \mu_{vu}$. On peut donc s'attendre \`a ce que le syst\`eme
lin\'eaire~\thetag{ \ref{syst-lin-E-F-G-lambda-mu}} ci-dessus, 
o\`u l'on envisage les quatre d\'eriv\'ees
partielles du premier ordre $\lambda_u$, $\lambda_v$, $\mu_u$ et
$\mu_v$ comme des coefficients, puisse \^etre invers\'e de telle sorte
que les six d\'eriv\'ees partielles $\lambda_{ uu}$, $\lambda_{ uv}$,
$\lambda_{ vv}$, $\mu_{ uu}$, $\mu_{ uv}$ et $\mu_{vv}$ s'expriment en
fonction des six d\'eriv\'ees partielles $E_u$, $E_v$, $F_u$, $F_v$,
$G_u$ et $G_v$. Cette intuition va se confirmer.
 
Pour commencer, on observe qu'il y a un certain 
d\'es\'equilibre\footnote{\,
%%%%%%%%%%%%%%%%%%%%%%%-------DEBUT--------%%%%%%%%%%%%%%%%%%%%%%%%%%%
Seuls des moments intermédiaires d'examen interrogatif
peuvent déclencher de nouveaux actes significatifs de
calcul.
} %%%%%%%%%%%%%%%%%%%%%%%%-----FIN-----%%%%%%%%%%%%%%%%%%%%%%%%%%%%%%%
dans ce système~\thetag{ 
\ref{syst-lin-E-F-G-lambda-mu}}: 
les troisi\`eme et quatri\`eme \'equations
comportent plus de termes que les premi\`ere, deuxi\`eme, cinqui\`eme
et sixi\`eme. Pour r\'etablir l'\'equilibre, il suffit de retrancher
la deuxi\`eme \'equation de la troisi\`eme et de retrancher la
cinqui\`eme \'equation de la quatri\`eme, ce qui donne, en
recopiant simultanément les quatre
\'equations non modifi\'ees:
\begin{equation}
\label{4-non-modifiees}
\aligned
E_u 
& \
=
\lambda_{uu} \, \mu_u + \lambda_u \, \mu_{uu}, \\
E_v 
& \
= 
\lambda_{uv} \, \mu_u + \lambda_u\, \mu_{uv}, \\
2 F_u -E_v 
& \
=
\lambda_{ uu} \, \mu_v + 
\lambda_v\, \mu_{uu}, \\
2F_v - G_u
& \
=
\lambda_u \, \mu_{ vv} + 
\lambda_{vv} \, \mu_u, \\
G_u 
& \
= 
\lambda_{vu} \,\mu_v +
\lambda_v \, \mu_{vu}, \\
G_v 
& \
=
\lambda_{vv} \, \mu_v + 
\lambda_v \, \mu_{vv}.
\endaligned
\end{equation}
Maintenant, il devient particuli\`erement ais\'e d'inverser ce nouveau
syst\`eme lin\'eaire. En effet, on voit qu'il se d\'ecompose en
trois sous-syst\`emes lin\'eaires \`a deux inconnues:

\begin{itemize}

\smallskip\item[{\bf (a)}]
la premi\`ere et la troisi\`eme \'equations, 
aux deux inconnues $\lambda_{ uu}$ et 
$\mu_{ uu}$;

\smallskip\item[{\bf (b)}]
la deuxi\`eme et la cinqui\`eme \'equations, aux
deux inconnues $\lambda_{ uv}$ et $\mu_{uv}$;

\smallskip\item[{\bf (c)}]
la quatri\`eme et la sixi\`eme
\'equations, aux deux inconnues
$\lambda_{ vv}$ et $\mu_{ vv}$.

\end{itemize}\smallskip

De plus, 
chacun de ces trois syst\`emes de deux \'equations
lin\'eaires \`a deux inconnues fait appara\^{\i}tre
le d\'eterminant
$2\times 2$
\[
\left\vert
\begin{array}{cc}
\lambda_u & \mu_u \\
\lambda_v & \mu_v 
\end{array}
\right\vert = 
\lambda_u \, \mu_v - \mu_u \, \lambda_v,
\]
car chacun d'entre eux s'\'ecrit sous la forme 
\[
\aligned
A 
& \
=
\lambda_u\, X + \mu_u \, Y, \\
B
& \
=
\lambda_v\, X + \mu_v \, Y,
\endaligned
\]
avec des constantes donn\'ees $A$, $B$ et des inconnues $X$, $Y$ qui
sont
différentes dans les trois cas. Or il existe des formules alg\'ebriques
\'el\'ementaires pour r\'esoudre de tels syst\`emes:
\[
\aligned
X 
& \
= 
\frac{\mu_v \, A - \mu_u \, B}{ 
\lambda_u \, \mu_v - \mu_u \, \lambda_v}, \\
Y
& \
= 
\frac{\lambda_u \, B - \lambda_v \, A}{ 
\lambda_u \, \mu_v - \mu_u \, \lambda_v},
\endaligned
\]
valables lorsque le d\'eterminant au dénominateur ne s'annule
pas. Heureusement par hypothèse, les diff\'erentielles de $\lambda$ et
de $\mu$ sont lin\'eairement ind\'ependantes en tout point de
coordonn\'ees $(u,\, v)$, d'o\`u il d\'ecoule que le d\'eteminant en
question ne s'annule en aucun point. On en d\'eduit imm\'ediatement
les formules de résolution pour les dérivées partielles du second
ordre:
\[
\footnotesize
\aligned
\lambda_{uu} 
& \
= 
\frac{1}{\lambda_u \mu_v - 
\mu_u \, \lambda_v} 
\left[
2\lambda_u \, F_u - 
\lambda_u \, E_v -
\lambda_v \, E_u
\right], \\
\mu_{uu} 
& \
= 
\frac{1}{\lambda_u \mu_v - 
\mu_u \, \lambda_v} 
\left[
\mu_v\, E_u - 2\mu_u \, F_u + 
\mu_u\, E_v 
\right], \\
\lambda_{uv} 
& \
= 
\frac{1}{\lambda_u \mu_v - 
\mu_u \, \lambda_v} 
\left[
\lambda_u \, G_u - \lambda_v \, E_v
\right], 
\endaligned
\]
\[
\footnotesize
\aligned
\mu_{uv} 
& \
= 
\frac{1}{\lambda_u \mu_v - 
\mu_u \, \lambda_v} 
\left[
\mu_v \, E_v - \mu_u \, G_u
\right], \\
\lambda_{vv} 
& \
= 
\frac{1}{\lambda_u \mu_v - 
\mu_u \, \lambda_v} 
\left[
\lambda_u \, G_v - 2 \lambda_v \, F_v + 
\lambda_v \, G_u
\right], \\
\mu_{vv} 
& \
= 
\frac{1}{\lambda_u \mu_v - 
\mu_u \, \lambda_v} 
\left[
2\mu_v \, F_v - \mu_v\, G_u - 
\mu_u \, G_v
\right].
\endaligned
\]
Ainsi, les six d\'eriv\'ees partielles d'ordre deux des fonctions
$\lambda$ et $\mu$ s'expriment-elles en fonction de leurs d\'eriv\'ees
partielles d'ordre $1$, et en fonction aussi des d\'eriv\'ees partielles
d'ordre $1$ des coefficients m\'etriques $E$, $F$, $G$.

\`A pr\'esent, on peut v\'erifier que l'\'equation
énoncée par le théorème est {\em
n\'ecessaire}\, pour l'existence de
ces deux fonctions $\lambda$ et $\mu$.

Pour cela,
on commence par rediff\'erentier les
\'equations~\thetag{ 
\ref{4-non-modifiees}}\,\,---\,\,m\'ecaniquement, une fois 
encore\,\,---\,\,par
rapport \`a $u$ et \`a $v$, ce qui donne dix
\'equations\footnote{\,
%%%%%%%%%%%%%%%%%%%%%%%-------DEBUT--------%%%%%%%%%%%%%%%%%%%%%%%%%%%
On v\'erifie que neuf d'entre elles seulement
sont lin\'eairement ind\'ependantes, 
exactement le même nombre, neuf, des 
d\'eriv\'ees
partielles d'ordre deux des trois fonctions $E$, $F$ et $G$, {\em
i.e.} $E_{uu}$, $E_{uv}$, $E_{vv}$, $F_{uu}$, $F_{uv}$, $F_{ vv}$,
$G_{ uu}$, $G_{uv}$ et $G_{vv}$. En effet, les \'equations num\'ero
$3$, $5$, $6$ et $8$ sont \'evidemment lin\'eairement
d\'ependantes, car la combinaison lin\'eaire $E_{ vv} + (2F_{uv} - E_{
vv}) - (2F_{ uv} - G_{uu}) - G_{ uu}$ donne trivialement un r\'esultat
nul.
}: %%%%%%%%%%%%%%%%%%%%%%%%-----FIN-----%%%%%%%%%%%%%%%%%%%%%%%%%%%%%%%
\[
\small
\aligned
\ \ \ \ \ \ \ \ \ \ \ \ \ 
E_{uu} 
& \
=
\lambda_{uuu} \, \mu_u + 2\lambda_{uu} \, \mu_{uu}
+ \lambda_u \, \mu_{ uuu}, 
\\
E_{uv}
& \
= 
\lambda_{uuv} \, \mu_u + \lambda_{uu} \, \mu_{uv}
+ \lambda_{uv} \, \mu_{ uu}+ \lambda_u \, \mu_{ uuv}, \\
E_{vv} 
& \
=
\lambda_{ uvv}\, \mu_u + 2\lambda_{ uv}\, \mu_{uv} +
\lambda_u \, \mu_{uvv}, 
\endaligned
\]
\[
\small
\aligned
2 F_{uu} -E_{uv} 
& \
=
\lambda_{ uuu} \, \mu_{v} + 
\lambda_{ uu} \, \mu_{uv} + 
\lambda_{uv}\, \mu_{uu} +
\lambda_{v}\, \mu_{uuu}, \\
2 F_{uv} -E_{vv} 
& \
=
\lambda_{ uuv} \, \mu_v +
\lambda_{ uu} \, \mu_{vv} + 
\lambda_{vv}\, \mu_{uu} +
\lambda_v \, \mu_{ uuv}, \\
2F_{uv} - G_{uu}
& \
=
\lambda_{uu} \, \mu_{ vv} + 
\lambda_{u} \, \mu_{ uvv} + 
\lambda_{uvv} \, \mu_u + 
\lambda_{vv} \, \mu_{uu}, \\
2F_{vv} - G_{uv}
& \
=
\lambda_{uv} \, \mu_{ vv} + 
\lambda_{u} \, \mu_{ vvv} + 
\lambda_{vvv} \, \mu_u+ 
\lambda_{vv} \, \mu_{uv},
\endaligned
\]
\[
\small
\aligned
\ \ \ \ \ \ \ \ \ \ \ \ 
G_{uu}
& \
= 
\lambda_{uuv} \,\mu_v + 2
\lambda_{uv}\, \mu_{uv} + \lambda_v\, \mu_{uuv}, 
\\
G_{uv} 
& \
=
\lambda_{uvv} \, \mu_v + 
\lambda_{vv} \, \mu_{uv} +
\lambda_{uv}\, \mu_{vv} +
\lambda_v \,\mu_{ uvv}, 
\\
G_{vv}
& \
=
\lambda_{vvv} \, \mu_v + 
2\lambda_{ vv} \, \mu_{ vv}+ 
\lambda_v \, \mu_{vvv}.
\endaligned
\]
Ce syst\`eme para\^{\i}t complexe au premier abord, mais pour y
discerner des structures, il suffit de se concentrer\footnote{\,
%%%%%%%%%%%%%%%%%%%%%%%-------DEBUT--------%%%%%%%%%%%%%%%%%%%%%%%%%%%
La phénoménologie étudie le `donné commun de la perception' pour les
structures transcendantales de l'intentionnalité et de la conscience.
La {\em lecture-en-recherche} des formules mathématiques ouvre vers
d'autres caractères constitués de la donation à la pensée.
} %%%%%%%%%%%%%%%%%%%%%%%%-----FIN-----%%%%%%%%%%%%%%%%%%%%%%%%%%%%%%%
seulement sur les d\'eriv\'ees partielles d'ordre trois des fonctions
$\lambda$ et $\mu$, en n\'egligeant\,\,---\,\,au moins pour
l'instant\,\,---\,\,la lecture des autres termes qui n'incorporent que
des dérivées d'ordre deux. On utilisera alors le signe
<<\,$\equiv$\,>> pour signifier que les d\'eriv\'ees partielles
d'ordre deux sont supprimées à la lecture:
\[
\small
\aligned
\ \ \ \ \ \ \ \ \ \ \ \ \ \ 
E_{uu} 
& \
\equiv
\lambda_{uuu} \, \mu_u
+ \lambda_u \, \mu_{ uuu}, 
\\
E_{uv}
& \
\equiv
\lambda_{uuv} \, \mu_u + \lambda_u \, \mu_{ uuv}, 
\\
E_{vv} 
& \
\equiv
\lambda_{ uvv}\, \mu_u +
\lambda_u \, \mu_{uvv},
\endaligned
\]
\[
\small
\aligned
2 F_{uu} -E_{uv} 
& \
\equiv
\lambda_{ uuu} \, \mu_{v} + 
\lambda_{v}\, \mu_{uuu}, 
\\
2 F_{uv} -E_{vv} 
& \
\equiv
\lambda_{ uuv} \, \mu_v +
\lambda_v \, \mu_{ uuv}, 
\\
2F_{uv} - G_{uu}
& \
\equiv
\lambda_{u} \, \mu_{ uvv} + 
\lambda_{uvv} \, \mu_u, 
\\
2F_{vv} - G_{uv}
& \
\equiv
\lambda_{u} \, \mu_{ vvv} + 
\lambda_{vvv} \, \mu_u,
\endaligned
\]
\[
\small
\aligned
\ \ \ \ \ \ \ \ \ \ \ \ \,
G_{uu}
& \
\equiv 
\lambda_{uuv} \,\mu_v + \lambda_v\, \mu_{uuv}, 
\\
G_{uv} 
& \
\equiv
\lambda_{uvv} \, \mu_v + 
\lambda_v \,\mu_{ uvv}, 
\\
G_{vv}
& \
\equiv
\lambda_{vvv} \, \mu_v + 
\lambda_v \, \mu_{vvv}.
\endaligned
\]
Alors, après quelque examen attentif, l'{\oe}il avisé du calculateur
verra appara\^{\i}tre une relation lin\'eaire non triviale: en
additionnant la cinqui\`eme et la sixi\`eme \'equations et en
soustrayant la troisi\`eme et la huiti\`eme \'equations, on constate
que:
\[
-\, 
2\,E_{vv} 
+ 
4\,F_{uv} 
-
2\,G_{uu} 
\equiv 
0.
\]
Ainsi cette combinaison lin\'eaire $-2 ( E_{ vv} - 2 F_{uv} + G_{
uu})$ de d\'eriv\'ees partielles d'ordre deux s'annule modulo les
d\'eriv\'ees partielles d'ordre inf\'erieur. Au facteur $-2$ pr\`es,
c'est la m\^eme combinaison lin\'eaire qui intervient dans la
derni\`ere ligne du th\'eor\`eme.

\medskip\noindent{\bf Complétion de la combinaison linéaire 
caractéristique.} \label{completion-combinaison-lineaire-caracteristique}
Il faut évidemment effectuer maintenant la m\^eme combinaison
lin\'eaire avec les \'equations compl\`etes, c'est-\`a-dire:
additionner la cinqui\`eme et la sixi\`eme \'equations, soustraire la
troisi\`eme et la huiti\`eme \'equations, et diviser le tout par $-2$,
ce qui donne:
\[
E_{ vv} -2 F_{ uv} + G_{uu} 
=
2\lambda_{ uv} \, \mu_{uv} - 
\lambda_{ uu} \, \mu_{ vv} - 
\lambda_{vv} \, \mu_{ uu}.
\]
Dans le membre de droite de cette \'equation, on doit à présent
remplacer les valeurs de $\lambda_{ uv}$, de $\mu_{ uv}$, de
$\lambda_{ uu}$, de $\mu_{ vv}$, de $\lambda_{ vv}$ et de $\mu_{ uu}$
trouv\'ees plus haut et d\'evelopper les produits:
\[
\small
\aligned
E_{vv}
-
2\,F_{uv}
+
G_{uu}
&
=
2\lambda_{uv}\,\mu_{uv}
-
\lambda_{uu}\,\mu_{vv}
-
\lambda_{vv}\,\mu_{uu}
\\
&
=
\frac{1}{[\lambda_u\,\mu_v-\mu_u\,\lambda_v]^2}
\big\{
2(\lambda_u\,G_u-\lambda_v\,E_v)
(\mu_v\,E_v-\mu_u\,G_u)
+
\\
&
\ \ \ \ \ 
+
(-2\lambda_u\,F_u+\lambda_u\,E_v+\lambda_v\,E_u)
(2\mu_v\,F_v-\mu_v\,G_u-\mu_u\,G_v)
+
\\
&
\ \ \ \ \
+
(-\lambda_u\,G_v+2\lambda_v\,F_v-\lambda_v\,G_u)
(\mu_v\,E_u-2\mu_u\,F_u+\mu_u\,E_v)
\big\}
\endaligned
\]
\[
\small
\aligned
\ \ \ \ \ \ \ \ \ \ \ \ \ \ \ \ \ \ \ \ \ \ \ \ \ \ \ \ \ \ \ \ \ 
&
=
\frac{1}{[\lambda_u\,\mu_v-\mu_u\,\lambda_v]^2}
\big\{
E_v\,G_u(\lambda_u\,\mu_v+\lambda_v\,\mu_u)
-
\\
&
\ \ \ \ \
-\,E_v\,E_v
(2\lambda_v\,\mu_v)-
G_u\,G_u(2\lambda_u\,\mu_u)-
\\
&
\ \ \ \ \
-
F_u\,F_v
(4\lambda_u\,\mu_v+4\lambda_v\,\mu_u)
+
F_u\,G_u
(2\lambda_u\,\mu_v+2\lambda_v\,\mu_u)
+
\\
&
\ \ \ \ \
+
F_u\,G_v(4\lambda_u\,\mu_u)
+
E_v\,F_v
(2\lambda_u\,\mu_v+2\lambda_v\,\mu_u)
-
\\
&
\ \ \ \ \
-
E_v\,G_v
(2\lambda_u\,\mu_u)
+
E_u\,F_v(4\lambda_v\,\mu_v)
-
\\
&
\ \ \ \ \
-
E_u\,G_u
(2\lambda_v\,\mu_v)
-
E_u\,G_v
(\lambda_u\,\mu_v+\lambda_v\,\mu_u)
\big\}.
\endaligned
\]
Pour achever d'\'etablir l'\'equation
du théorème, on observe que le
d\'enominateur ci-dessus s'exprime en fonction de $E$, $F$, $G$:
\[
\aligned
\big[
\lambda_u\,\mu_v-\mu_u\,\lambda_v\big]^2
&\
=
\lambda_u^2\,\mu_v^2-2\lambda_u\,\lambda_v\,\mu_u\,\mu_v+
\mu_u^2\,\lambda_v^2
\\
&\
=
-4\big[E\,G-F^2\big].
\endaligned
\]
En d\'efinitive, on obtient l'\'equation annoncée sous la forme:
\[
\aligned
E_{vv}
-
2\,F_{uv}
+
G_{uu}
=
\frac{1}{2\left[E\,G-F^2\right]}
\big\{
&
-E_v\,G_u(F)
+
\\
&
+E_v\,E_v(G)+G_u\,G_u(E)
+
\\
&
+F_u\,F_v(4F)-F_u\,G_u(2F) 
-
\\
&
- 
F_u\,G_v(2E)
- 
E_v\,F_v(2F) 
+
\\
&
+
E_v\,G_v(E)
-
E_u\,F_v(2G) 
+ 
\\
&
+ 
E_u\,G_u(G)
+ 
E_u\,G_v(F)
\big\},
\endaligned
\]
\'equivalente \`a celle annoncée par le théorème.
\qed

\Section{12.~Commentaire de la {\em computatio egregia}}
\label{commentaire-computatio-egregia}

\HEAD{12.~Commentaire de la {\em computatio egregia}}{
Jo\"el Merker, D\'partement de Math\'matiques d'Orsay}

\medskip\noindent{\bf Systématicité et complétude.}
\label{systematicite-et-completude}
Pour calculer la mesure de courbure d'une surface quelconque, Gauss
procède\,\,---\,\,on l'a vu\,\,---\,\,d'une manière systématique et
complète, traitant simultanément les trois représentations possibles
d'une surface: la représentation implicite, la représentation
paramétrée et la représentation graphée. Grâce à la totalisation
qu'offre une telle approche, il se donne la possibilité de {\em
circuler mathématiquement} d'une représentation à une autre, toujours
par les voies du calcul, 
à la fois en complexité, et en {\em explicite-té}.

\Scholie
En lui-même, le Calcul est non local, non sériel, non temporel, mais
il est un réseau de liens multiples qui expriment {\em tous} les chemins
qui peuvent être actualisables au sein de son graphe interne.
\qed

\smallskip
Gauss pense donc les mathématiques comme un tout de calcul que
l'esprit, enchaîné à la temporalité progressive des actes, 
doit\,\,---\,\,afin d'en synthétiser la totalité\,\,---\,\,parcourir 
et explorer jusqu'à en avoir compris (presque) tous les
délinéaments co-présents. La
charpente des {\em Disquisitiones} laisse en effet deviner une telle
intention, confirmée par l'élégance formelle des expressions
symboliques.

Spécialement pour l'exposition de la {\em computatio egregia}, la
disposition préliminaire des éléments géométriques va faire bénéficier
de la totalisation en parallèle qui avait été anticipée dans les
premiers paragraphes. Un lien constant de {\em traduction
fonctionnelle} depuis les formules graphées vers les formules
paramétrées va être mis en {\oe}uvre, d'où l'importance d'avoir
calculé d'abord la courbure $\kappa = \frac{ z_{xx} \, z_{ yy} - z_{
xy}^2}{ (1+ z_x^2 + z_y^2)^2}$ dans le cas où $z$ est une certaine
fonction de $x$ et de $y$.

Au contraire, la présentation immédiatement intrinsèque des tenseurs
de courbure que l'on apprend et que l'on reproduit dans les traités
modernes de géométrie différentielle {\em perd en complétude
d'explicitation} ce que Gauss mettait au fondement d'une saine et
adéquate {\em compréhesion-en-totalité} des concepts. C'est pourquoi,
dans ce commentaire\,\,---\,\,mathématique et spéculatif\,\,---\,\,de
la {\em computatio egregia}, nulle référence ne sera faite aux
concepts qui sont apparus ultérieurement pour encadrer le calcul et
pour lui conférer un support géométrique, concepts tels que, par
exemple: {\sl dérivées première et seconde du rayon-vecteur} ${\bf x}
(u, v) \in \R^3$; {\sl formules de produit scalaire}; {\sl première
forme fondamentale}; {\sl application de Weingarten}; {\sl
coefficients de Christoffel}; {\sl transport parallèle}; {\sl
connexion de Levi-Cività}; {\em etc.} 

\smallskip

L'objectif réflexif et métaphysique, ici, est surtout
de démontrer par l'Exemple 
en quoi les tensions du Calcul impliquent une {\em
morphogénèse} du réel mathématique qui est 
{\em antérieure à toute formation
de concepts}.

\medskip\noindent{\bf Composantes du vecteur normal.} 
\label{composantes-du-vecteur-normal}
En représentation graphée $z = z ( x, y)$, le champ de vecteurs ${\bf n}
\colon S \to T_\R^3 \vert_S$ normal à la surface possède trois composantes
$(X, Y, Z)$ qui sont des fonctions de $(x, y)$ satisfaisant:
\[
0
=
X\,dx+Y\,dy+Z\,dz,
\]
lors de tout déplacement infinitésimal $(dx, dy, dz)$ sur la surface,
où $dz = z_x\, dx + z_y\, dy$ s'exprime en fonction de $dx$ et
de $dy$.

\CITATION{
Dans la {\em seconde} méthode, on exprime les coordonnées en forme de
fonctions de deux variables $p$ et $q$. Supposons que par la
différentiation de ces fonctions, on ait:
\[
\aligned
dx
&
=
a\,dp+a'\,dq,
\\
dy
&
=
b\,dp+b'\,dq,
\\
dz
&
=
c\,dp+c'\,dq;
\endaligned
\]
par la substitution de ces valeurs dans une formule donnée
ci-dessus, nous obtenons:
\[
(a\,X+b\,Y+c\,Z)\,dp
+
(a'\,X+b'\,Y+c'\,Z)\,dq
=
0.
\vspace{-0.65cm}
\]
\REFERENCE{\cite{ga1828},~11}}

Ici donc, les trois composantes du vecteur normal
en représentation paramétrique:
\[
(p,q)
\longmapsto
\big(x(p,q),\,y(p,q),\,z(p,q)\big)
=:
{\bf x}(p,q)
\]
ne sont {\em pas} calculées en partant de la condition d'orthogonalité
de ${\bf n} (p, q)$ avec les deux vecteurs tangents indépendants ${\bf
x}_p ( p, q)$ et ${\bf x}_q ( p, q)$, mais elles sont calculées en se
ramenant au cas\,\,---\,\,exposé à l'instant par Gauss\,\,---\,\,d'un
graphe. {\em Une seule et même surface est vue sous trois angles
abstraitement équivalents mais calculatoirement différents}. La
nécessité interne des mathématiques commande en effet de tracer à tout
instant des voies de parcours pour passer d'une représentation à une
autre.

Observation sur les notations: la majeure partie des calculs de
courbure repose sur de l'algèbre différentielle dans des espaces de
jets de fonctions à plusieurs variables, laquelle algèbre
différentielle se ramène à de la pure algèbre commutative lorsqu'on
envisage les dérivées partielles impliquées comme de nouvelles
variables indépendantes. Par conséquent, la raison pour laquelle les
{\em Disquisitiones} introduisent systématiquement de nouvelles
lettres en lieu et place des dérivées partielles, telles que,
ici-même:
\[
a:=\frac{dx}{dp},\ \ 
a':=\frac{dx}{dq},\ \ \ \ \
b:=\frac{dy}{dp},\ \ 
b':=\frac{dy}{dq},\ \ \ \ \
c:=\frac{dz}{dp},\ \ 
c':=\frac{dz}{dq},
\]
ne tient donc pas seulement au fait qu'un gain en 
découle en termes d'abréviation, mais aussi à 
une réelle {\em conscience 
a posteriori} de l'essence algébrique des calculs.

\CITATION{
Comme cette équation doit avoir lieu indépendamment 
des valeurs des différentielles $dp$, $dq$, on aura
évidemment:
\[
a\,X+b\,Y+c\,Z=0,
\ \ \ \ \ \ \
a'\,X+b'\,Y+c'\,Z=0;
\]
d'où nous concluons que $X$, $Y$, $Z$ doivent être respectivement
proportionnels aux quantités:
\[
bc'-cb',
\ \ \ \ \
ca'-ac',
\ \ \ \ \
ab'-ba'.
\vspace{-0.5cm}
\]
\REFERENCE{\cite{ga1828},~11}}

La résolution formelle explicite d'un système linéaire de deux
équations indépendantes à trois inconnues doit faire ici partie de la
culture du lecteur. Point de calcul vectoriel, point de combinaisons
linéaires, point de rappel sur la non-annulation d'un déterminant $2
\times 2$, mais plutôt: {\em expression formelle pure} de la solution
(à un facteur multiplicatif près), symétrique, esthétique, rigoureuse.
Rien de plus satisfaisant pour l'esprit que de voir `jaillir' une
formule adéquate, vérité offerte de la chose à saisir, en 
particulier lorsqu'elle est déjà connue par ailleurs.

Pour terminer, il reste seulement à {\em normaliser} à $1$ la norme du
vecteur orthogonal.

\CITATION{
Posant donc, pour abréger:
\[
\sqrt{(bc'-cb')^2+(ca'-ac')^2+(ab'-ba')^2}
=
\Delta,
\]
on aura, soit:
\[
X
=
\frac{bc'-cb'}{\Delta},
\ \ \ \ \ \ \
Y
=
\frac{ca'-ac'}{\Delta},
\ \ \ \ \ \ \
Z
=
\frac{ab'-ba'}{\Delta},
\]
soit:
\[
X
=
\frac{cb'-bc'}{\Delta},
\ \ \ \ \ \ \
Y
=
\frac{ac'-ca'}{\Delta},
\ \ \ \ \ \ \
Z
=
\frac{ba'-ab'}{\Delta}.
\vspace{-0.5cm}
\]
\REFERENCE{\cite{ga1828},~11}}

\medskip\noindent{\bf Transfert à la représentation paramétrique.}
\label{transfert-representation-parametrique}
L'objectif du \S~X des {\em Disquisitiones} est de {\em 
transférer à la représentation paramétrique} 
la formule extrinsèque pour la courbure obtenue
précédemment au \S~VII:
\[
\kappa
=
\frac{TV-U^2}{(1+t^2+u^2)^2}.
\]
Si aucune relation 
d'inter-expression ne sera entreprise en 
liaison avec la formule obtenue
juste auparavant dans le cas d'une représentation implicite:
\[
\aligned
(P^2+Q^2+R^2)\,\kappa
&
=
P^2\,(Q'R'-{P''}^2)
+
Q^2(P'R'-{Q''}^2)
+
\\
&
\ \ \ \ \
+R^2(P'Q'-{R''}^2)
+
2\,QR\,(Q''R''-P'P'')
+
\\
&
\ \ \ \ \
+
2\,PR\,(P''R''-Q'Q'')
+
2\,PQ\,(P''Q''-R'R''),
\endaligned
\]
c'est probablement que Gauss estime inutile de {\em fermer un triangle
de calculs} dont deux côtés auront été explicités. En tout cas, 
il est suffisant pour satisfaire la nécessité de complétude 
d'avoir achevé les calculs dans les trois représentations,
lesquels sont détaillés dans l'ordre progressif
de leur complexité calculatoire.

\CITATION{
Nous obtiendrons une formule encore plus compliquée et renfermant
quinze éléments, si nous voulons suivre la seconde des méthodes
propres à l'étude des surfaces. Il est très-important cependant
d'arriver à cette formule.
\REFERENCE{\cite{ga1828},~22--23}}

\`A nouveau, des lettres algébriques se substitueront
aux dérivées partielles en tant que notations nouvelles:
\[
\aligned
\alpha
&
:=
\frac{d^2x}{dp^2},
\ \ \ \ \ \ \
\alpha'
:=
\frac{d^2x}{dp\,dq},
\ \ \ \ \ \ \
\alpha''
:=
\frac{d^2x}{dq^2},
\\
\beta
&
:=
\frac{d^2y}{dp^2},
\ \ \ \ \ \ \
\beta'
:=
\frac{d^2y}{dp\,dq},
\ \ \ \ \ \ \
\beta''
:=
\frac{d^2y}{dq^2},
\\
\gamma
&
:=
\frac{d^2z}{dp^2},
\ \ \ \ \ \ \
\gamma'
:=
\frac{d^2z}{dp\,dq},
\ \ \ \ \ \ \ \
\gamma''
:=
\frac{d^2z}{dq^2},
\endaligned
\]
propres à se combiner entre elles comme dans un simple anneau de 
polynômes.

\CITATION{
En outre nous ferons pour abréger:
\[
\aligned
bc'-cb'
&
=
A,
\\
ca'-ac'
&
=
B,
\\
ab'-ba'
&
=
C.
\endaligned
\]
Cela posé, nous observons d'abord que l'on a:
\[
A\,dx+B\,dy+C\,dz
=
0,
\]
ou bien:
\[
dz
=
-{\textstyle{\frac{A}{C}}}\,dx
-
{\textstyle{\frac{B}{C}}}\,dy;
\]
en sorte qu'en regardant $z$ comme une fonction de $x$, $y$, on
doit avoir:
\[
\aligned
\frac{dz}{dx}
&
=
t
=
-\frac{A}{C},
\\
\frac{dz}{dy}
&
=
u
=
-\frac{B}{C}.
\endaligned
\]
\REFERENCE{\cite{ga1828},~23}}

Ici, les composantes du champ de vecteurs normal à la surface entrent
nécessairement dans les calculs de différentielles, et elles vont
devoir être différentiées par rapport à $p$ et à $q$. Comme chacune de
leurs expressions incorpore déjà quatre quantités fondamentales, les
calculs ultérieurs, s'ils devaient être conduits de manière
complètement explicites, s'avèreraient prendre des proportions
inquiétantes, ce à cause de quoi l'{\oe}il du mathématicien risquerait
de ne plus être à même de {\em discerner les structures formelles
organisatrices}. L'acte d'introduire une notation temporaire pour une
quantité semi-complexe, {\em i.e.} d'introduire une {\em dénomination
résumée par lettre}, constitue un universel de la technique
mathématique, en tant qu'il subsume un multiple sous une unité
provisoire. Cette unité notationnelle entretiendra des relations
algébriques contractées avec les autres quantités formelles pendant
lesquelles certaines simplifications pourront être lisibles à un tel
niveau, {\em sans qu'il soit nécessaire de deviner les mêmes
simplifications formelles au niveau beaucoup plus foisonnant de
l'expansion complète des formules}. Autrement dit, d'un point de vue
strictement calculatoire qui se refuserait à céder à l'interprétation
de ces trois déterminants $2 \times 2$ comme composantes du {\sl
produit vectoriel hamiltonien} ${\bf x}_u \wedge {\bf x}_v$, Gauss
synthétise trois `molécules formelles de calcul' qu'il s'autorisera
ensuite, quand cela sera nécessaire,
à désintégrer en atomes fondamentaux. 

\CITATION{
Mais des équations $dx = a\, dp + a'\, dq$, $dy = b\, dp + b'\, dq$,
nous tirons:
\[
\aligned
C\,dp
&
=
b'\,dx-a'\,dy,
\\
C\,dq
&
=
-b\,dx+a\,dy.
\endaligned
\]
Par là, nous obtenons les différentielles complètes de $t$ et de $u$:
\[
\aligned
C^3\,dt
&
=
\bigg(
A\,\frac{dC}{dp}-C\,\frac{dA}{dp}
\bigg)\,
(b'\,dx-a'\,dy)
+
\bigg(C\,\frac{dA}{dq}-A\,\frac{dC}{dq}
\bigg)\,
(b\,dx-a\,dy),
\\
C^3\,du
&
=
\bigg(
B\,\frac{dC}{dp}-C\,\frac{dB}{dp}
\bigg)\,
(b'\,dx-a'\,dy)
+
\bigg(C\,\frac{dB}{dq}-B\,\frac{dC}{dq}
\bigg)\,
(b\,dx-a\,dy).
\endaligned
\]
\REFERENCE{\cite{ga1828},~23}}

Gauss ne signale pas que les valeurs du déterminant $C$ sont
automatiquement non nulles, puisque par hypothèse, $C$ exprime le
jacobien de l'application inversible (ou difféomorphisme local) $(p, q)
\mapsto (x, y)$: c'est clair; 
il n'y a pas d'erreur; au lecteur de s'en
rendre compte par lui-même. Cette hypothèse\,\,---\,\,toujours 
satisfiable\,\,---\,\,relative à la
disposition du système de coordonnées par rapport à la surface permet
alors non seulement de résoudre les deux différentielles $dp$ et $dq$
par rapport à $dx$ et $dy$ grâce à des formules élémentaires
considérées comme connues, mais encore de {\em diviser} par $C$, et de
raisonner ensuite dans un certain corps de fractions
rationnelles. Toutefois, Gauss préfère calculer {\em sans écrire aucun
dénominateur}, et c'est la raison pour laquelle il multiplie les
membres de gauche par $C$ ou par $C^3$, afin de gommer la lourdeur
d'une réduction au même dénominateur,
qui serait exigible après différentiation d'une
expression rationnelle. On peut voir en cela un signe de 
{\em distinction de calculateur}, si doué qu'il {\em connaît déjà}
toute l'élégance des raccourcis de ce type, d'autant
plus que les expressions offertes à la lecture sont
parfaitement polies et symétriques.

Par ailleurs, il faut aussi noter que l'évanouissement des
dénominateurs dans le calcul intermédiaire non présenté:
\[
\small
\aligned
C^3\,dt
&
=
C^3
\bigg\{
\bigg(
-\frac{\partial A}{\partial p}\,\frac{1}{C}
+
A\,\frac{\partial C}{\partial p}\,\frac{1}{C^2}
\bigg)\,dp
+
\bigg(
-\frac{\partial A}{\partial q}\,\frac{1}{C}
+
A\,\frac{\partial C}{\partial q}\,\frac{1}{C^2}
\bigg)\,dq
\bigg\}
\\
&
=
\bigg(
-\,C\,\frac{\partial A}{\partial p}
+
A\,\frac{\partial C}{\partial p}
\bigg)\,C\,dp
+
\bigg(
-\,C\,\frac{\partial A}{\partial q}
+
A\,\frac{\partial C}{\partial q}
\bigg)\,C\,dq
\\
&
=
\bigg(
A\,\frac{\partial C}{\partial p}-C\,\frac{\partial A}{\partial p}
\bigg)\,
(b'\,dx-a'\,dy)
+
\bigg(C\,\frac{\partial A}{\partial q}-A\,\frac{\partial C}{\partial q}
\bigg)\,
(b\,dx-a\,dy),
\endaligned
\] 
exige de {\em multiplier à l'avance} le membre de gauche par
$C^3$\,\,---\,\,et non pas par $C^2$\,\,---\,\,sachant qu'un premier
facteur $C^2$ sera absorbé par les dérivées partielles
entre parenthèses, et qu'un
second facteur $C^1$ sera nécessaire pour que les deux différentielles
$C \, dp$ et $C\, dq$ apparaissent {\em sans dénominateur}, 
exactement comme elles avaient été 
calculées juste auparavant.

\Scholie
En définitive, malgré une apparence de fluidité et de complétude,
certains calculs intermédiaires sont néanmoins régulièrement
sous-entendus; seules les expressions finales nettes et symétriques
sont offertes au lecteur. 

Toujours à ce sujet, une question de
philosophie générale des mathématiques demeure: quel statut conférer
aux calculs intermédiaires, en tant que traces effaçables d'un
parcours? Un doute s'installe, car il est tout à fait possible que
des causes importantes vis-à-vis du contenu se jouent et se
poursuivent même dans les calculs élémentaires\,\,---\,\,que
ces calculs 
soient sous-entendus
(laissés au lecteur), ou inconsciemment éludés, voire même cachés
intentionnellement.
\qed

\smallskip
\SHORTCITATION{
Maintenant, si dans ces formules nous substituons les valeurs
suivantes:}
\SHORTCITATION{
\[
\aligned
\frac{dA}{dp}
&
=
c'\,\beta+b\,\gamma'-c\,\beta'-b'\,\gamma,
\\
\frac{dA}{dq}
&
=
c'\,\beta'+b\,\gamma''-c\,\beta''-b'\,\gamma',
\\
\frac{dB}{dp}
&
=
a'\,\gamma+c\,\alpha'-a\,\gamma'-c'\,\alpha,
\\
\frac{dB}{dq}
&
=
a'\,\gamma'+c\,\alpha''-a\,\gamma''-c'\,\alpha',
\\
\frac{dC}{dp}
&
=
b'\,\alpha+a\,\beta'-b\,\alpha'-a'\,\beta,
\\
\frac{dC}{dq}
&
=
b'\,\alpha'+a\,\beta''-b\,\alpha''-a'\,\beta',
\endaligned
\]}
\SHORTCITATION{\!\!\!\!\!\!\!\!\!\!\!\!\!\!\!\!\!
et si nous remarquons que les valeurs des différentielles $dt$, $du$
ainsi obtenues doivent être égales, indépendamment des
différentielles $dx$, $dy$, aux quantités $T\, dx + U\, dy$, 
$U\, dx + V\, dy$ respectivement, nous trouverons, 
après quelques transformations assez aisées:}\smallskip
\SHORTCITATION{
\[
\aligned
C^3\,T
&
=
\alpha\,A\,{b'}^2+\beta\,B\,{b'}^2+\gamma\,C\,{b'}^2
-
\\
&
\ \ \ \ \
-
2\,\alpha'\,A\,b\,b'-2\,\beta'\,B\,b\,b'-2\,\gamma'\,C\,b\,b'
+
\\
&
\ \ \ \ \
+
\alpha''\,A\,b^2+\beta''\,B\,b^2+\gamma''\,C\,b^2,
\endaligned
\]}
\SHORTCITATION{
\[
\aligned
\ \ \ \ \ \ \ \ \ \ \ \ \ \ \ \ \ \ \ \ \ \ \ \
C^3\,U
&
=
-\,\alpha\,A\,a'\,b'-\beta\,B\,a'\,b'-\gamma\,C\,a'\,b'
+
\\
&
\ \ \ \ \
+
\alpha'\,A\,(ab'+ba')+\beta'\,B\,(ab'+ba')+\gamma'\,C\,(ab'+ba')
-
\\
&
\ \ \ \ \
-\alpha''\,A\,a\,b-\beta''\,B\,a\,b-\gamma''\,C\,a\,b,
\endaligned
\]}
\SHORTCITATION{
\[
\aligned
C^3\,V
&
=
\alpha\,A\,{a'}^2+\beta\,B\,{a'}^2+\gamma\,C\,{a'}^2
-
\\
&
\ \ \ \ \
-
2\,\alpha'\,A\,a\,a'-2\,\beta'\,B\,a\,a'-2\,\gamma'\,C\,a\,a'
+
\\
&
\ \ \ \ \
+\alpha''\,A\,a^2+\beta''\,B\,a^2+\gamma''\,C\,a^2.
\endaligned
\vspace{-0.5cm}
\]
\REFERENCE{\cite{ga1828},~24}}

\medskip
Ici, les calculs commencent à se déployer véritablement. Quelque peu
paradoxalement, la théorie des surfaces\,\,---\,\,objets 
géométriques à {\em deux} variables
seulement\,\,---\,\,implique 
des manipulations algébriques dans des anneaux de
polynômes à un nombre
beaucoup plus considérable de variables, puisque si l'on compte le
nombre des quantités (en partie redondantes et hybrides
à cause du choix des notations) qui
interviennent dans les expressions rationnelles de $U$, $V$, $T$
ci-dessus, à savoir: $a$, $a'$, $b$, $b'$, $c$, $c'$, $\alpha$,
$\alpha'$, $\alpha''$, $\beta$, $\beta'$, $\beta''$, $\gamma$,
$\gamma'$, $\gamma''$, $A$, $B$, $C$, on trouve au total {\bf 18}
lettres, c'est-à-dire nettement plus que la dimension {\bf 2} des 
surfaces étudiées. 

\Scholie
Le paradoxe de l'augmentation
(explosion?) de la dimension d'étude est ubiquitaire en
géométrie différentielle computationnelle. Question ouverte de
philosophie générale des mathématiques:
pourquoi en est-il ainsi? quelles
sont les causes profondes? {\em
A priori}, le synthétique mathématique doit-t-il nécessairement
s'exprimer dans un multiple qui le dépasse toujours
sur le plan numéraire?
\qed

\smallskip
Dans un premier moment, les six dérivées partielles $\frac{ dA}{ dp}$,
$\frac{ dA}{ dq}$, $\frac{ dB}{ dp}$, $\frac{ dB}{ dq}$, $\frac{ dC}{
dp}$, $\frac{ dC}{ dq}$ se calculent automatiquement en partant des
expressions quadratiques de $A$, $B$, $C$ et en appliquant la règle de
Leibniz, raison pour laquelle Gauss n'en dit rien, laissant à la
réactivité du lecteur le soin de deviner cela spontanément et de
vérifier au besoin toutes les étapes de calcul. 

\`A tout le moins
pour une simple lecture passive, les six expressions sont écrites en
totalisation, toutes leurs symétries symboliques étant
intentionnellement rendues visibles, comme si une {\em Gestalt} du
symbolique devait suffire pour confirmer la correction des
expressions. Ensuite, dans un second moment, ces six dérivées
partielles sont {\em insérées} dans les deux expressions précédentes
de $C^3\,dt$ et de $C^3\, du$. 

Actes non explicités: des expansions,
des réorganisations de termes, des ordonnancements de
monômes\,\,---\,\,la tâche est conséquente. Enfin, les quatre
coefficients de $dp$, $dq$ dans $C^3\, dt$ et de $dp$, $dq$ dans
$C^3\, du$ permettent d'extraire les valeurs désirées de $C^3\, T$,
$C^3\, U$, $C^3\, V$: cela, n'est pas non plus explicitement écrit,
exemple supplémentaire de `{\em non-littérarité}' du Calcul.

\CITATION{
Si maintenant nous posons, pour abréger:\smallskip
\\
(1)
\centerline{$A\,\alpha+B\,\beta+C\,\gamma=D$}\smallskip
\\
(2)
\centerline{$A\,\alpha'+B\,\beta'+C\,\gamma'=D'$}\smallskip
\\
(3)
\centerline{$A\,\alpha''+B\,\beta''+C\,\gamma''=D''$}\smallskip
\\
il viendra:
\[
\aligned
C^3\,T
&
=
D\,{b'}^2-2\,D'\,b\,b'+D''\,{b^2},
\\
C^3\,U
&
=
-D\,a'\,b'+D'\,(ab'+ba')-D''\,a\,b,
\\
C^3\,V
&
=
D\,{a'}^2-2\,D'\,a\,a'+D''\,a^2.
\endaligned
\vspace{-0.5cm}
\]
\REFERENCE{\cite{ga1828},~24--25}}

Certes, ces trois expressions $D$, $D'$, $D''$ ne sont autres que de
simples déterminants $3 \times 3$ de dérivées partielles du rayon
vecteur ${\bf x} ( u, v)$, mais il a été convenu ci-dessus de faire
abstraction de toute réinterprétation vectorielle des calculs de
Gauss. Davantage, il faut voir ici, dans une certaine opacité
des calculs due à leur
volume, un acte de synthèse visuelle par factorisation
inévidente. D'ailleurs, Gauss avait finement anticipé cette {\em
micro-découverte formelle locale} (que ses manuscrits ont dû scruter):
en effet, dans les expressions de $C^3\, T$, $C^3\, U$, $C^3\, V$, des
sommes identiques de trois termes se reconstituaient en facteurs,
puisque les équations brutes étaient {\em déjà
présentées dans un ordre
qui ferait facilement voir au lecteur après-coup l'émergence de $D$,
$D'$, $D''$}.

\CITATION{
Par là nous obtenons, tous calculs faits:
\[
C^6(TV-U^2)
=
(DD''-{D'}^2)(ab'-ba')^2
=
(DD''-{D'}^2)\,C^2,
\]
d'où résulte l'expression suivante, pour la
mesure de la courbure:
\[
k
=
\frac{DD''-{D'}^2}{(A^2+B^2+C^2)^2}.
\vspace{-0.5cm}
\]
\REFERENCE{\cite{ga1828},~25}}

La métamorphose importante et remarquable 
du terme carré au dénominateur:
\[
1+z_x^2+z_y^2
=
1+t^2+u^2
=
A^2+B^2+C^2
\]
passe ici dans un {\em tacet} implicite qui relègue une telle
vérification au statut de l'évidence. Au niveau où en sont
les calculs, l'auteur est en droit d'attendre que son lecteur soit
apte à résoudre les exercices implicites qu'il lui laisse. Cette
relation n'est pas explicitée par anticipation dans ce qui précède du
mémoire, et en fait, elle apparaît plus bas, sans calcul
détaillé. Toutefois elle est cruciale, puisque grâce à elle, le sort
de la métamorphose vers l'intrinsèque s'avère maintenant
définitivement réglé quant au dénominateur.

\Scholie
Voici donc un exemple de la {\sl pensée tue}, ce phénomène qui affecte
toutes les mathématiques aussi bien dans leur transmission didactique
que dans leur développement par la recherche. \`A aucun moment Gauss
n'écrit qu'il cherche à {\em transsubstantier vers l'intrinsèque}
la formule graphée pour la courbure. Un mouvement synthétique long et
irréversible s'engage sans rien en dire. Du point de vue
de l'exégète qui a pris connaissance des études sur le {\em Nachlass},
le paradoxe s'accentue quand on sait combien d'efforts Gauss a
dépensés afin de réaliser l'intrinsèque par une formule mathématique.
\qed

\medskip\noindent{\bf \'Elimination du dernier bastion d'extrinsèque.}
\label{elimination-dernier-bastion-extrinseque}
En résumé, la formule précédente pour la courbure:
\[
\aligned
\kappa
&
=
\frac{DD''-{D'}^2}{(A^2+B^2+C^2)^2},
\endaligned
\]
bien que relativement complexe lorsqu'on en explicite complètement les
composantes, demeure essentiellement {\em sans mystère}, puisqu'elle
découle de la formule graphée par simple transfert de ses termes
différentiels. Mais elle fait encore intervenir dans son numérateur
les trois dimensions de plongement: dernier bastion d'extrinsèque.

\CITATION{
\`A l'aide de la formule ci-dessus, nous allons
maintenant en obtenir une autre qui peut être
comptée au nombre des théorèmes les plus féconds
dans la théorie des surfaces courbes.
\REFERENCE{\cite{ga1828},~25}}

Le \S~XI des {\em Disquisitiones} entreprend alors sans le dire un
{\em calcul final}, lequel change de nature,
puisqu'il s'agit maintenant d'un véritable 
{\em calcul d'élimination}.

\Scholie
Gauss {\em sait} qu'il a découvert une grande vérité mathématique.
Pourquoi se dépenserait-il à en écrire les tenants et 
aboutissants\footnote{\,
%%%%%%%%%%%%%%%%%%%%%%%-------DEBUT--------%%%%%%%%%%%%%%%%%%%%%%%%%%%
D'après l'article ``{\em Gauss zum Gedächtniss}''
écrit en 1856 par Sartorius von Waltershausen:
<<\,Une fois qu'un bâtiment a été construit [$\cdots$]
il faut qu'on ne voie plus les échafaudages\,>>
répétait-il.
}? %%%%%%%%%%%%%%%%%%%%%%%%-----FIN-----%%%%%%%%%%%%%%%%%%%%%%%%%%%%%%%
Tel est son choix métaphysique: seule compte la vérité, dite par un
calcul à laquelle le résultat a conduit {\em irréversiblement}.
\qed

\smallskip
Sans aucune explication à cet endroit précis, 
trois notations par lettres\,\,---\,\,{\em retenues par la 
postérité}!\,\,---\,\,sont introduites pour les
coefficients métriques en représentation 
paramétrique:
\[
\aligned
a^2+b^2+c^2
&
=
E,
\\
aa'+bb'+cc'
&
=
F,
\\
{a'}^2+{b'}^2+{c'}^2
&
=
G.
\endaligned
\]
Systématiquement aussi, six autres notations sont introduites:
\[
\aligned
a\,\alpha+b\,\beta+c\,\gamma
&
=
m,
\\
a\,\alpha'+b\,\beta'+c\,\gamma'
&
=
m',
\\
a\,\alpha''+b\,\beta''+c\,\gamma''
&
=
m'',
\\
a'\alpha+b'\beta+c'\gamma
&
=
n,
\\
a'\alpha'+b'\beta'+c'\gamma'
&
=
n',
\\
a'\alpha''+b'\beta''+c'\gamma'
&
=
n'',
\endaligned
\]
et enfin, une dernière
notation est introduite, qui ré-effectue une synthèse du
dénominateur et qui explicite un acte passé sous silence
un peu plus haut:
\[
A^2+B^2+C^2=EG-F^2=\Delta.
\]

\CITATION{
\'Eliminons des équations (1), (4), (7) les quantités $\beta$,
$\gamma$, ce qui se fera en multipliant ces équations par $bc' - cb'$,
$b' C - c' B$, $cB - bC$, et ajoutant; nous aurons ainsi:
\[
\aligned
& \ \ \ \
\big[A\,(bc'-cb')+a\,(b'C-c'B)+a'\,(cB-bC)\big]\,\alpha
\\
&
=
D\,(bc'-cb')+m\,(b'C-c'B)+n\,(cB-bC),
\endaligned
\] 
équation que nous transformons facilement en celle-ci:
\[
A\,D
=
\alpha\,\Delta
+
a\,(nF-mG)
+
a'\,(mF-nE).
\vspace{-0.65cm}
\]
\REFERENCE{\cite{ga1828},~25--26}}

Maintenant, les calculs d'élimination commencent, sans explication
d'intention. Même en connaissance de l'objectif
visé (la {\em formula egregia}), il se peut que l'on
ne comprenne pas du premier coup ce que cherche à faire Gauss ici. En
effet, il a choisi de ne pas dire tout de suite que les six quantités
$m$, $m'$, $m''$, $n$, $n'$, $n''$ s'expriment en fonction des
premières dérivées partielles de $E$, $F$, $G$, donc qu'elles
conviennent pour la formule finale. On laissera donc temporairement
en suspens la question de savoir à quoi sert l'élimination des
quantités $\beta$, $\gamma$, et on modifiera légèrement l'ordre
d'apparition de l'information.

\Scholie
Au niveau de l'esprit, la {\sl compréhension mathématique} est une
synthèse qui embrasse une multiplicité d'actes de pensée suffisante
pour qu'il soit possible de {\em circuler librement} d'un moment à un
autre des raisonnements.  Un désordre, certes, s'installe dans le flux
des mentalisations, mais ce désordre d'`{\em après coup}' est le {\em
sang} d'une totalisation.
\qed

\Scholie
<<\,Suivez le guide\,>>; <<\,multipliez, additionnez, soustrayez
telles et telles équations\,>>; <<\,à la fin, vous obtiendrez le
résultat désiré\,>>: telles sont les
(seules) informations que fournit un
mémoire de mathématiques
publié, mais ce n'est jamais de cette manière présentationnelle {\em
a posteriori} que la pensée du calcul procède. {\sf
La {\em vraie pensée
causale des calculs} aura exploré le problème sous un {\em faisceau
d'angles interrogatifs 
beaucoup plus complexes}}. Elle aura développé
complètement le numérateur $DD'' - {D'}^2$, réorganisé ses termes,
examiné comment se combinent entre elles les trois dérivées partielles
d'ordre le plus élevé ${\bf x}_{ uu}$, ${\bf x}_{ uv}$, ${\bf
x}_{vv}$, du rayon vecteur ${\bf x} ( u, v)$ et cherché à éliminer
l'extrinsèque.
\qed

\CITATION{
Maintenant, il est évident que l'on a:
\[
\aligned
\frac{dE}{dp}=2m,
\ \ \ \ \ 
\frac{dE}{dq}
&
=
2\,m',
\ \ \ \ \
\frac{dF}{dp}=m'+n,
\ \ \ \ \
\frac{dF}{dq}=m''+n',
\\
\frac{dG}{dp}
&
=
2\,n',
\ \ \ \ \
\frac{dG}{dq}=2\,n'',
\endaligned
\]
ou bien:
\[
\aligned
m
&
=
\frac{1}{2}\,\frac{dE}{dp},
\ \ \ \ \
m'=\frac{1}{2}\,\frac{dE}{dq},
\ \ \ \ \
m''
=
\frac{dF}{dq}
-
\frac{1}{2}\,\frac{dG}{dp},
\\
n
&
=
\frac{dF}{dp}
-
\frac{1}{2}\,\frac{dE}{dq},
\ \ \ \ \
n'=\frac{1}{2}\,\frac{dG}{dp},
\ \ \ \ \
n''=\frac{1}{2}\,\frac{dG}{dq};
\endaligned
\]
de plus, il est facile de s'assurer que l'on a:
\[
\aligned
\alpha\alpha''+\beta\beta''+\gamma\gamma''
-{\alpha'}^2-{\beta'}^2-{\gamma'}^2
&
=
\frac{dn}{dq}
-
\frac{dn'}{dp}
=
\frac{dm''}{dp}
-
\frac{dm'}{dq}
\\
&
=
-\frac{1}{2}\cdot
\frac{d^2E}{dq^2}
+
\frac{d^2F}{dp\,dq}
-
\frac{1}{2}\cdot\frac{d^2G}{dp^2}.
\endaligned
\]
\REFERENCE{\cite{ga1828},~26--27}}

Le seul objectif de Gauss ici est de préparer une expression
appropriée de $DD'' - {D'}^2$ qui guidera à l'avance\,\,---\,\,sans
aucun moment d'exploration questionnante au sujet de son développement
brut\,\,---\,\,vers une réorganisation de ses termes où les neuf
dérivées partielles secondes apparaîtront seulement sous la forme
quadratique:
\[
\alpha\alpha''+\beta\beta''+\gamma\gamma''
-
{\alpha'}^2-{\beta'}^2-{\gamma'}^2.
\]
Telle est l'intention `technique' de ces calculs qui semblaient
mystérieux au premier abord: que seule cette combinaison quadratique
des dérivées partielles secondes apparaisse, puisqu'elle se ré-exprime
comme:
\[
-{\textstyle{\frac{1}{2}}}\,E_{qq} 
+
F_{pq} 
- 
{\textstyle{\frac{1}{2}}}\,G_{pp}.
\]

\Scholie
Cette explication a l'air satisfaisante au premier abord.  Toutefois,
la vérité de l'{\sl exploration-en-recherche} est beaucoup plus vaste
et complexe, parce que de nombreuses {\em questions-idées} auront été
éliminées à cause de l'accès irréversible à un but visé, et parce que
de nombreuses autres explications n'auront pas été apportées.

En effet, est-ce la seule combinaison linéaire de dérivées partielles
d'ordre deux qui s'exprime en fonction de $E$, $F$, $G$?  Comment les
restes d'ordre $1$ se réorganisent-ils?  Pourquoi cette combinaison
quadratique est-elle multipliée par le dénominateur $\Delta$?  Est-ce
un phénomène universel?  Que se passe-t-il en dimension supérieure?
Dans d'autres problèmes de géométrie différentielle, existe-t-il des
phénomènes semblables de `découverte' d'invariants par le biais de
purs calculs?

Toutes les questions de cette sorte, et d'autres encore situées plus
au c{\oe}ur des intentions techniques, sont irrémédiablement éliminées
du texte publié\,\,---\,\,{\sl pensée disparaissante}\,\,---\,\,bien
qu'elles structurent la {\sl pensée-productrice-en-recherche}.
\qed

\CITATION{
En éliminant des mêmes équations, soit $\alpha$ et $\gamma$,
soit $\alpha$ et $\beta$, on obtiendrait
de même:
\[
\aligned
BD
&
=
\beta\,\Delta
+
b\,(nF-mG)
+
b'(mF-nE),
\\
CD
&
=
\gamma\,\Delta
+
c\,(nF-mG)+c'(mF-nE).
\endaligned
\]
\REFERENCE{\cite{ga1828},~26}}

Ici comme auparavant, ces calculs élémentaires
d'élimination exigent une pensée systématique de la réorganisation
et de la symétrisation. Gauss place
à dessein les trois dérivées partielles secondes $\alpha$, 
$\beta$, $\gamma$ en première position à droite
du signe d'égalité, afin de reconstituer
comme suit le produit $DD''$:

\CITATION{
Multipliant ces trois équations respectivement par
$\alpha''$, $\beta''$, $\gamma''$ et ajoutant, nous
obtenons:
\def\theequation{10}\begin{equation}
\left\{
\aligned
DD''
&
=
(\alpha\alpha''+\beta\beta''+\gamma\gamma'')\,\Delta
+
\\
&
\ \ \ \ \
+
m''(nF-mG)+n''(mF-nE).
\endaligned\right.
\vspace{-0.65cm}
\end{equation}
\REFERENCE{\cite{ga1828},~26}}

Quasi `miraculeusement', lorsqu'on multiplie ainsi ces trois
équations par $\alpha''$, $\beta ''$, $\gamma''$, les deux
combinaisons linéaires $a\, \alpha'' + b \, \beta '' + c\, \gamma '' =
m''$ et $a' \alpha'' + b' \beta '' + c' \gamma'' = n''$ {\em
reconstituent des quantités intrinsèques}. Il est vrai que ces
calculs s'organiseraient plus économiquement si l'on s'autorisait à
utiliser la théorie des déterminants $3 \times 3$ et leur
développement le long de certaines colonnes, mais la systématicité
symétrique des calculs d'élimination de Gauss
perdrait de sa visibilité.

\CITATION{
Si nous traitons de même les équations (2), (5), (8), il vient:
\[
\aligned
AD'
&
=
\alpha'\Delta+a\,(n'F-m'G)+a'(m'F-n'E),
\\
BD'
&
=
\beta'\Delta+b\,(n'F-m'G)+b'(m'F-n'E),
\\
CD'
&
=
\gamma'\Delta+c\,(n'F-m'G)+c'(m'F-n'E);
\endaligned
\]
multipliant ces équations respectivement par $\alpha'$, 
$\beta'$, $\gamma'$ et ajoutant, on a:
\[
{D'}^2
=
({\alpha'}^2+{\beta'}^2+{\gamma'}^2)\,\Delta
+
m'(n'F-m'G)+n'(m'F-n'E).
\vspace{-0.65cm}
\]
\REFERENCE{\cite{ga1828},~26}}

Systématiques, symétriques\,\,---\,\,beaux et élaborés, tels sont les
calculs de Gauss. 

\CITATION{
La combinaison de cette équation avec l'équation~(10) donne:
\[
\aligned
DD''-{D'}^2
&
=
\big(\alpha\alpha''+\beta\beta''+\gamma\gamma''
-{\alpha'}^2-{\beta'}^2-{\gamma'}^2\big)\,\Delta
+
\\
&
\ \ \ \ \
+
E\,({n'}^2-nn'')
+
F\,(nm''-2\,m'n'+mn'')
+
G\,({m'}^2-mm'').
\endaligned
\]
\REFERENCE{\cite{ga1828},~26}}

Une fois ces calculs faits, le numérateur $DD'' - {D'}^2$ s'exprime
enfin sous la forme d'un polynôme différentiel en les dérivées
partielles d'ordre $\leqslant 2$ des trois coefficients métriques $E$,
$F$, $G$.

\CITATION{
Si maintenant nous substituons ces expressions diverses dans la
formule établie à la fin du \S précédent, pour la mesure de la
courbure, nous parvenons à la formule suivante, qui ne renferme que
les seules quantités $E$, $F$, $G$ avec leurs quotients différentiels
du premier et du second ordre:
\[
\boxed{
\aligned
4\big(EG-F^2\big)^2\,k
&
=
E\,\bigg[
\frac{dE}{dq}\cdot\frac{dG}{dq}-
2\,\frac{dF}{dp}\cdot\frac{dG}{dq}
+
\bigg(
\frac{dG}{dp}
\bigg)^2
\bigg]
+
\\
&
\ \ \ \ \
+
F\,\bigg(
\frac{dE}{dp}\cdot\frac{dG}{dq}
-
\frac{dE}{dq}\cdot\frac{dG}{dp}
-
2\,\frac{dE}{dq}\cdot\frac{dF}{dq}
+
\\
&
\ \ \ \ \ \ \ \ \ \ \ \ \ \ \ 
+
4\,\frac{dF}{dp}\cdot\frac{dF}{dq}
-
2\,\frac{dF}{dp}\cdot\frac{dG}{dp}
\bigg)
+
\\
&
\ \ \ \ \
+
G\,\bigg[
\frac{dE}{dp}\cdot\frac{dG}{dp}
-
2\,\frac{dE}{dp}\cdot\frac{dF}{dq}
+
\bigg(
\frac{dE}{dq}
\bigg)^2
\bigg]
-
\\
&
\ \ \ \ \
-
2\,\big(EG-F^2\big)\,
\bigg(
\frac{d^2E}{dq^2}
-
2\,\frac{d^2F}{dp\,dq}
+
\frac{d^2G}{dp^2}
\bigg).
\endaligned}
\vspace{-0.65cm}
\]
\REFERENCE{\cite{ga1828},~27}}

\medskip\noindent{\bf Theorema egregium.}
\label{theorema-egregium}
Sur toute surface $S$ plongée dans l'espace euclidien $\R^3$, la
longueur au carré d'un élément infinitésimal de coordonnées $(dx, dy,
dz)$ tangent à la surface s'exprime, {\em intrinsèquement} dans des
coordonnées paramétriques bidimensionnelles internes à la surface:
\[
(p,q)
\longmapsto
\big(x(p,q),\,y(p,q),\,z(p,q)\big)
\]
comme une certaine forme quadratique 
définie postive en le vecteur infinitésimal
de coordonnées $(dp, dq)$:
\[
ds^2
=
dx^2+dy^2+dz^2
=
E\,dp^2+2F\,dpdq+G\,dq^2,
\]
dont les trois coefficients s'expriment naturellement par:
\[
\aligned
E
&
=
x_u^2+y_u^2+z_u^2
=
a^2+b^2+c^2,
\\
F
&
=
x_ux_v+y_uy_v+z_uz_v
=
aa'+bb'+cc',
\\
G
&
=
x_v^2+y_v^2+z_v^2
=
{a'}^2+{b'}^2+{c'}^2.
\endaligned
\]
La {\em formula egregia} de Gauss apprend donc que pour 
trouver la mesure de la courbure d'une surface
en l'un de ses points quelconques
de coordonnées $(p, q)$, on n'a pas besoin de
formules finies donnant les coordonnées extrinsèques
$(x, y, z)$ en fonction de coordonnées paramétriques internes
$(p, q)$, mais qu'il suffit
de connaître l'expression générale de la grandeur de chaque
élément linéaire {\em dans les coordonnées bidimensionnelles
internes} $(p, q)$. 

L'application principale de cette {\em formula egregia}
est la suivante. 
Si, à chaque point $(x, y, z)$ d'une surface $S \subset \R^3$
correspond, de manière différentiable ou analytique,
un unique point $(x', y', z')$ d'une 
autre surface $S' \subset \R^3$, la correspondance étant
biunivoque et analytique réelle ou de classe au moins $\mathcal{ C}^2$, 
alors les coordonnées images $(x', y', z')$ pourront
aussi être considérées comme des fonctions
de $p$ et de $q$, d'où pour
l'élément linéaire pythagoricien induit 
$\sqrt{ d{x'}^2 + d {y'}^2 + d{ z'}^2}\big\vert_{ S'}$
une autre expression similaire telle que:
\[
\sqrt{
E'\,dp^2+2F'\,dpdq+G'\,dq^2},
\] 
dans laquelle
$E'$, $F'$, $G'$ sont aussi certaines fonctions
déterminées de $p$ et de $q$. 
Mais par la notion même d'{\em application}
entre les deux surfaces, les éléments
infinitésimaux qui se correspondent sur 
chaque surface seront alors nécessairement 
égaux l'un à l'autre point par point,
et l'on aura identiquement:
\[
E'=E,\ \ \ \ \ \ \
F'=F,\ \ \ \ \ \ \
G'=G;
\]

\CITATION{
\!\!\!\!\!\!\!\!\!\!\!
de sorte que la formule du \S~ précédent conduit spontanément à 
ce théorème remarquable \deutsch{egregium}:
\\
${}$ \ \ \ \ \
{\em Si une surface courbe est appliquée sur une autre
surface courbe quelconque, la mesure de courbure en chaque
point reste invariable}.
\\
${}$ \ \ \ \ \
Par suite, {\em La courbure intégrale d'une portion
finie quelconque de la surface ne changera pas}. 
\REFERENCE{\cite{ga1828},~26}}

\Section{13.~Leçons de m\'etaphysique gaussienne}
\label{lecons-metaphysique-gaussienne}

\HEAD{13.~Le\c cons de m\'etaphysique gaussienne}{
Jo\"el Merker, D\'partement de Math\'matiques d'Orsay}

\medskip\noindent{\bf Différer et mûrir, différer pour mûrir.} 
\label{differer-murir}
Comme la correspondance
de Gauss le montre\footnote{\,
%%%%%%%%%%%%%%%%%%%%%%%-------DEBUT--------%%%%%%%%%%%%%%%%%%%%%%%%%%%
Pour les extraits de correspondance, 
{\em cf.}~\cite{ stackel1890, do1979}. 
}, %%%%%%%%%%%%%%%%%%%%%%%%-----FIN-----%%%%%%%%%%%%%%%%%%%%%%%%%%%%%%%
l'élaboration de ses {\em Recherches générales sur les surfaces
courbes} est le fruit {\em mûri} d'une quinzaine d'années de
méditations continues, le produit d'efforts intellectuels soutenus et
répétés, et aussi le résultat de plus d'un an de travail intensif
durant la période qui précéda la publication effective du mémoire.  La
recherche du meilleur argument analytique qui démontrât l'invariance
de la mesure de courbure par isométries
infinitésimales\,\,---\,\,c'est-à-dire la quête renouvelée d'une
formule telle que {\em formula egregia} qui eût été valable dans un
système de coordonnées intrinsèques quelconques\,\,---\,\,aura
constitué le {\sl point de blocage principal} pour Gauss avant qu'il
ne s'autorisât à engager tout travail de rédaction définitive.

\medskip\noindent{\bf Maintenir en tension la recherche de vérités.}
\label{maintenir-tension-recherche-verites}
En effet, un principe métaphysique de raison déterminante commandait
de soupçonner qu'une formule centrale eût dû {\em expliquer à un
niveau supérieur} l'invariance de la mesure de courbure, découverte
aux alentours de 1816 grâce à la formule d'écart angulaire
$\widehat{A} + \widehat{B} +
\widehat{C} = \pi + \int\!\! \int_{ \mathcal{T}}\, \kappa\, d \sigma$
valable pour tout triangle géodésique $\mathcal{ T} = ABC$ tracé sur
une surface plongée (\deutschplain{schöne Theorem}). Un impératif
catégorique commandait donc à Gauss de poursuivre ses recherches
jusqu'à ce qu'il découvrît une vérité qui ne pouvait se réaliser que
dans et par un long calcul d'élimination.

Or, quant à la véritable nature de l'obstacle, ce calcul d'un genre
nouveau requérait une habileté très différente de celle que Gauss
avait développée au contact de l'arithmétique abstraite des nombres
entiers, puisqu'il s'agissait d'algèbre différentielle effective,
domaine qui ne s'est vraiment développé qu'après lui, d'abord avec les
travaux précurseurs de Riemann, Christoffel, Lipschitz sur la courbure
des variétés métriques, et surtout à la charnière du dix-neuvième et
vingtième siècles avec les {\oe}uvres de Lie et Cartan, poursuivies
aussi à la même époque par Engel, Scheffers, Kowalewski, Tresse,
Cartan, Drach, Janet, et d'autres.

\CITATION{
Il est intéressant de constater que les déclarations de Gauss les plus
connues et les plus profondes au sujet de son style propre de travail
proviennent de sa période d'activité en géométrie
différentielle.
\REFERENCE{\cite{do1979},~121}}

\medskip\noindent{\bf Calculer est une vérité de la chose mathématique
en elle-même.} \label{calcul-verite-chose-en-soi} Conséquence
incontournable quant à l'interprétation philosophique: du point de vue
supérieur auquel se place Gauss quant à l'exploration mathématique, le
calcul explicite ne constitue {\em pas}, par conséquent, une simple
réalisation en coordonnées d'un concept connu d'avance comme étant
intrinsèque, voire même `préempté' par une abstraction structuraliste
qui ré-organise {\em a posteriori} un champ théorique\,\,---\,\,non,
bien au contraire en fait: le calcul {\em est} vérité de la courbure
en elle-même, il {\em forme} et il {\em in-forme} l'essence de la
vérité mathématique concernée, c'est-à-dire qu'il confère
intrinsèquement à la courbure la forme générale explicite et complète
qu'elle doit posséder à l'intérieur même de son être véritable.

\medskip\noindent{\bf Voir émerger des entités autonomes organisées.}
\label{voir-emerger-entites}
Par ailleurs, dans une lettre à Olbers qu'il a écrite à Göttingen le
30 octobre 1825 (citée par Dombrovski), Gauss exprime une autre
exigence métaphysique qui rappelle un motif récurrent de la
philosophie hégélienne de la connaissance, en tant que ce n'est
qu'après avoir effectué un parcours dialectique complet que la pensée
peut se retourner et embrasser par compréhension la totalité d'une
chose méditée.

\CITATION{
Bien que les aspects mathématiques d'une recherche m'apparaissent
comme les plus intéressants, je ne peux pas contester d'un autre côté
qu'afin d'éprouver du plaisir avec un travail aussi long que celui-là,
je dois à la fin voir émerger une `entièreté' harmonieuse, qui ne sera
pas défigurée par une apparence bigarrée\footnotemark.
\REFERENCE{\cite{gauss1981},~vol.~8,~399}}

\footnotetext{\,
%%%%%%%%%%%%%%%%%%%%%%%-------DEBUT--------%%%%%%%%%%%%%%%%%%%%%%%%%%%
So wie die mathematische Seite einer Arbeit mir gewöhlich
die interessanteste ist, so kann ich auch von der andern
Seite nicht leugnen, dass ich, um an einer ausgedehnten Arbeit
Freude zu haben, doch
am Ende ein schön organisirtes Ganzes muss hervorgehen sehen,
was durch ain zu buntscheckiges Ansehen nicht verunstaltet wird.
} %%%%%%%%%%%%%%%%%%%%%%%%-----FIN-----%%%%%%%%%%%%%%%%%%%%%%%%%%%%%%%

\medskip\noindent{\bf Publier, ou ne pas publier.}
\label{publish-not-publish}
Schumacher et d'autres correspondants de Gauss l'ont régulièrement
encouragé, vu son âge grandissant, à rendre publiques ses nombreuses
idées pour les générations futures et à accepter de leur offrir du
travail d'exploration et d'approfondissement sur des sujets dans
lesquels il avait déjà considérablement progressé. \`A ce sujet,
Stäckel cite une lettre écrite par Gauss le 5 février 1850 à l'âge de
soixante-treize ans qui révèle à quel point ce dernier était
profondément investi par sa quête de la perfection des contenus et
dans laquelle s'exprime l'importance de maintenir ouverte toute
recherche tant qu'elle n'a pas abouti.

\medskip\noindent{\bf Nihil actum reputans si quid superesset agendum.}
\label{nihil-superesset-agendum}
Ainsi l'aspect le plus profond de la philosophie des mathématiques
de Gauss exprime-t-il que la quête de l'irréversible-synthétique
doit incontournablement 
s'étendre sur plusieurs années d'une vie de recherche, 
voire sur plusieurs générations de mathématiciens.

\CITATION{
Mes recherches sont devenues extrêmement difficiles pour moi à cause
du désir que j'ai toujours eu de leur conférer un degré de perfection
{\em ut nihil amplius desiderari possit} [{\em i.e.} telle que rien de
plus ample en compréhension ne puisse être désiré]\footnotemark.
\REFERENCE{\cite{gauss1981},~vol.~8,~400}}

\footnotetext{\,
%%%%%%%%%%%%%%%%%%%%%%%-------DEBUT--------%%%%%%%%%%%%%%%%%%%%%%%%%%%
Der Wunsch, den ich immer bei meinen Arbeiten gehabt habe, 
ihnen eine solche Vollendung zu geben, ut nihil amplius desidari
possit, erschwert sie mir freilich ausserdordentlich, eben
so wie die Notwendigkeit, heterogener Sachen
wegen oft davon abspringen zu müssen.
} %%%%%%%%%%%%%%%%%%%%%%%%-----FIN-----%%%%%%%%%%%%%%%%%%%%%%%%%%%%%%%

\medskip\noindent{\bf Le levier symbolique.}
\label{levier-symbolique} 
\`A cause de la grande difficult\'e qu'il y a \`a exposer pleinement
la signification des notions g\'eom\'etriques, on se contente
généralement de présenter seulement de manière formelle ou axiomatique
les concepts de connexion, de dérivée covariante, de courbure, de
torsion, ainsi que les identités de type Bianchi auxquelles satisfont
les dérivées covariantes des tenseurs fondamentaux. Les
pr\'esentations résumées des quantit\'es $g_{ kl}$, $g^{ pq}$,
$\Gamma_{ ij}^k$, $R_{ ijk}^l$, $R_{ ij}$, $R$ et $E_{ ij}$ que l'on
s'autorise en relativité générale souffrent parfois de cette
imperfection. \`A vrai dire, par l'effet d'un penchant au calcul
auquel la pens\'ee est entraînée d\`es que le langage symbolique
s'introduit dans les raisonnements, l'<<\,essence du tensoriel\,>>
incite \`a l'<<\,alg\'ebrisation des contenus\,>>.

\CITATION{
Il est dans la nature des mathématiques des temps modernes (en
contraste avec celles de l'Antiquité) que nous possédons un levier
sous la forme de notre langage symbolique et de notre terminologie,
gr\^ace à quoi les raisonnements les plus complexes sont réduits à un
certain mécanisme. De cette manière, la science a gagné infiniment en
richesse, mais elle a autant perdu en beauté et en caractère, comme
on le fait ordinairement dans la pratique.
Combien ce levier est-il fréquemment appliqué de manière purement
mécanique! bien que l'autorisation de l'employer implique, dans la
plupart des cas, certaines hypothèses passées sous
silence\footnotemark.
\REFERENCE{\cite{gauss1981},~{\bf 10},~1,~p.~434}}

\footnotetext{\,
<<\,Est ist der Character der Mathematik der neueren Zeit (im
Gegensatz gegen das Alterthum), dass durch unsere Zeichensprache und
Namengebungen wir eine Hebel besitzen, wodurch die verwickeltsten
Argumentationen auf einen gewissen Mechanismus reducirt werden. An
Reichthum hat dadurch die Wissenschaft unendlich gewonnen, an
Sch\"onheit und Solidit\"at aber, wie das Gesch\"aft gew\"ohnlich
betrieben wird, eben so sehr verloren. Wie oft wird jener Hebel eben
nur mecanisch angewandt, obgleich die Befugniss dazu in den meisten
F\"allen gewisse stillschweigende Voraussetzungen implicirt.\,>>}

\medskip\noindent{\bf Die Gaussche Strenge.}
\label{die-Gaussche-strenge} 
Cependant, comme l'exigeait Gauss lui-même dans cette lettre du
1\textsuperscript{er} septembre 1850 adress\'ee \`a son fidèle
interlocuteur Schumacher (où il critiquait l'absence de rigueur
que Prehn manifestait dans un mémoire sur les séries divergentes paru
la même année dans le {\em Journal de Crelle}):

\medskip
\hfill\fbox{
\begin{minipage}[t]{11cm}
\baselineskip=0.43cm
{\small{\sf\blue{ 
Chaque fois que l'on utilise le calcul et chaque fois qu'on emploie
des concepts, j'exige que l'on reste toujours conscient des
stipulations originelles [urspr\"unglichen Bedingungen], et que tous
les résultats du mécanisme ne soient jamais considérés comme des
propriétés en dehors d'une autorisation 
claire\footnotemark.}
\REFERENCE{\cite{gauss1981},~{\bf 10},~1,~p.~434}}}
\end{minipage}}

\footnotetext{\,
<<\,Ich fordere, man soll bei allem Gebrauch des Calculs, bei allen
Begriffsverwendungen sich immer der urspr\"unglichen Bedingungen
bewusst bleiben, und alle Producte des Mechanismus niemals \"uber die
klare Befugniss hinaus als Eigenthum betrachten.\,>>
}

\medskip
Ainsi s'exprime la rigoureuse et austère sévérité mathématique de
Gauss \deutsch{die Gau{\ss}che Strenge}: jamais acte de calcul
ne doit être engagé qui ne soit encadré au préalable par des
conditions précises quant à l'extension de ses significations. Par
exemple, quel que soit le sens qu'on donne aux séries divergentes en
choisissant un procédé de sommation déterminé, de type Cesàro, Fejér
ou Toeplitz, il faut impérativement justifier par des propositions
dûment établies à partir des définitions initiales que l'on peut
effectuer toutes les opérations algébriques élémentaires sur les
séries ainsi <<\,apprivoisées\,>>, et notamment la multiplication; on
ne pourrait certainement pas se contenter d'exécuter ces opérations au
prétexte qu'elles sont justifiées pour les séries convergentes. Tout
calcul que l'on transfère graduellement à des objets qui s'\'el\`event
en généralité expose en effet à des non-sens \'eventuels. Ainsi
l'addition, la soustraction et la multiplication entre séries
divergentes doivent-elle conserver la mémoire des processus sommatoires
qui ont \'et\'e choisis pour leur donner un sens. C'est donc la
première interprétation immédiate de cette citation de Gauss:
nécessité de se conformer au sens interne; n\'ecessit\'e de respecter
les bornes d\'efinitionnelles; sous peine d'incohérence.

\medskip\noindent{\bf Science et conscience copr\'esente.}
\label{science-conscience-copresente}
Mais au-delà, Gauss évoque sans la développer une pensée plus profonde
qui nous est suggérée par le membre de phrase <<\,j'exige que l'on
reste toujours conscient des stipulations originelles\,>>. On sait
combien il est délicat et malaisé de répondre à cet impératif de
<<\,conscientisation continue du sens\,>>, surtout lorsqu'il s'agit de
calculs de type tensoriel qui sont assez longs ou considérables pour
égarer l'intuition géométrique ou noyer l'esprit de synthèse, un peu
comme si les <<\,masses-pensées\,>> chères à Riemann devaient se
réserver à tout instant l'énergie de <<\,se rendre présentes les
totalités conceptuelles\,>>; toutefois, l'austérité ou la rigueur
gaussienne \deutsch{die Gau{\ss}che Strenge} sont là pour nous
rappeler les {\sl impératifs catégoriques de la pensée
mathématique}, similaires aux impératifs que Kant théorisa dans sa
{\em Critique de la raison pratique}, c'est-à-dire des règles qui sont
constituées par des principes scientifiques abstraits et qui sont
désignées par des devoirs.

\medskip\noindent{\bf L'imp\'eratif cat\'egorique de la morale kantienne.}
\label{imperatif-categorique-morale}
D'après Kant, l'impératif moral est un {\sl impératif
catégorique}: il commande absolument la poursuite d'une fin morale, en
elle-même et pour elle-même, et il détermine la volonté en indiquant
une loi objective de la raison: la {\sl loi morale}, valable
universellement pour tout être raisonnable en tant que tel. D'une
portée inférieure, les {\sl impératifs hypothétiques} ne se
rapportent qu'à la nécessité pratique d'une action considérée comme
{\em moyen}\, de parvenir à quelque chose, la {\em fin}\, vis\'ee
incitant \`a recourir à la technique, à l'habileté, à la prudence, à
la méthode et au pragmatisme. Quant à lui, l'impératif catégorique
représente une action comme nécessaire pour elle-même, cette action
n'étant subordonnée en tant que moyen à aucune fin déterminante
étrangère à son principe. Il commande de se conformer aux actions qui
sont bonnes en elles-mêmes, et par l'effet d'une nécessité
inconditionnée et véritablement objective, il soumet la volonté à une
{\em loi morale}\, interne et autonome que Kant énonce comme
suit: <<\,{\em Agis uniquement d'après la maxime qui fait que tu
peux vouloir en même temps qu'elle devienne une loi universelle parmi
les hommes}\,>>. Est donc {\em immorale}\, toute action qui, si
on la supposait perp\'etr\'ee par tous les hommes universellement en
même temps, conduirait \`a une perte dommageable d'\'equilibre de la
communaut\'e des hommes, \`a un chaos des m{\oe}urs
et des travaux, aux conflits, au
crime, \`a la mort. L'imp\'eratif cat\'egorique kantien \'enonce donc
une r\`egle absolue et universelle propre \`a d\'eterminer l'action
morale en toutes circonstances: il suffit de tester la moralité d'une
action en s'imaginant les cons\'equences d'une universalisation pour
se d\'ecider en cons\'equence. Ainsi, le cat\'egorique domine et se
ramifie dans l'hypoth\'etique.

\medskip\noindent{\bf Impératifs catégoriques de la pensée mathématique.}
\label{imperatifs-categoriques-mathematiques}
En math\'ematiques, l'imp\'eratif cat\'egorique\,\,---\,\,si tant est
qu'on puisse lui donner un sens tout analogue au sens
kantien\,\,---\,\,doit n\'ecessairement s'exprimer dans une pure
abstraction immanente, parce que la mati\`ere même des math\'ematiques
transcende les conditions biologiques ou neuronales de son exploration
effective\footnote{\, D'autres langages impr\'evisibles ou approches
dont la teneur reste insoupçonn\'ee se r\'ev\`eleront probablement
plus ad\'equats pour le traitement et la compr\'ehension de probl\`emes
encore ouverts \`a ce jour: tel est le {\em credo}\, fondamental en
l'être-disponible de la chose math\'ematique que partagent
tous les math\'ematiciens. }. De plus
les imp\'eratifs se d\'emultiplient et se pluralisent pour se disposer
en un tripôle fondamental abstrait qui n'influe pas aussi directement
sur la volont\'e que dans le champ de la raison pratique.

\bigskip$\square$\ \ \
D'abord l'exigence absolue de {\em coh\'erence}:
non-contradiction, vigilance architecturale, et v\'erit\'e
des raisonnements.

\bigskip$\square$\ \ \
Ensuite, le devoir p\'erenne et ind\'efectible de {\em
recherche}: ind\'efini potentiel, ouverture dynamique, et {\em in
situ}\, de l'Inconnu.

\bigskip$\square$\ \ \
Enfin, l'admission de la {\em nouveaut\'e}\, comme crit\`ere et
r\`egle de participation: absorption effective, progr\`es
(im)perceptibles, et enrichissement croissant des arborescences.

\bigskip
Trois imp\'eratifs, donc, inscrits dans une relation triangulaire o\`u
les hi\'erarchies sont interchangeables et o\`u les connexions sont
cycliques. Au centre du cercle circonscrit, l'exigence gaussienne de
pr\'esence permanente de la pens\'ee conceptuelle comme conscience:
li\'ee \`a chacun des trois pôles, elle exerce sa capacit\'e
d'activation en mesurant sa force au champ opaque de l'ouverture.
 
\medskip$\square$\ \ \
V\'erification constante de rigueur.

\medskip$\square$\ \ \
Maintien absolu des questions ind\'ecid\'ees.

\medskip$\square$\ \ \
\'Evaluation conceptuelle des apports.

\linestop

%%%%%%%%%%%%%%%%%%%%%%%%%%%%%%%%%%%%%%%%%%%%%%%%%%%%%%%%%%%%%%%%%%%%%%

%%%%%%%%%%%%%%%%%%%%%%%%%%%%%%%%%%%%%%%%%%%%%%%%%%%%%%%%%%%%%%%%%%%%%

\vfill\end{document}

%% file: macros.tex
\newcommand{\B}{\mathbb{B}}\newcommand{\C}{\mathbb{C}}
\newcommand{\E}{\mathbb{E}}

\newcommand{\R}{\mathbb{R}}

\definecolor{black}{cmyk}{1.,1.,1.,1.}
\definecolor{blue}{cmyk}{1.,1.,0.,0.43}
\definecolor{red}{cmyk}{0.,1.,1.,0.43}
\definecolor{green}{cmyk}{1.,0.,1.,0.43}

\newcommand{\blue}{\textcolor{blue}}
\newcommand{\green}{\textcolor{green}}
\newcommand{\red}{\textcolor{red}}

\makeatletter
\renewcommand{\@fnsymbol}[1]
{\ensuremath{\ifcase#1\or $*$\or $**$\or $***$\or $****$\or $*****$
\else\@ctrerr\fi}}
\makeatother

\newcommand{\Section}[1]{\bigskip\begin{center}{\bf #1}\end{center}}

\newcommand{\HEAD}[2]{%
\pagestyle{fancy}
\fancyhead[RO]{\scriptsize\sf\thepage}
\fancyhead[CO]{{\scriptsize\sf #1}}
\fancyhead[LE]{\scriptsize\sf\thepage}
\fancyhead[CE]{{\scriptsize\sf #2}}
\fancyfoot{}}

\newcommand{\CITATION}[1]{\smallskip\hfill
\begin{minipage}[t]{12cm}\baselineskip=0.43cm\parindent=0.91cm
\blue{\footnotesize{\sf \!\!\!\!\!\!\!#1}}\end{minipage}\medskip}
\newcommand{\REFERENCE}[1]{\hfill \green{#1}}

\newcommand{\SHORTCITATION}[1]{\hfill
\begin{minipage}[t]{12cm}\baselineskip=0.43cm\parindent=0.91cm
\blue{\footnotesize{\sf #1}}\end{minipage}}

\newcommand{\bignorm}{\big\vert\!\big\vert}

\newcommand{\deutsch}[1]{[{\sl \!\!\;#1}]}

\newcommand{\deutschplain}[1]{{\sl #1}}

\newcommand{\linestop}{\medskip\centerline{\bf 
-----------------}\medskip}

\newcommand{\norm}{\vert\!\vert}

\newcommand{\Scholie}{\medskip\noindent{\em Scholie.} }

\newcommand{\vf}{\vfill\end{document}}

\newcommand{\zero}[1]{\underline{#1}_{\red{\circ}}}

\newcommand{\zerozero}[1]{\underline{#1}_{\red{\circ\circ}}}

\def\boiteepaisseavecuntitre#1{%
  \def\thickhrulefill{\leavevmode \leaders \hrule height 1pt \hfill \kern \z@}%
  \def\bkvz@before@breakbox{\ifhmode\par\fi\vskip\breakboxskip\relax}%
  \fboxrule=1pt
  \def\bkvz@set@linewidth{\advance\linewidth -2\fboxrule
                          \advance\linewidth -2\fboxsep}%
  \def\bkvz@left{\vrule \@width\fboxrule\hskip\fboxsep}%
  \def\bkvz@right{\hskip\fboxsep\vrule \@width\fboxrule}%
  \def\bkvz@top{\hbox to \hsize{%
      \vrule\@width\fboxrule\@height 1.2pt %%% D'où vient ce 0.2pt ????
      \thickhrulefill{#1}\thickhrulefill
      \vrule\@width\fboxrule\@height 1.2pt}}%
  \def\bkvz@bottom{\hrule\@height\fboxrule}%
  \breakbox}

\newcommand{\encadre}{\smallskip\begin{boiteepaisseavecuntitre}{}
\normalsize\baselineskip=0.53cm\fboxrule=0.77pt}
\newcommand{\stopencadre}{\end{boiteepaisseavecuntitre}}

%% file: application-de-Gauss.pstex_t
\begin{picture}(0,0)%
\includegraphics{application-de-Gauss.pstex}%
\end{picture}%
\setlength{\unitlength}{4144sp}%
\begingroup\makeatletter\ifx\SetFigFont\undefined%
\gdef\SetFigFont#1#2#3#4#5{%
  \reset@font\fontsize{#1}{#2pt}%
  \fontfamily{#3}\fontseries{#4}\fontshape{#5}%
  \selectfont}%
\fi\endgroup%
\begin{picture}(5418,2047)(1022,-2040)
\put(1044,-1125){\makebox(0,0)[lb]{\smash{{\SetFigFont{10}{12.0}{\familydefault}{\mddefault}{\updefault}{\color[rgb]{0,0,0}$S$}%
}}}}
\put(2754,-842){\makebox(0,0)[lb]{\smash{{\SetFigFont{10}{12.0}{\familydefault}{\mddefault}{\updefault}{\color[rgb]{0,0,0}$\mathcal{A}_S$}%
}}}}
\put(2894,-516){\makebox(0,0)[lb]{\smash{{\SetFigFont{7}{8.4}{\familydefault}{\mddefault}{\updefault}{\color[rgb]{0,0,0}$p$}%
}}}}
\put(2999,-601){\makebox(0,0)[lb]{\smash{{\SetFigFont{6}{7.2}{\familydefault}{\mddefault}{\updefault}{\color[rgb]{0,0,0}$q$}%
}}}}
\put(3053,-689){\makebox(0,0)[lb]{\smash{{\SetFigFont{6}{7.2}{\familydefault}{\mddefault}{\updefault}{\color[rgb]{0,0,0}$q$}%
}}}}
\put(2649,-597){\makebox(0,0)[lb]{\smash{{\SetFigFont{6}{7.2}{\familydefault}{\mddefault}{\updefault}{\color[rgb]{0,0,0}$q$}%
}}}}
\put(2761,-554){\makebox(0,0)[lb]{\smash{{\SetFigFont{6}{7.2}{\familydefault}{\mddefault}{\updefault}{\color[rgb]{0,0,0}$q$}%
}}}}
\put(4487,-710){\makebox(0,0)[lb]{\smash{{\SetFigFont{8}{9.6}{\familydefault}{\mddefault}{\updefault}{\color[rgb]{0,0,0}Auxiliary unit sphere}%
}}}}
\put(4648,-1895){\makebox(0,0)[lb]{\smash{{\SetFigFont{10}{12.0}{\familydefault}{\mddefault}{\updefault}{\color[rgb]{0,0,0}$\S^2$}%
}}}}
\put(4420,-1236){\makebox(0,0)[lb]{\smash{{\SetFigFont{10}{12.0}{\familydefault}{\mddefault}{\updefault}{\color[rgb]{0,0,0}$\mathcal{A}_{\S^2}$}%
}}}}
\end{picture}%

%% file: plan.pstex_t
\begin{picture}(0,0)%
\includegraphics{plan.pstex}%
\end{picture}%
\setlength{\unitlength}{4144sp}%
\begingroup\makeatletter\ifx\SetFigFont\undefined%
\gdef\SetFigFont#1#2#3#4#5{%
  \reset@font\fontsize{#1}{#2pt}%
  \fontfamily{#3}\fontseries{#4}\fontshape{#5}%
  \selectfont}%
\fi\endgroup%
\begin{picture}(4704,1495)(248,-1427)
\put(1533,-1087){\makebox(0,0)[lb]{\smash{{\SetFigFont{8}{9.6}{\familydefault}{\mddefault}{\updefault}{\color[rgb]{0,0,0}$0$}%
}}}}
\put(2433,-1262){\makebox(0,0)[lb]{\smash{{\SetFigFont{8}{9.6}{\familydefault}{\mddefault}{\updefault}{\color[rgb]{0,0,0}$x$}%
}}}}
\put(1343,-642){\makebox(0,0)[lb]{\smash{{\SetFigFont{8}{9.6}{\familydefault}{\mddefault}{\updefault}{\color[rgb]{0,0,0}$y$}%
}}}}
\put(2563,-632){\makebox(0,0)[lb]{\smash{{\SetFigFont{8}{9.6}{\familydefault}{\mddefault}{\updefault}{\color[rgb]{0,0,0}$(x, \, y)$}%
}}}}
\put(4006,-16){\makebox(0,0)[lb]{\smash{{\SetFigFont{8}{9.6}{\familydefault}{\mddefault}{\updefault}{\color[rgb]{0,0,0}Ligne droite}%
}}}}
\put(4276,-1051){\makebox(0,0)[lb]{\smash{{\SetFigFont{8}{9.6}{\familydefault}{\mddefault}{\updefault}{\color[rgb]{0,0,0}Ligne courbe}%
}}}}
\put(946,-1387){\makebox(0,0)[lb]{\smash{{\SetFigFont{8}{9.6}{\familydefault}{\mddefault}{\updefault}{\color[rgb]{0,0,0}{\sc Plan math\'ematique id\'eal  et syst\`eme de coordonn\'ees}}%
}}}}
\end{picture}%

%% file: equation.pstex_t
\begin{picture}(0,0)%
\includegraphics{equation.pstex}%
\end{picture}%
\setlength{\unitlength}{4144sp}%
\begingroup\makeatletter\ifx\SetFigFont\undefined%
\gdef\SetFigFont#1#2#3#4#5{%
  \reset@font\fontsize{#1}{#2pt}%
  \fontfamily{#3}\fontseries{#4}\fontshape{#5}%
  \selectfont}%
\fi\endgroup%
\begin{picture}(3812,1662)(533,-1540)
\put(1038,-1179){\makebox(0,0)[lb]{\smash{{\SetFigFont{8}{9.6}{\familydefault}{\mddefault}{\updefault}{\color[rgb]{0,0,0}$0$}%
}}}}
\put(1053,-144){\makebox(0,0)[lb]{\smash{{\SetFigFont{8}{9.6}{\familydefault}{\mddefault}{\updefault}{\color[rgb]{0,0,0}$y$}%
}}}}
\put(2273,-1134){\makebox(0,0)[lb]{\smash{{\SetFigFont{8}{9.6}{\familydefault}{\mddefault}{\updefault}{\color[rgb]{0,0,0}$x$}%
}}}}
\put(533, 26){\makebox(0,0)[lb]{\smash{{\SetFigFont{8}{9.6}{\familydefault}{\mddefault}{\updefault}{\color[rgb]{0,0,0}{\bf I~: \'Equation}}%
}}}}
\put(1060,-1504){\makebox(0,0)[lb]{\smash{{\SetFigFont{8}{9.6}{\familydefault}{\mddefault}{\updefault}{\color[rgb]{0,0,0}{\sc Courbe d\'efinie par une \'equation $W(x,\, y)= 0$}}%
}}}}
\put(2937,-556){\makebox(0,0)[lb]{\smash{{\SetFigFont{8}{9.6}{\familydefault}{\mddefault}{\updefault}{\color[rgb]{0,0,0}$(x_1,\, y_1)$}%
}}}}
\put(2304,-221){\makebox(0,0)[lb]{\smash{{\SetFigFont{8}{9.6}{\familydefault}{\mddefault}{\updefault}{\color[rgb]{0,0,0}$(x_p,\, y_p)$ tel que}%
}}}}
\put(2303,-360){\makebox(0,0)[lb]{\smash{{\SetFigFont{8}{9.6}{\familydefault}{\mddefault}{\updefault}{\color[rgb]{0,0,0}$W(x_p,\, y_p)= 0$}%
}}}}
\put(3701,-497){\makebox(0,0)[lb]{\smash{{\SetFigFont{8}{9.6}{\familydefault}{\mddefault}{\updefault}{\color[rgb]{0,0,0}$\Delta_q$}%
}}}}
\put(3887,-703){\makebox(0,0)[lb]{\smash{{\SetFigFont{8}{9.6}{\familydefault}{\mddefault}{\updefault}{\color[rgb]{0,0,0}$q$}%
}}}}
\end{picture}%

%% file: graphe.pstex_t
\begin{picture}(0,0)%
\includegraphics{graphe.pstex}%
\end{picture}%
\setlength{\unitlength}{4144sp}%
\begingroup\makeatletter\ifx\SetFigFont\undefined%
\gdef\SetFigFont#1#2#3#4#5{%
  \reset@font\fontsize{#1}{#2pt}%
  \fontfamily{#3}\fontseries{#4}\fontshape{#5}%
  \selectfont}%
\fi\endgroup%
\begin{picture}(4649,1677)(452,-1555)
\put(1053,-144){\makebox(0,0)[lb]{\smash{{\SetFigFont{8}{9.6}{\familydefault}{\mddefault}{\updefault}{\color[rgb]{0,0,0}$y$}%
}}}}
\put(2273,-1134){\makebox(0,0)[lb]{\smash{{\SetFigFont{8}{9.6}{\familydefault}{\mddefault}{\updefault}{\color[rgb]{0,0,0}$x$}%
}}}}
\put(533, 26){\makebox(0,0)[lb]{\smash{{\SetFigFont{8}{9.6}{\familydefault}{\mddefault}{\updefault}{\color[rgb]{0,0,0}{\bf III~:  Graphe}}%
}}}}
\put(1190,-1519){\makebox(0,0)[lb]{\smash{{\SetFigFont{8}{9.6}{\familydefault}{\mddefault}{\updefault}{\color[rgb]{0,0,0}{\sc Courbe d\'efinie par un graphe $y= f(x)$}}%
}}}}
\put(3476,-806){\makebox(0,0)[lb]{\smash{{\SetFigFont{8}{9.6}{\familydefault}{\mddefault}{\updefault}{\color[rgb]{0,0,0}$p$}%
}}}}
\put(3573,-499){\makebox(0,0)[lb]{\smash{{\SetFigFont{8}{9.6}{\familydefault}{\mddefault}{\updefault}{\color[rgb]{0,0,0}$\Delta_p$}%
}}}}
\put(3190,-228){\makebox(0,0)[lb]{\smash{{\SetFigFont{8}{9.6}{\familydefault}{\mddefault}{\updefault}{\color[rgb]{0,0,0}$(x_p, f(x_p))$}%
}}}}
\put(3254,-1316){\makebox(0,0)[lb]{\smash{{\SetFigFont{8}{9.6}{\familydefault}{\mddefault}{\updefault}{\color[rgb]{0,0,0}$(x_p, 0)$}%
}}}}
\put(452,-843){\makebox(0,0)[lb]{\smash{{\SetFigFont{6}{7.2}{\familydefault}{\mddefault}{\updefault}{\color[rgb]{0,0,0}$(0,\, f(x_p))$}%
}}}}
\end{picture}%

%% file: parametrisation.pstex_t
\begin{picture}(0,0)%
\includegraphics{parametrisation.pstex}%
\end{picture}%
\setlength{\unitlength}{4144sp}%
\begingroup\makeatletter\ifx\SetFigFont\undefined%
\gdef\SetFigFont#1#2#3#4#5{%
  \reset@font\fontsize{#1}{#2pt}%
  \fontfamily{#3}\fontseries{#4}\fontshape{#5}%
  \selectfont}%
\fi\endgroup%
\begin{picture}(4451,1635)(533,-1525)
\put(533, 26){\makebox(0,0)[lb]{\smash{{\SetFigFont{8}{9.6}{\familydefault}{\mddefault}{\updefault}{\color[rgb]{0,0,0}{\bf II~: Param\'etrisation}}%
}}}}
\put(991,-133){\makebox(0,0)[lb]{\smash{{\SetFigFont{8}{9.6}{\familydefault}{\mddefault}{\updefault}{\color[rgb]{0,0,0}$y$}%
}}}}
\put(1741,-557){\makebox(0,0)[lb]{\smash{{\SetFigFont{8}{9.6}{\familydefault}{\mddefault}{\updefault}{\color[rgb]{0,0,0}$(x(t),\, y(t))$}%
}}}}
\put(3979,-595){\makebox(0,0)[lb]{\smash{{\SetFigFont{8}{9.6}{\familydefault}{\mddefault}{\updefault}{\color[rgb]{0,0,0}$\Delta_p$}%
}}}}
\put(4023,-856){\makebox(0,0)[lb]{\smash{{\SetFigFont{8}{9.6}{\familydefault}{\mddefault}{\updefault}{\color[rgb]{0,0,0}$p$}%
}}}}
\put(980,-1198){\makebox(0,0)[lb]{\smash{{\SetFigFont{8}{9.6}{\familydefault}{\mddefault}{\updefault}{\color[rgb]{0,0,0}$0$}%
}}}}
\put(2358,-1168){\makebox(0,0)[lb]{\smash{{\SetFigFont{8}{9.6}{\familydefault}{\mddefault}{\updefault}{\color[rgb]{0,0,0}$x$}%
}}}}
\put(1322,-1489){\makebox(0,0)[lb]{\smash{{\SetFigFont{8}{9.6}{\familydefault}{\mddefault}{\updefault}{\color[rgb]{0,0,0}{\sc Courbe d\'efinie par une param\'etrisation $(x(t),\, y(t))$}}%
}}}}
\end{picture}%

%% file: univers-courbes.pstex_t
\begin{picture}(0,0)%
\includegraphics{univers-courbes.pstex}%
\end{picture}%
\setlength{\unitlength}{4144sp}%
\begingroup\makeatletter\ifx\SetFigFont\undefined%
\gdef\SetFigFont#1#2#3#4#5{%
  \reset@font\fontsize{#1}{#2pt}%
  \fontfamily{#3}\fontseries{#4}\fontshape{#5}%
  \selectfont}%
\fi\endgroup%
\begin{picture}(4115,1758)(924,-1110)
\put(2206,552){\makebox(0,0)[lb]{\smash{{\SetFigFont{7}{8.4}{\familydefault}{\mddefault}{\updefault}{\color[rgb]{0,0,0}Ensemble des fonctions $W(x,\, y)$}%
}}}}
\put(2214,424){\makebox(0,0)[lb]{\smash{{\SetFigFont{7}{8.4}{\familydefault}{\mddefault}{\updefault}{\color[rgb]{0,0,0}satisfaisant $\frac{\partial W}{\partial Y} \neq 0$ ou $\frac{\partial W}{\partial x} \neq 0$}%
}}}}
\put(3901,-653){\makebox(0,0)[lb]{\smash{{\SetFigFont{7}{8.4}{\familydefault}{\mddefault}{\updefault}{\color[rgb]{0,0,0}Ensemble des couples}%
}}}}
\put(3896,-768){\makebox(0,0)[lb]{\smash{{\SetFigFont{7}{8.4}{\familydefault}{\mddefault}{\updefault}{\color[rgb]{0,0,0}de fonctions $(x(t),\, y(t))$}%
}}}}
\put(2186,-169){\makebox(0,0)[lb]{\smash{{\SetFigFont{7}{8.4}{\familydefault}{\mddefault}{\updefault}{\color[rgb]{0,0,0}{\bf Ensemble des courbes}}%
}}}}
\put(2187,-282){\makebox(0,0)[lb]{\smash{{\SetFigFont{7}{8.4}{\familydefault}{\mddefault}{\updefault}{\color[rgb]{0,0,0}{\bf diff\'erentiables du plan}}%
}}}}
\put(935,-660){\makebox(0,0)[lb]{\smash{{\SetFigFont{7}{8.4}{\familydefault}{\mddefault}{\updefault}{\color[rgb]{0,0,0}Ensemble des}%
}}}}
\put(924,-768){\makebox(0,0)[lb]{\smash{{\SetFigFont{7}{8.4}{\familydefault}{\mddefault}{\updefault}{\color[rgb]{0,0,0}fonctions $f(x)$ ou $g(y)$}%
}}}}
\put(1351,-1079){\makebox(0,0)[lb]{\smash{{\SetFigFont{8}{9.6}{\familydefault}{\mddefault}{\updefault}{\color[rgb]{0,0,0}{\sc Univers des courbes locales et univers de fonctions}}%
}}}}
\end{picture}%

%% file: echantillon-surfaces.pstex_t
\begin{picture}(0,0)%
\includegraphics{echantillon-surfaces.pstex}%
\end{picture}%
\setlength{\unitlength}{4144sp}%
\begingroup\makeatletter\ifx\SetFigFont\undefined%
\gdef\SetFigFont#1#2#3#4#5{%
  \reset@font\fontsize{#1}{#2pt}%
  \fontfamily{#3}\fontseries{#4}\fontshape{#5}%
  \selectfont}%
\fi\endgroup%
\begin{picture}(4297,1992)(624,-1742)
\put(1603,-1697){\makebox(0,0)[lb]{\smash{{\SetFigFont{10}{12.0}{\familydefault}{\mddefault}{\updefault}{\color[rgb]{0,0,0}{\sc Quelques surfaces math\'ematiques}}%
}}}}
\end{picture}%

%% file: ellipsoide.pstex_t
\begin{picture}(0,0)%
\includegraphics{ellipsoide.pstex}%
\end{picture}%
\setlength{\unitlength}{4144sp}%
\begingroup\makeatletter\ifx\SetFigFont\undefined%
\gdef\SetFigFont#1#2#3#4#5{%
  \reset@font\fontsize{#1}{#2pt}%
  \fontfamily{#3}\fontseries{#4}\fontshape{#5}%
  \selectfont}%
\fi\endgroup%
\begin{picture}(2573,1517)(1296,-1199)
\put(3670,-434){\makebox(0,0)[lb]{\smash{{\SetFigFont{8}{9.6}{\familydefault}{\mddefault}{\updefault}{\color[rgb]{0,0,0}$x$}%
}}}}
\put(2521,-489){\makebox(0,0)[lb]{\smash{{\SetFigFont{8}{9.6}{\familydefault}{\mddefault}{\updefault}{\color[rgb]{0,0,0}$y$}%
}}}}
\put(2682,209){\makebox(0,0)[lb]{\smash{{\SetFigFont{8}{9.6}{\familydefault}{\mddefault}{\updefault}{\color[rgb]{0,0,0}$z$}%
}}}}
\put(2682,-81){\makebox(0,0)[lb]{\smash{{\SetFigFont{8}{9.6}{\familydefault}{\mddefault}{\updefault}{\color[rgb]{0,0,0}$c$}%
}}}}
\put(2553,-346){\makebox(0,0)[lb]{\smash{{\SetFigFont{8}{9.6}{\familydefault}{\mddefault}{\updefault}{\color[rgb]{0,0,0}$b$}%
}}}}
\put(3048,-400){\makebox(0,0)[lb]{\smash{{\SetFigFont{8}{9.6}{\familydefault}{\mddefault}{\updefault}{\color[rgb]{0,0,0}$a$}%
}}}}
\put(2676,-391){\makebox(0,0)[lb]{\smash{{\SetFigFont{8}{9.6}{\familydefault}{\mddefault}{\updefault}{\color[rgb]{0,0,0}$p_1$}%
}}}}
\put(1296,-1156){\makebox(0,0)[lb]{\smash{{\SetFigFont{9}{10.8}{\familydefault}{\mddefault}{\updefault}{\color[rgb]{0,0,0}{\sc Axes principaux d'un ellipso{\i}de g\'en\'eral}}%
}}}}
\end{picture}%

%% file: u-v-x-y-z.pstex_t
\begin{picture}(0,0)%
\includegraphics{u-v-x-y-z.pstex}%
\end{picture}%
\setlength{\unitlength}{4144sp}%
\begingroup\makeatletter\ifx\SetFigFont\undefined%
\gdef\SetFigFont#1#2#3#4#5{%
  \reset@font\fontsize{#1}{#2pt}%
  \fontfamily{#3}\fontseries{#4}\fontshape{#5}%
  \selectfont}%
\fi\endgroup%
\begin{picture}(3488,985)(812,-700)
\put(1441,-654){\makebox(0,0)[lb]{\smash{{\SetFigFont{9}{10.8}{\familydefault}{\mddefault}{\updefault}{\color[rgb]{0,0,0}{\sc Coordonn\'ees curvilignes $(u, v)$}}%
}}}}
\end{picture}%

%% file: theoreme-pythagore.pstex_t
\begin{picture}(0,0)%
\includegraphics{theoreme-pythagore.pstex}%
\end{picture}%
\setlength{\unitlength}{4144sp}%
\begingroup\makeatletter\ifx\SetFigFont\undefined%
\gdef\SetFigFont#1#2#3#4#5{%
  \reset@font\fontsize{#1}{#2pt}%
  \fontfamily{#3}\fontseries{#4}\fontshape{#5}%
  \selectfont}%
\fi\endgroup%
\begin{picture}(4295,1500)(731,-1091)
\put(1415,-857){\makebox(0,0)[lb]{\smash{{\SetFigFont{9}{10.8}{\familydefault}{\mddefault}{\updefault}{\color[rgb]{0,0,0}$b$}%
}}}}
\put(1541,-197){\makebox(0,0)[lb]{\smash{{\SetFigFont{9}{10.8}{\familydefault}{\mddefault}{\updefault}{\color[rgb]{0,0,0}$\sqrt{ a^ 2+ b^ 2}$}%
}}}}
\put(4780,-853){\makebox(0,0)[lb]{\smash{{\SetFigFont{9}{10.8}{\familydefault}{\mddefault}{\updefault}{\color[rgb]{0,0,0}$x$}%
}}}}
\put(2976,272){\makebox(0,0)[lb]{\smash{{\SetFigFont{9}{10.8}{\familydefault}{\mddefault}{\updefault}{\color[rgb]{0,0,0}$y$}%
}}}}
\put(4570,292){\makebox(0,0)[lb]{\smash{{\SetFigFont{9}{10.8}{\familydefault}{\mddefault}{\updefault}{\color[rgb]{0,0,0}$C$}%
}}}}
\put(731,-324){\makebox(0,0)[lb]{\smash{{\SetFigFont{9}{10.8}{\familydefault}{\mddefault}{\updefault}{\color[rgb]{0,0,0}$a$}%
}}}}
\put(3983,-745){\makebox(0,0)[lb]{\smash{{\SetFigFont{6}{7.2}{\familydefault}{\mddefault}{\updefault}{\color[rgb]{0,0,0}$dx$}%
}}}}
\put(4446,-346){\makebox(0,0)[lb]{\smash{{\SetFigFont{5}{6.0}{\familydefault}{\mddefault}{\updefault}{\color[rgb]{0,0,0}$dy= f'(x) \, dx$}%
}}}}
\put(3211,-257){\makebox(0,0)[lb]{\smash{{\SetFigFont{5}{6.0}{\familydefault}{\mddefault}{\updefault}{\color[rgb]{0,0,0}$\sqrt{ 1+ (f'(x))^ 2} \, dx$}%
}}}}
\put(1261,-1051){\makebox(0,0)[lb]{\smash{{\SetFigFont{8}{9.6}{\familydefault}{\mddefault}{\updefault}{\color[rgb]{0,0,0}{\sc Th\'eor\`eme de Pythagore et sa version infinit\'esimale}}%
}}}}
\end{picture}%

%% file: courbure-cercle.pstex_t
\begin{picture}(0,0)%
\includegraphics{courbure-cercle.pstex}%
\end{picture}%
\setlength{\unitlength}{4144sp}%
\begingroup\makeatletter\ifx\SetFigFont\undefined%
\gdef\SetFigFont#1#2#3#4#5{%
  \reset@font\fontsize{#1}{#2pt}%
  \fontfamily{#3}\fontseries{#4}\fontshape{#5}%
  \selectfont}%
\fi\endgroup%
\begin{picture}(4774,1437)(500,-1081)
\put(3321,-800){\makebox(0,0)[lb]{\smash{{\SetFigFont{7}{8.4}{\familydefault}{\mddefault}{\updefault}{\color[rgb]{0,0,0}droite~: rayon infini~; courbure nulle}%
}}}}
\put(1222,236){\makebox(0,0)[lb]{\smash{{\SetFigFont{10}{12.0}{\familydefault}{\mddefault}{\updefault}{\color[rgb]{0,0,0}{\sl Courbure} $:=$ {\bf  inverse} du rayon = $\frac{ 1}{ R}$}%
}}}}
\put(500,-738){\makebox(0,0)[lb]{\smash{{\SetFigFont{8}{9.6}{\familydefault}{\mddefault}{\updefault}{\color[rgb]{0,0,0}courbure grande}%
}}}}
\put(508,-626){\makebox(0,0)[lb]{\smash{{\SetFigFont{8}{9.6}{\familydefault}{\mddefault}{\updefault}{\color[rgb]{0,0,0}rayon petit}%
}}}}
\put(3808,-281){\makebox(0,0)[lb]{\smash{{\SetFigFont{7}{8.4}{\familydefault}{\mddefault}{\updefault}{\color[rgb]{0,0,0}rayon grand}%
}}}}
\put(3800,-371){\makebox(0,0)[lb]{\smash{{\SetFigFont{7}{8.4}{\familydefault}{\mddefault}{\updefault}{\color[rgb]{0,0,0}courbure petite}%
}}}}
\put(2359,-302){\makebox(0,0)[lb]{\smash{{\SetFigFont{8}{9.6}{\familydefault}{\mddefault}{\updefault}{\color[rgb]{0,0,0}$R$}%
}}}}
\put(2115,  7){\makebox(0,0)[lb]{\smash{{\SetFigFont{8}{9.6}{\familydefault}{\mddefault}{\updefault}{\color[rgb]{0,0,0}$\mathcal{ C}$}%
}}}}
\put(2127,-430){\makebox(0,0)[lb]{\smash{{\SetFigFont{6}{7.2}{\familydefault}{\mddefault}{\updefault}{\color[rgb]{0,0,0}$O$}%
}}}}
\put(1707,-1045){\makebox(0,0)[lb]{\smash{{\SetFigFont{8}{9.6}{\familydefault}{\mddefault}{\updefault}{\color[rgb]{0,0,0}{\sc Courbure des cercles dans le plan}}%
}}}}
\end{picture}%

%% file: angle-courbure.pstex_t
\begin{picture}(0,0)%
\includegraphics{angle-courbure.pstex}%
\end{picture}%
\setlength{\unitlength}{4144sp}%
\begingroup\makeatletter\ifx\SetFigFont\undefined%
\gdef\SetFigFont#1#2#3#4#5{%
  \reset@font\fontsize{#1}{#2pt}%
  \fontfamily{#3}\fontseries{#4}\fontshape{#5}%
  \selectfont}%
\fi\endgroup%
\begin{picture}(4284,1378)(513,-1022)
\put(513,237){\makebox(0,0)[lb]{\smash{{\SetFigFont{8}{9.6}{\familydefault}{\mddefault}{\updefault}{\color[rgb]{0,0,0}Segments align\'es~:}%
}}}}
\put(517,111){\makebox(0,0)[lb]{\smash{{\SetFigFont{8}{9.6}{\familydefault}{\mddefault}{\updefault}{\color[rgb]{0,0,0}pas de courbure}%
}}}}
\put(1962,228){\makebox(0,0)[lb]{\smash{{\SetFigFont{8}{9.6}{\familydefault}{\mddefault}{\updefault}{\color[rgb]{0,0,0}Variation d'angle}%
}}}}
\put(3351,222){\makebox(0,0)[lb]{\smash{{\SetFigFont{8}{9.6}{\familydefault}{\mddefault}{\updefault}{\color[rgb]{0,0,0}Courbure}%
}}}}
\put(1081,-991){\makebox(0,0)[lb]{\smash{{\SetFigFont{8}{9.6}{\familydefault}{\mddefault}{\updefault}{\color[rgb]{0,0,0}{\sc Courbure de segments articul\'es et variation d'angle}}%
}}}}
\end{picture}%

%% file: normales-tangentes.pstex_t
\begin{picture}(0,0)%
\includegraphics{normales-tangentes.pstex}%
\end{picture}%
\setlength{\unitlength}{4144sp}%
\begingroup\makeatletter\ifx\SetFigFont\undefined%
\gdef\SetFigFont#1#2#3#4#5{%
  \reset@font\fontsize{#1}{#2pt}%
  \fontfamily{#3}\fontseries{#4}\fontshape{#5}%
  \selectfont}%
\fi\endgroup%
\begin{picture}(3993,1440)(768,-1207)
\put(3999,-795){\makebox(0,0)[lb]{\smash{{\SetFigFont{7}{8.4}{\familydefault}{\mddefault}{\updefault}{\color[rgb]{0,0,0}$p$}%
}}}}
\put(3388,-979){\makebox(0,0)[lb]{\smash{{\SetFigFont{7}{8.4}{\familydefault}{\mddefault}{\updefault}{\color[rgb]{0,0,0}$p$}%
}}}}
\put(1625,-431){\makebox(0,0)[lb]{\smash{{\SetFigFont{7}{8.4}{\familydefault}{\mddefault}{\updefault}{\color[rgb]{0,0,0}$p$}%
}}}}
\put(1296,-859){\makebox(0,0)[lb]{\smash{{\SetFigFont{7}{8.4}{\familydefault}{\mddefault}{\updefault}{\color[rgb]{0,0,0}$p$}%
}}}}
\put(2263,-357){\makebox(0,0)[lb]{\smash{{\SetFigFont{7}{8.4}{\familydefault}{\mddefault}{\updefault}{\color[rgb]{0,0,0}$p$}%
}}}}
\put(2713,-766){\makebox(0,0)[lb]{\smash{{\SetFigFont{7}{8.4}{\familydefault}{\mddefault}{\updefault}{\color[rgb]{0,0,0}$p$}%
}}}}
\put(806,-570){\makebox(0,0)[lb]{\smash{{\SetFigFont{7}{8.4}{\familydefault}{\mddefault}{\updefault}{\color[rgb]{0,0,0}$N(p)$}%
}}}}
\put(1368, 42){\makebox(0,0)[lb]{\smash{{\SetFigFont{7}{8.4}{\familydefault}{\mddefault}{\updefault}{\color[rgb]{0,0,0}$N(p)$}%
}}}}
\put(2523,113){\makebox(0,0)[lb]{\smash{{\SetFigFont{7}{8.4}{\familydefault}{\mddefault}{\updefault}{\color[rgb]{0,0,0}$N(p)$}%
}}}}
\put(3085,-412){\makebox(0,0)[lb]{\smash{{\SetFigFont{7}{8.4}{\familydefault}{\mddefault}{\updefault}{\color[rgb]{0,0,0}$N(p)$}%
}}}}
\put(3453,-468){\makebox(0,0)[lb]{\smash{{\SetFigFont{7}{8.4}{\familydefault}{\mddefault}{\updefault}{\color[rgb]{0,0,0}$N(p)$}%
}}}}
\put(3877,-280){\makebox(0,0)[lb]{\smash{{\SetFigFont{7}{8.4}{\familydefault}{\mddefault}{\updefault}{\color[rgb]{0,0,0}$N(p)$}%
}}}}
\put(3797,-962){\makebox(0,0)[lb]{\smash{{\SetFigFont{7}{8.4}{\familydefault}{\mddefault}{\updefault}{\color[rgb]{0,0,0}$T(p)$}%
}}}}
\put(1147,-411){\makebox(0,0)[lb]{\smash{{\SetFigFont{7}{8.4}{\familydefault}{\mddefault}{\updefault}{\color[rgb]{0,0,0}$T(p)$}%
}}}}
\put(2684,-344){\makebox(0,0)[lb]{\smash{{\SetFigFont{7}{8.4}{\familydefault}{\mddefault}{\updefault}{\color[rgb]{0,0,0}$T(p)$}%
}}}}
\put(2894,-1018){\makebox(0,0)[lb]{\smash{{\SetFigFont{7}{8.4}{\familydefault}{\mddefault}{\updefault}{\color[rgb]{0,0,0}$T(p)$}%
}}}}
\put(4408,-614){\makebox(0,0)[lb]{\smash{{\SetFigFont{7}{8.4}{\familydefault}{\mddefault}{\updefault}{\color[rgb]{0,0,0}$T(p)$}%
}}}}
\put(1817,-12){\makebox(0,0)[lb]{\smash{{\SetFigFont{7}{8.4}{\familydefault}{\mddefault}{\updefault}{\color[rgb]{0,0,0}$T(p)$}%
}}}}
\put(1311,-1176){\makebox(0,0)[lb]{\smash{{\SetFigFont{8}{9.6}{\familydefault}{\mddefault}{\updefault}{\color[rgb]{0,0,0}{\sc Vecteurs unitaires tangents et normaux}}%
}}}}
\end{picture}%

%% file: cercle-auxiliaire.pstex_t
\begin{picture}(0,0)%
\includegraphics{cercle-auxiliaire.pstex}%
\end{picture}%
\setlength{\unitlength}{4144sp}%
\begingroup\makeatletter\ifx\SetFigFont\undefined%
\gdef\SetFigFont#1#2#3#4#5{%
  \reset@font\fontsize{#1}{#2pt}%
  \fontfamily{#3}\fontseries{#4}\fontshape{#5}%
  \selectfont}%
\fi\endgroup%
\begin{picture}(4506,1978)(541,-1636)
\put(4171,248){\makebox(0,0)[lb]{\smash{{\SetFigFont{8}{9.6}{\familydefault}{\mddefault}{\updefault}{\color[rgb]{0,0,0}Cercle auxiliaire}%
}}}}
\put(541, 33){\makebox(0,0)[lb]{\smash{{\SetFigFont{8}{9.6}{\familydefault}{\mddefault}{\updefault}{\color[rgb]{0,0,0}Petit cercle}%
}}}}
\put(541,-102){\makebox(0,0)[lb]{\smash{{\SetFigFont{8}{9.6}{\familydefault}{\mddefault}{\updefault}{\color[rgb]{0,0,0}de rayon $r$}%
}}}}
\put(4473,-412){\makebox(0,0)[lb]{\smash{{\SetFigFont{5}{6.0}{\familydefault}{\mddefault}{\updefault}{\color[rgb]{0,0,0}$d \theta$}%
}}}}
\put(1038,-562){\makebox(0,0)[lb]{\smash{{\SetFigFont{5}{6.0}{\familydefault}{\mddefault}{\updefault}{\color[rgb]{0,0,0}$r d\theta$}%
}}}}
\put(4072,-788){\makebox(0,0)[lb]{\smash{{\SetFigFont{5}{6.0}{\familydefault}{\mddefault}{\updefault}{\color[rgb]{0,0,0}$O$}%
}}}}
\put(4168,123){\makebox(0,0)[lb]{\smash{{\SetFigFont{8}{9.6}{\familydefault}{\mddefault}{\updefault}{\color[rgb]{0,0,0}de rayon $1$}%
}}}}
\put(1774,246){\makebox(0,0)[lb]{\smash{{\SetFigFont{8}{9.6}{\familydefault}{\mddefault}{\updefault}{\color[rgb]{0,0,0}Grand cercle de rayon $R$}%
}}}}
\put(1468,-328){\makebox(0,0)[lb]{\smash{{\SetFigFont{5}{6.0}{\familydefault}{\mddefault}{\updefault}{\color[rgb]{0,0,0}$N(p)$}%
}}}}
\put(2200,-193){\makebox(0,0)[lb]{\smash{{\SetFigFont{5}{6.0}{\familydefault}{\mddefault}{\updefault}{\color[rgb]{0,0,0}$q$}%
}}}}
\put(3790,-727){\makebox(0,0)[lb]{\smash{{\SetFigFont{5}{6.0}{\familydefault}{\mddefault}{\updefault}{\color[rgb]{0,0,0}$N(p)$}%
}}}}
\put(3972,-405){\makebox(0,0)[lb]{\smash{{\SetFigFont{5}{6.0}{\familydefault}{\mddefault}{\updefault}{\color[rgb]{0,0,0}$N(q)$}%
}}}}
\put(1927,-498){\makebox(0,0)[lb]{\smash{{\SetFigFont{5}{6.0}{\familydefault}{\mddefault}{\updefault}{\color[rgb]{0,0,0}$p$}%
}}}}
\put(1770,-218){\makebox(0,0)[lb]{\smash{{\SetFigFont{5}{6.0}{\familydefault}{\mddefault}{\updefault}{\color[rgb]{0,0,0}$R d\varphi$}%
}}}}
\put(2016,110){\makebox(0,0)[lb]{\smash{{\SetFigFont{5}{6.0}{\familydefault}{\mddefault}{\updefault}{\color[rgb]{0,0,0}$N(q)$}%
}}}}
\put(3623,-420){\makebox(0,0)[lb]{\smash{{\SetFigFont{5}{6.0}{\familydefault}{\mddefault}{\updefault}{\color[rgb]{0,0,0}$d\varphi$}%
}}}}
\put(1201,-1605){\makebox(0,0)[lb]{\smash{{\SetFigFont{8}{9.6}{\familydefault}{\mddefault}{\updefault}{\color[rgb]{0,0,0}{\sc Courbure des cercles et vecteurs normaux unitaires}}%
}}}}
\end{picture}%

%% file: cercle-application-Gauss.pstex_t
\begin{picture}(0,0)%
\includegraphics{cercle-application-Gauss.pstex}%
\end{picture}%
\setlength{\unitlength}{4144sp}%
\begingroup\makeatletter\ifx\SetFigFont\undefined%
\gdef\SetFigFont#1#2#3#4#5{%
  \reset@font\fontsize{#1}{#2pt}%
  \fontfamily{#3}\fontseries{#4}\fontshape{#5}%
  \selectfont}%
\fi\endgroup%
\begin{picture}(4161,1677)(637,-1457)
\put(4253,-172){\makebox(0,0)[lb]{\smash{{\SetFigFont{8}{9.6}{\familydefault}{\mddefault}{\updefault}{\color[rgb]{0,0,0}$d\theta$}%
}}}}
\put(4317,-808){\makebox(0,0)[lb]{\smash{{\SetFigFont{8}{9.6}{\familydefault}{\mddefault}{\updefault}{\color[rgb]{0,0,0}$O$}%
}}}}
\put(1886,-437){\makebox(0,0)[lb]{\smash{{\SetFigFont{5}{6.0}{\familydefault}{\mddefault}{\updefault}{\color[rgb]{0,0,0}$p+dp$}%
}}}}
\put(1595,-441){\makebox(0,0)[lb]{\smash{{\SetFigFont{6}{7.2}{\familydefault}{\mddefault}{\updefault}{\color[rgb]{0,0,0}$p$}%
}}}}
\put(1442,112){\makebox(0,0)[lb]{\smash{{\SetFigFont{8}{9.6}{\familydefault}{\mddefault}{\updefault}{\color[rgb]{0,0,0}$N(p)$}%
}}}}
\put(1955,124){\makebox(0,0)[lb]{\smash{{\SetFigFont{8}{9.6}{\familydefault}{\mddefault}{\updefault}{\color[rgb]{0,0,0}$N(p+dp)$}%
}}}}
\put(1711,-235){\makebox(0,0)[lb]{\smash{{\SetFigFont{8}{9.6}{\familydefault}{\mddefault}{\updefault}{\color[rgb]{0,0,0}$d\ell$}%
}}}}
\put(2049,-950){\makebox(0,0)[lb]{\smash{{\SetFigFont{8}{9.6}{\familydefault}{\mddefault}{\updefault}{\color[rgb]{0,0,0}{\sl courbure en $p$} $=$ $\frac{d\theta}{d\ell}$}%
}}}}
\put(637,-1306){\makebox(0,0)[lb]{\smash{{\SetFigFont{8}{9.6}{\familydefault}{\mddefault}{\updefault}{\color[rgb]{0,0,0}$C$}%
}}}}
\put(3203,-235){\makebox(0,0)[lb]{\smash{{\SetFigFont{8}{9.6}{\familydefault}{\mddefault}{\updefault}{\color[rgb]{0,0,0}$C$}%
}}}}
\put(1392,-1421){\makebox(0,0)[lb]{\smash{{\SetFigFont{8}{9.6}{\familydefault}{\mddefault}{\updefault}{\color[rgb]{0,0,0}{\sc Application de Gauss et courbure des courbes}}%
}}}}
\end{picture}%

%% file: cercle-secant.pstex_t
\begin{picture}(0,0)%
\includegraphics{cercle-secant.pstex}%
\end{picture}%
\setlength{\unitlength}{4144sp}%
\begingroup\makeatletter\ifx\SetFigFont\undefined%
\gdef\SetFigFont#1#2#3#4#5{%
  \reset@font\fontsize{#1}{#2pt}%
  \fontfamily{#3}\fontseries{#4}\fontshape{#5}%
  \selectfont}%
\fi\endgroup%
\begin{picture}(4438,1732)(496,-1336)
\put(609,-1155){\makebox(0,0)[lb]{\smash{{\SetFigFont{8}{9.6}{\familydefault}{\mddefault}{\updefault}{\color[rgb]{0,0,0}$C$}%
}}}}
\put(2972,-322){\makebox(0,0)[lb]{\smash{{\SetFigFont{6}{7.2}{\familydefault}{\mddefault}{\updefault}{\color[rgb]{0,0,0}$O_p$}%
}}}}
\put(2905,300){\makebox(0,0)[lb]{\smash{{\SetFigFont{8}{9.6}{\familydefault}{\mddefault}{\updefault}{\color[rgb]{0,0,0}$p$}%
}}}}
\put(2597,242){\makebox(0,0)[lb]{\smash{{\SetFigFont{6}{7.2}{\familydefault}{\mddefault}{\updefault}{\color[rgb]{0,0,0}$q$}%
}}}}
\put(2222,-31){\makebox(0,0)[lb]{\smash{{\SetFigFont{6}{7.2}{\familydefault}{\mddefault}{\updefault}{\color[rgb]{0,0,0}$q$}%
}}}}
\put(2009,-430){\makebox(0,0)[lb]{\smash{{\SetFigFont{6}{7.2}{\familydefault}{\mddefault}{\updefault}{\color[rgb]{0,0,0}$q$}%
}}}}
\put(3284,203){\makebox(0,0)[lb]{\smash{{\SetFigFont{6}{7.2}{\familydefault}{\mddefault}{\updefault}{\color[rgb]{0,0,0}$r$}%
}}}}
\put(3551,-100){\makebox(0,0)[lb]{\smash{{\SetFigFont{6}{7.2}{\familydefault}{\mddefault}{\updefault}{\color[rgb]{0,0,0}$r$}%
}}}}
\put(3863,-457){\makebox(0,0)[lb]{\smash{{\SetFigFont{6}{7.2}{\familydefault}{\mddefault}{\updefault}{\color[rgb]{0,0,0}$r$}%
}}}}
\put(4474,-1044){\makebox(0,0)[lb]{\smash{{\SetFigFont{8}{9.6}{\familydefault}{\mddefault}{\updefault}{\color[rgb]{0,0,0}$\mathcal{ C}(p,\, q,\, r)$}%
}}}}
\put(496, 65){\makebox(0,0)[lb]{\smash{{\SetFigFont{7}{8.4}{\familydefault}{\mddefault}{\updefault}{\color[rgb]{0,0,0}l'unique cercle $\mathcal{ C}(p,\, q,\, r)$ passant par}%
}}}}
\put(496,-25){\makebox(0,0)[lb]{\smash{{\SetFigFont{7}{8.4}{\familydefault}{\mddefault}{\updefault}{\color[rgb]{0,0,0}les trois points $p$, $q$ et $r$ se rapproche}%
}}}}
\put(496,-115){\makebox(0,0)[lb]{\smash{{\SetFigFont{7}{8.4}{\familydefault}{\mddefault}{\updefault}{\color[rgb]{0,0,0}ind\'efiniment du cercle osculateur}%
}}}}
\put(496,-215){\makebox(0,0)[lb]{\smash{{\SetFigFont{7}{8.4}{\familydefault}{\mddefault}{\updefault}{\color[rgb]{0,0,0}\`a la courbe $C$  au point $p$}%
}}}}
\put(496,160){\makebox(0,0)[lb]{\smash{{\SetFigFont{7}{8.4}{\familydefault}{\mddefault}{\updefault}{\color[rgb]{0,0,0}de part et d'autre de $p$ tendent vers $p$,}%
}}}}
\put(496,260){\makebox(0,0)[lb]{\smash{{\SetFigFont{7}{8.4}{\familydefault}{\mddefault}{\updefault}{\color[rgb]{0,0,0}Lorsque les points distincts $q$ et $r$ situ\'es sur $C$}%
}}}}
\put(2957,-49){\makebox(0,0)[lb]{\smash{{\SetFigFont{6}{7.2}{\familydefault}{\mddefault}{\updefault}{\color[rgb]{0,0,0}$R_p$}%
}}}}
\put(1149,-1305){\makebox(0,0)[lb]{\smash{{\SetFigFont{8}{9.6}{\familydefault}{\mddefault}{\updefault}{\color[rgb]{0,0,0}{\sc Cercle osculateur comme limite de cercles s\'ecants}}%
}}}}
\end{picture}%

%% file: tp-np-projection.pstex_t
\begin{picture}(0,0)%
\includegraphics{tp-np-projection.pstex}%
\end{picture}%
\setlength{\unitlength}{4144sp}%
\begingroup\makeatletter\ifx\SetFigFont\undefined%
\gdef\SetFigFont#1#2#3#4#5{%
  \reset@font\fontsize{#1}{#2pt}%
  \fontfamily{#3}\fontseries{#4}\fontshape{#5}%
  \selectfont}%
\fi\endgroup%
\begin{picture}(4573,1956)(402,-1645)
\put(1477,215){\makebox(0,0)[lb]{\smash{{\SetFigFont{8}{9.6}{\familydefault}{\mddefault}{\updefault}{\color[rgb]{0,0,0}$y$}%
}}}}
\put(2026,-1028){\makebox(0,0)[lb]{\smash{{\SetFigFont{8}{9.6}{\familydefault}{\mddefault}{\updefault}{\color[rgb]{0,0,0}$T(x)$}%
}}}}
\put(1126,-543){\makebox(0,0)[lb]{\smash{{\SetFigFont{8}{9.6}{\familydefault}{\mddefault}{\updefault}{\color[rgb]{0,0,0}$N(x)$}%
}}}}
\put(2771,-1342){\makebox(0,0)[lb]{\smash{{\SetFigFont{8}{9.6}{\familydefault}{\mddefault}{\updefault}{\color[rgb]{0,0,0}$x$}%
}}}}
\put(2547,-944){\makebox(0,0)[lb]{\smash{{\SetFigFont{5}{6.0}{\familydefault}{\mddefault}{\updefault}{\color[rgb]{0,0,0}$\frac{f'}{\sqrt{1+f'{}^2}}$}%
}}}}
\put(402,-554){\makebox(0,0)[lb]{\smash{{\SetFigFont{5}{6.0}{\familydefault}{\mddefault}{\updefault}{\color[rgb]{0,0,0}$\frac{1}{\sqrt{1+f'{}^2}}$}%
}}}}
\put(951,-14){\makebox(0,0)[lb]{\smash{{\SetFigFont{5}{6.0}{\familydefault}{\mddefault}{\updefault}{\color[rgb]{0,0,0}$\frac{f'}{\sqrt{1+f'{}^2}}$}%
}}}}
\put(4256,-722){\makebox(0,0)[lb]{\smash{{\SetFigFont{6}{7.2}{\familydefault}{\mddefault}{\updefault}{\color[rgb]{0,0,0}$x+dx$}%
}}}}
\put(4115,-832){\makebox(0,0)[lb]{\smash{{\SetFigFont{6}{7.2}{\familydefault}{\mddefault}{\updefault}{\color[rgb]{0,0,0}$x$}%
}}}}
\put(4760,-623){\makebox(0,0)[lb]{\smash{{\SetFigFont{8}{9.6}{\familydefault}{\mddefault}{\updefault}{\color[rgb]{0,0,0}$C$}%
}}}}
\put(3262,-1137){\makebox(0,0)[lb]{\smash{{\SetFigFont{8}{9.6}{\familydefault}{\mddefault}{\updefault}{\color[rgb]{0,0,0}$C$}%
}}}}
\put(3470,103){\makebox(0,0)[lb]{\smash{{\SetFigFont{8}{9.6}{\familydefault}{\mddefault}{\updefault}{\color[rgb]{0,0,0}$dN(x)$}%
}}}}
\put(4122, 36){\makebox(0,0)[lb]{\smash{{\SetFigFont{8}{9.6}{\familydefault}{\mddefault}{\updefault}{\color[rgb]{0,0,0}$N(x+dx)$}%
}}}}
\put(3280,-501){\makebox(0,0)[lb]{\smash{{\SetFigFont{8}{9.6}{\familydefault}{\mddefault}{\updefault}{\color[rgb]{0,0,0}$N(x)$}%
}}}}
\put(1670,-491){\makebox(0,0)[lb]{\smash{{\SetFigFont{5}{6.0}{\familydefault}{\mddefault}{\updefault}{\color[rgb]{0,0,0}$\frac{1}{\sqrt{1+f'{}^2}}$}%
}}}}
\put(923,-1609){\makebox(0,0)[lb]{\smash{{\SetFigFont{8}{9.6}{\familydefault}{\mddefault}{\updefault}{\color[rgb]{0,0,0}{\sc Projections sur les axes de coordonn\'ees de $T(x)$ et de $N(x)$}}%
}}}}
\end{picture}%

%% file: articulation-chaine.pstex_t
\begin{picture}(0,0)%
\includegraphics{articulation-chaine.pstex}%
\end{picture}%
\setlength{\unitlength}{4144sp}%
\begingroup\makeatletter\ifx\SetFigFont\undefined%
\gdef\SetFigFont#1#2#3#4#5{%
  \reset@font\fontsize{#1}{#2pt}%
  \fontfamily{#3}\fontseries{#4}\fontshape{#5}%
  \selectfont}%
\fi\endgroup%
\begin{picture}(4087,1291)(671,-1277)
\put(1344,-1241){\makebox(0,0)[lb]{\smash{{\SetFigFont{8}{9.6}{\familydefault}{\mddefault}{\updefault}{\color[rgb]{0,0,0}{\sc Disparition de la courbure d'une cha\^{\i}ne articul\'ee}}%
}}}}
\end{picture}%